\documentclass[12pt,letterpaper]{amsart}

\usepackage{amsthm}

\usepackage{geometry}
\usepackage[english]{babel}
\usepackage[utf8]{inputenc}
\usepackage[T1]{fontenc}

\usepackage{amsmath,amssymb,amsthm,amscd}
\usepackage{mathtools}
\usepackage{mathrsfs}
\usepackage{calligra}

\usepackage{enumitem}
\usepackage{comment}
\usepackage{float}

\usepackage{tikz}
\usetikzlibrary{calc,arrows,matrix,shapes.geometric}

\usepackage{tikz-cd}

\usepackage[all]{xy}

\usepackage{hyperref}
\hypersetup{
  colorlinks=true,
  allcolors=cyan
}

\newcommand{\xr}[1]{\textcolor{red}{#1}} 

\newtheorem{thm}{Theorem}[section]
\newtheorem{prop}[thm]{Proposition}
\newtheorem{lemm}[thm]{Lemma}
\newtheorem{cor}[thm]{Corollary}
\newtheorem{claim}[thm]{Claim}

\newtheorem{conj}[thm]{Conjecture}

\numberwithin{equation}{section}

\theoremstyle{definition}
\newtheorem{defi}[thm]{Definition}
\newtheorem{ex}[thm]{Example}
\newtheorem{step}{Step}
\newtheorem{setup}[thm]{Setting}

\theoremstyle{remark}
\newtheorem{rem}[thm]{Remark}

\newcommand{\Ker}[0]{\operatorname{Ker}}

\newcommand{\Alb}[0]{\operatorname{Alb}}
\newcommand{\codim}[0]{\operatorname{codim}}
\newcommand{\Aut}[0]{\operatorname{Aut}}
\newcommand{\pr}[0]{\operatorname{pr}}
\newcommand{\Image}[0]{\operatorname{Im}}
\newcommand{\Supp}[0]{\operatorname{Supp}}
\newcommand{\End}[0]{\operatorname{End}}
\newcommand{\QQ}{\mathbb{Q}}
\newcommand{\Exc}{\mathrm{Exc}}

\newcommand{\ddbar}{dd^c}

\newcommand{\dbar}{\overline{\partial}}
\newcommand{\e}{\varepsilon}

\newcommand{\OX}{\mathcal{O}}

\newcommand{\wlim}{\mathop{\mathrm{w\text{-}lim}}\limits}
\newcommand{\Unv}[1]{{#1}^{\rm{univ}}}

\newcommand{\cal}[1]{\mathcal{#1}}
\newcommand{\bb}[1]{\mathbb{#1}}

\def \NN {\mathbb N} %
\def \ZZ {\mathbb Z} %
\def \CC {\mathbb C} %
\def \RR {\mathbb R} %
\def \QQ {\mathbb Q} %
\def \PP {\mathbb P}  %

\newcommand{\inv}{^{-1}} 
\def \ZZ {\mathbb Z} %
\def \scrO {\mathscr O} %

\def \scrH {\mathscr H} %
\def \scrI {\mathscr I} %
\def \scrK {\mathscr K} %
\def \scrL {\mathscr L} %
\def \calF {\mathcal F} %
\def \calG {\mathcal G} %
\def \calH {\mathcal H} %
\def \calN {\mathcal N} %
\def \calU {\mathcal U} %
\def \calV {\mathcal V} %
\def \calW {\mathcal W} %

\DeclareMathOperator{\SheafHom}{\mathscr{H}\!\!\!\text{\calligra om}} %
\DeclareMathOperator{\SheafTor}{\mathscr{T}\!\!\text{\calligra or}} %
\DeclareMathOperator{\id}{id} %
\DeclareMathOperator{\Coh}{H} %
\DeclareMathOperator{\dimcoh}{h} %
\DeclareMathOperator{\alb}{alb} %

\DeclareMathOperator{\rank}{rk} %
\DeclareMathOperator{\reg}{reg} %
\DeclareMathOperator{\sing}{sing} %
\DeclareMathOperator{\Sing}{Sing} %
\DeclareMathOperator{\tor}{tor} %
\DeclareMathOperator{\Sym}{Sym} %

\DeclareMathOperator{\multiplicity}{mult}
\DeclareMathOperator{\Ramification}{Ram} %

\newcommand{\qf}{\text{\rm\tiny qf}}
\DeclareMathOperator{\free}{free} %
\DeclareMathOperator{\base}{base} %
\DeclareMathOperator{\GL}{GL}

\DeclareMathOperator{\BottChern}{BC}
\DeclareMathOperator{\depth}{depth}
\let\Gamma\varGamma
\let\Delta\varDelta
\let\Theta\varTheta
\let\Lambda\varLambda
\let\Xi\varXi
\let\Pi\varPi
\let\Sigma\varSigma
\let\Upsilon\varUpsilon
\let\Phi\varPhi
\let\Psi\varPsi
\let\Omega\varOmega

\begin{document}

\date{\today, version 0.01}

\title[BBY decomposition for generalized pairs]
{Beauville--Bogomolov--Yau decomposition \\for K\"ahler generalized pairs}

\author[S. MATSUMURA]{Shin-ichi MATSUMURA}
\address{Mathematical Institute 
$\&$ Division for the Establishment of Frontier Science of Organization for Advanced Studies, 
Tohoku University, 
6-3, Aramaki Aza-Aoba, Aoba-ku, Sendai 980-8578, Japan.}
\email{{\tt mshinichi-math@tohoku.ac.jp}}
\email{{\tt mshinichi0@gmail.com}}

\author[J. WANG]{Juanyong WANG}
\address{State Key Laboratory of Mathematical Sciences, Academy of Mathematics and Systems Science, Chinese Academy of Sciences, 
55 Zhongguancun East Road, 100190 Beijing, China} 
\email{{\tt juanyong.wang@amss.ac.cn}}
\email{{\tt serredeciel@gmail.com}}

\author[X. WU]{Xiaojun WU}
\address{Laboratoire J.\,A. Dieudonn\'e, Universit\'e C\^ote d'Azur
CNRS UMR 7351
Parc Valrose 06108 NICE CEDEX 2, 
France} 
\email{{\tt xiaojun.wu@univ-cotedazur.fr}}

\author[Q. Zhang]{Qimin ZHANG}
\address{Laboratoire J.\,A. Dieudonn\'e, Universit\'e C\^ote d'Azur
CNRS UMR 7351
Parc Valrose 06108 NICE CEDEX 2, 
France} 
\email{{\tt Qimin.ZHANG@univ-cotedazur.fr}}
\email{{\tt qiminzhang28@gmail.com}}


\renewcommand{\subjclassname}{%
\textup{2020} Mathematics Subject Classification}


\maketitle



\begin{abstract}
In this paper, we establish the Beauville--Bogomolov--Yau decomposition for generalized pairs in the K\"ahler setting, thereby extending this structure theorem to the natural framework of the K\"ahler generalized Minimal Model Program. More precisely, we prove that, after passing to a finite quasi-\'etale cover, a K\"ahler generalized klt pair $(X,\Delta,\boldsymbol{\beta})$ of Calabi--Yau type admits a locally constant fibration $(X,\Delta)\to Y$ over a Calabi--Yau variety $Y$ whose fiber $(F,\Delta|_{F})$ is rationally connected. Equivalently, after base change to the universal cover of $Y$, the fibration becomes a product, and its global structure is determined by a monodromy action preserving the pair $(F,\Delta|_{F})$. 

As a principal application, when the nef b-part $\boldsymbol{\beta}$ vanishes, we show that the monodromy can be eliminated after a further finite quasi-\'etale cover, yielding a product decomposition into a rationally connected pair, strict Calabi--Yau varieties, irreducible holomorphic symplectic varieties, and complex tori. The proof introduces new analytic methods involving relative projectivity, localized positivity and flatness of direct image sheaves, generalized pairs in the analytic Minimal Model Program, and foliations on singular K\"ahler varieties.
\end{abstract}

\tableofcontents

\newpage

\section{Overview and results}\label{sec-intro}

\subsection{The main result}\label{subsec-result}

The Beauville--Bogomolov--Yau decomposition,
established in \cite{Bea83} using Yau's theorem
on the existence of Ricci-flat K\"ahler metrics,
is a fundamental structure theorem for Calabi--Yau manifolds.
It asserts that, after passing to a finite \'etale cover,
such a manifold decomposes as a product of complex tori,
strict Calabi--Yau manifolds, and holomorphic symplectic manifolds.
Throughout this paper, a \textit{Calabi--Yau variety} means a compact K\"ahler klt variety $X$ whose canonical divisor is numerically trivial. 
The notions of \textit{holomorphic symplectic varieties} and \textit{strict Calabi--Yau varieties} 
are defined as in \cite[Definition~1.1]{BGL22} (see Subsection~\ref{subsec-notation}).

The philosophy of the Minimal Model Program and the Abundance Conjecture
suggests that varieties should be understood in terms of three basic
geometric regimes: varieties of general type, varieties of Fano type,
and varieties of Calabi--Yau type \textit{with mild singularities}.
From this perspective, it is natural to seek an extension
of the classical Beauville--Bogomolov--Yau decomposition to singular varieties.
The classical decomposition was first extended
from smooth manifolds to Calabi--Yau varieties with klt singularities,
without any projectivity assumption, through a series of works
\cite{GKP16a,Druel18a,GGK19,HP19,Cam20}.
We refer to these results collectively as the
\textit{singular Beauville--Bogomolov--Yau decomposition}.
In the projective setting, the work \cite{MW} further extends this structure theory
to klt pairs of Calabi--Yau type within the framework of the log Minimal Model Program,
building on the work of several authors
(see \cite{Amb05,CH19,CCM19,Zhang05} and the references therein).
Such a decomposition cannot be expected for lc pairs in general
(see \cite{BFPT24}).
Thus, the klt assumption is  optimal with respect
to singularities.

The aim of the present paper is to establish the extension of the
Beauville--Bogomolov--Yau decomposition to the framework 
of the generalized Minimal Model Program in the K\"ahler setting, 
that is, to the remaining case combining these two directions:
generalized pairs and non-projective setting.
More precisely, we consider a K\"ahler generalized klt pair
$(X,\Delta,\boldsymbol{\beta})$ of Calabi--Yau type,
that is, a triple consisting of a compact K\"ahler variety $X$,
an effective $\mathbb{R}$-divisor $\Delta$, and a nef
b-$(1,1)$-current $\boldsymbol{\beta}$ such that
$(X,\Delta,\boldsymbol{\beta})$ is klt and
\[
K_X+\Delta+\boldsymbol{\beta}_X
\equiv 0
\quad\text{in } H_{\BottChern}^{1,1}(X).
\]
See Subsection~\ref{subsec-gpairs} for the precise definitions.

Our main result, Theorem~\ref{thm-main}, shows that, after passing to a finite quasi-\'etale cover, 
a K\"ahler generalized klt pair
$(X,\Delta,\boldsymbol{\beta})$ of Calabi--Yau type
decomposes into a Calabi--Yau base variety and rationally connected fibers 
in the sense of a locally constant fibration, 
thereby admitting an almost product structure.
See Definition~\ref{def-constant} for locally constant fibrations.

\begin{thm}\label{thm-main}
Let $(X,\Delta,\boldsymbol{\beta})$ be a K\"ahler generalized klt pair
of Calabi--Yau type.
Then, there exist a finite quasi-\'etale cover
$\nu\colon X'\to X$ and a fibration
$f\colon X'\to Y'$ satisfying the following properties:
\begin{itemize}
\item[$\mathrm{(1)}$]
the fibration $f$ is locally constant
with respect to the pair $(X',\Delta')$, where $\Delta':=\nu^{*}\Delta$;
\item[$\mathrm{(2)}$]
the fiber $F'$ of $f$ is a rationally connected variety;
\item[$\mathrm{(3)}$]
the base variety $Y'$ is a compact Calabi--Yau klt variety.
\end{itemize}

In particular, the fibration $f\colon (X',\Delta')\to Y'$ is locally trivial; that is, for every sufficiently small open subset $B\subset Y'$, 
there exists an isomorphism 
$$ 
\bigl(f^{-1}(B),\Delta'|_{f^{-1}(B)}\bigr) \cong B\times(F',\Delta'_{F'}) $$ 
over $B$, where $\Delta'_{F'}:=\Delta'|_{F'}$.
\end{thm}

Theorem~\ref{thm-main} should be viewed as the
Beauville--Bogomolov--Yau decomposition in the framework 
of the generalized Minimal Model Program.
In contrast to the case where $\boldsymbol{\beta}=0$,
a generalized pair need not admit a global product decomposition.
Nevertheless, by the definition of locally constant fibrations,
the pullback of the fibration to the universal cover of the base is a product, 
and the original fibration is recovered as the quotient of this product
by a monodromy action preserving the fiber and its boundary.
Thus, the obstruction to a global product decomposition is encoded
by the monodromy representation.
Moreover, when $\boldsymbol{\beta}=0$, the monodromy can be eliminated
after passing to a further finite quasi-\'etale cover,
and we obtain an actual product decomposition
(see Corollary~\ref{cor-bby}), thereby recovering the result of \cite{MW}.
In this sense, a locally constant fibration is the precise substitute
for a product decomposition in the setting of generalized pairs.

Compared to previous work, the main novel contributions of this paper are twofold:
\begin{itemize}
\item[$\bullet$] developing the analytic theory required to treat non-projective K\"ahler varieties;
\item[$\bullet$] incorporating the theory of nefness into the generalized Minimal Model Program. 
\end{itemize}
The first advance is the extension from projective varieties to possibly
non-projective K\"ahler varieties. 
Although such an extension is naturally suggested by the
Beauville--Bogomolov--Yau philosophy, it entails substantial technical
difficulties.
The standard tools from the projective setting do not extend directly
to the singular K\"ahler setting.
In particular, the arguments relying on ample line bundles,
hyperplane sections, algebraic covering constructions,
or projective results on foliations must be replaced by analytic methods.
One of the main contributions of this paper is to develop
the analytic techniques required to overcome these difficulties.

The second advance is the incorporation of generalized pairs,
which is new even in the projective setting.
Beginning with \cite{DPS94},
a substantial body of work has studied
varieties with various forms of ``non-negative curvature.''
A central problem in this direction concerns varieties
with nef anti-canonical divisor and is formulated as follows.

\begin{conj}\label{main-conj}
Let $X$ be a compact K\"ahler manifold with the nef anti-canonical bundle
$-K_X$.
Then, there exists a fibration $f\colon X\to Y$ satisfying
the following properties:
\begin{itemize}
\item[$\mathrm{(1)}$]
the fibration $f\colon X\to Y$ is a locally constant fibration;
\item[$\mathrm{(2)}$]
the fiber $F$ is a rationally connected manifold;
\item[$\mathrm{(3)}$]
the base variety  $Y$ is a Calabi--Yau manifold.
\end{itemize}
\end{conj}

The previous works \cite{Amb05,CH19,CCM19,Zhang05} also study nef
anti-canonical divisors in accordance with this conjectural framework.
However, the nefness of the anti-canonical divisor is not preserved
under birational contractions or strict transforms. 
By contrast, its nef contribution can be recorded on higher birational models
as a nef b-object, which naturally leads to the notion
of generalized pairs.
More generally, if $(X,\Delta)$ is a log pair such that
$-(K_X+\Delta)$ is nef, then $(X,\Delta)$ can be viewed
as a generalized pair of Calabi--Yau type by taking
$\boldsymbol{\beta}$ to be the nef b-$(1,1)$-current induced by
$-(K_X+\Delta)$.
Thus, from a birational point of view, the generalized-pair formulation provides a natural
framework that both encompasses and extends the previous setting.
Identifying this framework as the appropriate setting for the structure
theorem is another contribution of this paper.

\subsection{Applications of the main result}\label{subsec-cor}

The main result, Theorem~\ref{thm-main}, has several important applications,
parallel to those obtained in the projective setting.

The first application concerns the uniformization problem.
By the singular Beauville--Bogomolov--Yau decomposition,
after passing to a further finite quasi-\'etale cover,
the base variety $Y'$ in Theorem~\ref{thm-main} decomposes as a product
of complex tori, strict Calabi--Yau varieties,
and holomorphic symplectic varieties.
Moreover, the rationally connected fiber $F'$ in Theorem~\ref{thm-main}
is simply connected, whereas factors of the latter two types are conjectured
to have finite fundamental group.

Therefore, pulling back the fibration along the universal cover of the torus factor 
suggests that the total space should admit a product uniformization.
Interestingly, even without assuming this conjecture, we establish the expected
uniformization by applying the finiteness theorem \cite[Theorem~13.6]{GGK19} and \cite[Theorem~7.8]{Cam04} for monodromy representations
of the fundamental groups of the latter two types of factors. 

\begin{cor}\label{cor-main1}
In the setting of Theorem~\ref{thm-main},
there exist a finite quasi-\'etale cover
$\nu\colon X'\to X$ and a possibly infinite \'etale cover
$p\colon X''\to X'$ such that
\[
(X'',\Delta'')
\cong
(F',\Delta'_{F'})\times Z\times \mathbb{C}^{d}
\]
where $\Delta'':=(\nu\circ p)^{*}\Delta$.
Here $\mathbb{C}^{d}$ is the universal cover of a complex torus,
the variety $Z$ is a product of strict Calabi--Yau varieties
and holomorphic symplectic varieties,
and $(F',\Delta'_{F'})$ is the fiber pair obtained in
Theorem~\ref{thm-main}.
\end{cor}

The second application gives a product decomposition
for K\"ahler klt pairs $(X,\Delta)$ with numerically trivial
log-canonical divisor, that is, $K_X+\Delta\equiv 0$.
Equivalently, in the generalized-pair formulation,
this is precisely the case in which the nef b-part
$\boldsymbol{\beta}$ vanishes.
Combining Theorem~\ref{thm-main} with \cite{Amb05},
we show that the monodromy of locally constant fibrations becomes trivial 
after passing to a further finite quasi-\'etale cover.
This yields an actual product decomposition for klt pairs
of Calabi--Yau type in the K\"ahler setting.

\begin{cor}\label{cor-bby}
In the setting of Theorem~\ref{thm-main},
assume further that the nef b-part $\boldsymbol{\beta}$ vanishes.
Then, there exists a finite quasi-\'etale cover
$\nu\colon X'\to X$ such that 
the pair $(X',\Delta')$ admits a decomposition
\[
(X',\Delta')
\cong
(F',\Delta'_{F'})
\times \prod_i Y_i
\times \prod_j Z_j
\times T,
\]
where $\Delta':=\nu^{*}\Delta$, 
the pair $(F',\Delta'_{F'})$ is the rationally connected fiber pair
obtained in Theorem~\ref{thm-main},
each $Y_i$ is a strict Calabi--Yau variety,
each $Z_j$ is an irreducible holomorphic symplectic variety,
and $T$ is a complex torus.
\end{cor}

The third application concerns a question posed by
Hacon--M$^{\mathrm{c}}$Kernan~\cite[Question~3.1]{HM07},
which asks for a generalization of the classical theorem \cite{KoMM92,Cam92}
that Fano manifolds are rationally connected. 
This question was resolved in the projective setting
in \cite{EG19,CCM19,EIM20}, but has remained open
in the K\"ahler setting.
We generalize the results of \cite{EG19,CCM19} to the K\"ahler setting.

\begin{cor}\label{cor-HM}
Let $(X,\Delta)$ be a K\"ahler klt pair such that
$-(K_X+\Delta)$ is a nef $\mathbb{Q}$-Cartier divisor.
After replacing $(X,\Delta)$ with its pullback to a finite
quasi-\'etale cover, let $f\colon X\to Y$ be the fibration
obtained in Theorem~\ref{thm-main}, and let $F$ be a fiber of $f$.
Then, the following assertions hold:
\begin{itemize}
\item[$\mathrm{(1)}$]
$\kappa\bigl(-(K_X+\Delta)\bigr)
\leq
\kappa\bigl(-(K_F+\Delta_F)\bigr)$.

\item[$\mathrm{(2)}$]
$
\nu\bigl(-(K_X+\Delta)\bigr)
=
\nu\bigl(-(K_F+\Delta_F)\bigr), 
$
\end{itemize}
where $\Delta_{F}:=\Delta|_{F}$. 
Here $\kappa(\cdot)$ and $\nu(\cdot)$ denote the Kodaira dimension
and the numerical dimension, respectively.
\end{cor}

Even in the special case of compact K\"ahler manifolds
with nef anti-canonical bundle, Theorem~\ref{thm-main}
yields a new structure theorem.
The fourth application gives an affirmative answer
to Conjecture~\ref{main-conj}.
Cao--H\"oring~\cite{CH19} proved the conjecture
for smooth projective varieties.
Their argument, however, relies crucially on the existence
of ample line bundles.
Consequently, the general K\"ahler case has remained largely open,
apart from a few partial results, for instance when
$\dim X=3$~\cite{MW25} or when $-K_X$ is semipositive
\cite{CDP15,DPS96}.
Theorem~\ref{thm-main} proves the conjecture in full generality.

\begin{cor}\label{cor-conj}
Conjecture~\ref{main-conj} holds true.

\end{cor}

As a further consequence of Corollary~\ref{cor-conj},
we obtain an affirmative answer to the Kodaira problem 
for compact K\"ahler manifolds with nef anti-canonical bundle,
which asks whether a compact K\"ahler manifold can be approximated
by projective manifolds through deformations
(see, for instance, \cite{Huy06}).
More strongly, the nefness of the anti-canonical bundle is preserved
throughout the deformation, and the deformation is realized
by a K\"ahler fibration.

\begin{cor}\label{cor-kodaira}
Let $X$ be a compact K\"ahler manifold with nef anti-canonical bundle.
Then there exists a smooth fibration
$\mathfrak{X}\to\mathbb{D}$ from a complex manifold $\mathfrak{X}$
onto an open disc $\mathbb{D}\subset\mathbb{C}$ containing $0$
with the following properties:
\begin{itemize}
\item[$\mathrm{(1)}$]
The fibration $\mathfrak{X}\to\mathbb{D}$ is a K\"ahler fibration.
\item[$\mathrm{(2)}$]
The central fiber $\mathfrak{X}_0$ is isomorphic to $X$.
\item[$\mathrm{(3)}$]
For every $t\in\mathbb{D}$, the fiber $\mathfrak{X}_t$
has nef anti-canonical bundle.
\item[$\mathrm{(4)}$]
There exists a sequence
$\{t_k\}_{k=1}^{\infty}\subset\mathbb{D}$ converging to $0$
such that $\mathfrak{X}_{t_k}$ is projective.
\end{itemize}
\end{cor}

Furthermore, we prove the generic nefness of the tangent sheaf $T_X$
(see Theorem~\ref{nefness}) and establish a Miyaoka--Yau type inequality
involving the second Chern class
(see Theorem~\ref{thm-second}).
The generic-nefness statement extends
\cite[Theorem~1.4]{Ou23} to K\"ahler generalized pairs.
Combining this generic nefness with the solution to the Kodaira problem,
we also characterize, in terms of the second Chern class $c_2(X)$,
when $X$ admits the structure of a $\PP^1$-bundle over a complex torus
(see Theorem~\ref{thm-second}).
This generalizes
\cite[Corollary~1.7]{IM22},
\cite[Theorem~1.8]{Ou23}, and
\cite[Theorem~1.4]{IMM}.

\subsection{Strategy of the proof of the main result}
\label{subsec-strategy}

We outline the proof of Theorem~\ref{thm-main},
explain the organization of the paper,
and describe the main new ideas introduced in this paper.
The proof consists of two conceptually distinct parts:
\begin{itemize}
\item[$\mathrm{(A)}$]
We first prove Theorem~\ref{thm-main} under the additional technical
assumptions that $X$ is strongly $\mathbb{Q}$-factorial
and maximally quasi-\'etale. 
This part is carried out mainly in
Sections~\ref{sec-main} and~\ref{sec-mainthm}.

\item[$\mathrm{(B)}$]
We then remove these additional assumptions and prove
Theorem~\ref{thm-main} in full generality.
The part is carried out mainly
in Section~\ref{sec-red}.
\end{itemize}

We briefly explain the contents of the other sections.

In Section~\ref{sec-pre}, we develop the foundations of the generalized
Minimal Model Program in the analytic setting, together with the analytic
tools needed in the proof.
In particular, in Subsection~\ref{subsec-Q-fact},
we prove the existence of small $\mathbb{Q}$-factorializations
for generalized klt pairs.
We also develop the theory of singular Hermitian metrics
and their positivity in a form that is flexible under
bimeromorphic modifications.

In Section~\ref{sec-druel}, we prove an analytic version of Druel's theorem
on weakly regular foliations \cite{Druel21}, which plays a crucial role in the construction of an everywhere-defined MRC fibration.

In Section~\ref{sec-application}, we derive
Corollaries~\ref{cor-main1}, \ref{cor-bby}, \ref{cor-HM}, 
\ref{cor-conj}, and \ref{cor-kodaira} 
as applications of Theorem~\ref{thm-main}.
We also establish the generic nefness result stated in Theorem~\ref{nefness}.

In Section~\ref{sec-app}, we collect technical material
that might be familiar to experts, with the exception of the results
developed in Subsections~\ref{subsec-psef} and~\ref{subsec-segre}.
To the best of our knowledge, several of these results
have not been explicitly formulated in the existing literature.
Therefore, we include complete proofs for future reference.
In Subsection~\ref{subsec-psef},
we establish a characterization of pseudo-effective sheaves,
which is a key ingredient in the proof of the flatness criterion in \cite{CDM26}.

\smallskip

Part~$\mathrm{(A)}$ takes the approach developed in
\cite{CH19,CCM19} as its starting point.
However, the essential steps of this approach rely on algebraic tools
that are unavailable in the K\"ahler setting.
Replacing these tools requires new analytic arguments involving
relative projectivity, the positivity and flatness of direct image sheaves,
and foliations on singular K\"ahler varieties.

Part~$\mathrm{(B)}$ requires an idea completely different from that used in the previous work.
The approach of the previous work \cite{MW} for the projective case does not
extend to the K\"ahler setting.
Thus, based on the idea of \cite[Lemma~7.5]{BGL22}, 
we introduce a new descent method based on the locally constant
structure.
This method is new even in the projective setting.

Our proof draws on several recent major advances in K\"ahler geometry:
the analytic Minimal Model Program
\cite{Fuj22,DHP24,DHY26},
the projectivity criterion \cite{CH24},
the flatness criterion \cite{CDM26},
and the uniruledness criterion \cite{Ou25}.
None of these results alone yields the structure theorem.
A central contribution of the present paper is to connect them
through a collection of new results on generalized pairs,
MRC fibrations, direct image sheaves, and foliations.
These results are formulated in a way that is compatible with
the K\"ahler generalized Minimal Model Program
and should also be useful for other problems involving
K\"ahler and projective varieties.

\smallskip 
We first explain the details of Part~$\mathrm{(A)}$ in a simplified setting.
Assume that $X$ is smooth, the boundary divisor $\Delta$ is zero, the trace
$\boldsymbol{\beta}_X$ is nef, so that $-K_X$ is nef,
and that the MRC fibration is represented
by an everywhere-defined fibration
\[
\phi\colon X\to Y
\]
onto a compact K\"ahler manifold $Y$.
Even under these assumptions, the proof in the previous work
does not extend formally to the K\"ahler setting.
The argument proceeds in three steps.

\smallskip

\noindent
\textbf{Step~1 of Part~$\mathrm{(A)}$: construction of a relative polarization.}
In the projective setting, we start with an ample line bundle on $X$
and construct a relatively $\phi$-ample line bundle $G$ such that,
for a suitable integer $m>0$, the direct image sheaf
\[
\mathcal{V}:=\phi_{*}\mathcal{O}_X(mG)
\]
is pseudo-effective and has vanishing first Chern class.

The first obstruction in the K\"ahler setting is the construction
of a relatively $\phi$-ample line bundle $G$, owing to the possible absence
of ample line bundles on $X$.
However, a general fiber $F$ of the MRC fibration is
a rationally connected compact K\"ahler variety and thus it is projective.
Using the projectivity criterion of \cite{CH24},
we prove that a suitable model of the MRC fibration is projective over the base variety $Y$.
This allows us to construct a relatively $\phi$-ample line bundle $G$
with the required properties.

In the actual proof, the MRC fibration is only a rational map.
We must resolve its indeterminacy while retaining precise control
of the exceptional divisors, the boundary divisor $\Delta$, and the nef b-part $\boldsymbol{\beta}$.
The strong $\mathbb{Q}$-factoriality of $X$ is essential at this stage.
The resulting construction of the relative polarization is carried out
in Subsection~\ref{subsec-setup}.

We might instead attempt to use the Albanese map,
as in \cite{DPS94,Cao13,NW23,Wu23a}.
However, for the present argument, this does not provide an adequate
substitute for the MRC fibration, since the Albanese map need not be
a projective morphism.
Therefore, the use of the MRC fibration is essential to our construction
of the direct image sheaves.

\smallskip

\noindent
\textbf{Step~2 of Part~$\mathrm{(A)}$: positivity and flatness
of direct image sheaves.}
The next step is to prove that the reflexive hull $\mathcal{V}^{**}$ of
$\mathcal{V}=\phi_{*}\mathcal{O}_X(mG)$
is locally free and numerically flat.
These arguments are developed in
Subsections~\ref{subsec-mrc} and~\ref{subsec-direct}.
The previous approach proves pseudo-effectivity by combining
an extension theorem with an ample line bundle on the base variety $Y$
and Kawamata's covering trick.
These tools are unavailable when the base variety $Y$ is non-projective.

To overcome the absence of an ample line bundle on $Y$,
we replace these global algebraic ingredients
with a localized positivity theory for direct image sheaves
based on Bergman kernel techniques.
This localization is one of the key ideas 
of the paper, allowing us to establish the required positivity
without using a global ample line bundle or an algebraic covering
construction.
Kawamata's covering trick was also a key tool in studying
the bimeromorphic geometry of the MRC fibration.
We prove the required birational properties
in Proposition~\ref{prop_bir-geometry-psi},
using the generic nefness established in
Subsection~\ref{subsec-gen-new}.
This argument avoids
the use of Kawamata's covering trick in the earlier approach \cite{Zhang05}.

A further difficulty is to establish the flatness
of pseudo-effective sheaves with vanishing first Chern class.
In the projective setting, this is achieved using the flatness criterion
of \cite{HP19}, whose proof relies on hyperplane sections
and intersection-theoretic arguments involving the second Chern class.
Such arguments are not available on compact K\"ahler varieties.
We instead apply the K\"ahler flatness criterion of \cite{CDM26}.
The assumption that $X$ is maximally quasi-\'etale is essential at this point.
We emphasize that the flatness criterion in \cite{CDM26} relies on
a characterization of pseudo-effective sheaves
developed in Subsection~\ref{subsec-psef}.

\smallskip

\noindent
\textbf{Step~3 of Part~$\mathrm{(A)}$: from flatness to a locally
constant fibration.}
We next translate the flatness of the reflexive hull $\mathcal{V}^{**}$ of
$\mathcal{V}=\phi_{*}\mathcal{O}_X(mG)$
into geometric information about the MRC fibration.
The induced flat structure makes the projective bundle
$\mathbb{P}(\mathcal{V}^{**})\to Y$ locally constant.
By comparing $\phi\colon X\to Y$ with this projective bundle
through the relative embedding defined by $mG$,
we obtain a splitting of the tangent sheaf
\[
T_X\cong\calF\oplus\calG
\]
into two foliations.
Here $\calF$ is an algebraically integrable foliation whose general leaves are the fibers of the MRC fibration, while $\calG$ satisfies
$K_{\calG}\sim_{\QQ}0$.

To promote this infinitesimal splitting to an everywhere-defined
equidimensional fibration, we require an analytic analogue
of Druel's theorem on weakly regular foliations
\cite{Druel21}.
Druel's theorem was originally proved in the projective setting,
and its extension to singular compact K\"ahler varieties 
is not formal.
We establish the required analytic version in
Section~\ref{sec-druel}.
Together with the flat direct image sheaves and the relative embedding,
this yields the locally constant fibration in Theorem~\ref{thm-main} under the technical assumptions of Part~$\mathrm{(A)}$.

\smallskip

\noindent
\textbf{Step~1 of Part~$\mathrm{(B)}$: passage to a strongly
$\mathbb{Q}$-factorial model.}
We now turn to the general case.
After replacing $(X,\Delta,\boldsymbol{\beta})$ with a maximally quasi-\'etale cover and retaining the same notation,
we construct a small bimeromorphic morphism
\[
q\colon
(X^{\qf},\Delta^{\qf},\boldsymbol{\beta}^{\qf})
\longrightarrow
(X,\Delta, \boldsymbol{\beta})
\]
such that $X^{\qf}$ is strongly $\mathbb{Q}$-factorial.
The variety $X^{\qf}$ remains maximally quasi-\'etale.
The existence of strongly $\mathbb{Q}$-factorializations and maximally
quasi-\'etale covers in the framework of generalized pairs is itself
nontrivial and relies on the analytic Minimal Model Program. 
The construction is given in Subsection~\ref{subsec-Q-fact}.

Applying Part~$\mathrm{(A)}$ to
$(X^{\qf},\Delta^{\qf},\boldsymbol{\beta}^{\qf})$, we obtain a locally constant fibration
\[
f_{\qf}\colon
(X^{\qf},\Delta^{\qf})
\longrightarrow
Y^{\qf}
\]
with rationally connected fibers and a Calabi--Yau base variety.
It remains to descend this fibration through $q$
to the original pair $(X,\Delta)$.

\smallskip

\noindent
\textbf{Step~2 of Part~$\mathrm{(B)}$: descent of the locally
constant fibration.}
The descent is carried out in Section~\ref{sec-red}.
In the projective setting, the corresponding descent  was established in \cite[Section~4]{MW}, 
but heavily depends on the projectivity of $X$. 
Our argument avoids the use of projectivity by exploiting the locally constant structure itself.

Although the fibration
$f_{\qf}\colon(X^{\qf},\Delta^{\qf})\to Y^{\qf}$
need not be a product globally, its pullback to the universal cover
of the base variety is a product.
The contraction $q$ can therefore be studied equivariantly
on this product.
Extending the idea of \cite[Lemma~7.5]{BGL22},
which treats the special case of an actual product,
we prove that the contraction is compatible with the locally constant
structure and hence descends fiberwise.
This produces the desired fibration on $(X,\Delta)$
and completes the proof of Theorem~\ref{thm-main} (see Theorem~\ref{thm-reduction}).

The descent method is new even in the projective setting.
There are several obstacles to adapting the previous argument.
One of obstacles was the transcendental base-point-free conjecture for Calabi--Yau varieties.
Although this conjecture has now been resolved in \cite{HX26}, our proof
does not rely on it; instead, it uses a structural argument based on
monodromy and equivariance.
This method may therefore be applicable to other descent problems in the
K\"ahler Minimal Model Program.

\medskip

The preceding arguments illustrate a recurring feature of the proof:
the algebraic tools used in the projective setting are replaced not by
a single transcendental analogue, but by an interplay of relative
projectivity, localized positivity, flatness, pseudo-effectivity, and
analytic foliation theory.
This framework appears well suited to other problems arising in the
generalized Minimal Model Program in the K\"ahler setting.

Several results obtained in the course of the proof have applications beyond
Theorem~\ref{thm-main}, including the construction of strongly
$\mathbb{Q}$-factorial models for generalized K\"ahler pairs,
birational properties of MRC fibrations, positivity and flatness criteria
for direct image sheaves, a characterization of pseudo-effectivity,
and an analytic version of Druel's theorem.
Taken together, these results provide a useful collection of tools for
the study of the K\"ahler Minimal Model Program and are also of
independent interest.

\subsection*{Acknowledgments}\label{subsec-ack}
The authors would like to thank Professors Junyan Cao and Andreas H\"oring
for fruitful discussions and continued encouragement.
In particular, the investigation leading to Theorems~\ref{nefness}
and~\ref{thm-second} was initiated at the suggestion of Professor Junyan Cao.
They are grateful to Professor Mihai P{\u a}un
for kindly answering their questions.
They also thank Professor Masataka Iwai
for pointing out the reference \cite{CW24}.
The first author thanks Professor Kenta Hashizume
for providing a proof of Theorem~\ref{thm-Q-factorialization}.
The authors used generative AI tools solely for grammatical correction and to improve the fluency of the English; all mathematical content was written and verified by the authors.

The first author was partially supported by JSPS KAKENHI
Grant-in-Aid for Scientific Research (B), Grant Number 21H00976,
and by the JST FOREST Program, Grant Number JPMJFR2368.
The second author was supported by the NSFC
(Grant Numbers 12301060 and 12288201)
and the National Key R\&D Program of China
(Grant Number 2021YFA1003100).
The third author was supported by a 2023 Young Researcher Excellence Fellowship
from the IdEx Universit\'e C\^ote d'Azur.

\section{Preliminaries and technical results}\label{sec-pre}

\subsection{Notation and conventions}\label{subsec-notation}

Throughout this paper, we work over the field of complex numbers.

We use the terms ``Cartier divisors,'' ``invertible sheaves,'' and
``line bundles'' interchangeably, and adopt the additive notation for tensor
products (e.g.,\,$L+M:=L\otimes M$ for line bundles $L$ and $M$).
Similarly, we use the terms ``locally free sheaves'' and
``vector bundles'' interchangeably.
A \textit{fibration} denotes a proper surjective morphism with connected fibers.
An \textit{analytic variety} denotes an irreducible and reduced complex analytic
space.
A \textit{K\"ahler variety} denotes a normal analytic variety that admits
a K\"ahler form
(i.e.,\,a smooth positive $(1,1)$-form with local potentials). 
Note that the normality is included in our definition of K\"ahler varieties.
A \textit{Calabi--Yau variety} is a K\"ahler klt variety
$X$ whose canonical divisor $K_X$ is $\mathbb{Q}$-Cartier and
numerically trivial
(i.e.,\,$K_X$ determines the trivial Bott--Chern class).
The terms ``\textit{strict Calabi--Yau variety}'' and
``\textit{holomorphic symplectic variety}'' 
are used in the sense of \cite[Definition~1.1]{BGL22}
(see also \cite[Definition~8.16]{GKP16b}).
In the proof, we do not use their precise definitions, but only the
fact that such varieties have vanishing augmented irregularity.

Unless otherwise stated, all sheaves are assumed to be coherent.
Let $\mathcal{E}$ be a torsion-free sheaf on a normal analytic variety $X$.
The reflexive hull $\mathcal{E}^{**}$ denotes its double dual, where
the dual sheaf is defined by
$
\mathcal{E}^*
\coloneqq
\SheafHom(\mathcal{E},\mathcal{O}_X).
$
For a positive integer $m\in\mathbb{Z}_{+}$, we define
\[
\mathcal{E}^{[\otimes m]}
:=
(\mathcal{E}^{\otimes m})^{**},
\qquad
\operatorname{Sym}^{[m]}\mathcal{E}
:=
(\operatorname{Sym}^m\mathcal{E})^{**},
\qquad
\wedge^{[m]}\mathcal{E}
:=
(\wedge^m\mathcal{E})^{**}.
\]
The \textit{determinant sheaf} $\det\mathcal{E}$ is defined by
\[
\det\mathcal{E}
:=
\wedge^{[\mathrm{rk}\,\mathcal{E}]}\mathcal{E}.
\]
A torsion-free sheaf $\mathcal{E}$ of rank one is said to be
\textit{$\mathbb{Q}$-Cartier} if there exists an integer
$m\in\mathbb{Z}_{+}$ such that
$\mathcal{E}^{[\otimes m]}$ is invertible.

For torsion-free sheaves $\mathcal{E}$ and $\mathcal{F}$ on $X$
and a (not necessarily surjective) morphism
$f\colon Y\to X$ between normal analytic varieties, we define
\[
\mathcal{E}[\otimes]\mathcal{F}
:=
(\mathcal{E}\otimes\mathcal{F})^{**},
\qquad
f^{[*]}\mathcal{E}
:=
(f^*\mathcal{E})^{**}.
\]
For a subset $U\subset X$, we often say that a property~\textup{(P)}
holds \textit{over} $U$ if it holds on the inverse image
$f^{-1}(U)\subset Y$.

For a normal analytic variety $X$, we define the canonical sheaf by
\[
K_X
:=
j_{[*]}\mathcal{O}_{X_{\reg}}(K_{X_{\reg}})
:=
\left(
j_*\mathcal{O}_{X_{\reg}}(K_{X_{\reg}})
\right)^{**},
\]
where $j\colon X_{\reg}\hookrightarrow X$ is the natural inclusion of the smooth locus $X_{\reg}$
and $K_{X_{\reg}}$ is the canonical bundle defined on $X_{\reg}$.
In contrast to the projective case, the sheaf $K_X$ need not be represented
by a Weil divisor, but we often refer to $K_{X}$ as the
\textit{canonical divisor}.
Similarly, for a generalized pair
$(X,\Delta,\boldsymbol{\beta})$, we refer to
$K_X+\Delta+\boldsymbol{\beta}_X$ as the
\textit{log-canonical divisor}, although $K_X+\Delta+\boldsymbol{\beta}_X$ need not be a divisor
(even when $\boldsymbol{\beta}_X=0$).

Following \cite[Definition~4.6.2]{BEG13}, we define the Bott--Chern
cohomology group $H_{\BottChern}^{1,1}(X)$ by
\[
H_{\BottChern}^{1,1}(X)
:=
\frac{
\left\{
\alpha \,\middle|\,
\begin{array}{c}
\alpha \text{ is a real $d$-closed smooth $(1,1)$-form on $X$}\\
\text{admitting local potentials}
\end{array}
\right\}
}{
\left\{
dd^c\varphi \mid \varphi\in C^\infty(X,\mathbb{R})
\right\}
}.
\]
The same cohomology group can equivalently be defined using real closed
$(1,1)$-currents with local potentials, modulo currents of the form
$dd^c u$, where $u$ is a globally defined real-valued distribution on $X$.
The symbol $\equiv$ denotes equality in Bott--Chern cohomology.
For a $(1,1)$-current $T$ with local potentials
(resp.\,an $\mathbb{R}$-line bundle $L$) on $X$,
we denote by $c_1(T)$ (resp.\,$c_1(L)$) its associated Bott--Chern class.

The following proposition will be used repeatedly, often without explicit
mention, to extend semipositive currents representing Bott--Chern
cohomology classes (see \cite{Dem85} for currents on analytic varieties).

\begin{prop}[{\cite[Proposition~4.6.3]{BEG13}}]\label{BEG}
Let
$\alpha\in H^{1,1}_{\BottChern}(X)$
be a Bott--Chern cohomology class on a normal analytic variety $X$,
and let $T$ be a semipositive current on $X_{\rm{reg}}$ representing
the restricted class
$
\alpha|_{X_{\rm{reg}}}
\in
H^{1,1}_{\BottChern}(X_{\rm{reg}}).
$
Then, the current $T$ uniquely extends to a semipositive current on $X$
with local potentials representing
$\alpha\in H^{1,1}_{\BottChern}(X)$.
\end{prop}

\subsection{Locally constant fibrations}

In this subsection, following \cite{MW}, we introduce the notion of
locally constant fibrations with respect to Weil $\RR$-divisors $\Delta$.
We also provide a characterization of locally constant fibrations in
terms of local systems.

\begin{defi}\label{def-constant}
Let $\phi\colon X\to Y$ be a fibration between normal analytic varieties,
and let $\Delta$ be a Weil $\RR$-divisor on $X$.
The morphism $\phi\colon (X,\Delta)\to Y$ is called a
\textit{locally constant fibration}
if the following conditions hold:
\begin{itemize}
\item[$\bullet$]
The morphism $\phi\colon X\to Y$ is an analytic fiber bundle
with fiber $F$.
\item[$\bullet$]
Every component $\Delta_i$ of $\Delta$ is horizontal
(i.e.,\,$\phi(\Delta_i)=Y$).
\item[$\bullet$]
There exist a Weil $\RR$-divisor $\Delta_F$ on $F$ and a representation
\[
\rho\colon\pi_1(Y)\to\Aut(F)
\]
of the (topological) fundamental group $\pi_1(Y)$ into the automorphism group
$\Aut(F)$ satisfying the following conditions:
\begin{itemize}
\item[$\bullet$]
The divisor $\Delta_F$ is invariant under the action of $\pi_1(Y)$.
\item[$\bullet$]
The pair $(X,\Delta)$ is isomorphic over $Y$ to the quotient
\[
\bigl(\Unv{Y}\times(F,\Delta_F)\bigr)/\pi_1(Y).
\]
Here $\Unv{Y}$ denotes the universal cover of $Y$, and
$\gamma\in\pi_1(Y)$ acts diagonally on $\Unv{Y}\times F$ by
\[
\gamma\cdot(y,z)
:=
\bigl(\gamma\cdot y,\rho(\gamma)(z)\bigr)
\quad\text{for }(y,z)\in\Unv{Y}\times F.
\]
\end{itemize}
\end{itemize}
\end{defi}

For a locally constant fibration $\phi\colon X\to Y$,
the fiber product $X\times_Y\Unv{Y}$ is isomorphic to
$\Unv{Y}\times F$.
Under this identification, the induced morphism
$X\times_Y\Unv{Y}\to\Unv{Y}$ coincides with the first projection
$\Unv{Y}\times F\to\Unv{Y}$.
Thus, we obtain the following diagram:
\begin{equation*}\label{locally-constant}
\begin{gathered}
\xymatrix@C=40pt@R=30pt{
X \ar[d]_{\phi}
&
X\times_Y\Unv{Y}\cong\Unv{Y}\times F
\ar[d]_{\pr_1}
\ar[l]^{p_X \qquad\qquad}
\ar[r]^{\qquad\qquad \pr_2}
&
F.
\\
Y
&
\Unv{Y}\ar[l]^{p_Y}
&
}
\end{gathered}
\end{equation*}

A locally constant fibration
$\phi\colon (X,\Delta)\to Y$ is locally trivial.
More precisely, it is an analytic fiber bundle with fiber $F$, and
there exists a Weil $\RR$-divisor $\Delta_F$ on $F$ such that
\[
\bigl(\phi^{-1}(B),\Delta|_{\phi^{-1}(B)}\bigr)
\cong
B\times(F,\Delta_F)
\]
over every sufficiently small open subset $B\subset Y$.
Nevertheless, the local constancy is substantially stronger than the local
triviality, and some arguments in this paper fail for fibrations that are merely
locally trivial.
The formulation in terms of local constancy was shown in \cite{MW}
to be essential, and we follow the same principle throughout this paper.

The following lemma is a reformulation of
\cite[Lemma~2.4]{MW}.
Conclusion~\textup{(1)} follows from its proof, whereas
Conclusion~\textup{(2)} is an immediate consequence of its statement.

\begin{lemm}
\label{lemma_lc_line}
Let $\phi\colon X\to Y$ be a projective locally constant fibration
such that the fiber $F$ has vanishing irregularity.
Let $L$ be a relatively $\phi$-ample line bundle on $X$.
Then, the following assertions hold:
\begin{itemize}
\item[\textup{(1)}]
For every sufficiently large and divisible integer
$m\in\mathbb{Z}_{+}$, the restriction $mL|_{F}$ to the fiber $F$ is
$\rho$-equivariant and $\rho$-linearizable.

\item[\textup{(2)}]
There exist a line bundle $L_F$ on $F$ and a
$\mathbb{Q}$-line bundle $L_Y$ on $Y$ such that
\[
p_X^*(L-\phi^*L_Y)
\sim_{\QQ}
\pr_2^*L_F. 
\]
Moreover, for every sufficiently divisible integer $m\in\mathbb{Z}_{+}$, 
the isomorphism 
$$
p_X^*(m(L-\phi^*L_Y))
\cong
\pr_2^*(m L_F). 
$$
is $\pi_1(Y)$-equivariant. 
In particular, the direct image sheaf
$\phi_*(\mathcal{O}_{X}(m(L-\phi^*L_Y)))$ is a flat vector bundle on $Y$
for every sufficiently divisible integer
$m\in\mathbb{Z}_{+}$.
\end{itemize}
\end{lemm}

We reformulate \cite[Proposition~2.5]{MW}, which gives a criterion
for locally constant fibrations.
The only difference from the cited proposition is that 
$D$ is not assumed to be a $\QQ$-divisor on $X$ in this paper.

\begin{prop}
\label{prop_flat-lcf}
Let $\phi\colon X\to Y$ be a fibration between normal analytic varieties

and let $D$ be an effective Weil $\RR$-divisor on $X$.
Assume that $\phi$ is a flat projective morphism and that there exists
a relatively $\phi$-ample line bundle $L$ on $X$ satisfying the
following conditions:
\begin{itemize}
\item[\rm(1)]
For every $m\in\ZZ_{+}$, the sheaf
$E_m:=\phi_*(\mathcal{O}_{X}(mL))$ is a local system $($equivalently, a flat vector bundle$)$.

\item[\rm(2)]
For every $m\in\ZZ_{+}$, the natural morphism
$\Sym^m\!E_1\to E_m$ is a morphism of local systems
$($i.e.,\,it is induced by a morphism between the corresponding
representations of $\pi_1(Y)$$)$.

\item[\rm(3)]
For every $m\in\ZZ_{+}$, the subsheaf
$F_m:=\phi_*(\mathcal{O}_{X}(mL-\lfloor D\rfloor))$ is a local subsystem of $E_m$, 
where $\lfloor D\rfloor$ denotes the round-down of $D$.
\end{itemize}
Then, the fibration 
$
\phi\colon(X,\lfloor D\rfloor)\to Y
$
is a locally constant fibration. 
\end{prop}

\begin{rem}
\label{rmk_prop_flat-lcf}
If $Y$ is a compact K\"ahler manifold, then conditions~(1), (2), and~(3) can be replaced by the requirement that both $E_m$ and $F_m$ be numerically flat vector bundles. 
Indeed, the compatibility condition~(2) is automatically satisfied by \cite[Lemma 4.3.3]{Cao13}.
\end{rem}

\subsection{Generalized pairs in the analytic setting}\label{subsec-gpairs}

In this subsection, following \cite[\S\S~2.1--2.2]{DHY26},
we recall the definition of generalized pairs and fix the notation used
throughout the paper.
The formulation in terms of nef b-$(1,1)$-currents is a natural extension
of the corresponding definition in terms of nef b-divisors
and is well suited to the non-projective setting.

\begin{defi}\label{def-gpairs}
\textup{(1) (Generalized pairs.)}
A \textit{generalized pair} $(X,\Delta,\boldsymbol{\beta})$ is a triple
consisting of the following data:
\begin{itemize}
\item[$\bullet$] a normal analytic variety $X$;
\item[$\bullet$] an effective $\mathbb{R}$-divisor $\Delta$ on $X$;
\item[$\bullet$] a nef b-$(1,1)$-current $\boldsymbol{\beta}$ over $X$,
\end{itemize}
such that the Bott--Chern class
\[
c_1(K_X+\Delta+\boldsymbol{\beta}_X)
\in H_{\BottChern}^{1,1}(X)
\]
is well-defined.
More precisely, this means that the Bott--Chern class on $X_{\reg}$
determined by $K_X+\Delta+\boldsymbol{\beta}_X$ extends to
a Bott--Chern class on $X$.
A generalized pair $(X,\Delta,\boldsymbol{\beta})$ is called a
\textit{K\"ahler generalized pair of Calabi--Yau type}
if $X$ is a compact K\"ahler variety and
\[
K_X+\Delta+\boldsymbol{\beta}_X
\equiv 0
\quad\text{in }
H_{\BottChern}^{1,1}(X).
\]

\smallskip

\textup{(2) (Nef b-$(1,1)$-currents.)}
A closed \textit{b-$(1,1)$-current} $\boldsymbol{\beta}$ over $X$
is a collection
\[
\boldsymbol{\beta}
:=
\{\boldsymbol{\beta}_{X'}\}_{X'}
\]
of closed $(1,1)$-currents $\boldsymbol{\beta}_{X'}$ on $X'$,
where $X'$ runs through all proper bimeromorphic models $X'\to X$,
satisfying the following compatibility condition:
for every bimeromorphic morphism $g\colon X''\to X'$ over $X$, we have 
\[
g_*\boldsymbol{\beta}_{X''}
=
\boldsymbol{\beta}_{X'}.
\]
The $(1,1)$-current $\boldsymbol{\beta}_{X'}$ is called the
\textit{trace} of $\boldsymbol{\beta}$ on $X'$.
Note that $\boldsymbol{\beta}_{X'}$ need not admit local potentials, 
and thus, its pullback cannot be defined in general.

A b-$(1,1)$-current $\boldsymbol{\beta}$ is said to
\textit{descend to a model $M\to X$}
if $\boldsymbol{\beta}_M$ admits local potentials and,
for every higher model $g\colon X''\to M$ over $X$,
the current $\boldsymbol{\beta}_{X''}$ also admits local potentials and satisfies 
\[
c_{1}(\boldsymbol{\beta}_{X''})
\equiv 
g^* c_{1}(\boldsymbol{\beta}_M)
\quad\text{in }
H_{\BottChern}^{1,1}(X'').
\]
A b-$(1,1)$-current $\boldsymbol{\beta}$ is said to be
\textit{nef} if it descends to a model $M\to X$ such that
$
c_{1}(\boldsymbol{\beta}_M)
\in H_{\BottChern}^{1,1}(M)
$
is nef.

\smallskip

\textup{(3) (Singularities of generalized pairs.)} 
A bimeromorphic morphism $\pi\colon M\to X$ is called a
\textit{log resolution} of a generalized pair
$(X,\Delta,\boldsymbol{\beta})$ if
$\pi\colon M\to X$ is a resolution of singularities of $X$ such that
$\boldsymbol{\beta}$ descends to $M$ and
$\pi^{-1}_*\Delta+E$ has SNC support,
where $E$ is the reduced divisor supported on the
$\pi$-exceptional locus.
Since the pullback of the Bott--Chern class
$c_{1}(K_X+\Delta+\boldsymbol{\beta}_X)$ is well-defined, 
we may write, by abuse of notation,
\begin{align}
\label{eq-gpair}
K_M+\Delta_M+\boldsymbol{\beta}_M
\equiv
\pi^*(K_X+\Delta+\boldsymbol{\beta}_X)+F_+
\quad\text{in}\quad
H_{\BottChern}^{1,1}(M),
\end{align}
where $F_+$ is an effective $\pi$-exceptional $\mathbb{R}$-divisor
and $\Delta_M$ is an effective $\mathbb{R}$-divisor such that
$\Delta_M$ and $F_+$ have no common irreducible components.  
Here, we implicitly use the second theorem of support
\cite[III.2.13, pp.~142--143]{agbook} (see also Lemma~\ref{topo-lemma}.

A generalized pair $(X,\Delta,\boldsymbol{\beta})$ is called
\textit{klt} (resp.\,\textit{lc}) if
$\Delta_M<1$ (resp.\,$\Delta_M\leq 1$)
for some, or equivalently every, log resolution
$\pi\colon M\to X$.
In this case, we also say that $(X,\Delta,\boldsymbol{\beta})$ is a
generalized klt pair (resp.\,a generalized lc pair).
\end{defi}

We provide four examples of generalized pairs of Calabi--Yau type.
These examples illustrate both the breadth of the extension to
generalized pairs and its natural compatibility with birational geometry.
In particular, allowing anti-log-canonical divisors that admit
birational Zariski decompositions considerably broadens the scope of
the theory and may lead to applications.

\begin{ex}\label{example-compare}
\textup{(1)}
A basic source of generalized klt pairs of Calabi--Yau type is given by
ordinary klt pairs with nef anti-log-canonical divisor.
This example shows that the objects studied in this paper naturally
generalize those considered in the previous work \cite{MW}.

Let $(X,\Delta)$ be a klt pair such that $-(K_X+\Delta)$ is a nef
$\QQ$-Cartier divisor.
Let $\boldsymbol{\beta}$ be the nef b-$(1,1)$-current induced by $-(K_{X}+\Delta)$, 
which is constructed as follows: 
Choose a smooth $(1,1)$-form $\theta$ representing 
$-c_1(K_X+\Delta)\in H_{\BottChern}^{1,1}(X)$, and define a
b-$(1,1)$-current
$\boldsymbol{\beta}:=\{\boldsymbol{\beta}_{X'}\}_{X'}$ by
$\boldsymbol{\beta}_{X'}:=f^*\theta$ for every bimeromorphic morphism
$f\colon X'\to X$.
Then $\boldsymbol{\beta}$ descends to $X$ itself, and
$c_1(\boldsymbol{\beta}_X)=-c_1(K_X+\Delta)$ is nef.
Moreover, we see that 
$K_X+\Delta+\boldsymbol{\beta}_X\equiv 0$ in
$H_{\BottChern}^{1,1}(X)$ (i.e.,\,$(X,\Delta,\boldsymbol{\beta})$ is of Calabi--Yau type). 
Further, we have 
$c_{1}(\boldsymbol{\beta}_M)=\pi^* c_{1}(\boldsymbol{\beta}_X) $
for every log resolution $\pi\colon M\to X$.
Consequently, after cancelling the b-$(1,1)$-current terms
in the pullback formula \eqref{eq-gpair}, the generalized discrepancies of
$(X,\Delta,\boldsymbol{\beta})$ coincide with the ordinary discrepancies
of $(X,\Delta)$.
Hence, the triple $(X,\Delta,\boldsymbol{\beta})$ is a generalized klt pair.

\smallskip

\textup{(2)}
The following proposition gives a general construction based on
Zariski-type decompositions.

\begin{prop}
\label{prop-zariski}
Let $(X,\Delta)$ be an ordinary K\"ahler klt pair.
Assume that there exist  a resolution $\pi\colon M\to X$ of singularities,
a nef $\mathbb{R}$-divisor $P$, and an effective
$\pi$-exceptional $\mathbb{R}$-divisor $N$ such that
\[
-\pi^*(K_X+\Delta)
\equiv
P+N
\quad\text{in}\quad
H_{\BottChern}^{1,1}(M).
\]
Assume further that $(M,\Delta'_M+N)$ is klt.
Here $\Delta'_M$ is the effective $\mathbb{R}$-divisor defined by
\[
K_M+\Delta'_M
=
\pi^*(K_X+\Delta)+F'_+,
\]
where $F'_+$ is an effective $\pi$-exceptional
$\mathbb{R}$-divisor having no common components with $\Delta'_M$.
Let $\boldsymbol{\beta}$ be the nef b-$(1,1)$-current induced by $P$,
constructed as follows.
Choose a smooth $(1,1)$-form $\theta$ representing
$c_1(P)\in H_{\BottChern}^{1,1}(M)$.
For every higher model $f\colon X'\to M$, set
$
\boldsymbol{\beta}_{X'}:=f^*\theta.
$
For an arbitrary bimeromorphic model $X''\to X$, choose a common
resolution
\[
\begin{tikzcd}
& X' \arrow[ld, "f"'] \arrow[rd, "g"] & \\
M & & X''
\end{tikzcd}
\]
and define
$
\boldsymbol{\beta}_{X''}
:=
g_*\boldsymbol{\beta}_{X'}.
$
This definition is independent of the choice of the common resolution
and determines a nef b-$(1,1)$-current over $X$ descending to $M$.

Then $(X,\Delta,\boldsymbol{\beta})$ is a generalized klt pair
of Calabi--Yau type.
\end{prop}

\begin{proof}
By construction, 
we see that $\boldsymbol{\beta}$ is a nef b-$(1,1)$-current,
since $P$ is nef.
Since $N$ is $\pi$-exceptional, we have
\[
c_1(\boldsymbol{\beta}_X)
=
c_1(\pi_*\boldsymbol{\beta}_M)
=
-c_1(K_X+\Delta).
\]
Moreover, by the definition of the generalized boundary $\Delta_M$ on $M$
in \eqref{eq-gpair}, we have
\[
\Delta'_{M}-F'_+
\equiv
\Delta_M-F_++N.
\]
By Lemma~\ref{topo-lemma}, the associated integration currents are in fact equal.
Although $F_+$ and $N$ may have common components, the pair
$(M,\Delta_M+N)$ is klt and $N$ is effective.
Therefore, removing from $\Delta_M$ the irreducible components appearing
in $F_+$,  we deduce that  $(M,\Delta_M')$ remains klt.
Hence, the triple $(X,\Delta,\boldsymbol{\beta})$ is a generalized klt pair
of Calabi--Yau type.

\end{proof}

\smallskip

\textup{(3)}
We next give an example that does not necessarily arise from the
construction in \textup{(1)}.
We construct a bimeromorphic morphism $\pi\colon M\to X$ such that
$-K_X$ is not nef, whereas $-K_M$ is nef.
We then consider the nef b-$(1,1)$-current
$\boldsymbol{\beta}$ induced by the nef divisor $-K_M$ 
(see Proposition~\ref{prop-zariski} for the construction).
Since $c_1(\boldsymbol{\beta}_X)=-c_1(K_X)$, we have
$K_X+\boldsymbol{\beta}_X\equiv0$.
This shows that $(X,0,\boldsymbol{\beta})$ is a generalized klt pair of
Calabi--Yau type.

To construct such an example, consider the projective bundle
\[
p\colon
X
=
\mathbb{P}_{\mathbb{P}^1}
\left(
\mathcal{O}_{\mathbb{P}^1}
\oplus
\mathcal{O}_{\mathbb{P}^1}(2)
\oplus
\mathcal{O}_{\mathbb{P}^1}(2)
\right)
\longrightarrow
\mathbb{P}^1,
\]
and set $G:=\mathcal{O}_X(1)$ and
$H:=p^*\mathcal{O}_{\mathbb{P}^1}(1)$.
Let $C\cong\mathbb{P}^1$ be the section corresponding to the quotient
\[
\mathcal{O}_{\mathbb{P}^1}
\oplus
\mathcal{O}_{\mathbb{P}^1}(2)
\oplus
\mathcal{O}_{\mathbb{P}^1}(2)
\longrightarrow
\mathcal{O}_{\mathbb{P}^1}.
\]
Then $G\cdot C=0$ and $H\cdot C=1$.
The canonical bundle formula gives $K_X=-3G+2H$.
This shows that $-K_X$ is not nef, since $(-K_X)\cdot C=-2$.

On the other hand, let
$\pi\colon M:=\operatorname{Bl}_C X\to X$
be the blow-up of $X$ along $C$, and let $F$ be its exceptional divisor.
Then, we have 
\[
-K_M
=
-\pi^*K_X-F
=
2\pi^*G
+
\left(
\pi^*(G-2H)-F
\right).
\]
The divisor $\pi^*G$ is nef, since $G$ is globally generated.
Moreover,
\[
p_*\mathcal{O}_X(G-2H)
\cong
\mathcal{O}_{\mathbb{P}^1}(-2)
\oplus
\mathcal{O}_{\mathbb{P}^1}
\oplus
\mathcal{O}_{\mathbb{P}^1},
\]
and the base ideal sheaf of $|G-2H|$ coincides with
$\mathcal{I}_C$.
The line bundle $\pi^*(G-2H)-F$ is base-point-free.
Thus, we see that $-K_M$ is nef.

\smallskip

\textup{(4)}
We give an example arising from the Atiyah flop. 
We consider the projective space bundle 
\[
p\colon
X
=
\mathbb{P}_{\mathbb{P}^1}
\left(
\mathcal{O}_{\mathbb{P}^1}
\oplus
\mathcal{O}_{\mathbb{P}^1}(1)
\oplus
\mathcal{O}_{\mathbb{P}^1}(1)
\right)
\longrightarrow
\mathbb{P}^1,
\]
and set $G:=\mathcal{O}_X(1)$ and
$H:=p^*\mathcal{O}_{\mathbb{P}^1}(1)$.
Let $C\cong\mathbb{P}^1$ be the section corresponding to the quotient
\[
\mathcal{O}_{\mathbb{P}^1}
\oplus
\mathcal{O}_{\mathbb{P}^1}(1)
\oplus
\mathcal{O}_{\mathbb{P}^1}(1)
\longrightarrow
\mathcal{O}_{\mathbb{P}^1}.
\]
Then, we see that 
\[
G\cdot C=0,
\qquad
H\cdot C=1,
\qquad
N_{C/X}\cong\mathcal{O}_{\mathbb{P}^1}(-1)^{\oplus 2}.
\]
The linear system $|G|$ defines a small contraction
$s\colon X\to Z$ of $C$.
The standard Atiyah flop of $C$ gives a small birational map
$X\dashrightarrow X'$ over $Z$.
Let $s'\colon X'\to Z$ be the induced small contraction with the following diagram: 
\[
\begin{tikzcd}
X \arrow[dr, "s:=\Phi_{|G|}"'] \arrow[rr, dashed, "\text{flop}"] & & X' \arrow[dl, "s'"] \\
& Z. &
\end{tikzcd}
\]

Fix a rational number $0<\varepsilon\ll 1$, and set
$A:=G-\varepsilon H$.
Let $G'$, $H'$, and $A'$ denote the strict transforms of
$G$, $H$, and $A$ on $X'$, respectively.
Then $A$ is not nef since $A\cdot C=-\varepsilon<0$. 
Nevertheless $A'$ is ample.
Indeed, there exists an ample divisor $L$ on $Z$ such that
$G=s^*L$ and $G'={s'}^*L$.
Since $-H'$ is relatively $s'$-ample, the divisor
\[
A'
=
G'-\varepsilon H'
=
{s'}^*L+\varepsilon(-H')
\]
is ample for $0<\varepsilon\ll 1$.

The canonical divisor formula yields $K_X=-3G$, 
showing that $-(K_X+A)=2G+\varepsilon H$ is ample.
Therefore, we may  choose an effective $\mathbb{Q}$-divisor
$\Delta\sim_{\mathbb{Q}}-(K_X+A)$ such that $(X,\Delta)$ is klt.
Then, by $-(K_X+\Delta)\sim_{\mathbb{Q}}A$, we see that 
$-(K_X+\Delta)$ is not nef.

Let $\boldsymbol{\beta}$ be the nef b-$(1,1)$-current determined by
the ample $\mathbb{Q}$-divisor $A'$ on $X'$.
Since $X\dashrightarrow X'$ is small, the trace of
$\boldsymbol{\beta}$ on $X$ satisfies that $\boldsymbol{\beta}_X \equiv A$.
Therefore, we have 
\[
K_X+\Delta+\boldsymbol{\beta}_X
\equiv
K_X+\Delta+A
\equiv
0,
\]
and hence $(X,\Delta,\boldsymbol{\beta})$ is a generalized pair
of Calabi--Yau type.

It remains to verify that this generalized pair is klt.
Let $\pi \colon M:=\operatorname{Bl}_C X\to X$ be the blow-up,
let $E$ be its exceptional divisor, and let
$\pi'\colon M\to X'$ be the contraction of the other ruling of $E$.
By choosing $\Delta$ sufficiently general, we may assume that
$C\not\subset\operatorname{Supp}\Delta$ and that the union of its strict transform
and the exceptional divisor $E$
has SNC supprt.
Note that  $M$ is a common resolution of the Atiyah flop
$X\dashrightarrow X'$.

Set $\tilde{\Delta}:=\pi_*^{-1}\Delta$.
Since $C\not\subset\operatorname{Supp}\Delta$, we have
$\pi^*\Delta=\tilde{\Delta}$.
Moreover, the standard calculation for the Atiyah flop gives
$\boldsymbol{\beta}_M\equiv{\pi'}^*A'\equiv\pi^*A-\varepsilon E$.
Since $\pi$ is the blow-up of a smooth curve in a smooth threefold,
we also have $K_M=\pi^*K_X+E$.
Consequently, we obtain 
\[
\begin{aligned}
K_M+\tilde{\Delta}+\boldsymbol{\beta}_M
&\equiv
\pi^*K_X+E+\tilde{\Delta}+\pi^*A-\varepsilon E \\
&=
\pi^*(K_X+\Delta+A)+(1-\varepsilon)E \\
&\equiv
\pi^*(K_X+\Delta+\boldsymbol{\beta}_X)
+(1-\varepsilon)E.
\end{aligned}
\]
Thus, in the notation of \eqref{eq-gpair}, we have
$\Delta_M=\tilde{\Delta}$ and $F_+=(1-\varepsilon)E$.
The divisors $\Delta_M$ and $F_+$ have no common irreducible components
by $C\not\subset\operatorname{Supp}\Delta$.
Finally, since $(X,\Delta)$ is klt, we can conclude that 
$\Delta_M=\tilde{\Delta}<1$.
Therefore $(X,\Delta,\boldsymbol{\beta})$ is a generalized klt pair.
\end{ex}

We record the behavior of generalized pairs under finite quasi-\'etale
covers and small bimeromorphic morphisms.
In particular, the Calabi--Yau condition and generalized klt
(resp.\,lc) singularities are preserved under these operations.

\begin{prop}\label{prop-pullback}
Let $(X,\Delta,\boldsymbol{\beta})$ be a generalized pair.
Then, the following assertions hold.
\begin{enumerate}
\item[\textup{(1)}]
Let $\nu\colon X'\to X$ be a finite quasi-\'etale cover.
Set $\Delta':=\nu^*\Delta$, where the pullback is taken in the sense
of Weil divisors.
Then, there exists a nef b-$(1,1)$-current
$\boldsymbol{\beta}'$ over $X'$ such that
$(X',\Delta',\boldsymbol{\beta}')$ is a generalized pair and
\[
K_{X'}+\Delta'+\boldsymbol{\beta}'_{X'}
\equiv 
\nu^*(K_X+\Delta+\boldsymbol{\beta}_X)
\quad\text{in}\quad
H_{\BottChern}^{1,1}(X').
\]

\item[\textup{(2)}]
Let $q\colon\bar X\to X$ be a small bimeromorphic morphism.
Set $\bar\Delta:=q^{-1}_*\Delta$, and let
$\bar{\boldsymbol{\beta}}$ be the b-$(1,1)$-current over $\bar X$
induced by $\boldsymbol{\beta}$.
Then $(\bar X,\bar\Delta,\bar{\boldsymbol{\beta}})$ is a generalized pair
and
\[
K_{\bar X}+\bar\Delta+\bar{\boldsymbol{\beta}}_{\bar X}
\equiv 
q^* (K_X+\Delta+\boldsymbol{\beta}_X )
\quad\text{in}\quad
H_{\BottChern}^{1,1}(\bar X).
\]
\end{enumerate}

In particular, in either case, if $(X,\Delta,\boldsymbol{\beta})$ is of
Calabi--Yau type, then so is the induced generalized pair.
The same assertion holds for generalized klt and generalized lc
singularities.
\end{prop}

\begin{proof}
Set
$\alpha:=c_1(K_X+\Delta+\boldsymbol{\beta}_X) \in H_{\BottChern}^{1,1}(X)$.
Let $\pi\colon M\to X$ be a log resolution of
$(X,\Delta,\boldsymbol{\beta})$.

\textup{(1)}
Let $Z$ be a resolution of singularities of the normalization of the
main component of the fiber product $M\times_X X'$, and consider the
commutative diagram
\[
\begin{tikzcd}
Z \arrow[r, "h"] \arrow[d, "g"'] &
M \arrow[d, "\pi"] \\
X' \arrow[r, "\nu"'] &
X.
\end{tikzcd}
\]
Define $\boldsymbol{\beta}'$ to be the b-$(1,1)$-current induced 
by a representative of $h^* c_{1}(\boldsymbol{\beta}_M)$.
Since $[\boldsymbol{\beta}_M]$ is nef, we see that 
$c_{1}(\boldsymbol{\beta}'_Z)=h^* (c_{1}\boldsymbol{\beta}_M)$ is also nef.

Pulling back the discrepancy formula for
$(X,\Delta,\boldsymbol{\beta})$ to $Z$,
we obtain effective $\mathbb{R}$-divisors $\Delta'_Z$ and $F'_+$,
with no common irreducible components, such that
\[
K_Z+\Delta'_Z+\boldsymbol{\beta}'_Z
\equiv
g^*\nu^*\alpha+F'_+,
\]
where $F'_+$ is $g$-exceptional. 
Since $\nu \colon X' \to X$ is quasi-\'etale, 
we deduce that $g_*\Delta'_Z=\nu^*\Delta=\Delta'$.
Pushing forward by $g$, we obtain
\[
K_{X'}+\Delta'+\boldsymbol{\beta}'_{X'}
\equiv
\nu^*\alpha.
\]

\smallskip

\textup{(2)}
After replacing $M$ by a higher resolution, we may assume that there
exists a morphism $p\colon M\to\bar X$ such that $\pi=q\circ p$
and that $M$ is a common log resolution:
\[
\begin{tikzcd}
& M \arrow[ld, "\pi"'] \arrow[rd, "p"] & \\
X & & \bar X. \arrow[ll, "q"']
\end{tikzcd}
\]
We regard $\boldsymbol{\beta}$ as a b-$(1,1)$-current over $\bar X$
and denote it by $\bar{\boldsymbol{\beta}}$.
By \eqref{eq-gpair}, there exist effective $\mathbb{R}$-divisors
$\Delta_M$ and $F_+$, with no common irreducible components, such that
\[
K_M+\Delta_M+\boldsymbol{\beta}_M
\equiv
\pi^*\alpha+F_+,
\]
where $F_+$ is $\pi$-exceptional.

Since $q$ is small, every $\pi$-exceptional divisor on $M$ is
$p$-exceptional.
Thus, we have $p_*\Delta_M=q^{-1}_*\Delta=\bar\Delta$.
Pushing forward by $p$, we obtain
\[
K_{\bar X}+\bar\Delta+\bar{\boldsymbol{\beta}}_{\bar X}
\equiv
q^*\alpha.
\]
\end{proof}

At the end of this subsection, we confirm \cite[Theorem~2.19]{DHY26} and \cite{Tak03} in our formulation. 
We begin by recalling the notion of potentially klt singularities.

\begin{defi}[{cf.~\cite[Lemmas~2.37 and~3.5]{Fuj22}}]
Let $X$ be a normal analytic variety.
We say that $X$ is \textit{locally potentially klt}
if it locally satisfies the following equivalent conditions:
\begin{itemize}
\item[\textup{(1)}]
There exists an effective $\mathbb{R}$-divisor
$\Delta_{\mathbb{R}}$ such that
$(X,\Delta_{\mathbb{R}})$ is a klt pair.

\item[\textup{(2)}]
There exists an effective $\mathbb{Q}$-divisor
$\Delta_{\mathbb{Q}}$ such that
$(X,\Delta_{\mathbb{Q}})$ is a klt pair.
\end{itemize}
\end{defi}

\begin{proof}[Proof of the equivalence]
It suffices to prove that \textup{(1)} implies \textup{(2)}.

Assume that there exists an effective $\mathbb{R}$-divisor
$\Delta_{\mathbb{R}}$ such that
$(X,\Delta_{\mathbb{R}})$ is klt.
Let $V$ be the finite-dimensional $\mathbb{R}$-vector space generated
by the irreducible components of $\Delta_{\mathbb{R}}$,
and let $V_{\mathbb{Q}}$ be the $\mathbb{Q}$-vector space generated
by the same components.
Since the subset
\[
V_K
:=
\left\{
B\in V
\,\middle|\,
K_X+B \text{ is $\mathbb{R}$-Cartier}
\right\}
\]
is an affine subspace defined over $\mathbb{Q}$, 
the rational points $V_K\cap V_{\mathbb{Q}}$ are dense in $V_K$.
Moreover, for every $B\in V_K\cap V_{\mathbb{Q}}$,
the divisor $K_X+B$ is $\mathbb{Q}$-Cartier
(see \cite[Lemma~2.37 and the proof of Lemma~3.5]{Fuj22}).

By the standard rational-polytope argument, the klt locus is open
in $V_K$ 
(see \cite[Lemma~3.5]{Fuj22} or \cite[Lemma~3.7.2]{BCHM10}.)
Indeed, after taking a log resolution, the relevant discrepancies
are affine-linear functions of the coefficients of $B$.
Thus, we can choose an effective rational point
$\Delta_{\mathbb{Q}}\in V_K\cap V_{\mathbb{Q}}$
 such that
$(X,\Delta_{\mathbb{Q}})$ remains klt.
By the property of $V_K$ recalled above,
we see that $K_X+\Delta_{\mathbb{Q}}$ is $\mathbb{Q}$-Cartier.
\end{proof}

The first proposition below is a direct consequence of
\cite[Theorem~2.19]{DHY26}, while the second follows from the first
proposition and \cite{Tak03}.
Note that Takayama's invariance theorem for fundamental groups applies
to varieties with locally potentially klt singularities.
The first proposition asserts that the underlying space is locally
potentially klt. However, it is not known whether there exists a global divisor
$\Delta$ such that $(X,\Delta)$ is klt.

\begin{prop}[{\cite[Theorem~2.19]{DHY26}}]
\label{prop-local-klt}
Let $(X,\Delta,\boldsymbol{\beta})$ be a generalized klt pair.
Then, the variety $X$ is locally potentially klt in the analytic topology.
\end{prop}

\begin{prop}[{\cite{Tak03}}]
\label{prop-funda}
Let $(X,\Delta,\boldsymbol{\beta})$ be a generalized klt pair,
and let $\pi\colon M\to X$ be a resolution of singularities of $X$.
Then, the induced homomorphism
\[
\pi_*\colon\pi_1(M)\to\pi_1(X)
\]
is an isomorphism.
\end{prop}

\subsection{Linear systems and Bertini type theorems on complex analytic spaces}
\label{subsec-Bertini}
In this subsection, we collect some Bertini-type theorems for
(finite-dimensional) linear systems on (possibly non-compact) complex analytic
spaces that will be used in the subsequent subsections. Before stating the
main results, we first recall the notion of an {\itshape analytically
sufficiently general} point, introduced by Fujino in
\cite[Definitions 2.49 and 2.50]{Fuj17}:
\begin{defi}
Let $Y$ be an analytic variety. A subset $S\subseteq Y$ is called analytically
meagre if
\[
S\subseteq\bigcup_{n\in\NN} Y_n,
\]
where each $Y_n$ is a locally analytic subset of codimension $\ge 1$ in $Y$
(\cite[Definition (0.1)]{Man82}). We say that a property~\textup{(P)} holds
for an {\itshape analytically sufficiently general} point $y\in Y$ if
property~\textup{(P)} holds for every $y\in Y\backslash S$ for some
analytically meagre subset $S\subseteq Y$. For a morphism $f: X\to Y$, we
say that a property~\textup{(P)} holds for an {\itshape analytically
sufficiently general} fibre of $f$ if property~\textup{(P)} holds for every
fibre $X_y$ with $y$ contained in $Y\backslash S$ for some analytically
meagre subset $S\subseteq Y$.
\end{defi}

\begin{rem}
In algebraic geometry, we have the notion of a {\itshape general} point,
and we can similarly define this notion in the category of analytic spaces:
we say that a property~\textup{(P)} holds for a general point of $Y$ if
property~\textup{(P)} holds for every $y$ contained in a Zariski open subset
of $Y$. However, Zariski open subsets do not behave well in complex analytic
geometry:
\begin{enumerate}
\item As explained in \cite[Remark 2]{Fuj17}, a Zariski open subset of a
Zariski open subset is usually no longer Zariski open in the ambient space.
\item A key ingredient in the proof of Bertini-type theorems is Frisch's
generic flatness theorem \cite[Th{\'e}or{\`e}me (IV,9)]{Fri67} (see also
\cite[\S V.4, Theorem V.4.5, p.~189]{BS76}), according to which the non-flat
locus of a morphism between analytic varieties is an analytic subset. We then
need to consider its image, outside of which the Bertini-type theorems hold.
However, unless the morphism is proper, the image of the non-flat locus is
in general no longer an analytic subset, but is only analytically meagre
(as shown in the theorem below). In particular, `analytically sufficiently
general', rather than `general', is the appropriate notion for stating
Bertini-type theorems for non-compact analytic spaces.
\end{enumerate}
\end{rem}

\begin{thm}[{\cite[Th{\'e}or{\`e}me (IV,9), Proposition (IV,14)]{Fri67}}]
\label{thm-generic-flatness}
Let $f:X\to Y$ be a morphism between analytic varieties, and let $\mathcal{F}$ be
a coherent sheaf on $X$. Then the non-$f$-flat locus of $\mathcal{F}$
\[
S:=\{x\in X\;|\;\mathcal{F}_x\text{{ is not flat over }}\OX_{Y,f(x)}\}
\]
is an analytic subset of $X$. Moreover, if $Y$ is reduced, then $f(S)$ is
analytically meagre in $Y$.
\end{thm}

We can now state the Bertini-type theorems:
\begin{thm}
\label{thm-Bertini}
Let $X$ be an analytic variety, let $L$ be a line bundle on $X$, and let
$V\subseteq H^0(X,L)$ be a finite-dimensional vector subspace that generates
$L$ at every point of $X$; in other words, $V$ is a basepoint-free linear
system on $X$. Assume that $\mathcal{F}$ is a coherent sheaf on $X$
satisfying property~\textup{(P)}, where property~\textup{(P)} is one of the
following:
\begin{itemize}
\item $\mathcal{F}$ is locally free;
\item $\mathcal{F}$ satisfies $(S_k)$;
\item $\mathcal{F}$ is reflexive;
\item $\mathcal{F}=\OX_X$, and $X$ satisfies $(R_k)$.
\end{itemize}
Then, for an analytically sufficiently general $v\in V$,
$\mathcal{F}|_{H_v}$ satisfies property~\textup{(P)}, where $H_v$ denotes
the effective divisor defined by $v$. In particular, if $X$ is smooth
(resp.\ normal, resp.\ reduced), then so is $H_v$ for an analytically
sufficiently general $v\in V$.
\end{thm}

Here, for a Noetherian ring $A$, a finite $A$-module $M$ is said to satisfy
Serre's condition $(S_k)$ if
$\depth_{A_\mathfrak{p}}(M_{\mathfrak{p}})
\ge \min(\dim M_{\mathfrak{p}},k)$
for every prime ideal $\mathfrak{p}$ of $A$. A coherent sheaf $\mathcal{F}$
on an analytic variety $X$ is said to satisfy Serre's condition $(S_k)$ if
$\mathcal{F}_x$ satisfies $(S_k)$ as an $\OX_{X,x}$-module for every
$x\in X$. Note that the ``in particular'' part of
Theorem~\ref{thm-Bertini} was essentially established in
\cite[Theorem II.5]{Man82}. Here, we give a uniform treatment using the more
precise notion of an analytically sufficiently general point.

For a linear system $V$ on $X$, let
$\mathscr{H}_V\subseteq V\times X$ be its universal family. More concretely,
if we choose a basis $\{v_i\}$ of $V$, then $\mathscr{H}_V$ is defined by
the equation $\sum_i a_i v_i(x)=0$, where the $a_i$ are the coordinates on
$V$ with respect to the basis $\{v_i\}$. Clearly, the fibre of
$\pi:=\pr_1|_{\scrH_V}:\mathscr{H}_V\to V$ over $v\in V$ is $H_v$, while
the fibre of $e:=\pr_2|_{\scrH_V}:\mathscr{H}_V\to X$ over $x\in X$ is the
subspace of $V$ consisting of those $v\in V$ such that $x\in H_v$.
Moreover, if $V$ is basepoint-free, then
$e:\mathscr{H}_V\to X$ is an affine bundle
(see \cite[Lemma 3.4.4, p.~108]{FOV99}). Hence, to prove
Theorem~\ref{thm-Bertini}, it suffices to establish the following generic
principle:
\begin{prop}
\label{prop-generic}
Let $f:X\to Y$ be a morphism between analytic varieties, and assume that $Y$ is
smooth. Let $\mathcal{F}$ be a coherent sheaf on $X$ satisfying
property~\textup{(P)}, where property~\textup{(P)} is one of those listed
in Theorem~\ref{thm-Bertini}. Then $\mathcal{F}|_{X_y}$ satisfies
property~\textup{(P)} for an analytically sufficiently general $y\in Y$,
where $X_y$ denotes the fibre of $f$ over $y$.
\end{prop}

The proof of Proposition~\ref{prop-generic} is essentially the same as that
of \cite[Theorem 3.3.10]{FOV99}.

As a corollary of Theorem~\ref{thm-Bertini}, we obtain the following result,
which will be used in the proof of the Bertini-type theorem for foliations:
\begin{cor}
\label{cor-Bertini}
Let $X$ be an analytic variety, let $L$ be a line bundle on $X$, and
let $V\subseteq H^0(X,L)$ be a finite-dimensional vector subspace that
generates $L$ at every point of $X$. Assume that $\mathcal{G}$ is a
reflexive sheaf on $X$ and that $\mathcal{F}$ is a saturated subsheaf of
$\mathcal{G}$, i.e., $\mathcal{G}/\mathcal{F}$ is torsion free. Then, for
an analytically sufficiently general $v\in V$, $\mathcal{G}|_{H_v}$ is
reflexive and
$(\calG/\calF)|_{H_v}\simeq\calG|_{H_v}/\calF|_{H_v}$ is torsion free.
In particular, $\mathcal{F}|_{H_v}$ is a saturated subsheaf of
$\mathcal{G}|_{H_v}$.
\end{cor}

\begin{proof}
For every $v\in V$, the short exact sequence
\[
0\to \OX_X(-H_v)\to \OX_X\to\OX_{H_v}\to 0
\]
induces the exact sequence
\[
\SheafTor_{\,1}^{\OX_X}(\mathcal{G}/\mathcal{F},\OX_{H_v})
\to \mathcal{G}/\mathcal{F}\otimes\OX_X(-H_v)
\to \mathcal{G}/\mathcal{F}
\to(\mathcal{G}/\mathcal{F})|_{H_v}\to 0.
\]
Since $\SheafTor_{\,1}^{\OX_X}(\mathcal{G}/\mathcal{F},\OX_{H_v})$ is a
torsion sheaf, whereas
$\mathcal{G}/\mathcal{F}\otimes\OX_X(-H_v)$ is torsion free, the image in
$\mathcal{G}/\mathcal{F}\otimes\OX_X(-H_v)$ must be zero. We thus
obtain a short exact sequence

\[
0\to \mathcal{F}|_{H_v}\to \mathcal{G}|_{H_v}
\to(\mathcal{G}/\mathcal{F})|_{H_v}\to 0.
\]
By Theorem~\ref{thm-Bertini}, for an analytically sufficiently general
$v\in V$, $H_v$ is reduced, $\mathcal{G}|_{H_v}$ is reflexive, and
$(\mathcal{G}/\mathcal{F})|_{H_v}$ satisfies $(S_1)$ and is therefore
torsion free. This proves the corollary.
\end{proof}

\begin{rem}
\label{rmk-Bertini}
Finally, let us point out that, unlike in algebraic geometry, the Bertini
irreducibility theorem fails for analytic varieties, as illustrated by the
following example. Let
$X=D\times\CC_{z_3}\subseteq\CC^3_{z_1,z_2,z_3}$, where $D$ is a small
polydisc in $\CC^2_{z_1,z_2}$ containing no integral point of $\CC^2$, and
let
\[
V:=\text{span}\{e^{2\pi iz_1}-1,e^{2\pi iz_2}-1,e^{2\pi iz_3}-1\}
\subseteq H^0(X,\OX_X)
\]
be a linear system. Clearly, $V$ is basepoint-free: regarded as a linear
system on $\CC^3$, its base locus is $\ZZ^3$. Moreover, it is not composed
of a pencil, since the three exponential functions are algebraically
independent. However, a general member of $V$ has infinitely many connected
components. An example in dimension $2$ can be obtained by slightly
modifying the example above.
\end{rem}

\subsection{Strong \texorpdfstring{$\mathbb{Q}$}{text}-factorializations}
\label{subsec-Q-fact}

In this subsection, we discuss the existence of small
$\mathbb Q$-factorializations for generalized pairs in the analytic
setting.

\begin{defi}\label{def:strongly_Q_fac}
Let $X$ be a normal analytic variety. We say that $X$ is strongly $\mathbb Q$-factorial if every Weil divisor on $X$ is $\mathbb Q$-Cartier, i.e. for every Weil divisor $D$ on $X$, there exists an integer $m>0$ such that $mD$ is Cartier.
\end{defi}

The following proposition is a generalization of \cite[Theorem 2.19(2)]{DHY26}.

The existence of small
$\mathbb{Q}$-factorializations for generalized pairs has been proved locally in \cite[Theorem 2.19]{DHY26}, 
but we need the global version.

\begin{thm}[see {\cite[Lemma 2.11, Lemma 2.12]{CH24}}]
\label{thm-Q-factorialization}

Let $X$ be a compact normal analytic variety, and let $\Delta$ be an effective divisor on $X$ such that $\lfloor \Delta\rfloor=0$.

Then there exists a strong $\QQ$-factorialization of $(X,\Delta)$, namely, a compact strongly $\QQ$-factorial klt pair $(X^{\qf},\Delta^{\qf})$ equipped with a projective modification $q\colon X^{\qf}\to X$ satisfying the following properties:
\begin{itemize}
\item[\rm(a)] $q_\ast \Delta^{\qf}=\Delta$;
\item[\rm(b)] $K_{X^\qf}+\Delta^{\qf}$ is $q$-nef; in other words, $(X^\qf,\Delta^\qf)$ is a relative minimal model over $X$;
\item[\rm(c)] Suppose that there exist an open subset $X_0$ of $X$ and a nef b-$(1,1)$-current $\boldsymbol{\beta}$ on $X$ such that the generalized pair $(X,\Delta,\boldsymbol{\beta})$ has klt singularities on $X_0$. Set $\Delta_0:=\Delta|_{X_0}$, $\boldsymbol{\beta}_0:=\boldsymbol{\beta}|_{X_0}$, $X_0^{\qf}:=q^{-1}(X_0)$, $\Delta^{\qf}_0:=\Delta^{\qf}|_{X_0^{\qf}}$, and $q_0:=q|_{X_0^{\qf}}$. Then $q_0$ is a small contraction, and we have
\[
K_{X^{\qf}_0}+\Delta^{\qf}_0+\boldsymbol{\beta}_{X^{\qf}}|_{X^{\qf}_0}\equiv q_0^\ast(K_{X_0}+\Delta_0+\boldsymbol{\beta}_0);
\]
\item[\rm(d)] In the setting of \textup{(c)}, if, moreover, $X_0$ is strongly $\QQ$-factorial, then $q_0$ is an isomorphism;
\item[\rm(e)] If $X$ is K\"ahler, then so is $X^\qf$.
\end{itemize}

In particular, if $(X,\Delta,\boldsymbol{\beta})$ is a g-klt pair (i.e., if we can take $X_0=X$ in \textup{(c)}), then $q$ is a small strong $\QQ$-factorialization of $(X,\Delta,\boldsymbol{\beta})$.
\end{thm}

\begin{proof}
Take a log resolution $f\colon W\to X$ of $(X,\Delta)$, and let $\Phi=\sum_i\Phi_i$ be the exceptional divisor of $f$.
Then, for any sufficiently small $\varepsilon>0$, we can run the relative log-MMP for $K_W+f^{-1}_\ast\Delta+(1-\varepsilon)\Phi$ over $X$, as in \cite[Theorem~1.4]{DHP24} (see also \cite[Theorem~1.7]{Fuj22}), to obtain a bimeromorphic map $g\colon W\dashrightarrow X^{\qf}$. The pair
\[
(X^{\qf},\Delta^{\qf}:=g_\ast f^{-1}_\ast\Delta+(1-\varepsilon)g_*\Phi)
\]
is a strongly $\mathbb{Q}$-factorial klt pair and is equipped with a bimeromorphic morphism $q\colon X^\qf\to X$ such that $K_{X^{\qf}}+\Delta^{\qf}$ is $q$-nef. Clearly, $q_\ast\Delta^\qf=\Delta$. Thus, \textup{(a)} and \textup{(b)} follow.

We now assume the setting of \textup{(c)}. After reordering the components if necessary, we may assume that $\Phi_1,\ldots,\Phi_k$ are precisely the components of $\Phi$ intersecting $W_0:=f^{-1}(X_0)$, and we set $f_0:=f|_{W_0}$. Moreover, we may assume that $f$ is a log resolution of the generalized pair $(X,\Delta,\boldsymbol{\beta})$. Since this pair has klt singularities on $X_0$, we can write
\begin{equation}
\label{eq-X_0-klt}
K_{W_0}+(f_0)^{-1}_\ast\Delta_0+\boldsymbol{\beta}_W|_{W_0}
\equiv
f_0^\ast(K_{X_0}+\Delta_0+\boldsymbol{\beta}_0)
+\sum_{i=1}^k a_i\Phi_i|_{W_0}
\end{equation}
with $a_i>-1$ for all $i$.
Choosing $\varepsilon$ sufficiently small so that
$\varepsilon<\min_i\{1+a_i\}$, we see that
\[
\Psi^\circ:=\sum_{i=1}^k(1+a_i-\varepsilon)
g_*\Phi_i|_{X^\qf_0}\geq 0
\]
is $q_0$-exceptional and $q_0$-nef. Pushing forward
\eqref{eq-X_0-klt}, we obtain
\[
K_{X_0^\qf}+\Delta_0^{\qf}
+\boldsymbol{\beta}_{X^{\qf}}|_{X_0^{\qf}}
\equiv
q_0^\ast(K_{X_0}+\Delta_0+\boldsymbol{\beta}_0)+\Psi^\circ.
\]
By the generalized negativity lemma below, we conclude that
$\Psi^\circ=0$. In other words, every component of
$\Phi=\mathrm{Exc}(f)$ intersecting $W_0$ is contracted by $g$.
Therefore, the restriction $q_0$ is a small contraction. Moreover, pushing
forward \eqref{eq-X_0-klt} by $g_\ast$, we obtain
\[
K_{X^{\qf}_0}+\Delta^{\qf}_0
+\boldsymbol{\beta}_{X^{\qf}}|_{X^{\qf}_0}
\equiv
q_0^\ast(K_{X_0}+\Delta_0+\boldsymbol{\beta}_0).
\]
This proves \textup{(c)}. The ``in particular'' assertion follows by taking
$X_0=X$.

If, moreover, $X_0$ is strongly $\mathbb{Q}$-factorial, then
\cite[Lemma~2.4]{DH25} implies that $q_0$ is an isomorphism. This proves
\textup{(d)}.
\end{proof}

\begin{lemm}[Negativity lemma with a b-nef term]
\label{lem-negativity-bnef}
Let $h\colon Y\to X$ be a projective bimeromorphic morphism between normal
analytic varieties. Let $F\geq 0$ be an $h$-exceptional $\mathbb R$-divisor
on $Y$. Let $\boldsymbol{\beta}$ be a nef b-$(1,1)$-current on $X$, and let
$\boldsymbol{\beta}_Y$ be its trace on $Y$.
Assume that $Y$ is strongly $\mathbb Q$-factorial and that
$F-\boldsymbol{\beta}_Y$ admits local potentials and is $h$-nef. Then
$F=0$.
\end{lemm}

\begin{proof}
Suppose, for contradiction, that $F\neq 0$.

We first show, following the proof of the classical negativity lemma, that
there exists an $h$-contracted curve $C\subset\Supp(F)$ passing through a
general point of some component of $\Supp(F)$ such that $F\cdot C<0$.

Indeed, by the proof of \cite[III.5.1, pp.~103--104]{Nak04}, we can find an
$h$-contracted curve $C\subseteq\Supp(F)$ passing through a general point of
some component of $\Supp(F)$ such that $F\cdot C<0$. For the convenience of
the reader, we briefly recall the argument. Since the problem is local, we
may assume that $X$ is a Stein variety. Set
$n:=\dim Y$ and $d:=\dim h(\Supp(F))$. Cut $X$ by $d$ analytically
sufficiently general hypersurfaces to obtain a general complete intersection
subvariety $V$ such that $V\cap h(\Supp(F))$ is $0$-dimensional and both
$V$ and $h^{-1}(V)$ are normal by Theorem~\ref{thm-Bertini}
(see also \cite[Theorem II.5]{Man82}). We then cut $h^{-1}(V)$ by
$n-d-2$ relatively ample hypersurfaces on $Y$ to obtain a surface $T'$,
which is again normal, and every component of $T'$ maps bimeromorphically
onto a surface in $X$. Choose a component $T$ of $T'$ that meets
$\Supp(F)$, and set $S:=h(T)$. Applying the classical negativity lemma
\cite[Lemma 1.3]{JWang19} to $h|_T:T\to S$ and $-F|_T$, we obtain an
$h$-contracted curve $C\subseteq T$ such that
$F|_T\cdot C=F\cdot C<0$. By construction, $C$ moves in a family covering
a component of $\Supp(F)$.

Take a higher model $\alpha:M\to Y$ such that the trace
$\boldsymbol{\beta}_M$ on $M$ is nef. Since $C$ passes through a general
point of a component of $\Supp(F)$, it is not contained in the exceptional
locus of $\alpha$. Let $C_M$ be the strict transform of $C$ on $M$.
Since $Y$ is strongly $\QQ$-factorial and
$F-\boldsymbol{\beta}_Y$ admits local potentials,
$\boldsymbol{\beta}_Y$ also admits local potentials. Since
$[\boldsymbol{\beta}_M]$ is nef, we can choose a closed semipositive
$(1,1)$-current $\gamma_M\in[\boldsymbol{\beta}_M]$, which induces a closed
semipositive $(1,1)$-current
$\gamma_Y\in[\boldsymbol{\beta}_Y]$ by
\cite[Proposition 4.6.3, pp.~230--231]{BEG13}. Applying the second theorem
of support \cite[(III.2.14), p.~143]{agbook} to
$\gamma_M-\alpha^\ast\gamma_Y$, which has order $0$ because both currents
have order $0$, we can write
\[
\boldsymbol{\beta}_M
\equiv
\alpha^\ast\boldsymbol{\beta}_Y+G
\]
for some $\alpha$-exceptional $\RR$-divisor $G$. Since
$\boldsymbol{\beta}_M$ is nef, $G$ is $\alpha$-nef. The negativity lemma
\cite[Lemma 1.3]{JWang19} therefore implies that $G\leq 0$. Since $C$
passes through a general point of a component of $\Supp(F)$, the curve
$C_M$ is not contained in $\Supp(G)$, and hence $G\cdot C_M\leq 0$.
Consequently,
\[
F\cdot C
\geq
\boldsymbol{\beta}_Y\cdot C
\geq
\boldsymbol{\beta}_M\cdot C_M
\geq 0,
\]
where the first and last inequalities follow from the $h$-nefness of
$F-\boldsymbol{\beta}_Y$ and the nefness of $\boldsymbol{\beta}_M$,
respectively. This contradicts $F\cdot C<0$.
\end{proof}

\begin{rem}\label{rem-klt-type}
The following simple observation will be used repeatedly in this paper when
applying results concerning klt singularities. Let
$(X,\Delta,\boldsymbol{\beta})$ be a generalized klt pair. If $X$ is
strongly $\mathbb{Q}$-factorial, then $X$ itself has klt singularities.

Indeed, we can write
\[
K_M+\Delta'_M
\equiv
\pi^*(K_X+\Delta)+F'_+,
\]
where $F'_+$ is an effective $\pi$-exceptional $\mathbb{R}$-divisor and
$\Delta'_M$ is an effective $\mathbb{R}$-divisor such that $\Delta'_M$ and
$F'_+$ have no common components. Moreover, since $X$ is strongly
$\mathbb{Q}$-factorial, the log-canonical divisor, and hence the trace
$\boldsymbol{\beta}_X$, determines a Bott--Chern class on $X$. In
particular, the pullback of $\boldsymbol{\beta}_X$ is well-defined.
Comparing the above relation with \eqref{eq-gpair}, we obtain
\[
\Delta'_M-F'_+
\equiv
\Delta_M
-F_+-
\bigl(\pi^*\boldsymbol{\beta}_X-\boldsymbol{\beta}_M\bigr).
\]

By the second theorem of support and the classical negativity lemma, the
class
$c_{1}(\pi^*\boldsymbol{\beta}_X-\boldsymbol{\beta}_M)$
is represented by an effective $\pi$-exceptional $\mathbb{R}$-divisor.

Hence,
\[
\Delta'_M\leq\Delta_M.
\]
The generalized klt condition $\Delta_M<1$ implies that
$\Delta'_M<1$, and hence $(X,\Delta)$ is a klt pair. Since the boundary
divisor $\Delta$ is $\mathbb{Q}$-Cartier, it follows that $X$ has klt
singularities.

\end{rem}

\subsection{Maximally quasi-\'etale covers}\label{subsec-quasi-etale}

In this subsection, we discuss maximally quasi-\'etale covers.
Combining the result on small
$\mathbb{Q}$-factorializations established in
Theorem~\ref{thm-Q-factorialization}, we prove the existence of maximally
quasi-\'etale covers for generalized klt pairs.

\begin{defi}
A normal analytic variety $X$ is said to be \textit{maximally quasi-\'etale} if the natural
homomorphism of profinite fundamental groups 
\[
\iota_*\colon
\widehat{\pi}_1(X_{\reg})
\longrightarrow
\widehat{\pi}_1(X),
\]
induced by the inclusion $\iota\colon X_{\reg}\hookrightarrow X$ is an isomorphism.
A finite quasi-\'etale cover $X'\to X$ is called a
\textit{maximally quasi-\'etale cover} of $X$ if $X'$ is maximally
quasi-\'etale.
\end{defi}

\begin{thm}[{cf.~\cite[Theorem~1.5]{GKP16b}, \cite[Theorem~2.5]{FO25},
and \cite[Theorem~2.12]{MW}}]
\label{thm-maxqetale}
Let $(X,\Delta,\boldsymbol{\beta})$ be a generalized klt pair, where $X$ is compact.
Then, the following assertions hold:

\begin{itemize}
\item[$(1)$] There exists a maximally quasi-\'etale cover
$\nu\colon X'\to X$.
\item[$(2)$] Suppose that $X$ is maximally quasi-\'etale.
Then, for every linear representation
$\rho_0\colon\pi_1(X_{\reg})\to\GL(r,\mathbb{C})$, 
there exists a representation
$\rho\colon\pi_1(X)\to\GL(r,\mathbb{C})$ such that
$\rho_0=\rho\circ\iota_*$:
\[
\xymatrix@C=40pt@R=30pt{
\pi_1(X_{\reg})
\ar@/^18pt/[rr]^{\rho_0}
\ar@{->>}[r]_{\iota_*}
&
\pi_1(X)
\ar[r]_{\rho}
&
\GL(r,\mathbb{C}).
}
\]
\item[$(3)$] Let $\pi\colon\bar X\to X$ be a bimeromorphic morphism, 
where $\bar X$ is the underlying space of a generalized klt pair.
If $X$ is maximally quasi-\'etale, then so is $\bar X$.
\end{itemize}
\end{thm}

\begin{proof}
\textup{(1)}
If $X$ were globally potentially klt, then Conclusion~\textup{(1)}
would follow directly from \cite[Subsection~2.B]{FO25}.
By Proposition~\ref{prop-local-klt}, the variety $X$ is locally
potentially klt in the analytic topology, and $X$ need not be globally
potentially klt.
To overcome this difficulty, by Theorem~\ref{thm-Q-factorialization}, 
we take a small $\mathbb{Q}$-factorialization
$
q\colon X^{\qf}\to X.
$
Then, by Proposition \ref{prop-pullback} and Remark \ref{rem-klt-type}, the variety $X^{\qf}$ has globally potentially klt singularities.
Hence, there exists a maximally
quasi-\'etale Galois cover
$\nu\colon Y\to X^{\qf}$.

Let
$
Y\xrightarrow{p}X'\xrightarrow{\tau}X
$
be the Stein factorization of the composition
$Y\xrightarrow{\nu}X^{\qf}\xrightarrow{q}X$.
We have the following commutative diagram:
\[
\xymatrix@C=4em@R=3em{
Y \ar[r]^{p} \ar[d]_{\nu}
&
X' \ar[d]^{\tau}
\\
X^{\qf} \ar[r]_{q}
&
X.
}
\]

We show that $\tau$ is quasi-\'etale.
Since $q$ is small, there exists a Zariski-open subset $U\subset X$
such that
$\codim(X\setminus U)\geq 2$,
$\codim(X^{\qf}\setminus q^{-1}(U))\geq 2$, and
$q^{-1}(U)\to U$ is an isomorphism.
After shrinking $U$ further by removing a subset of codimension at
least two, we may also assume that $\nu$ is \'etale over $q^{-1}(U)$.
Over $U$, the morphism $\tau$ is identified with $\nu$ via the
isomorphism $q^{-1}(U)\cong U$.
This shows that $\tau$ is \'etale over $U$ and thus $\tau$ is quasi-\'etale.
The same argument shows that $p\colon Y\to X'$ is a small bimeromorphic morphism.

It remains to prove that $X'$ is maximally quasi-\'etale.
Take a common resolution $W$ with the following diagram:
\[
\begin{tikzcd}
& W \arrow[ld, "a"'] \arrow[rd, "b"] & \\
Y \arrow[rr, "p"'] & & X'. 
\end{tikzcd}
\]
By Proposition~\ref{prop-pullback}, the varieties
$Y$ and $X'$ are underlying varieties of generalized klt pairs.
Thus, by Proposition~\ref{prop-funda}, both
$a_*\colon\pi_1(W)\to\pi_1(Y)$ and
$b_*\colon\pi_1(W)\to\pi_1(X')$ are isomorphisms.
It follows that
$
p_*\colon\pi_1(Y)\longrightarrow\pi_1(X')
$
is an isomorphism.

Since $p$ is small, there exists a Zariski-open subset
$U'\subset X'_{\reg}$ such that
$\codim(X'\setminus U')\geq 2$ and $p$ is an isomorphism over $U'$. 
After shrinking $U'$ further, 
we may assume that $V':=p^{-1}(U')$ satisfies 
$$
V'\subset Y_{\reg}, \quad 
\codim (Y\setminus V')\geq 2,  \quad 
p|_{V'}\colon V'\cong U'.
$$
By the condition of codimension,  the natural inclusions induce
isomorphisms
\[
\pi_1(V')
\cong
\pi_1(Y_{\reg})
\qquad\text{and}\qquad
\pi_1(U')
\cong
\pi_1(X'_{\reg}).
\]
Consequently, we obtain natural identifications
\[
\widehat{\pi}_1(X'_{\reg})
\cong
\widehat{\pi}_1(Y_{\reg})
\qquad\text{and}\qquad
\widehat{\pi}_1(X')
\cong
\widehat{\pi}_1(Y).
\]
Under these identifications, the natural homomorphism
$\widehat{\pi}_1(X'_{\reg})\to\widehat{\pi}_1(X')$
corresponds to
$\widehat{\pi}_1(Y_{\reg})\to\widehat{\pi}_1(Y)$.
The latter is an isomorphism since $Y$ is maximally quasi-\'etale.
This shows that $X'$ is maximally quasi-\'etale.
This proves~\textup{(1)}.

\textup{(2)}
Conclusion~\textup{(2)} follows from Malcev's theorem by the argument
in \cite{GKP16b} (see also \cite[Theorem~2.9]{CDM26}.)

\textup{(3)}
Conclusion~\textup{(3)} follows from the proof of
\cite[Theorem~2.12]{MW}.
\end{proof}

\subsection{Singular Hermitian metrics on torsion-free sheaves}
\label{subsec-pre}

In this subsection, we prepare several results on singular Hermitian
metrics for use in the proof of Theorem~\ref{thm-faked-flat}.
We begin by briefly recalling the definitions of singular Hermitian
metrics and their positivity (see \cite{Mat23} for details).

Let $\mathcal{E}$ be a torsion-free sheaf on a normal analytic variety
$X$.
By a \textit{singular Hermitian metric} on $\mathcal{E}$, we mean a
singular Hermitian metric on the vector bundle
$\mathcal{E}|_{X_{\mathcal{E}}\cap X_{\reg}}$, where
$X_{\mathcal{E}}\subset X$ is the largest Zariski-open subset on which
$\mathcal{E}$ is locally free.
We often refer to singular Hermitian metrics simply as ``metrics.''

Let $\theta$ be a smooth $(1,1)$-form on $X$ with local potentials.
We write
\[
\sqrt{-1}\Theta_h\geq\theta\otimes\id
\quad\text{on }X
\]
if, for every open subset $U\subset X$ and every local section
$e\in H^0(U,\mathcal{E}^*)$, the function
$\log|e|_{h^*}-\psi$ is plurisubharmonic (psh) on
$U\cap X_{\mathcal{E}}\cap X_{\reg}$, where $\psi$ is a local potential
of $\theta$ and $h^*$ is the metric induced on the dual sheaf
$\mathcal{E}^*$.
Note that this function extends to a psh
function on $U$ by the Hartogs-type extension theorem.

\begin{defi}\label{def-weak}
Let $\mathcal{E}$ be a torsion-free sheaf on a normal analytic variety
$X$.

\textup{(1)}
Let $\theta$ be a smooth $(1,1)$-form on $X$ with local potentials,
and let $\omega$ be a K\"ahler form on $X$. 
The sheaf $\mathcal{E}$ is said to be
\textit{$\theta$-weakly positively curved}
if, for every $\varepsilon>0$, there exist singular Hermitian metric
$\{g_{\varepsilon}\}_{\e>0}$  on $\mathcal{E}$ such that
\[
\sqrt{-1}\Theta_{g_{\varepsilon}}
\geq
(\theta - \varepsilon\omega) \otimes\id
\quad\text{on }X.
\]
We say that $\mathcal{E}$ is \textit{weakly positively curved}
when $\theta=0$.

\textup{(2)}
Suppose that $X$ is a compact K\"ahler variety and that $\omega$ is
a K\"ahler form on $X$.
The sheaf $\mathcal{E}$ is said to be \textit{pseudo-effective}
if there exists singular
Hermitian metrics $\{h_m\}_{m=1}^{\infty}$ on $\operatorname{Sym}^{[m]}\mathcal{E}$
such that
\[
\sqrt{-1}\Theta_{h_m}
\geq
-\omega\otimes\id
\quad\text{on }X.
\]
\end{defi}

By applying the positivity results for direct image sheaves established in
\cite{PT18,HPS18}, we obtain a positivity result in a form suited to
our later applications.

\begin{prop}\label{prop-psef}
Let $f\colon X\to Y$ be a fibration between compact K\"ahler manifolds.
Consider the following data:
\begin{itemize}
\item[$\bullet$]
$G$ is a relatively $f$-big line bundle on $X$;
\item[$\bullet$]
$N$ is a nef $\mathbb{R}$-line bundle on $X$;
\item[$\bullet$]
$L$ is an $\mathbb{R}$-line bundle admitting a singular Hermitian
metric $h_L$ with semipositive curvature such that
$\mathcal{I}(h_L)|_{X_y}=\mathcal{O}_{X_y}$ for a general fiber $X_y$.
\end{itemize}
For every sufficiently large integer $m\gg 1 $, let $P=P_m$ be a line bundle
on $X$ such that
\[
P\equiv G+mN+mL
\quad\text{in}\quad
H_{\BottChern}^{1,1}(X)
\]
and  $f_*\OX_X(mK_{X/Y}+P)$ is nonzero.
Then, the following assertions hold:
\begin{itemize}
\item[\textup{(1)}]
The $\mathbb{R}$-line bundle $K_{X/Y}+N+L$ is pseudo-effective.

\item[\textup{(2)}]
Let $\theta$ be a smooth $(1,1)$-form on $Y$ with local potentials.
Assume further that $G$ is $f^*\theta$-pseudo-effective, that is,
there exists  a singular Hermitian metric $h_{\theta}$  on $G$ such that
\[
\sqrt{-1}\Theta_{h_{\theta}}
\geq
f^*\theta.
\]
Then, the direct image sheaf
$f_*\OX_X(mK_{X/Y}+P)$ is $\theta$-weakly positively curved.
\end{itemize}
\end{prop}

\begin{proof}
We first remark that $f\colon X\to Y$ is a projective fibration 
by \cite[Theorem~1.1]{CH24}.

\textup{(1)}
By assumption, there exists a singular Hermitian metric $h$ on $G$
such that
\[
\sqrt{-1}\Theta_h+f^*\omega_Y
\geq
\omega_X,
\]
and there exist smooth Hermitian metrics
$\{g_{\varepsilon}\}_{\varepsilon>0}$ on $mN$ such that
$\sqrt{-1}\Theta_{g_{\varepsilon}}\geq-\varepsilon\omega_X$.
Here $\omega_X$ and $\omega_Y$ are K\"ahler forms on $X$ and $Y$,
respectively.

The singular Hermitian metric
$H:=h\cdot g_{\varepsilon}\cdot h_L^m$ on
$G+mN+mL$ satisfies
$\sqrt{-1}\Theta_H\geq-f^*\omega_Y$.
Using the numerical equivalence $P\equiv G+mN+mL$, we regard $H$ as a singular
Hermitian metric on $P$ satisfying the same curvature inequality.
Thus, we obtain 
\[
\sqrt{-1}\Theta_{H\cdot e^{-f^*\psi}}
\geq
0,
\]
where $\psi$ is a local potential of $\omega_Y$.

Consider the $m$-Bergman kernel metric $B_m^{-1}$ on
$mK_{X/Y}+P$, which is numerically equivalent to
$mK_{X/Y}+G+mN+mL$, induced by $H$
(see \cite[\S~4.2]{PT18} for the definition).
Let $B_m^{\prime-1}$ be the $m$-Bergman kernel metric induced by
$H\cdot e^{-f^*\psi}$, which is defined locally over $Y$.
By the assumption on $h_L$, the $L^{2/m}$-multiplier ideal sheaf
$\mathcal{J}_m(H^{1/m})$ satisfies
\begin{align}\label{eq-Lm}
\mathcal{J}_m\bigl(H^{1/m}|_{X_y}\bigr)
=
\mathcal{J}_m\bigl(h^{1/m}h_L|_{X_y}\bigr)
=
\mathcal{O}_{X_y}
\end{align}
for a general fiber $X_y$ and $m\gg 1$.
Here we use the solution of the strong openness conjecture
and H\"older's inequality.
Thus, by \cite[Theorem~0.1]{BP10}
(see also \cite[Theorem~4.2.7]{PT18} and
\cite[Theorem~3.5]{Cao17}), we deduce that the curvature of
$B_m^{\prime-1}$ is semipositive.
By construction, we see that 
\[
B_m^{\prime-1}
=
B_m^{-1}\cdot e^{-f^*\psi}.
\]
This shows that
\[
\sqrt{-1}\Theta_{B_m^{-1}}(mK_{X/Y}+P)
\geq
-f^*\omega_Y. 
\]
Dividing the corresponding class by $m$ and letting $m\to\infty$,
we conclude that $K_{X/Y}+N+L$ is pseudo-effective.

\textup{(2)}
After replacing $\omega_Y$ by a sufficiently large multiple,
we may assume that $\omega_Y\geq\theta$.
For $0<\varepsilon<\delta$, the singular Hermitian metric
\[
H_{\delta}
:=
h^{\delta}\cdot h_{\theta}^{1-\delta}
\cdot g_{\varepsilon}\cdot h_L^m
\]
on $G+mN+mL$ satisfies
\[
\begin{aligned}
\sqrt{-1}\Theta_{H_{\delta}}
&\geq
-\delta f^*\omega_Y
+\delta\omega_X
-\varepsilon\omega_X
+(1-\delta)f^*\theta \\
&\geq
f^*(\theta-2\delta\omega_Y).
\end{aligned}
\]
As above, using $P\equiv G+mN+mL$ again, we regard $H_{\delta}$
as a singular Hermitian metric on $P$ with the same curvature
lower bound.
The $L^{2/m}$-multiplier ideal sheaf
$\mathcal{J}_m(H_{\delta}^{1/m})$ satisfies the same property
as in \eqref{eq-Lm}.
Hence, \cite[Theorem~21.1]{HPS18}
(see also \cite[Theorem~5.3]{DWZZ18}) shows the curvature of the induced singular Hermitian
metrics on $f_*\OX_X(mK_{X/Y}+P)$ is bounded from below by  $(\theta-2\delta\omega_Y)\otimes\id$.
Letting $\delta\to 0$, we obtain the desired
the desired $\theta$-weak positivity.
\end{proof}

In the following lemmas, we examine how weak positive curvature is
preserved under pullbacks by fibrations and pushforwards by
bimeromorphic morphisms.
The same arguments also apply to pseudo-effectivity and may be useful
more generally for studying the bimeromorphic behavior of positivity
properties.

\begin{lemm}\label{lem-push}
Let $\pi\colon M\to X$ be a bimeromorphic morphism between normal
compact K\"ahler varieties, and let $\mathcal{E}$ be a torsion-free sheaf
on $M$ that is weakly positively curved on $M\setminus E_{\pi}$
with respect to a K\"ahler form $\omega_M$ defined on $M$,
where $E_{\pi}$ is the exceptional locus of $\pi$.
Then the pushforward $\pi_*\mathcal{E}$ is weakly positively curved
on $X$.
\end{lemm}

\begin{proof}

Since $\mathcal{E}$ is weakly positively curved on
$M\setminus E_{\pi}$ with respect to $\omega_M$,
there exist singular Hermitian metrics
$\{g_{\varepsilon}\}_{\e>0}$ on $\mathcal{E}$ such that
\[
\sqrt{-1}\Theta_{g_{\varepsilon}}
\geq
-\varepsilon\omega_M\otimes\id
\quad\text{on }M\setminus E_{\pi}.
\]
We would like to push forward the metrics and the preceding curvature
inequality to obtain the desired curvature lower bound.
The main difficulty is that the pushforward $\pi_*\omega_M$ of the
K\"ahler form $\omega_M$ need not define a Bott--Chern class on $X$
(i.e.,\,$\pi_*\omega_M$ need not admit local potentials).
For this reason, instead of using $\omega_M$ directly, we construct
a current $T$ representing $c_{1}(\pi^*\omega_X)-c_1(G)$ whose
pushforward defines a Bott--Chern class on $X$ by
Proposition~\ref{BEG}.
Here $\omega_X$ is a K\"ahler form on $X$, and $G$ is an effective
$\pi$-exceptional divisor.

After replacing $M$ by a further blow-up, we may assume that $M$ is smooth. 
Take a sequence of blow-ups $p\colon\tilde{X}\to X$ with smooth centers
and a projective bimeromorphic morphism
$q\colon\tilde{X}\to M$ fitting into the following commutative
diagram:
\begin{equation}\label{eq:hironaka}
\xymatrix{
\tilde{X} \ar[d]_q \ar[rd]^p & \\
M \ar[r]_\pi & X.
}
\end{equation}
Since $p$ is projective, there exists an effective
$p$-exceptional $\mathbb{Q}$-divisor $E$ such that
$p^*c_{1}(\omega_X)-c_1(E)$ is a K\"ahler class on $\tilde{X}$.
Choose a K\"ahler form $\omega_{\tilde{X}}$ 
such that $\omega_{\tilde{X}}\geq q^*\omega_M$ 
and $\omega_{\tilde{X}} \in p^*c_{1}(\omega_X)-c_1(E)$.

Set $T:=q_*\omega_{\tilde{X}}$.
Then $T$ is a closed current on $M$ such that
$T\geq\omega_M$ and 
\[
\begin{aligned}
c_{1}(T)
=
q_*
\bigl(
p^* c_{1}(\omega_X)-c_1(E)
\bigr) 
=
\pi^* c_{1}(\omega_X)-c_1(G),
\end{aligned}
\]
where $G:=q_*E$ is an effective $\pi$-exceptional
$\mathbb{Q}$-divisor.

Hence, there exists a quasi-psh function $\psi$ on $M$ such that
\[
T
=
\pi^*\omega_X-[G]+\ddbar\psi.
\]

Over $X\setminus\pi(E_{\pi})$, the morphism $\pi$ is an isomorphism,
and hence $\psi$ induces a quasi-psh function
$\varphi^{\circ}$ on $X\setminus\pi(E_{\pi})$.
By the Hartogs extension theorem, the function $\varphi^{\circ}$
extends to a quasi-psh function $\varphi$ on $X$.
Since $\psi$ and $\pi^*\varphi$ are quasi-psh functions on $M$
and coincide outside $E_{\pi}$, we have $\psi=\pi^*\varphi$ on $M$.
Thus, we obtain
\begin{equation}\label{eq:T}
T
=
\pi^*\omega_X-[G]+\ddbar\pi^*\varphi.
\end{equation}

Define a singular Hermitian metric $h_{\varepsilon}$ on
$\pi_*\mathcal{E}|_{X\setminus\pi(E_{\pi})}$ by
$h_{\varepsilon}:=\pi_*g_{\varepsilon}$.
We construct weakly positively curved metrics on $\pi_*\mathcal{E}$
by modifying $h_{\varepsilon}$.
For simplicity, we use the same notation $X$ for a sufficiently small
open subset of $X$ and $M$ for its inverse image under $\pi$.
Let $e$ be a local section of $(\pi_*\mathcal{E})^*$ and consider
$\log|e|_{h_{\varepsilon}^*}$.
By the reflexivity of $(\pi_*\mathcal{E})^*$, we have the natural
identifications
\[
\begin{aligned}
H^0\bigl(X,(\pi_*\mathcal{E})^*\bigr)
=
H^0\bigl(
X\setminus\pi(E_{\pi}),
(\pi_*\mathcal{E})^*
\bigr) 
=
H^0\bigl(
M\setminus E_{\pi},
\mathcal{E}^*
\bigr).
\end{aligned}
\]
Thus, we may regard $e$ as a local section of $\mathcal{E}^*$ on
$M\setminus E_{\pi}$.
Then, we have 
\[
\pi^*\log|e|_{h_{\varepsilon}^*}
=
\log|\pi^*e|_{\pi^*h_{\varepsilon}^*}
=
\log|\pi^*e|_{g_{\varepsilon}^*}
\quad\text{on }M\setminus E_{\pi}.
\]
Combining this with \eqref{eq:T} and
$\operatorname{Supp}(G)\subset E_{\pi}$, we obtain
\[
\begin{aligned}
\ddbar (\pi^*\log|e|_{h_{\varepsilon}^*})
=
\ddbar (\log|\pi^*e|_{g_{\varepsilon}^*} )
\geq
-\varepsilon\omega_M 
\geq
-\varepsilon T 
=
-\varepsilon
\bigl(
\pi^*\omega_X+\ddbar\pi^*\varphi
\bigr)
\end{aligned}
\]
on $M\setminus E_{\pi}$.
Therefore, we obtain 
\[
\ddbar
\bigl(
\pi^*\log|e|_{h_{\varepsilon}^*}
+\varepsilon\pi^*\varphi
\bigr)
\geq
-\varepsilon\pi^*\omega_X 
\quad\text{on }M\setminus E_{\pi}. 
\]
 Equivalently,
\[
\ddbar
\bigl(
\log|e|_{h_{\varepsilon}^*}
+\varepsilon\varphi
\bigr)
\geq
-\varepsilon\omega_X
\quad\text{on }X\setminus\pi(E_{\pi}).
\]
By the Hartogs extension theorem, the function
$\log|e|_{h_{\varepsilon}^*}+\varepsilon\varphi$ extends to an
$\varepsilon\omega_X$-psh function on $X$.
Consequently, the metric
$
H_{\varepsilon}
:=
h_{\varepsilon}e^{-\varepsilon\varphi}
$
on $\pi_*\mathcal{E}$ 
satisfies the desired curvature condition.  
\end{proof}

\begin{lemm}\label{lem-pull}
Let $\phi\colon M\to Y$ be a fibration between compact K\"ahler varieties, and assume that $\phi$ is flat over a Zariski-open subset $Y_1\subset Y$.
Let $\mathcal{E}$ be a torsion-free sheaf on $Y$.
Assume that $\mathcal{E}$ is weakly positively curved on $Y_1$
with respect to a K\"ahler form $\omega_Y$ on $Y$.
Then, the following assertions hold:
\begin{itemize}
\item[\textup{(1)}]
The reflexive pullback $\phi^{[*]}\mathcal{E}$ is weakly positively
curved on $M_1:=\phi^{-1}(Y_1)$ with respect to a K\"ahler form
$\omega_M$ on $M$.

\item[\textup{(2)}]
Let $\pi\colon M\to X$ be a bimeromorphic morphism onto a normal
compact K\"ahler variety $X$ such that
$M\setminus M_1\subset E_{\pi}$.

Then the reflexive pushforward
$\pi_{[*]}\phi^{[*]}\mathcal{E}$ is weakly positively curved on $X$.
\end{itemize}

\end{lemm}

\begin{proof}
\textup{(1)}

We explain only the strategy of the proof
(see \cite[Lemma~2.8]{MW} for details).
By assumption, there exists singular
Hermitian metrics $\{g_\varepsilon\}_{\e>0}$ on $\mathcal{E}|_{Y_1}$ such that
\[
\sqrt{-1}\Theta_{g_\varepsilon}
\geq
-\varepsilon\omega_Y\otimes\id
\quad\text{on }Y_1.
\]
Note that $g_\varepsilon$ is well-defined only on $Y_{1} \cap Y_{\mathcal{E}}\cap Y_{\reg}$. 
After replacing $\omega_M$ by a sufficiently large multiple, we may
assume that $\omega_M\geq\phi^*\omega_Y$.
On the locus where the induced pullback metric is well-defined, the
metric $\phi^*g_\varepsilon$ on $\phi^{[*]}\mathcal{E}$ satisfies
\[
\sqrt{-1}\Theta_{\phi^*g_\varepsilon}
\geq
-\varepsilon\phi^*\omega_Y\otimes\id
\geq
-\varepsilon\omega_M\otimes\id.
\]
By the extension argument in \cite[Lemma~2.8]{MW}, these metrics extend
to singular Hermitian metrics on $\phi^{[*]}\mathcal{E}|_{M_1}$ with
the same curvature lower bound.
Hence, $\phi^{[*]}\mathcal{E}$ is weakly positively curved on $M_1$.

\textup{(2)}
Since $M\setminus M_1\subset E_{\pi}$, Conclusion~\textup{(1)} shows
that $\phi^{[*]}\mathcal{E}$ is weakly positively curved on
$M\setminus E_{\pi}$.
Therefore, Conclusion~\textup{(2)} follows from
Lemma~\ref{lem-push}.
\end{proof}

The following result is straightforward; however, we state it as a lemma
since it will be used repeatedly.

\begin{lemm}\label{lem-Qfact}
Let $\pi\colon M\to X$ be a bimeromorphic morphism between strongly
$\mathbb{Q}$-factorial compact K\"ahler varieties, and let $L$ be a line
bundle on $M$.
For a sufficiently divisible integer $m\in\mathbb{Z}_{+}$, consider
the $\mathbb{Q}$-line bundle
\[
B
:=
\frac{1}{m}
\pi_{[*]}\mathcal{O}_M(mL).
\]
Then, the following assertions hold:
\begin{itemize}
\item[\textup{(1)}]
There exists a $($not necessarily effective$)$ $\pi$-exceptional
$\mathbb{Q}$-divisor $F$ such that
\[
L
\sim_{\mathbb{Q}}
\pi^*B+F.
\]

\item[\textup{(2)}]
If the sheaf $\pi_{[*]}\mathcal{O}_M(L)$ is weakly positively curved
on $X$, then $L+dE_{\pi}$ is pseudo-effective for $d\gg 1$, where
$E_{\pi}$ is the $\pi$-exceptional divisor.

\item[\textup{(3)}]
Assume further that $L$ is weakly positively curved on
$M\setminus E_{\pi}$ with respect to a K\"ahler form $\omega_M$ on $M$.
Then the $\mathbb{Q}$-line bundle $B$ is pseudo-effective on $X$.
\end{itemize}

\end{lemm}

\begin{proof}
\textup{(1)}
The natural sheaf morphisms
\[
\pi^*\pi_*\mathcal{O}_M(mL)
\longrightarrow
\mathcal{O}_M(mL)
\]
and
\[
\pi^*\pi_*\mathcal{O}_M(mL)
\longrightarrow
\pi^*\pi_{[*]}\mathcal{O}_M(mL)
=
\mathcal{O}_M(m\pi^*B)
\]
are isomorphisms on $M\setminus E_{\pi}$.
Hence their kernels and cokernels are supported in $E_{\pi}$.
By Proposition~\ref{det-tor}, their determinant sheaves are
$\QQ$-linearly equivalent to divisors supported on $E_{\pi}$.
Moreover, Proposition~\ref{det-exact} shows that taking determinants
is additive in short exact sequences.
Then, since $M$ is strongly $\QQ$-factorial, we see that
\[
L
\sim_{\QQ}
\pi^*B+F
\]
for some $\pi$-exceptional $\QQ$-divisor $F$.

We next show that $B$ is independent of the choice of a divisible integer $m$.
Let $m$ be an integer such that 
$
A_m:=\pi_{[*]}\mathcal{O}_M(mL)
$
is invertible. 
Choose a Zariski-open subset $U\subset X$ such that
$\pi^{-1}(U)\to U$ is an isomorphism and $\codim (X \setminus U) \geq 2$.
Then, for every $k\in\mathbb{Z}_{+}$, we have
$
A_{km}|_U
=
(A_m|_U)^{\otimes k}
$
and this equality extends to $X$ by reflexivity. 
Consequently, we obtain 
\[
\frac{1}{km}A_{km}
=
\frac{1}{km}A_m^{\otimes k}
=
\frac{1}{m}A_m
\]
as $\mathbb{Q}$-line bundles.
This argument shows that
$
(1/m_1)\mathcal{A}_{m_1}
=
(1/m_2)\mathcal{A}_{m_2}.
$
Thus, we see that $B$ is independent of the choice of $m$.

\textup{(2)}
By assumption, the $\mathbb{Q}$-line bundle $B$ is weakly
positively curved and hence pseudo-effective.
Choose $k\gg 1$ such that $kE_{\pi}+F$ is effective.
Then \textup{(1)} implies that $L+kE_{\pi}$ is pseudo-effective.

\textup{(3)}
Lemma~\ref{lem-push} shows that
$\pi_*\mathcal{O}_M(mL)$ is weakly positively curved on $X$.
Hence its reflexive hull
$\pi_{[*]}\mathcal{O}_M(mL)=mB$ is weakly positively curved by $\codim(X\setminus\pi(E_{\pi}))\geq 2$.
It follows that $B$ is pseudo-effective.
\end{proof}

\subsection{Foliations}
\label{subsec-foliations}
In this subsection, we collect several basic results concerning foliations.
The main reference for this subsection is \cite[\S 3]{Druel21}. Throughout
this subsection, let $X$ be a normal analytic variety. Recall that a
foliation $\mathcal{F}$ on $X$ is a saturated subsheaf of $T_X$ that is
closed under the Lie bracket. Note that it suffices to check the latter
condition on any Zariski open subset of $X$. Moreover, the normal
(resp.\ conormal) sheaf of $\calF$ is defined by
$\calN_{\calF}:=T_X/\calF$ (resp.\ by the dual sheaf
$\calN^\ast_{\calF}$ of $\calN_{\calF}$). The canonical divisor of
$\mathcal{F}$ is defined as the rank-one reflexive sheaf
$K_{\calF}:=\det(\calF)^\ast$.

\begin{defi}[projectable foliations, see {\cite[\S 3.3]{Druel21}}]
\label{def-projectable}
Let $p:X\to Y$ be a morphism between normal analytic varieties such that
$p(X)$ contains a Zariski open subset of $Y$, and let $\calF$ be a
foliation on $X$. We say that $\calF$ is projectable under $p$ if there
exists a saturated distribution $\calF_Y\subseteq T_Y$ such that the
restriction to $\calF|_{p^{-1}(Y_{\reg})}$ of the tangent map
\[
Tp|_{p^{-1}(Y_{\reg})}:
T_X|_{p^{-1}(Y_{\reg})}\to
(p|_{p^{-1}(Y_{\reg})})^\ast T_{Y_{\reg}}
\]
induces an isomorphism
\[
\calF|_{p^{-1}(Y_{\reg})}
\xrightarrow{\simeq}
\left((p|_{p^{-1}(Y_{\reg})})^\ast
\calF_Y|_{Y_{\reg}}\right)^{\text{sat}},
\]
where the right-hand side denotes the saturation of
$(p|_{p^{-1}(Y_{\reg})})^\ast\calF_Y|_{Y_{\reg}}$ in
\[
(p^\ast T_Y)|_{p^{-1}(Y_{\reg})}
\simeq
(p|_{p^{-1}(Y_{\reg})})^\ast T_{Y_{\reg}}.
\]

As shown in \cite[\S 2.7]{Druel17b}, a foliation $\calF$ of rank $r$ is
projectable if and only if, for a general $y\in Y$, $T_xp(\calF_x)$ is
independent of the choice of $x\in p^{-1}(y)$ and has dimension $r$, where
$T_xp:T_{X,x}\to T_{Y,y}$ is the tangent map.
\end{defi}

\subsubsection{Foliations and twisted reflexive differential forms}
\label{subsub-foliations-forms}
If $\mathcal{F}$ is a foliation of corank $q$ on $X$, then taking the wedge
product of the natural injection
$\calN_{\calF}^\ast\hookrightarrow\Omega_X$ yields a section
$\omega\in H^0(X,\Omega_X^q[\otimes]\det\calN_{\calF})$. Here,
$[\otimes]$ denotes the reflexive hull of a tensor product. One can check
that $\omega$ satisfies the following two conditions:
\begin{itemize}
\item $\omega$ is locally decomposable near a general point of $X$; that is,
near a general point of $X$, we can write
$\omega=\omega_1\wedge\omega_2\wedge\dots\wedge\omega_q$ for some
$1$-forms $\omega_i$.
\item $\omega$ is integrable; that is, for the local decomposition above,
we have $d\omega_i\wedge\omega=0$ for every $i$.
\end{itemize}
Conversely, if
$\omega\in H^0(X,\Omega_X^q[\otimes]\scrL)$ for some rank-one reflexive
sheaf $\scrL$ satisfies the two conditions above, then
\[
\calF:=\Ker(\lrcorner\omega:
T_X\to\Omega_X^{q-1}[\otimes]\scrL)
\]
is a foliation. See \cite[\S 3.2]{AD14}; further details can be found in
\cite[Proposition-Definition 2.4.3]{Wang-thesis}.

Using this construction, we can pull back foliations via dominant
meromorphic maps. More precisely, let $\phi:Y\dashrightarrow X$ be a
meromorphic map between normal analytic varieties such that the image of
$\phi$ contains a Zariski open subset $X^\circ$ of $X_{\reg}$. Then
$\phi$ restricts to a morphism
$\phi^\circ:Y^\circ\to X^\circ$, where $Y^\circ$ is a Zariski open subset
of $Y_{\reg}$. A foliation $\calF$ of corank $q$ on $X$ gives rise to a
twisted reflexive $q$-form
\[
\omega\in H^0(X,\Omega_X^q[\otimes]\det\calN_{\calF}).
\]
The pullback $(\phi^\circ)^\ast(\omega|_{X^\circ})$ is locally decomposable
and integrable and hence gives rise to a foliation on $Y^\circ$. By taking
the saturation in $T_Y$, this foliation extends uniquely to a foliation on
$Y$, denoted by $\phi^{-1}\calF$. Moreover, by
\cite[\S 2.5]{Druel17b}, $\phi^{-1}\calF$ can be characterized as the
unique foliation on $Y$ whose restriction to $Y^\circ$ is the saturation
of
\[
(T\phi^\circ)^{-1}
\bigl((\phi^\circ)^\ast\calF|_{X^\circ}\bigr)
\]
in $T_{Y^\circ}$, where
$T\phi^\circ:T_{Y^\circ}\to(\phi^\circ)^\ast T_{X^\circ}$ is the tangent
map.

\subsubsection{Invariant subvarieties}
\label{subsub-foliations-invariant}
Let $Y\subseteq X$ be an analytic subvariety. We say that $Y$ is invariant
under a local derivation $\partial$ on $X$ if
$\partial(\scrI_Y)\subseteq\scrI_Y$.

Let $\calF$ be a foliation on $X$. We say that $Y$ is invariant under
$\calF$ if, for every local section $\partial$ of $\calF$ over an open
subset $U$ of $X$, the subvariety $Y\cap U$ is invariant under $\partial$.
To prove that $Y$ is invariant under $\calF$, it is enough to show that
$Y\cap X^\circ$ is invariant under $\calF|_{X^\circ}$ for some Zariski
open subset $X^\circ\subseteq X$ such that $Y\cap X^\circ$ is dense in
$Y$.

\paragraph{(A)} Under some additional assumptions, we have the following
equivalent characterizations:

\paragraph{(A1)} Suppose that $X$ and $Y$ are smooth and that
$\calF\subseteq T_X$ is a subbundle. Then $Y$ is invariant under $\calF$
if and only if
$\calF|_Y\subseteq T_Y\subseteq T_X|_Y$. To see this, consider the exact
sequence
\begin{equation}
\label{seq-conormal}
0\to
\calN_{Y/X}^\ast=\scrI_Y/\scrI_Y^2
\xrightarrow{\alpha}
\Omega_X|_Y\to\Omega_Y\to 0,
\end{equation}
where the morphism $\alpha$ is given by
\[
f\text{ mod }\scrI_Y^2
\longmapsto
df\text{ mod }\scrI_Y\cdot\Omega_X
\]
(see \cite[II.1.2, p.~101]{CDGPR94}). Moreover, for every local section
$v$ of $\calF$ and every local holomorphic function $f\in\scrI_Y$, we have
$v(f)=\langle v,df\rangle$. The conclusion is then immediate.

\paragraph{(A2)} Recall that a {\itshape Pfaff field} of rank $r$ on $X$
is a nonzero morphism $\eta:\Omega_X^r\to\scrL$, where $\scrL$ is a
rank-one reflexive sheaf. The singular locus of $\eta$ is the analytic
subset of $X$ defined by the ideal sheaf
\[
\Image\bigl(
\eta[\otimes]\id_{\scrL^\ast}:
\Omega_X^r[\otimes]\scrL^\ast\to\OX_X
\bigr),
\]
and is denoted by $\Sing(\eta)$
(see \cite[Definition 3.4]{AD14} and
\cite[Definition 2.4.9]{Wang-thesis}). Note that our definition of Pfaff
fields differs slightly from that in \cite[\S 4.1]{DGP24}. In particular,
a foliation $\calF$ of rank $r$ induces a Pfaff field of rank $r$ by taking
the $r$-th wedge product of the dual of the natural inclusion
$\calF\hookrightarrow T_X$:
\[
\eta_{\calF}:
\Omega_X^r\to\bigwedge^r\calF^\ast
\to\det\calF^\ast=\OX(K_{\calF}).
\]
The singular locus of $\calF$, denoted by $\Sing(\calF)$, is defined as
$\Sing(\eta_{\calF})$, and $\calF$ is called {\itshape weakly regular} if
$\Sing(\calF)=\emptyset$.

\begin{defi}
\label{def-Pfaff_invariant}
Let $\eta:\Omega_X^r\to\scrL$ be a Pfaff field of rank $r$ on $X$, and let
$Y$ be an analytic subset of $X$. Assume that $\scrL$ is $\QQ$-Cartier;
that is, some reflexive power of $\scrL$ is invertible. We say that $Y$ is
invariant under $\eta$ if:
\begin{itemize}
\item none of the irreducible components of $Y$ is contained in
$\Sing(\eta)$;
\item for some $m\in\ZZ_{>0}$ such that $\scrL^{[\otimes]m}$ is invertible,
the restriction map
\[
(\Omega_X^r)^{\otimes m}|_Y\to\scrL^{[\otimes]m}|_Y
\]
induced by $\eta^{\otimes m}$ factors through the natural map
\[
(\Omega_X^r)^{\otimes m}|_Y\to(\Omega_Y^r)^{\otimes m}.
\]
\end{itemize}
\end{defi}

In particular, we have the following:
\begin{lemm}
\label{lemma-invariant=Pfaff}
Let $\calF$ be a foliation on $X$ such that $K_{\calF}$ is $\QQ$-Cartier,
and let $Y$ be an analytic subset of $X$ such that none of the irreducible
components of $Y$ is contained in $\Sing(\calF)$. Then $Y$ is invariant
under $\calF$ if and only if it is invariant under the Pfaff field
$\eta_{\calF}$ induced by $\calF$.
\end{lemm}

\begin{proof}
Choose $m\in\ZZ_{>0}$ such that $\OX_X(mK_{\calF})$ is invertible. By
\cite[Lemma 2.4.14]{Wang-thesis}, $Y$ is invariant under $\eta_{\calF}$ if
and only if
\[
\scrK:=
\Ker\bigl(
(\Omega_X^r)^{\otimes m}|_Y
\to(\Omega_Y^r)^{\otimes m}
\bigr)
\]
is annihilated by the morphism
\[
(\Omega_X^r)^{\otimes m}|_Y\to\OX_Y(mK_{\calF}).
\]
Moreover, by \cite[Remark 2.4.13]{Wang-thesis}, we may replace $X$ by a
small neighbourhood of a point $x\in Y_{\reg}\cap X_{\reg}$. It then
follows from \eqref{seq-conormal} that a section of $\scrK$ is of the form
\[
\big((df_1\wedge\varphi_1)\otimes\beta_1
+\alpha_2\otimes(df_2\wedge\varphi_2)\otimes\beta_2
+\cdots
+\alpha_{m-1}\otimes(df_{m-1}\wedge\varphi_{m-1})\otimes\beta_{m-1}
+\alpha_m\otimes(df_m\wedge\varphi_m)\big)\big|_Y\,,
\]
where the $f_i$ are local holomorphic functions vanishing along $Y$, the
$\varphi_i$ are local $(r-1)$-forms on $X$, the $\alpha_i$ are local
sections of $(\Omega_X^r)^{\otimes(i-1)}$, and the $\beta_i$ are local
sections of $(\Omega_X^r)^{\otimes(m-i)}$. Hence, $\scrK$ is annihilated
by
\[
(\Omega_X^r)^{\otimes m}|_Y\to\OX_Y(mK_{\calF})
\]
if and only if, for any local sections $v_1,\cdots,v_r$ of $\calF$, we have
\[
(df_i\wedge\varphi_i)(v_1,\cdots,v_r)|_Y=0.
\]
The calculation in \cite[Proof of Lemma 2.4.14]{Wang-thesis} gives
\[
(df_i\wedge\varphi_i)(v_1,\cdots,v_r)|_Y
=
\frac{1}{r}\sum_{j=1}^r(-1)^{j-1}v_j(f_i)|_Y
\cdot
\varphi_i(v_1,\cdots,\hat v_j,\cdots,v_r)|_Y.
\]
Hence, $Y$ is invariant under $\eta_{\calF}$ if and only if, for every
holomorphic function $f$ vanishing along $Y$ and every local section
$\partial$ of $\calF$, we have $\partial(f)|_Y=0$. This condition is
precisely $\partial(\scrI_Y)\subseteq\scrI_Y$, and the lemma follows.
\end{proof}

\paragraph{(B)}
We next give some examples of invariant subvarieties. A well-known theorem
of Seidenberg \cite[Theorem 5]{Sei67} states that the singular locus of a
normal complex algebraic variety is invariant under every derivation. The
same result holds for normal analytic varieties:
\begin{thm}
\label{thm-Seidenberg}
Let $X$ be a normal complex analytic variety. Then $X_{\sing}$ is invariant
under every derivation on $X$.
\end{thm}

\begin{proof}
The same calculation as in \cite[Theorem 3]{Har74} shows that every Fitting
ideal on $X$ is invariant under every derivation. Here, one should note that
the calculation can be carried out in the local ring $\OX_{X,x}$, which
can be written as $\OX_{\CC^N,0}/I$ for some ideal
$I\subseteq\OX_{\CC^N,0}=\CC\{z_1,\cdots,z_N\}$. In particular, the
Jacobian ideal defining the singular locus of $X$ is invariant under every
derivation.
\end{proof}

Moreover, the following result shows that the singular locus of a foliation
$\calF$ on a complex manifold $X$ is invariant under $\calF$.
\begin{lemm}[{\cite[Lemma 3.5]{Druel21}}]
\label{lemma-singular_invariant}
Let $X$ be a complex manifold, and let $\calF$ be a foliation of rank
$r\geq 1$ on $X$. Then every component of the singular locus of $\calF$ is
invariant under $\calF$.
\end{lemm}

\begin{proof}
The proof is essentially the same as that of
\cite[Lemma 3.5]{Druel21}. For the convenience of the reader, we include
the proof to verify that the argument also applies to possibly
non-algebraic analytic varieties. We argue by induction on the rank $r$.

Suppose first that $r=1$. Then $\calF$ is a line bundle and
$\OX_X(K_{\calF})=\calF^\ast$. In this case, the singular locus of $\calF$
is defined by the ideal sheaf
\[
\scrI:=\Image(\Omega_X\otimes\calF\to\OX_X).
\]
Let $S$ be an irreducible component of $\Sing(\calF)$. Near a general point
of $S$, the tangent bundle $T_X$ is freely generated by
$\partial_1,\cdots,\partial_n$, where $n=\dim X$, and
$\calF=\OX_X\cdot\partial$. We can therefore write
\[
\partial=\sum_{i=1}^nf_i\partial_i
\]
for some local holomorphic functions $f_i$. By construction, near a general
point of $S$, the ideal sheaf $\scrI$ is generated by
$f_1,\cdots,f_n$. Moreover,
\[
\partial(f_j)=\sum_i f_i\partial_i(f_j)\in(f_1,\cdots,f_n).
\]
Thus, $\scrI$ is invariant under $\partial$, and hence $S$ is invariant
under $\calF$.

Suppose now that $r\geq 2$. Let $S$ be an irreducible component of
$\Sing(\calF)$, and consider a general point $x\in S$. We may assume that
there exists a local section $\partial$ of $\calF$ such that
$\partial(\mathfrak{m}_x)\not\subseteq\mathfrak{m}_x$; otherwise, we are
done. This means that $\OX_X\cdot\partial$ defines a regular foliation near
$x$. By the Frobenius theorem
(see \cite[\S II.4, Theorem of Frobenius, p.~36]{CLN85}), there exist an
open neighbourhood $U$ of $x$ in $X$, a complex manifold $W$, and an
isomorphism $U\simeq\mathbb{D}\times W$ such that
$\OX_U\cdot\partial$ is induced by the projection
$\pr_2:U\simeq\mathbb{D}\times W\to W$, where $\mathbb{D}$ denotes the
open unit disc in $\CC$. The image of $\calF|_U$ under the surjective
bundle map
\[
T_U\twoheadrightarrow\pr_2^\ast T_W
\]
is of the form $\pr_2^\ast\calG$ for a saturated subsheaf

$\calG$ of
$T_W$. By the discussion in
Subsubsection~\ref{subsub-foliations-forms}, $\calG$ is a foliation of rank
$r-1$ on $W$, and
$\calF|_U=\pr_2^{-1}\calG$. 

In particular, there exists a component $T$ of $\Sing(\calG)$ such that
$S\cap U=\pr_2^{-1}(T)$. By the induction hypothesis, $T$ is invariant
under $\calG$. Hence, $S\cap U$ is invariant under $\calF|_U$, and the
lemma follows.
\end{proof}

\begin{lemm}[{\cite[Lemma 4.19]{Druel21}}]
\label{lemma-exc_invariant_rank1}
Let $q:Z\to X$ be a bimeromorphic morphism between normal analytic
varieties, and let $\calF$ be a foliation of rank $1$ on $X$ such that
$K_{\calF}$ is a line bundle. Write
\[
K_{q^{-1}\calF}\sim q^\ast K_{\calF}+B
\]
for some possibly non-effective $q$-exceptional divisor $B$ on $Z$. Let
$E$ be a prime divisor on $Z$ that is not contained in $\Supp(B)$. If
$q(E)$ is contained in a proper analytic $\calF$-invariant subvariety
$Y\subsetneqq X$, then $E$ is $q^{-1}\calF$-invariant.
\end{lemm}

\paragraph{(D)}
Finally, we have the following observation, which follows from a direct
computation:
\begin{lemm}[{\cite[Lemma 3.4]{Druel21}}]
\label{lemma-inv_divisor}
Let $f:Y\to X$ be a dominant quasi-finite morphism between normal analytic
varieties, and let $\calF$ be a foliation on $X$ of rank
$1\leq r\leq\dim X-1$. Let $B\subseteq X$ be a divisorial component of the
branch locus of $f$.
\begin{itemize}
\item Suppose that $B$ is $\calF$-invariant. Let $D$ be an irreducible
component of $f^{-1}(B)$, and let $m$ denote the ramification index of $f$
along $D$. Then the natural map
\[
\det f^{[\ast]}\calN_{\calF}^\ast
\to
\det\calN_{f^{-1}\calF}^\ast
\]
vanishes to order $m-1$ along $D$.
\item If $B$ is not $\calF$-invariant, then the map
\[
\det f^{[\ast]}\calN_{\calF}^\ast
\to
\det\calN_{f^{-1}\calF}^\ast
\]
is an isomorphism at a general point of every component of $f^{-1}(B)$.
\end{itemize}
In particular, we have
\[
K_{f^{-1}\calF}
\sim
f^\ast K_{\calF}
+\Ramification(f)^{{\rm non}-\calF-{\rm inv}},
\]
where $\Ramification(f)^{{\rm non}-\calF-{\rm inv}}$ denotes the
non-$\calF$-invariant part of the ramification divisor $\Ramification(f)$
of $f$. In other words, $\calF$-invariant divisors do not contribute to the
canonical class of $f^{-1}\calF$.
\end{lemm}

\subsubsection{Bertini-type theorem for foliations}
\label{subsub-foliations-Bertini}
We now consider the restriction of a foliation to a hypersurface. In
Subsubsection~\ref{subsub-foliations-forms}, we defined the pullback of a
foliation $\calF$ via dominant maps. We define the restriction of $\calF$
to a hypersurface $H$ in the same way, namely, by pulling back the twisted
reflexive differential form $\omega$ induced by $\calF$. However, some
care is required: we need to ensure that $\omega$ does not vanish
identically along $H$. This amounts to checking that $\calF$ is generically
transverse to $H$, or equivalently, that $\calF$ is not tangent to $H$ at
a general point of $H$. More precisely, we have the following observation:
\begin{lemm}
\label{lemma-restriction_foliation}
Let $H\subseteq X$ be a prime divisor on a normal analytic variety $X$,
and assume that $H$ is itself normal. Let $\calF\neq\{0\}$ be a foliation
of corank $q\leq\dim X-2$ on $X$ that induces a twisted reflexive $q$-form
\[
\omega\in
H^0(X,\Omega_X^q[\otimes]\det\calN_{\calF}).
\]
Let $X^\circ\subseteq X_{\reg}$ be the Zariski open subset where
$\calF|_{X_{\reg}}$ is a subbundle of $T_{X_{\reg}}$, and assume that
$H\backslash(H\cap X^\circ)$ has codimension $\geq 2$ in $H$. Let
\[
\omega_H\in
H^0(H,\Omega_H^q[\otimes]\det\calN_{\calF}|_H)
\]
be the twisted reflexive $q$-form induced by $\omega$ on $H$.

Let $B$ be a subvariety of $H$ such that
$B\cap X^\circ\neq \emptyset$. Then $\omega_H$ vanishes along $B$ if and
only if $\calF$ is tangent to $H$ at a general point $x\in B$, that is,
$\calF_x\subseteq T_{H,x}$. In particular, $\omega_H\not\equiv 0$ if and
only if $\calF$ is generically transverse to $H$.
\end{lemm}

\begin{proof}
By the discussion in Subsubsection~\ref{subsub-foliations-forms}
(see also \cite[Proposition-Definition 2.4.3]{Wang-thesis}), near a general
point $x\in B\cap X^\circ$, we have a local decomposition
$\omega=\omega_1\wedge\cdots\wedge\omega_q$ for some local holomorphic
$1$-forms $\omega_1,\cdots,\omega_q$, and
\[
\calF\overset{\mathrm{loc}}{=}
\left\{
v\text{ a local holomorphic vector field }
\mathrel{\big|}
\omega_i(v)=0\text{ for every }i
\right\}.
\]
Then
\[
(\omega_H)_x=\omega_x|_{T_{H,x}}=0
\]
if and only if
\[
\calF_x=\bigcap_i\Ker((\omega_i)_x)\subseteq T_{H,x},
\]
by the following elementary lemma from linear algebra.
\end{proof}

\begin{lemm}
\label{lemma-lin_algebra2}
Let $V$ be an $n$-dimensional vector space, and let
$0\neq\omega\in\bigwedge^qV^\ast$ be an alternating $q$-form on $V$ such
that
$\omega=\omega_1\wedge\cdots\wedge\omega_q$ for some
$\omega_i\in V^\ast$. Set
$F:=\bigcap_i\Ker(\omega_i)\subseteq V$. Then, for every $m$-dimensional
subspace $W\subseteq V$ with $m\geq q$, the following conditions are
equivalent:
\begin{itemize}
\item[\fbox{a}] $\omega|_W=0$;
\item[\fbox{b}] every $q$-dimensional subspace of $W$ has nonzero
intersection with $F$;
\item[\fbox{c}] $\dim(W\cap F)\geq m-q+1$.
\end{itemize}
In particular, if $m=n-1$, these conditions are equivalent to
\begin{itemize}
\item[\fbox{d}] $F\subseteq W$.
\end{itemize}
\end{lemm}

Lemma~\ref{lemma-lin_algebra2} follows immediately from the following
elementary observation:
\begin{lemm}
\label{lemma-lin_algebra1}
Let $V$, $\omega$, and $F$ be as above. Then, for every $q$-dimensional
subspace $W\subseteq V$, the following conditions are equivalent:
\begin{itemize}
\item[\fbox{a}] $\omega|_W=0$;
\item[\fbox{b}] $W\cap F\neq\{0\}$.
\end{itemize}
\end{lemm}

We can now state the Bertini-type theorem for foliations:
\begin{thm}
\label{thm-Bertini-foliation}
Let $X$ be a normal analytic variety, and let $\calF\neq\{0\}$ be a
foliation on $X$ of corank $q\leq\dim X-2$. Let $L$ be a line bundle on
$X$, and let $V\subseteq H^0(X,L)$ be a finite-dimensional subspace that
generates $L$ at every point of $X$; in other words, $V$ is a
basepoint-free linear system on $X$. Suppose that there exists
$0\neq v_0\in V$ such that $\calF$ is generically transverse to a
component of $H_{v_0}$. Then, for an analytically sufficiently general
$v\in V$ and every component $H$ of $H_v$, the space $H$ is a normal
analytic variety, $H\cap X_{\reg}$ is smooth, and
\[
\calF_H:=\calF|_H\cap T_H
\]
is a foliation of corank $q$ on $H$.
\begin{itemize}
\item[\textup{(a)}] There exists an effective divisor $\Delta_H$ on $H$
such that
\[
K_{\calF_H}\sim(K_{\calF}+H)|_H-\Delta_H.
\]
Here, $K_{\calF}|_H$ is the rank-one reflexive sheaf on $H$ that uniquely
extends the line bundle $K_{\calF}|_{H\cap X_{\reg}}$ on
$H\cap X_{\reg}\subseteq H_{\reg}$. Moreover, at a general point $x$ of
every component of $\Delta_H$, the foliation $\calF$ is tangent to $H$.

\item[\textup{(b)}] Suppose that there exists a Zariski open subset
$X_0\subseteq X_{\reg}$ such that $\codim(X\backslash X_0)\geq 2$ and, for
every $x\in X_0$, $\calF$ is a subbundle of $T_{X_{\reg}}$ at $x$ and
there exists a $2$-dimensional subspace $W_x\subseteq V_x$ such that
$\calF$ is not tangent to $H_w$ at $x$ for every
$0\neq w\in W_x$. Here, $\scrH_V\subseteq V\times X$ denotes the universal
family over $V$, with natural maps
$\pi:\scrH_V\to V$ and $e:\scrH_V\to X$, and
$V_x:=e^{-1}(x)$ is a hyperplane in
$V\simeq V\times\{x\}$ (see the related discussion in
\S\,\ref{subsec-Bertini}). Then
\[
K_{\calF_H}\sim(K_{\calF}+H)|_H.
\]

\item[\textup{(c)}] Suppose that $X$ is a locally analytic subset of
$\CC^{N-1}$ and that $V\subseteq H^0(X,\OX_X)$ is the linear system
consisting of the restrictions to $X$ of affine-linear functions on
$\CC^{N-1}$. Thus, the divisor defined by every nonconstant element of $V$
is the intersection of $X$ with an affine hyperplane in $\CC^{N-1}$. Then
$V$ is a basepoint-free linear system satisfying the condition in
\textup{(b)}.
\end{itemize}
\end{thm}

\begin{proof}
We first introduce some notation:
\begin{itemize}
\item Let $X^\circ\subseteq X_{\reg}$ be the locus where $\calF$ is a
subbundle of $T_{X_{\reg}}$. Then $X^\circ$ is a Zariski open subset of
$X$, and $\codim(X\backslash X^\circ)\geq 2$.

\item Recall that $\scrH_V\subseteq V\times X$ is the universal family,
with natural morphisms
\[
\pi:=\pr_1|_{\scrH_V}:\scrH_V\to V
\quad\text{and}\quad
e:=\pr_2|_{\scrH_V}:\scrH_V\to X.
\]
Since $V$ is basepoint-free, $e$ is a locally trivial
$\CC^{N-1}$-bundle over $X$ by the discussion in
\S\,\ref{subsec-Bertini}, where $N:=\dim V$.

\item Set $\scrH_V^\circ:=e^{-1}(X^\circ)$. Then $\scrH_V^\circ$ is a
Zariski open subset of $\scrH_V$ whose complement has codimension
$\geq 2$. Consider the subset
$\scrH_V^{\text{tan},\calF}\subseteq\scrH_V^\circ$ defined by
\[
\scrH_V^{\text{tan},\calF}
:=
\left\{
(x,v)\in\scrH_V^\circ
\mathrel{\big|}
\calF\text{ is tangent to }H_v\text{ at }x
\right\}.
\]
The subset $\scrH_V^{\text{tan},\calF}$ is defined by the condition
$dv_x(\calF_x)=0$ and is therefore an analytic subset of
$\scrH_V^\circ$. Moreover, since this condition is linear when $x$ is
fixed,
\[
V_x^{\text{tan},\calF}
:=
V_x\cap\scrH_V^{\text{tan},\calF}
\]
is a linear subspace of $V_x\subseteq V$.

\item Let
\[
\omega\in
H^0(X,\Omega_X^q[\otimes]\det\calN_{\calF})
\]
be the twisted reflexive $q$-form induced by $\calF$. By
Subsubsection~\ref{subsub-foliations-forms}, near every
$x\in X^\circ$, we can write
$\omega=\omega_1\wedge\cdots\wedge\omega_q$ for some local holomorphic
$1$-forms such that
\[
\calF\overset{\mathrm{loc}}{=}
\left\{
v\text{ a local holomorphic vector field }
\mathrel{\big|}
\omega_i(v)=0\text{ for every }i
\right\}.
\]
\end{itemize}

First, Theorem~\ref{thm-Bertini} implies that, for an analytically
sufficiently general $v\in V$, $H_v$ is normal and
$H_v\cap X_{\reg}$ is smooth. Moreover,
$H_v\backslash(H_v\cap X^\circ)$ has codimension $\geq 2$ in $H_v$.
Thus, for every component $H$ of $H_v$, $\omega$ induces a twisted
reflexive $q$-form
\[
\omega_H\in
H^0(H,\Omega_H^q[\otimes]\det\calN_{\calF}|_H).
\]
By assumption,
$\scrH_V^{\text{tan},\calF}\neq\scrH_V^\circ$. Hence, every component of
$\scrH_V^{\text{tan},\calF}$ has dimension at most
\[
\dim\scrH_V-1=\dim X+N-2.
\]
Let $S$ be a component of $\scrH_V^{\text{tan},\calF}$. By Remmert's
mapping theorem, Theorem~\ref{thm-Remmert-mapping}, $\pi(S)$ is a
countable union of locally analytic subsets of $V$.
\begin{itemize}
\item If $\pi(S)$ has dimension $\leq N-1$, then $\pi(S)$ is analytically
meagre. Thus, for an analytically sufficiently general $v\in V$,
$\pi^{-1}(v)\cap S=\emptyset$.

\item Suppose that $\pi(S)$ has dimension $N$, and let $U$ be the union of
the $N$-dimensional components of $\pi(S)$. Then $U$ is an open subset of
$V$, and $\pi(S)\backslash U$ is analytically meagre. Applying Frisch's
generic flatness theorem, Theorem~\ref{thm-generic-flatness}, to
$\pi|_{\pi^{-1}(U)\cap S}$ and using the dimension formula for flat
morphisms \cite[\S V.2, Proposition V.2.3, p.~178]{BS76}, we see that, for
an analytically sufficiently general $v\in V$, either
$\pi^{-1}(v)\cap S=\emptyset$ or
\[
\dim(\pi^{-1}(v)\cap S)
\leq
\dim X-2
<
\dim H_v.
\]
\end{itemize}
It follows that, for an analytically sufficiently general $v\in V$,
$\calF$ is not generically tangent to any component $H$ of $H_v$. Hence,
by Lemma~\ref{lemma-restriction_foliation}, $\omega_H\not\equiv 0$ and
therefore defines a foliation $\calG$ on $H$.

On the other hand, consider
$\calF_H:=\calF|_H\cap T_H$. We verify the following:
\begin{itemize}
\item[(1)] The sheaf $\calF_H$ is saturated in $T_H$. To see this, consider
the following commutative diagram:
\begin{center}
\begin{tikzpicture}[scale=1.2]
\node (A0) at (-3.8,0) {$0$};
\node (A1) at (-2,0) {$\calF|_H$};
\node (A2) at (0,0) {$T_X|_H$};
\node (A3) at (2.5,0) {$(T_X|_H)/(\calF|_H)$};
\node (A4) at (4.5,0) {$0$,};
\node (B0) at (-3.8,2) {$0$};
\node (B1) at (-2,2) {$\calF_H$};
\node (B2) at (0,2) {$T_H$};
\node (B3) at (2.5,2) {$T_H/\calF_H$};
\node (B4) at (4.5,2) {$0$};
\path[->,font=\scriptsize,>=angle 90]
(A0) edge (A1)
(A1) edge (A2)
(A2) edge (A3)
(A3) edge (A4)
(B0) edge (B1)
(B1) edge (B2)
(B2) edge (B3)
(B3) edge (B4)
(B1) edge (A1)
(B2) edge (A2)
(B3) edge node[right]{$\alpha$} (A3);
\end{tikzpicture}
\end{center}
Here, the rows are exact, while the left and central columns are natural
injections. The snake lemma yields an exact sequence
\[
0\to\Ker(\alpha)\to
(\calF|_H)/\calF_H
\xrightarrow{\beta}
(T_X|_H)/T_H.
\]
Since $\calF_H=\calF|_H\cap T_H$, the morphism $\beta$ is injective.
Therefore, $\Ker(\alpha)=0$, and we obtain an injection
\[
T_H/\calF_H
\hookrightarrow
(T_X|_H)/(\calF|_H).
\]
The latter sheaf is torsion free by
Corollary~\ref{cor-Bertini}. Hence, $T_H/\calF_H$ is also torsion free,
and therefore $\calF_H$ is saturated in $T_H$.

\item[(2)] Near a point $x\in X^\circ$, the sheaf $\calF_H$ is defined by
$\omega_{i,H}(v)=0$, and hence $\calF_H=\calG$, where $\omega_{i,H}$ is
the local holomorphic $1$-form on $H$ induced by $\omega_i$. Indeed, if we
regard $(\omega_i)_x$ as a linear form on $T_{X,x}$, then
$(\omega_{i,H})_x$ is its restriction to
$T_{H,x}\subseteq T_{X,x}$. Thus,
\[
(\calF_H)_x
=
T_{H,x}\cap\calF_x
=
\bigcap_{i=1}^q\Ker(\omega_{i,H})_x,
\]
and the assertion follows.
\end{itemize}

We now prove \textup{(a)}. We have seen that $\calF_H$ is defined by the
twisted reflexive $q$-form
\[
\omega_H\in
H^0(H,\Omega_H^q[\otimes]\det\calN_{\calF}|_H).
\]
By Lemma~\ref{lemma-restriction_foliation}, if $\calF$ is tangent to $H$ at
a general point of a prime divisor on $H$, then $\omega_H$ vanishes along
that divisor. Hence, there exists an effective divisor $\Delta_H\geq 0$
such that $\omega_H$, viewed as a section of
\[
H^0\bigl(
H,\Omega_H^q[\otimes]\det\calN_{\calF}|_H
[\otimes]\OX_H(-\Delta_H)
\bigr),
\]
has a vanishing locus of codimension $\geq 2$. Moreover, $\calF$ is tangent
to $H$ at a general point of every component of $\Delta_H$. By reflexivity,
we then have
\[
\det\calN_{\calF_H}
\sim
\det\calN_{\calF}|_H-\Delta_H,
\]
since the left- and right-hand sides coincide outside an analytic subset
of codimension $\geq 2$. A direct computation, together with the adjunction
formula for $H$, now gives
\[
K_{\calF_H}
\sim
(K_{\calF}+H)|_H-\Delta_H,
\]
where $K_{\calF}|_H$ is the rank-one reflexive sheaf that uniquely extends
the line bundle $K_{\calF}|_{X_{\reg}\cap H}$.

We next prove \textup{(b)}. After shrinking $X^\circ$, we may assume that
$X^\circ=X_0$. For every $x\in X^\circ$, we have
\[
V_x^{\text{tan},\calF}\cap W_x=\{0\}.
\]
Hence, by Frisch's generic flatness theorem,
Theorem~\ref{thm-generic-flatness}, and the dimension formula
\cite[\S V.2, Proposition V.2.3, p.~178]{BS76}, we obtain
\[
\dim\scrH_V^{\text{tan},\calF}
\leq
\dim X+\dim V_x^{\text{tan},\calF}
\leq
\dim X+N-3.
\]
Let $S$ be a component of $\scrH_V^{\text{tan},\calF}$. By Remmert's
mapping theorem, Theorem~\ref{thm-Remmert-mapping}, $\pi(S)$ is a
countable union of locally analytic subsets of $V$.
\begin{itemize}
\item If $\pi(S)$ has dimension $\leq N-1$, then $\pi(S)$ is analytically
meagre. Thus, for an analytically sufficiently general $v\in V$,
$\pi^{-1}(v)\cap S=\emptyset$.

\item Suppose that $\pi(S)$ has dimension $N$, and let $U$ be the union of
the $N$-dimensional components of $\pi(S)$. Then $U$ is an open subset of
$V$, and $\pi(S)\backslash U$ is analytically meagre. Applying Frisch's
generic flatness theorem, Theorem~\ref{thm-generic-flatness}, to
$\pi|_{\pi^{-1}(U)\cap S}$ and using the dimension formula for flat
morphisms \cite[\S V.2, Proposition V.2.3, p.~178]{BS76}, we see that, for
an analytically sufficiently general $v\in V$, either
$\pi^{-1}(v)\cap S=\emptyset$ or
\[
\dim(\pi^{-1}(v)\cap S)
=
\dim S-N
\leq
\dim X-3.
\]
\end{itemize}
The discussion above and \textup{(a)} imply that $\Delta_H=0$. This proves
\textup{(b)}.

Finally, we prove \textup{(c)}. Let
$x\in X\subseteq\CC^{N-1}$. There exists an affine hyperplane in
$\CC^{N-1}$ that does not contain $x$, and hence
$V\subseteq H^0(X,\OX_X)$ is basepoint-free. Moreover, for
$x\in X^\circ$, the affine-linear functions vanishing at $x$ whose
differentials vanish on
$\calF_x\subseteq T_{X,x}\subseteq T_{\CC^{N-1},x}\simeq\CC^{N-1}$
form a subspace $E_x$ of
$V_x\simeq(\CC^{N-1})^\ast$ of dimension
\[
\dim E_x=N-1-\rank\calF.
\]
Since $\rank\calF\geq 2$, there exists a $2$-dimensional subspace
$W_x\subseteq V_x$ such that $W_x\cap E_x=\{0\}$. It follows that $W_x$
satisfies the condition in \textup{(b)}.
\end{proof}

\subsubsection{Singularity and regularity}\label{subsub-sing}
In Subsubsection~\ref{subsub-foliations-invariant}\textup{(A2)}, we defined
the notion of weak regularity for foliations: a foliation $\calF$ on $X$ is
said to be weakly regular if $\Sing(\calF)=\emptyset$. As pointed out in
\cite[Remark 3.12]{AD13}, this is the appropriate notion of regularity for
foliations on singular varieties.

We next recall the notion of canonical singularities for foliations:
\begin{defi}[{\cite[Definition 4.1]{Druel21}}]
\label{def-canonical}
Let $\calF$ be a foliation on a normal analytic variety $X$ such that
$K_{\calF}$ is $\QQ$-Cartier. Let $\beta:Z\to X$ be a projective
modification. Then there are uniquely defined rational numbers
$a(E,X,\calF)$ such that
\[
K_{\beta^{-1}\calF}
\sim_{\QQ}
\beta^\ast K_{\calF}
+\sum_Ea(E,X,\calF)E,
\]
where $E$ runs through all exceptional prime divisors of $\beta$. The
rational number $a(E,X,\calF)$ depends only on the valuation associated
with $E$ and not on the modification $\beta$. We say that $\calF$ has
{\itshape canonical singularities}, or simply that $\calF$ is
{\itshape canonical}, if $a(E,X,\calF)\geq 0$ for every $E$ exceptional
over $X$.

\end{defi}

An example similar to that in \cite[Remark 3.12]{AD13} is given in
\cite[Remark 4.16]{Druel21}, illustrating that the singularities of a
foliation do not control those of the ambient variety except along the
leaves.

In the sequel, we collect some results concerning the behaviour of weakly
regular foliations and foliations with canonical singularities under finite
covers and bimeromorphic morphisms.

\begin{lemm}[{\cite[Lemma 4.3 \& Lemma 5.13]{Druel21}}]
\label{lemma-foliation_finite}
Let $f:X_1\to X$ be a quasi-finite dominant morphism between normal
analytic varieties, and let $\calF$ be a foliation on $X$. Assume that every
divisorial component of the branch locus of $f$ is $\calF$-invariant. This
condition holds automatically if $f$ is quasi-\'etale. Then
\[
K_{f^{-1}\calF}\sim f^\ast K_{\calF}.
\]
Moreover, the following assertions hold:
\begin{itemize}
\item[\textup{(a)}] If $\calF$ is canonical, then so is
$f^{-1}\calF$. The converse fails in general, even if $f$ is a finite
quasi-\'etale cover, as illustrated in \cite[Example 4.4]{Druel21}.
\item[\textup{(b)}] If $\calF$ is weakly regular, then so is
$f^{-1}\calF$. The converse holds if $f$ is a finite cover.
\end{itemize}
\end{lemm}

\begin{lemm}[slightly generalized version of {\cite[Lemma 4.2 \& Lemma 5.15]{Druel21}}]
\label{lemma-foliation_bimero}
Let $\beta:Z\to X$ be a proper modification between normal analytic
varieties, and let $\calF$ be a foliation on $X$ such that $K_{\calF}$ is
$\QQ$-Cartier. Assume that
\[
K_{\beta^{-1}\calF}
\sim_{\QQ}
\beta^\ast K_{\calF}+E
\]
for some effective $\beta$-exceptional divisor $E$ on $Z$.
\begin{itemize}
\item[\textup{(a)}] If $\beta^{-1}\calF$ is canonical, then so is
$\calF$. The converse holds if $E=0$.
\item[\textup{(b)}] Assume that $X$ has klt singularities and that $E=0$.
If $\calF$ is weakly regular, then so is $\beta^{-1}\calF$.
\end{itemize}
\end{lemm}

\begin{proof}
Let $m\in\ZZ_{>0}$ be an integer such that $mK_{\calF}$ is a line bundle.
Assertion \textup{(a)} is proved in \cite[Lemma 4.2]{Druel21}. When $m=1$,
assertion \textup{(b)} is proved in \cite[Lemma 5.15]{Druel21}. We explain
how to deduce \textup{(b)} from the case $m=1$, using an argument similar
to that of \cite[Lemma 5.16(3)]{Druel21}.

Since the problem is local on $X$, we may assume that $X$ is a small
neighbourhood of a point $x\in X$ such that
$\OX_X(mK_{\calF})\simeq\OX_X$. By assumption,
\[
\OX_Z(mK_{\beta^{-1}\calF})
\simeq
\beta^\ast\OX_X(mK_{\calF})
\simeq
\OX_Z.
\]
We can therefore take the associated cyclic cover $f:Z_1\to Z$, which is
quasi-\'etale. Let $\beta_1:Z_1\to X_1$ be the Stein factorization of
$\beta\circ f$. Then $X_1$ is normal, and we obtain a finite cover
$g:X_1\to X$. These maps are summarized in the following commutative
diagram:
\begin{center}
\begin{tikzpicture}[scale=2.0]
\node (A) at (0,0) {$Z$};
\node (A1) at (0,1) {$Z_1$};
\node (B) at (1,0) {$X$};
\node (B1) at (1,1) {$X_1$};
\path[->,font=\scriptsize,>=angle 90]
(A1) edge node[left]{$f$} (A)
(B1) edge node[right]{$g$} (B)
(A1) edge node[above]{$\beta_1$} (B1)
(A) edge node[below]{$\beta$} (B);
\end{tikzpicture}
\end{center}
Since $f$ is quasi-\'etale, so is $g$. Moreover, since $X$ has klt
singularities and $g$ is finite and quasi-\'etale, $X_1$ also has klt
singularities. By Lemma~\ref{lemma-foliation_finite}, we have
\[
K_{g^{-1}\calF}\sim g^\ast K_{\calF},
\]
and $g^{-1}\calF$ is weakly regular. Furthermore,
\[
K_{f^{-1}(\beta^{-1}\calF)}
\sim
f^\ast K_{\beta^{-1}\calF}
\sim
f^\ast\beta^\ast K_{\calF}
=
\beta_1^\ast g^\ast K_{\calF}
\sim
\beta_1^\ast K_{g^{-1}\calF}.
\]
By construction,
\[
\OX_{Z_1}\bigl(K_{f^{-1}(\beta^{-1}\calF)}\bigr)
\simeq
\OX_{Z_1},
\]
and hence $K_{g^{-1}\calF}$ is also Cartier. We can therefore apply
\cite[Lemma 5.15]{Druel21} to $g^{-1}\calF$ and $\beta_1$ to conclude that
\[
\beta_1^{-1}(g^{-1}\calF)
=
f^{-1}(\beta^{-1}\calF)
\]
is weakly regular. Since $f$ is a finite cover,
Lemma~\ref{lemma-foliation_finite}\textup{(b)} implies that
$\beta^{-1}\calF$ is weakly regular.
\end{proof}

Weak regularity and canonicity are related by the following fact: a weakly
regular foliation on a klt variety is canonical provided that its canonical
divisor is Cartier.
\begin{lemm}[{\cite[Lemma 5.9]{Druel21}}]
\label{lemma-foliation_Gorenstein_wr_canonical}
Let $X$ be a normal analytic variety with klt singularities, and let
$\calF$ be a foliation on $X$. Assume that $K_{\calF}$ is a line bundle.
If $\calF$ is weakly regular, then it has canonical singularities.
\end{lemm}

Finally, we recall the following result of Druel
\cite[Lemma 5.4]{Druel17a}, which controls the singularities of the ambient
variety in terms of those of an algebraically integrable foliation.
\begin{lemm}[{simplified version of \cite[Lemma 5.4]{Druel17a}}]
\label{lemma-canonical_foliation_variety}
Let $p:Z\to Y$ be a dominant equidimensional morphism between normal
analytic varieties, and let $\calH$ denote the foliation induced by $p$.
Suppose that $Y$ is smooth and that $\calH$ is canonical. Then $Z$ has
canonical singularities.
\end{lemm}

\subsection{Miscellaneous}
To conclude this preliminary section, we collect some miscellaneous
results.

\paragraph{(A) Hartshorne's connectedness theorem}
The following is the analytic version of Hartshorne's connectedness theorem
\cite[Theorem 2.2]{Har62}:

\begin{thm}[{\cite[\S II.3, Corollary II.3.11, p.~68]{BS76}}]
\label{thm-connectedness}
Let $X$ be a connected analytic variety, and let $Y$ be an analytic subset of
$X$. Assume that $\depth\OX_{X,x}\geq 2$ for every $x\in Y$. Then
$X\backslash Y$ is connected.
\end{thm}

\paragraph{(B) Functorial pullback for reflexive differentials}
\begin{thm}[{\cite[Theorem 10.1]{KS21}}]
\label{thm-pullback-reflexive}
Let $f:X\to Y$ be a morphism between normal analytic varieties with
rational singularities. For every $p\in\ZZ_{\geq 0}$, write
$\Omega_X^{[p]}=(\Omega_X^p)^{\ast\ast}$, and similarly for
$\Omega_Y^{[p]}$. Then there exists a pullback morphism
\[
d_{\text{refl}}f:\Omega_Y^{[p]}\to\Omega_X^{[p]}
\]
uniquely determined by its natural universal properties. In particular,
$d_{\text{refl}}f$ is compatible with the usual pullback of K\"ahler
differentials
$df:\Omega_Y^p\to\Omega_X^p$.
\end{thm}

\paragraph{(C) Local structure of klt singularities}
\begin{prop}[{\cite[Lemma 5.8]{GK20}}]
\label{prop-local_structure_klt_codim2}
Let $X$ be an analytic variety with klt singularities. Then there exists an
analytic subset $Z\subseteq X$ of codimension $\geq 3$ such that, for every
$x\in X\backslash Z$, there exists a klt surface singularity $(S,o)$ with
\[
(X,x)
\simeq
(S,o)\times(\CC^{\dim X-2},0)
\]
as germs of analytic varieties. In particular, $X\backslash Z$ has quotient
singularities.

\end{prop}

\section{Flatness of certain direct image sheaves}\label{sec-main}

The purpose of this section is to prove Theorem~\ref{thm-faked-flat} and Corollary~\ref{thm-flat}.
To this end, we construct a suitable relatively big line bundle for an MRC fibration in Subsection~\ref{subsec-setup},
study certain bimeromorphic properties of MRC fibrations in Subsection~\ref{subsec-mrc},
and finally prove both results in Subsection~\ref{subsec-direct}.
As explained in Subsection~\ref{subsec-strategy},
our approach is inspired by the previous work \cite{CH19,CCM19}.
However, the methods developed there rely heavily on projectivity and therefore cannot be applied directly
in the analytic setting.
To overcome this difficulty, we develop new arguments adapted to the compact K\"ahler setting,
particularly for constructing relatively big line bundles and localizing positivity arguments
for direct image sheaves associated with MRC fibrations.

\subsection{MRC fibrations and relatively big line bundles}\label{subsec-setup}

We begin by fixing the setting and notation used throughout this section.

\begin{setup}\label{setting}
Let $(X,\Delta,\boldsymbol{\beta})$ be a K\"ahler generalized klt pair of Calabi--Yau type
in the sense of Definition~\ref{def-gpairs}, where $X$ is a compact K\"ahler variety of dimension $n$.
Throughout this section, we further assume that $X$ is strongly $\QQ$-factorial.

Let $\psi\colon X\dashrightarrow Y$ be an MRC fibration onto a compact K\"ahler manifold $Y$.
Let $\pi\colon M\to X$ be a log resolution of the generalized pair
$(X,\Delta,\boldsymbol{\beta})$ such that the induced rational map
$M\dashrightarrow Y$ is resolved by a morphism $\phi\colon M\to Y$, as in the following diagram:
\begin{align}\label{comm-start}
\xymatrix{
M \ar[rr]^{\pi} \ar[rd]_{\phi} && X \ar@{.>}[ld]^{\psi} \\
& Y. &
}
\end{align}
After replacing $M$ and $Y$ with higher models if necessary, we assume that the discriminant locus
$\Sigma$ of $\phi \colon M \to Y$ is an SNC divisor on $Y$ and that
its pullback $\phi^{*}\Sigma$ has SNC support.
Let
\[
E := \Exc(\pi)
\]
be the reduced $\pi$-exceptional divisor on $M$.

Since $(X,\Delta,\boldsymbol{\beta})$ is of Calabi--Yau type,
the formula~\eqref{eq-gpair} of canonical divisors yields
\begin{align}
\label{eq-gpair2}
-(K_M+\Delta_M)+F_+
\equiv
\boldsymbol{\beta}_M.
\end{align}
In particular, the left-hand side represents a nef Bott--Chern class.
We emphasize that, even when $\Delta$ is a $\mathbb{Q}$-divisor,
the effective divisors $\Delta_M$ and $F_+$ are $\mathbb{R}$-divisors
and need not be $\mathbb{Q}$-divisors.
This is one of the distinctions between our setting and those considered in previous works.
\end{setup}

For later use, we choose a relatively $\phi$-big line bundle $G$ on $M$
satisfying the properties established in the following lemma.
When $X$ is projective, it suffices to take $G$ to be the pullback of an ample line bundle on $X$.
However, in the K\"ahler setting, the construction of $G$ is considerably more delicate.
In particular, properties~\textup{(1)} and~\textup{(2)} below, which will play a crucial role
in the proof of Theorem~\ref{thm_splitting}, require special care.

Henceforth, let $m\in\mathbb{Z}_{+}$ denote a sufficiently large and divisible integer.

\begin{lemm}\label{lemm-good}
Under the assumptions of Setting~\ref{setting}, 
there exists a relatively $\phi$-big line bundle $G$ on $M$ satisfying the following properties:
\begin{itemize}
\item[$(1)$] There exist a line bundle $A$ on $X$ and an effective $\pi$-exceptional divisor $F$ on $M$
such that
\[
G := \pi^{*}A + F,
\]
is relatively $\phi$-big 
and the restriction of $A$ to a general fiber of $\psi$ is ample.
Note that this restriction is well-defined since $\psi$ is almost holomorphic.

\item[$(2)$] Let $S$ be an irreducible component of a fiber of $\phi$
such that $S$ is not contained in $E$. Then, the restriction $G|_{S}$ is big.

\item[$(3)$] For every positive integer $c\in\mathbb{Z}_{+}$ and a general point $y\in Y$,
the natural multication morphism 
\[
\Sym^{m} H^{0}\bigl(M_{y}, \OX_{M_y}(G+cE)\bigr)
\longrightarrow
H^{0}\bigl(M_{y}, \OX_{M_y}(m(G+cE))\bigr)
\]
is surjective, where $M_{y}:=\phi^{-1}(y)$ denotes the fiber of
$\phi \colon M \to Y$ at $y$.
\end{itemize}
\end{lemm}

\begin{proof}
We would like to apply \cite[Corollary 4.2]{CH24} to the morphism $M \to Y$
in order to construct a relatively ample line bundle with the desired properties.
However, the $\pi$-exceptional locus $E$ dominates $Y$, 
preventing the resulting line bundle from having the desired properties. 
Therefore, we consider the normalization $\Gamma'$ of the graph of $\psi$,
with the natural morphisms $\varphi' \colon \Gamma' \to Y$ and
$\beta' \colon \Gamma' \to X$.
An advantage of considering $\Gamma'$ is that the $\beta'$-exceptional locus does not dominate $Y$.
On the other hand, the variety $\Gamma'$ need not be strongly $\mathbb{Q}$-factorial,
and thus, \cite[Corollary 4.2]{CH24} cannot be applied directly to $\Gamma'$.
To overcome this difficulty, we take a strong $\QQ$-factorialization
$\mu\colon \Gamma\to \Gamma'$ of $\Gamma'$,
as constructed in Theorem~\ref{thm-Q-factorialization},
and set $\beta := \beta' \circ \mu$ and $\varphi := \varphi' \circ \mu$.
After replacing $M$ with a higher model if necessary,
we may assume that $\pi \colon M \to X$ factors through $\Gamma$,
as in the following diagram:

\[
\xymatrix@C=4.8em@R=3.4em{
M
\ar[r]_{\alpha}
\ar[dr]_{\phi}
\ar@/^3.0pc/[rrr]^{\pi}
&
\Gamma
\ar[r]_{\mu}
\ar[d]^{\varphi}
\ar@/^1.0pc/[rr]^{\beta}
&
\Gamma'
\ar[r]_{\beta'}
\ar[dl]^{\varphi'}
&
X
\ar@{.>}@/^1.0pc/[dll]^{\psi}
\\
&
Y & &
}
\]

Since $\psi \colon X \dashrightarrow Y$ is almost holomorphic, there exists a Zariski-open subset
$U_Y \subset Y$ such that
$\beta' \colon \Gamma' \to X$ is an isomorphism over $U_Y$.
In particular, the inverse image $(\varphi')^{-1}(U_Y)$ is strongly $\mathbb{Q}$-factorial and has klt singularities.
Then, it follows from Theorem~\ref{thm-Q-factorialization} that $\mu \colon \Gamma \to \Gamma'$
is an isomorphism over $U_{Y}$.
Consequently, the morphism $\beta \colon \Gamma \to X$ is also an isomorphism over $U_{Y}$.
In particular, the $\beta$-exceptional locus
$E_\beta$  does not dominate $Y$, 
where $E_\beta$ denotes the effective reduced divisor whose support coincides with the $\beta$-exceptional locus.

By construction, the variety $\Gamma$ is strongly $\mathbb{Q}$-factorial and has klt singularities.
A general fiber $\Gamma_y$ is a rationally connected klt variety, and hence $H^i(\Gamma_y,\mathcal{O}_{\Gamma_y})=0$ for every $i>0$.
Therefore, by \cite[Corollary 4.2]{CH24}, there exists a relatively $\varphi$-ample line bundle $G'$.
We may assume that the resolution $\alpha\colon M\to\Gamma$ is a projective morphism.
Consequently, both $\varphi$ and $\phi=\varphi\circ\alpha$ are projective.
After replacing $G'$ with a sufficiently divisible multiple,
Lemma~\ref{lem-Qfact} allows us to decompose $G'$ as 
\[
G'=\beta^{*}A+F', 
\]
where $A$ is a line bundle on $X$ and  $F'$ is a $\beta$-exceptional divisor.
By the construction of $\Gamma$, the morphism $\beta \colon \Gamma \to X$ is an isomorphism over $U_Y \subset Y$.
This shows that the restriction of $A$ to a general fiber of $\psi \colon X \dashrightarrow Y$ is also ample.

Fix $d \gg 1$ such that $dE_{\beta}+F'$ is effective, and define $G$ and $F$ by
\[
G:=\pi^{*}A+F
:=\pi^{*}A+\alpha^{*}(dE_{\beta}+F')
=\alpha^{*}(G'+dE_{\beta}).
\]
Then $F$ is an effective $\pi$-exceptional divisor, and hence property~\textup{(1)}
follows directly from the construction.

We next verify property~\textup{(2)}.
Let $S$ be an irreducible component of $M_y$ that is not contained in $E$.
Then $\alpha|_S \colon S \to \alpha(S) \subseteq \Gamma_y$ is bimeromorphic 
and the image $\alpha(S)$ is not contained in $E_\beta$.
Since $G'$ is relatively $\varphi$-ample, the restriction $G'|_{\alpha(S)}$ is ample.
Moreover, the divisor $E_\beta|_{\alpha(S)}$ is effective, and hence
$(G'+dE_\beta)|_{\alpha(S)}$ is big.
This shows that 
$G|_S = (\alpha|_S)^* \bigl((G'+dE_\beta)|_{\alpha(S)}\bigr)$ is big,
which proves property~\textup{(2)}. 

Finally, we verify property~\textup{(3)}.
Let  $E_\alpha$ be the effective reduced divisor whose support coincides with the $\alpha$-exceptional locus.
Then, we can easily see that 
\[
E_{\alpha}\leq
E=E_{\alpha}+\alpha^{-1}_{*}E_{\beta}
\leq E_{\alpha}+\alpha^{*}E_{\beta}.
\]
For a general point $y\in Y$, we have $\Gamma_y\cap E_\beta=\emptyset$, 
where $\Gamma_{y} $ is the fiber of $\Gamma \to Y$ at $y \in Y$. 
Moreover, the induced morphism
$\alpha|_{M_y}\colon M_y\to\Gamma_y$ is bimeromorphic and
$E_{\alpha}|_{M_y}$ is $\alpha|_{M_y}$-exceptional.

Hence, by the projection formula and $E\leq E_{\alpha}+\alpha^{*}E_{\beta}$, we obtain
\begin{align*}
&H^{0}\bigl(M_y,\OX_{M_y}(m(G+cE))\bigr) \\
\quad\subset
&H^{0}\Bigl(
    M_y,
    \OX_{M_y}\bigl(
        mG+mc(E_{\alpha}+\alpha^{*}E_{\beta})
    \bigr)
\Bigr) \\
\quad\cong
&H^{0}\Bigl(
    \Gamma_y,
    (\alpha|_{M_y})_{*}
    \OX_{M_y}\bigl(
        mG+mc(E_{\alpha}+\alpha^{*}E_{\beta})
    \bigr)
\Bigr) \\
\quad\cong
&H^{0}\Bigl(
    \Gamma_y,
    \OX_{\Gamma_y}\bigl(m(G'+dE_{\beta})\bigr)
    \otimes
    (\alpha|_{M_y})_{*}
    \OX_{M_y}\bigl(
        mc(E_{\alpha}+\alpha^{*}E_{\beta})
    \bigr)
\Bigr) \\
\quad\cong
&H^{0}\bigl(
    \Gamma_y,
    \OX_{\Gamma_y}(mG')
\bigr).
\end{align*}
Here, to obtain the last isomorphism, we used the facts that
\[
(\alpha|_{M_y})_{*}
\OX_{M_y}(mcE_{\alpha})
=
\OX_{\Gamma_y}
\quad\text{and}\quad
\Gamma_y\cap E_{\beta}=\emptyset.
\]
Applying the same argument and using $E_{\alpha}\leq E$, we obtain
\begin{align}
\label{eq-rank}
H^{0}\bigl(
    M_y,
    \OX_{M_y}(m(G+cE))
\bigr)
\cong
H^{0}\bigl(
    \Gamma_y,
    \OX_{\Gamma_y}(mG')
\bigr).
\end{align}

Since $G'$ is relatively $\varphi$-ample, 
after replacing $G'$ with a sufficiently divisible multiple, 
we may assume that, for a general point $y\in Y$, the section ring
$
\oplus_{k\in\mathbb{Z}_{\geq 0}}
H^{0}(\Gamma_y,kG')
$
is generated in degree one.
Combining this with \eqref{eq-rank}, we obtain the surjectivity of the natural multiplication map
\[
\Sym^{m}H^{0}\bigl(M_y,\OX_{M_y}(G+cE)\bigr)
\longrightarrow
H^{0}\bigl(M_y,\OX_{M_y}(m(G+cE))\bigr),
\]
and hence property~\textup{(3)}.
\end{proof}

Throughout this section, we freely use the notation introduced in Lemma~\ref{lemm-good}.
Fix an integer $c\in\mathbb{Z}_{+}$, which will be chosen sufficiently large later, and consider the determinant sheaf
\[
\det \phi_*\mathcal{O}_{M}(G+cE)
:=
\wedge^{[r]}\bigl(\phi_*\mathcal{O}_{M}(G+cE)\bigr),
\]
where
$
r:=\rank\phi_*\mathcal{O}_{M}(G+cE).
$
Note that the rank $r$ is independent of $c$ by \eqref{eq-rank}. 
Moreover, since $Y$ is smooth, this determinant sheaf is invertible.
We define the $\phi$-big line bundle $L_m$ on $M$ by
\begin{align}\label{eq-line}
L_m:=L_{c,m}:=m(G+cE)
-\frac{m}{r}\phi^* \det \phi_* \mathcal{O}_{M}(G+cE),
\end{align}
where $m$ is a positive integer such that $m/r \in \mathbb{Z}_{+}$.
We also define the direct image sheaf $\calV_m$ on $Y$ by
\begin{align}\label{eq-line2}
\calV_m:=\calV_{c,m}:=\phi_*\OX_M(L_m).
\end{align}
Let $Y_0 \subset Y$ be the maximal Zariski-open subset with the following properties:
\begin{itemize}
\item The morphism $\phi \colon M\to Y$ is flat over $Y_0$.
\item For every prime divisor $P$ on $Y_0$, the divisor $\phi^*P$ is not $\pi$-exceptional.
\end{itemize}

We can now state the main results of this section.

\begin{thm}\label{thm-faked-flat}
Under the assumptions of Setting~\ref{setting},
let $c\in\mathbb{Z}_{+}$ be a fixed sufficiently large integer and
$m\in\mathbb{Z}_{+}$ be a sufficiently large and divisible integer.
Then, the reflexive sheaf
\[
\mathcal{E}_{m} := \pi_{[*]}(\phi^{*}\calV_{m})
\]
is weakly positively curved and satisfies
$
c_{1}(\mathcal{E}_{m})=0.
$
\end{thm}

\begin{cor}\label{thm-flat}
Under the assumptions of Setting~\ref{setting},
assume in addition that $X$ is maximally quasi-\'etale.
Then, the sheaf $\calV_{m}$ is locally free on $Y_{0}$ and admits a flat connection over $Y_{0}$.
\end{cor}

\subsection{Bimeromorphic properties of MRC fibrations}\label{subsec-mrc}

In this subsection, we study bimeromorphic properties of the MRC fibration
$\psi\colon X\dashrightarrow Y$, extending the original result of \cite{Zhang05}
to the setting of K\"ahler generalized pairs.

\begin{prop}
\label{prop_bir-geometry-psi}
Under the assumptions of Setting~\ref{setting},
the following assertions hold:
\begin{itemize}
\item[\rm(a)]
The canonical divisor $K_Y$ is $\QQ$-linearly equivalent to an effective
$\QQ$-divisor $N_Y$ such that $\phi^{*}N_Y$ is $\pi$-exceptional.
In particular, we have $\kappa(Y)=0$.

\item[\rm(b)]
Every irreducible component of $\Delta$ is horizontal with respect to $\psi$;
that is,
$
Y=\phi(\pi^{-1}_{*}\Delta_i)
$
for every irreducible component $\Delta_i$ of $\Delta$.

\item[\rm(c)]
The subset $\pi(\phi\inv(Y\backslash Y_0))$ has codimension at least two in $X$.
In particular, every $\phi$-exceptional divisor on $M$ is $\pi$-exceptional.

\item[\rm(d)]
The variety $Y_0$ satisfies the generalized Liouville property in the following sense:
for every flat vector bundle $(\calH_0,\nabla_0)$ on $Y_0$ satisfying
condition \rm($\bullet$) below, every global section of $\calH_0$ is parallel
with respect to $\nabla_0$.
\begin{itemize}
\item[\rm($\bullet$)]
There exists a numerically flat vector bundle $\calH$ on $X$ such that
\[
(\phi^\ast\calH_0,\phi^\ast\nabla_0)
\cong
(\pi^\ast\calH,\nabla)
\quad\text{on } M_0:=\phi\inv(Y_0),
\]
where $\phi^\ast\nabla_0$ denotes the pullback connection on
$\phi^\ast\calH_0$, and $\nabla$ is the (unique) flat connection on
$\pi^\ast\calH$ given by
\cite[Corollary 3.10 and the discussion thereafter]{Sim92},
which is compatible with the filtration provided by
\cite[Theorem 1.18]{DPS94}.
\end{itemize}

\item[\rm(e)]
The MRC fibration $\psi$ is semistable in codimension one; that is, for every
prime divisor $P$ on $Y_0$, if
$
\phi^\ast P=\sum_i c_iP_i
$
is its decomposition into prime divisors on $\phi\inv(Y_0)$, then every component
$P_i$ with $c_i>1$ is $\pi$-exceptional.
\end{itemize}
\end{prop}

\begin{rem}\label{rem-difference1}
As in the projective case, the K\"ahler version of the uniruledness criterion
of \cite{BDPP13}, recently established in \cite{Ou25}, plays an important role in the proof.
The main differences between our proof and those in previous works are as follows: 
First, in previous works, a suitable boundary divisor is constructed
using an ample line bundle on $X$.
In our proof, the positivity of direct image sheaves plays the corresponding role.
Second, we must work with $\mathbb{R}$-divisors, which introduces additional
technical difficulties, since the divisors $\Delta_M$ and $F_{+}$ appearing
in \eqref{eq-gpair2} are not necessarily $\mathbb{Q}$-divisors.
Third, the previous proof of \textup{(e)} relies on the existence of global
Kawamata covers, a technique that is not available in the K\"ahler setting.
Therefore, we give a genuinely new proof based on the generic nefness.
To this end, we extend the generic nefness theorem to generalized pairs
(see Theorem~\ref{nefness}).
\end{rem}

\begin{proof}
Conclusions~\textup{(a)} and~\textup{(e)} constitute the core of the proposition,
and their proofs are new.
Conclusion~\textup{(b)} follows from a slight modification of the new proof of
Conclusion~\textup{(a)}.
Conclusions~\textup{(c)} and~\textup{(d)} follow directly from the other
conclusions, as explained in the previous work, but we include the details
for completeness.

\begin{proof}[Proof of {\rm(a)}]
We divide the proof of \textup{(a)} into four steps.

\begin{step}
In this step, we show that $F_{+}-\phi^{*}K_{Y}$ is pseudo-effective.
To this end, for a sufficiently large and divisible integer $m$, we consider
the direct image sheaf
\[
\phi_{*}\left(
\mathcal{O}_{M}\left(
mK_{M/Y} 
+G
+(m\Delta_{M}+\lceil mF_{+}\rceil-mF_{+})
+m(F_{+}-K_{M}-\Delta_M)
\right)
\right).
\]
This direct image sheaf satisfies the assumptions of
Proposition~\ref{prop-psef} (1).
Indeed, the line bundle $G$ is $\phi$-big by construction, and
\[
N:=F_{+}-K_{M}-\Delta_M
\]
is nef by \eqref{eq-gpair2}.
Moreover, since $\Delta_{M}$ is a klt boundary and
\[
\frac{1}{m}(\lceil mF_{+}\rceil-mF_{+})\longrightarrow 0 \text{ as } m\to\infty, 
\]
the boundary
\[
L:=\Delta_{M}
+\frac{1}{m}(\lceil mF_{+}\rceil-mF_{+})
\]
satisfies the klt condtion for $m \gg 1$.

On the other hand, 
a direct computation using $K_{M/Y}=K_M-\phi^{*}K_Y$ gives 
\[
\begin{aligned}
m(K_{M/Y}+\Delta_M)
+G
+(\lceil mF_{+}\rceil-mF_{+})
+m(F_{+}-K_M-\Delta_M) 
=\lceil mF_{+}\rceil-m\phi^{*}K_Y+G.
\end{aligned}
\]
Hence, by the projection formula, the above direct image sheaf is
isomorphic to
\[
\begin{aligned}
\phi_{*}\left(
\mathcal{O}_{M}
\bigl(\lceil mF_{+}\rceil-m\phi^{*}K_Y+G\bigr)
\right) \cong
\phi_{*}\left(
\mathcal{O}_{M}
\bigl(\lceil mF_{+}\rceil+G\bigr)
\right)
\otimes\mathcal{O}_{Y}(-mK_Y).
\end{aligned}
\]
The relevant sheaf is nonzero since $F_{+}$ is effective and $G$ is $\phi$-big.
By applying Proposition~\ref{prop-psef} (1), we conclude that 
\[
\frac{1}{m}\lceil mF_{+}\rceil-\phi^{*}K_Y + \frac{1}{m} G
\]
is pseudo-effective.
Letting $m\to\infty$, we deduce that 
$F_{+}-\phi^{*}K_Y$ is pseudo-effective.
\end{step}

\begin{step}
In this step, we show that $\phi^{*}K_Y$ is numerically equivalent to a
$\pi$-exceptional $\mathbb{Q}$-divisor.
By Lemma~\ref{lem-Qfact}, we have
\[
\phi^{*}K_Y\sim_{\mathbb{Q}}\pi^{*}L+N
\]
for some $\mathbb{Q}$-line bundle $L$ on $X$ and some
$\pi$-exceptional $\mathbb{Q}$-divisor $N$.
Thus, it suffices to show that $L$ is numerically trivial.

By Step~1, the class
$c_1(F_{+}-\phi^{*}K_Y)$ is represented by a closed semipositive
$(1,1)$-current.
Pushing this current forward by $\pi$, we obtain a closed semipositive
$(1,1)$-current representing the class
\[
\pi_{*}c_1(F_{+}-\phi^{*}K_Y)
=
-c_1(L)
\in H^{1,1}_{\BottChern}(X).
\]
Here we used Proposition~\ref{BEG}, together with the facts that
$F_{+}$ and $N$ are $\pi$-exceptional and that
$\pi_*c_1(\pi^*L)=c_1(L)$.
Therefore, we see that $-L$ is pseudo-effective.

On the other hand, although $M$ need not be projective, 
we can apply \cite[Theorem~1.1]{GHS03} to conclude that 
the base $Y$ of the MRC fibration is not uniruled.
Indeed, suppose that there exists a rational curve $R\subseteq Y$
belonging to a covering family, and let
$\nu\colon\mathbb{P}^{1}\to R$ be its normalization.
Let $M_R$ be the main component of the fiber product 
$
M\times_Y\mathbb{P}^{1}
$
that dominates $\mathbb{P}^{1}$.
Since $\phi \colon M \to Y$ is a projective fibration, the variety $M_R$ is projective.
After taking a projective resolution of $M_R$, we can apply
\cite[Theorem~1.1]{GHS03} to the induced fibration over
$\mathbb{P}^{1}$, by noting that a general fiber is rationally connected.
Therefore, we can obtain a lift of $\mathbb{P}^{1}$.
Together with the rational connectedness of the general fibers, 
this lift shows that $M_R$ is rationally chain connected.
Its image in $X$ is rationally chain connected and dominates $R$.
Since the curves $R$ form a covering family, this contradicts the maximality of MRC fibrations.

By the K\"ahler analogue of the uniruledness criterion of
\cite{BDPP13}, recently established in \cite{Ou25}, it follows that
$K_Y$ is pseudo-effective.
Hence $\phi^*K_Y$ is pseudo-effective.
Applying Proposition~\ref{BEG} to the pushforward by $\pi$ again and using
$
\phi^*K_Y\sim_{\mathbb{Q}}\pi^*L+N, 
$
we deduce that $c_1(L)$ is pseudo-effective.
Combined with the preceding conclusion, we see that $L$ is numerically trivial.

\end{step}

\begin{step}
In this step, we show that there exists an effective divisor $N_Y$ on $Y$
such that $N=\phi^{*}N_Y$.
Let $T$ be a closed semipositive $(1,1)$-current representing $c_1(K_Y)$.
The pullback $\phi^{*}T$ is a closed semipositive $(1,1)$-current representing
$
c_1(\phi^{*}K_Y)=c_1(N).
$
Since $N$ is $\pi$-exceptional, the pushforward
$\pi_{*}\phi^{*}T$ represents the zero Bott--Chern class.
Since $\pi_{*}\phi^{*}T$ is a closed semipositive current, we have
$
\pi_{*}\phi^{*}T=0.
$
On the other hand, the integration current $[N]$ represents the same
Bott--Chern class $c_1(\phi^{*}K_Y)$, and $\pi_{*}[N]$ is also the zero current.
Furthermore, the integration current $[N]$, and hence
$\phi^{*}T-[N]$, can be written as the difference of semipositive currents.
Thus, by Lemma~\ref{topo-lemma}, we conclude that
$\phi^{*}T=[N]$.
In particular, it follows that $N$ is an effective divisor.

Let $y\in Y$ be a general point outside the polar set of a local
potential of $T$.
Then, the restriction $(\phi^*T)|_{M_y}$ to the fiber $M_y$ is
well-defined, and we have
\[
[N|_{M_y}]
=
(\phi^*T)|_{M_y}
=
0.
\]
Since $N$ is effective, this implies that $N$ has no
$\phi$-horizontal component.
Moreover, since
$
\Supp{T}\subset\phi(\operatorname{Supp}N),
$
the support of $T$ is contained in a proper Zariski-closed subset of
$Y$.
By the support theorem for closed semipositive currents
(see \cite[III.2.14, p.~143]{agbook}), there exists an effective
$\mathbb{R}$-divisor $N_Y$ on $Y$ such that
$
T=[N_Y].
$
The equality of currents
\[
[N]=\phi^*T=[\phi^*N_Y]
\]
gives
$
N=\phi^*N_Y.
$
Since $N$ is a $\mathbb{Q}$-divisor, so is $N_Y$.
Thus $N_Y$ is an effective $\mathbb{Q}$-divisor such that 
$N=\phi^{*}N_Y$ and 
$
c_1(N_Y)=c_1(K_Y).
$
\end{step}

\begin{step}
In this step, we show that $K_Y\sim_{\QQ}N_Y$ and $\kappa(K_Y)=0$.
By Steps~2 and~3, the $\mathbb{Q}$-line bundle
$
L_0:=N_Y-K_Y
$
is numerically trivial.
Since $N_Y$ is effective, it follows from
\cite[Corollary 5.8(a)]{JWang19} that
$$
\kappa(K_Y)
\geq
\kappa(K_Y+L_0)
=
\kappa(N_Y)
\geq 0.
$$
On the other hand, for every sufficiently divisible positive integer
$m\in\mathbb{Z}_+$, we have
\[
\begin{aligned}
h^0(Y,mK_Y)
=h^0(M,m\phi^*K_Y) 
=h^0(M,m(\pi^*L+N)) 
=h^0(X,mL) 
\leq 1.
\end{aligned}
\]
Here, the third equality follows from the fact that $N$ is an effective
$\pi$-exceptional divisor, and the last inequality follows from the
numerical triviality of $L$.
Therefore, we obtain $\kappa(K_Y)=0$.
Then, the Campana--Koziarz--P{\u a}un-type result
\cite[Corollary 5.8(b)]{JWang19} shows that $L_0\sim_{\QQ}0$.
This means that 
$
K_Y\sim_{\QQ}N_Y.
$
\end{step}

\vspace{-0.5cm}
\end{proof}
\begin{proof}[Proof of \textup{(b)}]
Consider the decomposition of $\Delta$ into its horizontal and vertical
parts with respect to $\psi\colon X\dashrightarrow Y$:
\[
\Delta=\Delta^v+\Delta^h.
\]
As in Step~1 of the proof of Conclusion~\textup{(a)}, we apply
Proposition~\ref{prop-psef}\textup{(1)} to the direct image sheaf
\begin{align*}
\phi_*\mathcal{O}_M\Bigl(
mK_{M/Y}+G
+\bigl(m\Delta_M-\lfloor m\pi^{-1}_*\Delta^v\rfloor\bigr)
+\bigl(\lceil mF_+\rceil-mF_+\bigr)
+m(F_+-K_M-\Delta_M)
\Bigr).
\end{align*}
Indeed $G$ is $\phi$-big and $F_+-K_M-\Delta_M$ is nef.
Moreover, the $\mathbb{R}$-divisor
\[
\frac{1}{m}\left(
m\Delta_M-\lfloor m\pi^{-1}_*\Delta^v\rfloor
+\lceil mF_+\rceil-mF_+
\right)
\]
is effective and satisfies the klt condition for $m\gg1$.
By a direct computation and the projection formula, the above sheaf is
isomorphic to
\[
\phi_*\mathcal{O}_M\left(
\lceil mF_+\rceil
-\lfloor m\pi^{-1}_*\Delta^v\rfloor
+G
\right)
\otimes\mathcal{O}_Y(-mK_Y).
\]
Since $\pi^{-1}_*\Delta^v$ is $\phi$-vertical and $F_+$ is effective, this direct image sheaf is nonzero.

Proposition~\ref{prop-psef}(1) shows that
$$
\lceil m F_{+} \rceil
-\lfloor m\pi^{-1}_{*}\Delta^{v}\rfloor
+G
-m\phi^{*}K_{Y}
$$
is pseudo-effective.
Passing to the limit as $m\to\infty$, we find that
\[
F_{+}-\pi^{-1}_{*}\Delta^{v}-\phi^{*}K_{Y}
\]
is pseudo-effective.
The pullback $\phi^{*}K_{Y}$ is $\pi$-exceptional by Conclusion~\textup{(a)}.
Pushing forward via $\pi$, we conclude that $-\Delta^{v}$ is pseudo-effective. 
Since $\Delta^{v}$ is effective, it follows that $\Delta^{v}=0$, which proves Conclusion~\textup{(b)}.
\end{proof}

\begin{proof}[Proof of {\rm (c)}]
Let \(V\subset M\) be a prime divisor such that
$
\phi(V)\subset Y\setminus Y_{0}.
$
By the definition of \(Y_{0}\), if \(\phi(V)\) has codimension one in \(Y\), then \(V\) is already \(\pi\)-exceptional. Hence, we may assume that
$
\codim \phi(V)\geq 2,
$ (i.e.,\,\(V\) is \(\phi\)-exceptional.)

Let \(\beta_Y\colon Y_1\to Y\) be a desingularization of the blow-up of \(Y\) along \(\phi(V)\).
Since \(Y\) is smooth (in particular, $Y$ has terminal singularities), we have
\[
K_{Y_1}=\beta_Y^*K_Y+F_Y,
\]
where \(F_Y\) is an effective \(\beta_Y\)-exceptional divisor such that
$
\Supp(F_Y)=\Exc(\beta_Y).
$
Let \(M_1\) be a desingularization of the component of the fiber product
\(M\times_Y Y_1\) that dominates \(Y_1\). We then have the following commutative diagram:
\begin{center}
\begin{tikzpicture}[scale=2.5]
\node (A) at (0,0) {$Y$};
\node (B) at (0,1) {$M$};
\node (C) at (1,1) {$X$};
\node (A1) at (-1,0) {$Y_1$};
\node (B1) at (-1,1) {$M_1$};
\path[->,font=\scriptsize,>=angle 90]
(B) edge node[left]{$\phi$} (A)
(B) edge node[above]{$\pi$} (C)
(B1) edge node[left]{$\phi_1$} (A1)
(B1) edge node[above]{$\beta_M$} (B)
(A1) edge node[below]{$\beta_Y$} (A);
\path[dashed, ->,font=\scriptsize,>=angle 90]
(C) edge node[below right]{$\psi$} (A);
\end{tikzpicture}
\end{center}

Let \(V_1\subset M_1\) be the strict transform of \(V\). 
Then, by construction, we have 
$
\phi_1(V_1)\subset \Exc(\beta_Y).
$
By Conclusion~{\rm (a)}, there exists an effective \(\mathbb Q\)-divisor \(N_Y\) on \(Y\) such that
$
N_Y\sim_{\mathbb Q} K_Y.
$
The $\mathbb{Q}$-divisor \(\beta_Y^*N_Y+F_Y\)
is an effective \(\mathbb Q\)-divisor on \(Y_1\) that is \(\mathbb Q\)-linearly equivalent to \(K_{Y_1}\).
Applying Conclusion~{\rm (a)} to the MRC fibration \(X\dashrightarrow Y_1\) instead of $\phi \colon X\dashrightarrow Y$, we find that
$
\phi_1^*(F_Y+\beta_Y^*N_Y)
$
is \(\mathbb Q\)-linearly equivalent to a \((\pi\circ\beta_M)\)-exceptional divisor.
Since \(\phi_1^*(F_Y+\beta_Y^*N_Y)\) has Kodaira dimension zero, we deduce that
$
\phi_1^*(F_Y+\beta_Y^*N_Y)
$
is itself \((\pi\circ\beta_M)\)-exceptional.

On the other hand, we have 
\[
\phi_1(V_1)\subset \Exc(\beta_Y)
=\Supp(F_Y)
\subset \Supp(\beta_Y^*N_Y+F_Y).
\]
This shows that 
$
V_1\subset \Supp\phi_1^*(\beta_Y^*N_Y+F_Y),
$
and thus \(V_1\) is \((\pi\circ\beta_M)\)-exceptional. Consequently,
the prime divisor $
V=\beta_M(V_1)
$
is \(\pi\)-exceptional.
\end{proof}

\begin{proof}[Proof of \textup{(e)}]

The proof of Conclusion~\textup{(e)} relies on the generic nefness and the theory of foliations; this argument is new even in the projective setting.

Consider the exact sequence
\[
0 \to \mathcal{F} \to T_{M} \to \mathcal{Q}:=T_{M}/\mathcal{F} \to 0,
\]
where \(\mathcal{F} \subset T_{M}\) is the foliation induced by
\(\phi \colon M \to Y\).
Taking reflexive pushforwards, we obtain a sheaf morphism on $X$
\[
T_{X}=\pi_{[*]}T_{M} \to \pi_{[*]}\mathcal{Q},
\]
which is surjective over \(X \setminus \pi(E)\).
The generic nefness theorem, Theorem~\ref{nefness}, shows that 
$
c_{1}(\pi_{[*]}\mathcal{Q}) \cdot \alpha \geq 0,
$
where \(\alpha:=\{\omega_{X}\}^{n-1}\) and \(\omega_{X}\) is a K\"ahler form on \(X\).
Since \(\codim \pi(E) \geq 2\) and determinant sheaves are reflexive, we have
$
\pi_{[*]}\det \mathcal{Q}
=
\det \pi_{[*]}\mathcal{Q}.
$
Consequently, we obtain 
\begin{equation}\label{eq-generic}
c_{1}(\pi_{[*]}\det \mathcal{Q}) \cdot \alpha \geq 0.
\end{equation}
Here, we also use the fact that \(X\) is strongly \(\QQ\)-factorial.

We next compute \(\pi_{[*]}\det \mathcal{Q}\) using the theory of foliations.
By the construction of \(Y_{0}\), the morphism
\(\phi \colon M \to Y\) is flat, and hence equidimensional, over \(Y_{0}\).
Therefore, over \(\phi^{-1}(Y_{0})\), we have
\[
\det \mathcal{F}
\cong
\mathcal{O}_{M}(-K_{M/Y}+R(\phi))
\quad\text{and hence}\quad
\det \mathcal{Q}
\cong
\mathcal{O}_{M}(-\phi^*K_Y-R(\phi)),
\]
where \(R(\phi)\) denotes the ramification divisor of \(\phi\) 
(see, for example, \cite[Lemma~2.31]{CKT16}).
Although this formula is established only over \(\phi^{-1}(Y_{0})\), 
every divisorial component of \(\phi^{-1}(Y \setminus Y_{0})\) is \(\pi\)-exceptional by
Conclusion~\textup{(c)}.
Hence, we see that 
\[
\pi_{[*]}\det \mathcal{Q}
=
\pi_{[*]}\mathcal{O}_{M}(-\phi^*K_Y-R(\phi))
=
\mathcal{O}_{X}(-\pi_{*}R(\phi)).
\]
The last equality follows from Conclusion~\textup{(a)}. 
It follows from \eqref{eq-generic} that
\[
-c_{1}\bigl(\pi_{*}R(\phi)\bigr)\cdot\alpha
\geq 0.
\]
Since \(\pi_{*}R(\phi)\) is effective and every nonzero effective divisor has positive intersection with \(\alpha\), 
we conclude that
$
\pi_{*}R(\phi)=0.
$
\end{proof}

\begin{proof}[Proof of \textup{(d)}]
The generalized Liouville property follows essentially from Conclusion~\textup{(c)}.
Let \(s\in\Coh^0(Y_0,\calH_0)\).
The section \(s\) is parallel with respect to \(\nabla_0\) if and only if
\(\phi^\ast s\) is parallel with respect to
$
\phi^\ast\nabla_0=\nabla|_{\phi^{-1}(Y_0)}.
$

By Conclusion~\textup{(c)}, the complement of
$
\pi\bigl(\phi^{-1}(Y_0)\setminus E\bigr)
$
in \(X\) has codimension at least two.
Via the isomorphism
\[
\pi|_{\phi^{-1}(Y_0)\setminus E}\colon
\phi^{-1}(Y_0)\setminus E
\xrightarrow{\ \cong\ }
\pi\bigl(\phi^{-1}(Y_0)\setminus E\bigr),
\]
the section
$
\phi^\ast s|_{\phi^{-1}(Y_0)\setminus E}
$
extends to a section
$
\sigma\in\Coh^0(X,\calH)
$
by the reflexivity of \(\calH\).
Then \(\pi^\ast\sigma\) is parallel with respect to \(\nabla\) by
\cite[Theorem~4.3.3]{Cao13}, and hence so is \(\phi^\ast s\).
\end{proof}
The above arguments complete the proof.
\end{proof}

As in \cite{Zhang05}, Proposition~\ref{prop_bir-geometry-psi} implies that
the Albanese map is a fibration without passing to a quasi-\'etale cover.
This extends \cite[Corollary~2]{Zhang05} to K\"ahler generalized pairs.

\begin{cor}
Let $(X,\Delta,\boldsymbol{\beta})$ be a K\"ahler generalized klt pair
of Calabi--Yau type.
Note that we do not assume that $X$ is strongly $\mathbb{Q}$-factorial.
Then, the Albanese map
$\alb_X\colon X\to\Alb_X$ is a fibration, that is, it is a proper
surjective morphism with connected fibers.
\end{cor}

\begin{proof}
Since $(X,\Delta,\boldsymbol{\beta})$ is generalized klt,
the variety $X$ is locally potentially klt by
Proposition~\ref{prop-local-klt}.
Since a complex torus contains no rational curves,
the Albanese map $\alb_X\colon X\to\Alb_X$ is a morphism by
\cite[Corollary~1.4]{Fuj26}
(see also \cite[Corollary~1.7]{HM07}).

Let $q\colon X^{\qf}\to X$ be a small $\QQ$-factorialization provided by
Theorem~\ref{thm-Q-factorialization}.
The generalized pair induced on $X^{\qf}$ (see Theorem~\ref{thm-Q-factorialization}) 
is again generalized klt and of
Calabi--Yau type, and $X^{\qf}$ is strongly $\QQ$-factorial.
Take a desingularization $\pi'\colon M\to X^{\qf}$, 
an MRC fibration of $X^{\qf} \dashrightarrow Y$, 
a fibration $\phi \colon M \to Y$ as in Setting \ref{setting}. 
By Proposition~\ref{prop_bir-geometry-psi}, 
we obtain $\kappa(Y)=0$.

Set $\pi:=q\circ\pi'\colon M\to X$.
By \cite[Theorem~C]{JWang19}, the Albanese map
$\alb_Y\colon Y\to\Alb_Y$ is a fibration.
Consequently, the composition
$\alb_Y\circ\phi\colon M\to\Alb_Y$ is also a fibration.
Since $X$ is locally potentially klt, it has rational singularities.
It follows that $\pi$ induces an isomorphism
$\Alb_M\cong\Alb_X$.
After choosing compatible base points, we may therefore identify
$\alb_X\circ\pi\colon M\to\Alb_X$ with the Albanese map of $M$.
By the universal property of the Albanese map, there exists a
homomorphism of complex tori
$f\colon\Alb_X\to\Alb_Y$ such that
\[
\alb_Y\circ\phi=f\circ\alb_X\circ\pi.
\]
Since $\alb_Y\circ\phi$ is surjective, so is $f$.
\begin{equation*}
\xymatrix@C=5em@R=3.5em{
M
\ar[r]^{\pi'}
\ar[d]_{\phi}
&
X^{\qf}
\ar[r]^{q}
\ar@{-->}[dl]
&
X
\ar[d]^{\alb_X}
\\
Y
\ar[r]_{\alb_Y}
&
\Alb_Y
&
\Alb_X
\ar[l]^{f}.
}
\end{equation*}

By the properties of MRC fibrations, we have
$
H^1(M,\OX_M)\cong H^1(Y,\OX_Y).
$
Hence, we obtain 
$
\dim\Alb_X=\dim\Alb_M=\dim\Alb_Y.
$
Thus, the surjective homomorphism $f$ is an isogeny 
(see \cite[\S\S~2.3--2.4, pp.~7--8]{Deb99}).

Since $f$ is finite and
\[
f\bigl((\alb_X\circ\pi)(M)\bigr)=\Alb_Y,
\]
the image $(\alb_X\circ\pi)(M)$ has dimension $\dim\Alb_X$ and thus is equal to $\Alb_X$.
Hence, we see that $\alb_X\circ\pi$ is surjective.
If $f$ were not an isomorphism, a fiber of $f$ would contain at least
two points. Since $\alb_X\circ\pi$ is surjective, the corresponding
fiber of $\alb_Y\circ\phi$ would then be a disjoint union of at least
two nonempty fibers of $\alb_X\circ\pi$, contradicting the
connectedness of the fibers of $\alb_Y\circ\phi$.
Therefore $f$ is an isomorphism, and hence
$\alb_X\circ\pi$ is a fibration.

Finally, since $\pi$ is a proper bimeromorphic morphism onto the normal
variety $X$, it is surjective and has connected fibers.
It follows that $\alb_X$ is surjective and has connected fibers.
Thus $\alb_X\colon X\to\Alb_X$ is a fibration.
\end{proof}

\subsection{Proof of flatness}\label{subsec-direct}
In this subsection, we prove Theorem~\ref{thm-faked-flat} and Corollary~\ref{thm-flat}.
To this end, we fix an effective $\mathbb{Q}$-divisor $N_Y$ such that
$K_Y\sim_{\mathbb{Q}}N_Y$, and define the $\mathbb{R}$-divisor
\begin{align}\label{eq-e'}
E' := F_{+} - \phi^{*} N_Y
\end{align}
on $M$. 
The divisor $E'$ is $\pi$-exceptional by Proposition~\ref{prop_bir-geometry-psi}\textup{(a)}.
By \eqref{eq-gpair2}, we obtain

\begin{align}
\label{eq_rel-canonical-pi-exc}
-(K_{M/Y}+\Delta_{M})+E' \equiv
\boldsymbol{\beta}_{M}. 
\end{align}
In particular, the left-hand side represents a nef Bott--Chern class.
We define the $\mathbb{Q}$-line bundle $D_{c,m}$ on $Y$ by
\begin{align*}
D_{c,m} := \frac{1}{r_m}\det \phi_{*} \mathcal{O}_M(m(G + cE)),
\end{align*}
where $r_m := \rank \phi_{*} \mathcal{O}_M(m(G + cE))$.

With these preparations in place, we first prove the following proposition.

\begin{prop}
\label{prop_more-psef-det}
Retain the notation and assumptions introduced above.
Then, the reflexive pushforward
\[
\pi_{[*]}\left(
\mathcal{O}_M\left(r_m(mG-\phi^*D_{c,m})\right)
\right)
\]
is weakly positively curved.
In particular, for a sufficiently large integer $d\gg1$, the $\mathbb{Q}$-line bundle
\[
r_m(mG-\phi^*D_{c,m})+dE
\]
is pseudo-effective on $M$ $($see Lemma~\ref{lem-Qfact}$)$.
\end{prop}

\begin{rem}\label{rem-more-psef-det}
This proposition corresponds to \cite[Proposition~3.3]{Wang20}, whose proof is divided into Steps~(A)--(C).
Step~(C) relies on the existence of an ample line bundle on $Y$.
However, the base variety $Y$ in our setting is non-projective.
Our strategy is to \textit{localize} the projective argument by using the positivity of direct image sheaves.
The new ingredients are contained in Steps~(C) and~(D).
Steps~(A) and~(B) remain essentially unchanged, 
but we will check several new properties for Steps~(C) and~(D).
\end{rem}

\begin{proof}
A natural approach would be to apply Proposition~\ref{prop-psef}\textup{(2)}
directly to $mG$.
However, this would require the lower bound of the curvature of $mG$ ,
which is precisely what we need to prove.
To overcome this difficulty, we use Viehweg's fiber-product trick.

For simplicity of notation, we consider only the case $m=1$ and set
$r:=r_1$.
The same argument applies to the general case.

\paragraph{\textbf{(A) Construction of the fiber product and the canonical section.}}

Let $M^r$ be the main component of the fiber product of $r$ copies of
$M$ over $Y$ 
\[
M\times_Y M\times_Y\cdots\times_Y M
\]
that dominates $Y$, and let $\pr_i\colon M^r\to M$ denote the natural
projections. Let $\mu\colon M^{(r)}\to M^r$ be a desingularization
such that $\mu$ is an isomorphism over the smooth locus $(M^r)_{\reg}$.
By replacing $M^{(r)}$ with a higher resolution, we may further assume
that the strict transform $\bar D\subset M^{(r)}$ of the diagonal
subvariety $D\subset M^r$ is smooth.
Set $p_i:=\pr_i\circ\mu$ and
$\phi^{(r)}:=\phi^r\circ\mu$, as in the following diagram:
\[
\xymatrix@C=5em@R=4em{
M^{(r)} \ar[r]^{\mu} \ar[dr]_{p_i}
\ar@/_4pc/[rrd]_{\phi^{(r)}} &
M^r \ar[r]^{\pr_i} \ar[d]_{\pr_j} \ar[dr]^{\phi^r} &
M \ar[d]^{\phi} \\
& M \ar[r]_{\phi} & Y.
}
\]
We further define
\begin{align*}
G^r&:=\sum_{i=1}^r\pr_i^*G,
&\quad
G^{(r)}&:=\mu^*G^r=\sum_{i=1}^rp_i^*G,\\
E^r&:=\sum_{i=1}^r\pr_i^*E,
&\quad
E^{(r)}&:=\mu^*E^r=\sum_{i=1}^rp_i^*E,\\
\Delta_{M^r}&:=\sum_{i=1}^r\pr_i^*\Delta_M,
&\quad
\Delta_{M^{(r)}}&:=\mu^*\Delta_{M^r}
=\sum_{i=1}^rp_i^*\Delta_M.
\end{align*}

Let $Y_{\mathrm{ff}}\subset Y$ be the maximal Zariski-open subset
over which $\phi$ is flat and
$\phi_*\mathcal{O}_M(G+cE)$ is locally free.

Although $\phi^{-1}(Y\setminus Y_{\mathrm{ff}})$ may have divisorial
components lying over $Y\setminus Y_0$, we have
\begin{equation}\label{eq-codim-Yff}
\codim 
\pi\bigl(\phi^{-1}(Y\setminus Y_{\mathrm{ff}})\bigr)\geq2.
\end{equation}
Indeed, Proposition~\ref{prop_bir-geometry-psi}\textup{(c)} gives
$
\codim 
\pi\bigl(\phi^{-1}(Y\setminus Y_0)\bigr)\geq2.
$
Since $\phi$ is flat over $Y_0$ and flat pullback preserves
codimension, we have
$
\codim \phi^{-1}(Y_0\setminus Y_{\mathrm{ff}})\geq2.
$
These inequalities prove \eqref{eq-codim-Yff}.

We prove that there exists an effective divisor $B_1$ on $M^{(r)}$
with the following properties:
\begin{itemize}
\item[$\bullet$]
$\Supp(B_1)\subset
M^{(r)}\setminus(\phi^{(r)})^{-1}(Y_{\mathrm{ff}})$;

\item[$\bullet$]
$L_0$ is pseudo-effective and relatively $\phi^{(r)}$-big,
\end{itemize}
where $L_0$ is the line bundle on $M^{(r)}$ defined by
\[
L_0
:=
G^{(r)}+cE^{(r)}
-r(\phi^{(r)})^*D_{c,1}+B_1.
\]
Since $M$ and $Y$ are smooth and $\phi$ is flat over
$Y_{\mathrm{ff}}$, every fiber $M_y$ at a point
$y\in Y_{\mathrm{ff}}$ is Gorenstein
(see \cite[\S~23, Theorem~23.4, p.~181]{Mat89}).
Over $Y_{\mathrm{ff}}$, there is a natural inclusion
\begin{equation}\label{eq_incl-det-tensor}
\det\phi_*\mathcal{O}_M(G+cE)|_{Y_{\mathrm{ff}}}
\cong
\mathcal{O}_{Y_{\mathrm{ff}}}
(rD_{c,1}|_{Y_{\mathrm{ff}}})
\hookrightarrow
\bigotimes^r
\phi_*\mathcal{O}_M(G+cE)|_{Y_{\mathrm{ff}}}.
\end{equation}
By flat base change and the projection formula applied to the full
$r$-fold fiber product, the right-hand side of
\eqref{eq_incl-det-tensor} is identified with the direct image of
the external tensor product of $\mathcal{O}_M(G+cE)$.
Restricting the resulting section to the component $M^r$ dominating
$Y$ and pulling it back by $\mu$, we obtain a nonzero section
\[
s_0\in
H^0\left(
(\phi^{(r)})^{-1}(Y_{\mathrm{ff}}),
\mathcal{O}_{M^{(r)}}\left(
G^{(r)}+cE^{(r)}
-r(\phi^{(r)})^*D_{c,1}
\right)
\right).
\]
The section $s_0$ may not extend to a global section of
$
G^{(r)}+cE^{(r)}-r(\phi^{(r)})^*D_{c,1}.
$
However, by
\cite[\S~III.5, Lemma~5.10, pp.~107--108]{Nak04}
(see also \cite[Theorem~1.14]{JWang19}), we can find an effective
divisor $B_1$ satisfying the required support condition such that
$s_0$ extends to a section of $L_0$.
In particular $L_0$ is pseudo-effective.
Moreover $G^{(r)}$ is relatively $\phi^{(r)}$-big,
$cE^{(r)}+B_1$ is effective, and
$(\phi^{(r)})^*D_{c,1}$ is trivial on every fiber.
It follows that $L_0$ is relatively $\phi^{(r)}$-big.

\smallskip
\paragraph{\textbf{(B) Comparison of the relative canonical divisors.}}

By induction and the base-change formula for relative canonical
bundles (see \cite[Proposition~(9)]{Kle80}), the inverse image
$
M^r_{\mathrm{ff}}:=(\phi^r)^{-1}(Y_{\mathrm{ff}})
$
is Gorenstein, and its relative dualizing sheaf satisfies
\[
\omega_{M^r_{\mathrm{ff}}/Y}
\cong
\mathcal{O}_{M^r_{\mathrm{ff}}}
\left(\sum_{i=1}^r\pr_i^*K_{M/Y}\right).
\]
Moreover, there is a natural morphism
\[
\mu_*
\mathcal{O}_{M^{(r)}}(K_{M^{(r)}/Y})
|_{M^r_{\mathrm{ff}}}
\longrightarrow
\omega_{M^r_{\mathrm{ff}}/Y},
\]
which is an isomorphism over the locus where $M^r$ has rational
singularities
(see \cite[Construction~3.19 and Definition~5.16]{Hor10}).

By assumption, the discriminant locus $\Sigma$ of
$\phi\colon M\to Y$ is an SNC divisor on $Y$, and the pullback
$\phi^*\Sigma$ has SNC support.
Write
\[
\phi^*\Sigma
=
\sum_\lambda W_\lambda+\sum_\mu a_\mu V_\mu,
\]
where $a_\mu>1$ for every $\mu$, and set
\[
W:=\sum_\lambda W_\lambda
\quad\text{and}\quad
V:=\sum_\mu V_\mu.
\]
Set $M_{\mathrm{ff}}:=\phi^{-1}(Y_{\mathrm{ff}})$.
Then, by \cite[Lemma~3.13]{Hor10}, the variety
$M^r_{\mathrm{ff}}$ has rational singularities along
\[
\underbrace{
\bigl(
M_{\mathrm{ff}}
\setminus(V\cup\phi^{-1}(\Sigma_{\sing}))
\bigr)
\times_{Y_{\mathrm{ff}}\setminus\Sigma_{\sing}}
\cdots
\times_{Y_{\mathrm{ff}}\setminus\Sigma_{\sing}}
\bigl(
M_{\mathrm{ff}}
\setminus(V\cup\phi^{-1}(\Sigma_{\sing}))
\bigr)
}_{r\ \text{factors}}.
\]

We define the (not necessarily effective) $\mathbb{R}$-divisor $B_2$ by
\begin{equation}\label{eq-B2-nef}
-(K_{M^{(r)}/Y}+\Delta_{M^{(r)}})+B_2
=
\sum_{i=1}^r
p_i^*\left(
-(K_{M/Y}+\Delta_M)+E'
\right).
\end{equation}
Then, the following properties hold:
\begin{itemize}
\item[$\bullet$]
The right-hand side of \eqref{eq-B2-nef}, and hence the left-hand
side, is nef.

\item[$\bullet$]
The support of $B_2$ satisfies
\[
\operatorname{Supp}(B_2)
\subset
E^{(r)}
\cup\Exc(\mu)
\cup\left(M^{(r)}\setminus\mu^{-1}(M^r_{\mathrm{ff}})\right)
\cup\operatorname{Supp}\left(\sum_{i=1}^rp_i^*V\right)
\cup(\phi^{(r)})^{-1}(\Sigma_{\sing}).
\]

\item[$\bullet$]
On the locus where $\phi\colon M\to Y$ is smooth on each factor and
$\mu\colon M^{(r)}\to M^r$ is an isomorphism, we have
\begin{equation}\label{eq-B2-smooth}
B_2=\sum_{i=1}^rp_i^*E'.
\end{equation}
\end{itemize}

We now verify these assertions.
The first follows directly from
\eqref{eq_rel-canonical-pi-exc}.
For the second, fix a point of $M^{(r)}$ outside the union appearing
on the right-hand side of the asserted inclusion.
At such a point, the desingularization $\mu\colon M^{(r)}\to M^r$ is an isomorphism, and the
point lies in the flat Gorenstein locus of the fiber product described
above.
Consequently, in a neighborhood of this point, the relative canonical
divisor satisfies
\[
K_{M^{(r)}/Y}
=
\sum_{i=1}^rp_i^*K_{M/Y}.
\]
Moreover, by the definition of $\Delta_{M^{(r)}}$, we have
\[
\Delta_{M^{(r)}}=\sum_{i=1}^rp_i^*\Delta_M.
\]
Since
$
\operatorname{Supp}(E')\subset E
$
by the construction of $E'$, we also have
$
\sum_{i=1}^rp_i^*E'=0
$
in a neighborhood of this point.
Substituting these equalities into \eqref{eq-B2-nef}, we conclude that $B_2$ vanishes in a neighborhood of
the point.
This proves the second assertion.

Finally, consider the locus where $\phi$ is smooth on each factor and
$\mu$ is an isomorphism.
On this locus, the fiber product is smooth over $Y$, and the relative
canonical bundle formula gives
$
K_{M^{(r)}/Y}=\sum_{i=1}^rp_i^*K_{M/Y}.
$
By definition, we also have
$
\Delta_{M^{(r)}}=\sum_{i=1}^rp_i^*\Delta_M.
$
Substituting these equalities into \eqref{eq-B2-nef} gives \eqref{eq-B2-smooth}.
This proves the third assertion.

\smallskip

\paragraph{\textbf{(C) Positivity of direct images.}}

In this step, we apply Proposition~\ref{prop-psef}\textup{(2)} to the
direct image sheaf
\begin{align*}
\phi^{(r)}_*
\mathcal{O}_{M^{(r)}}\Big(
mK_{M^{(r)}/Y}+L_0
+m\bigl(-(K_{M^{(r)}/Y}+\Delta_{M^{(r)}})+B_2\bigr)
+m\Delta_{M^{(r)}}
+\lceil mB_2\rceil-mB_2
\Big)
\end{align*}
for $m\gg 1$.
Here $m$ denotes the auxiliary integer occurring in
Proposition~\ref{prop-psef}, rather than the integer appearing in
Proposition~\ref{prop_more-psef-det}.
The assumptions of Proposition~\ref{prop-psef}\textup{(2)} are
satisfied.
Indeed, Step~\textup{(A)} shows that $L_0$ is pseudo-effective and
relatively $\phi^{(r)}$-big, while \eqref{eq-B2-nef} shows that
\[
-(K_{M^{(r)}/Y}+\Delta_{M^{(r)}})+B_2
\]
is nef.
It remains to verify the klt condition along a general fiber.
Over a general point of $Y$, the morphism $\phi$ is smooth, $\mu$ is
an isomorphism, and the support of $\Delta_{M^{(r)}}$ is SNC.
Since $(M,\Delta_M)$ is klt, every coefficient of
$\Delta_{M^{(r)}}$ along a general fiber satisfies the klt condition.
On the other hand, all the coefficients of
$\lceil mB_2\rceil-mB_2$ converge to zero.
Consequently, for $m\gg 1$, the divisor
\[
\Delta_{M^{(r)}}
+\frac{1}{m}\bigl(\lceil mB_2\rceil-mB_2\bigr)
\]
satisfies the klt condition along a general fiber.
Hence, Proposition~\ref{prop-psef}\textup{(2)} shows that
\[
\phi^{(r)}_*
\mathcal{O}_{M^{(r)}}(L_0+\lceil mB_2\rceil)
\]
is weakly positively curved on $Y$.

Let $C\subset Y_0$ be a Zariski-closed subset such that
$\phi^{(r)}$ is flat and
$
\phi^{(r)}_*
\mathcal{O}_{M^{(r)}}(L_0+\lceil mB_2\rceil)
$
is locally free over $Y_0\setminus C$.
We may assume that $\codim C\geq2$ and
$Y_0\setminus Y_{\mathrm{ff}}\subset C$.

Consider the evaluation morphism
\begin{equation}\label{eq-evaluation-map}
(\phi^{(r)})^*
\phi^{(r)}_*
\mathcal{O}_{M^{(r)}}(L_0+\lceil mB_2\rceil)
\longrightarrow
\mathcal{O}_{M^{(r)}}(L_0+\lceil mB_2\rceil)
\end{equation}
over $(\phi^{(r)})^{-1}(Y_0\setminus C)$.
In Step~\textup{(D)}, we will prove that this morphism is surjective at a
general point of $\bar D$.
Recall that $\bar D\subset M^{(r)}$ is the strict transform of the
diagonal subvariety $D\subset M^r$.
Assuming this assertion, we complete the proof in this step.

By Lemma~\ref{lem-pull}, the sheaf on the left-hand side of
\eqref{eq-evaluation-map} is weakly positively curved on
$(\phi^{(r)})^{-1}(Y_0\setminus C)$ with respect to a K\"ahler form
on $M^{(r)}$.
Let $g_\varepsilon$ denote the singular Hermitian metrics on
$\mathcal{O}_{M^{(r)}}(L_0+\lceil mB_2\rceil)$ induced by the
evaluation morphism.
The metrics on the direct image sheaf are nonsingular at a general
point of $Y$, and Step~\textup{(D)} shows that the evaluation morphism
is surjective at a general point of $\bar D$.
Consequently, the metric $g_\varepsilon$ is not identically singular along
$\bar D$, and its restriction to $\bar D$ is well-defined.
It follows that
\[
\mathcal{O}_{\bar D}(L_0+\lceil mB_2\rceil)
\]
is weakly positively curved on
$
(\phi^{(r)})^{-1}(Y_0\setminus C)\cap\bar D.
$

The diagonal $D\subset M^r$ is isomorphic to $M$.
Hence, there is an induced bimeromorphic morphism
$f\colon\bar D\to D\cong M$.
We claim the following assertions:
\begin{itemize}
\item[$\bullet$]
The complement of
\[
(\pi\circ f)
\bigl(
(\phi^{(r)})^{-1}(Y_0\setminus C)\cap\bar D
\bigr)
=
\pi\bigl(\phi^{-1}(Y_0\setminus C)\bigr)
\]
has codimension at least two.

\item[$\bullet$]
The divisor
\[
f^*(rcE)+B_1|_{\bar D}
+\lceil mB_2\rceil|_{\bar D}
\]
is $(\pi\circ f)$-exceptional.
\end{itemize}

We first complete the proof assuming these assertions.
By these assertions and reflexivity, we obtain
\[
\begin{aligned}
&(\pi\circ f)_{[*]}
\mathcal{O}_{\bar D}(L_0+\lceil mB_2\rceil)\\
\cong &
(\pi\circ f)_{[*]}
\mathcal{O}_{\bar D}\left(
f^*\bigl(rG+rcE-r\phi^*D_{c,1}\bigr)
+B_1|_{\bar D}
+\lceil mB_2\rceil|_{\bar D}
\right)\\
\cong &
\pi_{[*]}
\mathcal{O}_M\bigl(r(G-\phi^*D_{c,1})\bigr).
\end{aligned}
\]
Therefore, together with the support condition, 
Lemma~\ref{lem-push} implies that this sheaf is weakly
positively curved on $X$.

We now prove the first assertion.
The complement in question is contained in
\[
\pi\bigl(\phi^{-1}(Y\setminus Y_0)\bigr)
\cup
\pi\bigl(\phi^{-1}(C)\bigr).
\]
The first subset has codimension at least two by
Proposition~\ref{prop_bir-geometry-psi}\textup{(c)}, while the second
has codimension at least two since $\phi$ is flat over $Y_0$ and
$\codim C\geq2$.
This proves the first assertion.

We next prove the second assertion.
The divisor $f^*(rcE)$ is $(\pi\circ f)$-exceptional
since $E$ is $\pi$-exceptional.
Moreover, $B_1|_{\bar D}$ is $(\pi\circ f)$-exceptional by
\eqref{eq-codim-Yff} and the support condition on $B_1$.

It remains to consider $B_2|_{\bar D}$.
Suppose that $R$ is a divisorial component of $B_2|_{\bar D}$ that
is not $(\pi\circ f)$-exceptional.
Then $f(R)$ is a prime divisor on $M$ that is not $\pi$-exceptional,
and its generic point lies in $\phi^{-1}(Y_0)$ by
Proposition~\ref{prop_bir-geometry-psi}\textup{(c)}.
We claim that $\phi$ is smooth at the generic point of $f(R)$.
This follows from generic smoothness if $f(R)$ dominates $Y$.
Otherwise, the flatness of $\phi$ over $Y_0$ implies that
$\phi(f(R))$ is a prime divisor on $Y_0$.
By Proposition~\ref{prop_bir-geometry-psi}\textup{(e)}, the
multiplicity of $f(R)$ in $\phi^*\phi(f(R))$ is one since $f(R)$ is not $\pi$-exceptional.
This implies that $\phi$ is smooth at the generic
point of $f(R)$.
Hence $M^r$ is smooth at the corresponding point of the diagonal,
and $\mu$ is an isomorphism there.
Moreover, since $\operatorname{Supp}(E')\subset E$ and $f(R)$ is not
$\pi$-exceptional, the $\mathbb{Q}$-divisor $E'$ vanishes at the generic point of $f(R)$.
Since all the morphisms $p_i|_{\bar D}$ coincide with $f$,
\eqref{eq-B2-smooth} shows that $B_2$ vanishes at the generic point
of $R$, which is a contradiction.
Thus, every divisorial component of
$\lceil mB_2\rceil|_{\bar D}$ is $(\pi\circ f)$-exceptional.
This proves the second assertion.

\smallskip

\paragraph{\textbf{(D) Non-vanishing of the evaluation morphism along the diagonal.}}

We finally prove that \eqref{eq-evaluation-map} does not vanish at the
generic point of $\bar D$.
Choose a general point
$
y\in Y_0\setminus C
$
such that $\phi$ is smooth over $y$,
$y\notin\operatorname{Supp}(N_Y)$, and cohomology and base change hold
at $y$ for
$
\phi^{(r)}_*
\mathcal{O}_{M^{(r)}}(L_0+\lceil mB_2\rceil).
$
Since $Y_0\setminus C\subset Y_{\mathrm{ff}}$ and $\phi$ is smooth
over $y$, the fiber product $M^r$ is smooth along $M_y^r$.
The morphism $\mu$ is an isomorphism over this fiber, and
hence
$
M^{(r)}_y\cong M_y^r.
$
Moreover, the support condition on $B_1$ gives
$B_1|_{M_y^r}=0$.
By the definition of $E'$ and the choice of
$y\notin\operatorname{Supp}(N_Y)$, we have
$
E'|_{M_y}=F_+|_{M_y}.
$
Consequently, \eqref{eq-B2-smooth} gives
\[
B_2|_{M_y^r}
=
\sum_{i=1}^r\pr_i^*(E'|_{M_y})
=
\sum_{i=1}^r\pr_i^*(F_+|_{M_y})\geq0.
\]

Since $r=\operatorname{rank}\phi_*\mathcal{O}_M(G+cE)>0$ and
cohomology and base change hold at $y$, we may choose a nonzero
section
\[
u\in H^0(M_y,\mathcal{O}_{M_y}(G+cE)).
\]
Choose a general point $x\in M_y$ at which $u$ does not vanish and
which does not belong to $\operatorname{Supp}(F_+|_{M_y})$.
After choosing a nonzero element of
$\mathcal{O}_Y(-rD_{c,1})\otimes\mathbb{C}(y)$, the product of
$u^{\boxtimes r}$ and the canonical section of
$\lceil mB_2\rceil|_{M_y^r}$ gives a section of
\[
\mathcal{O}_{M_y^r}
\bigl((L_0+\lceil mB_2\rceil)|_{M_y^r}\bigr)
\]
that does not vanish at $(x,\ldots,x)$.
The base-change isomorphism
\[
\phi^{(r)}_*
\mathcal{O}_{M^{(r)}}(L_0+\lceil mB_2\rceil)
\otimes\mathbb{C}(y)
\longrightarrow
H^0\left(
M_y^r,
\mathcal{O}_{M_y^r}
\bigl((L_0+\lceil mB_2\rceil)|_{M_y^r}\bigr)
\right)
\]
shows that this section belongs to the image of
\eqref{eq-evaluation-map}.
Hence the evaluation morphism is nonzero at $(x,\ldots,x)$.
Therefore, the restriction to $\bar D$ is not identically zero, so it
does not vanish at the generic point of $\bar D$.
Since its target has rank one, it is surjective at a general point of
$\bar D$.
\end{proof}

We now study the effect of adding the $\pi$-exceptional divisor $E$ on the direct image
$\phi_* \mathcal{O}_M(G)$ by means of Propositions~\ref{prop_det-stablize} and~\ref{prop_iso-direct-image-Y0}.

\begin{prop}

\label{prop_det-stablize}
Let everything be as above.
Then, there exists an integer $c_0\in\ZZ_{+}$ such that, for every $c\geqslant c_0$, the natural inclusion
\[
\det \phi_* \mathcal{O}_M(G + cE) \to
\det \phi_* \mathcal{O}_M(G + (c+1)E)
\]
is an isomorphism over $Y_0$.
\end{prop}

\begin{rem}\label{rem-det-stablize}
The proof in the projective setting relies on the existence of an ample line bundle on $Y$, which is used to apply an extension theorem.
Therefore, we carefully localize the argument to the K\"ahler setting by using pseudo-effectivity and intersection numbers.
\end{rem}

\begin{proof}
Recall that
\[
r := \rank \phi_* \mathcal{O}_M(G + cE)
= \dimcoh^0(M_y, \mathcal{O}_{M_y}(G + cE))
\]
is independent of $c$ by \eqref{eq-rank}.
Assume, for contradiction, that the natural inclusion
\[
\det \phi_* \mathcal{O}_M(G + cE)
\to
\det \phi_* \mathcal{O}_M(G + (c+1)E)
\]
fails to be an isomorphism over $Y_0$ for infinitely many positive integers $c>0$.
Then, there exist a strictly increasing sequence $\{c_k\}_{k=1}^\infty$ and 
effective $\mathbb{Z}$-divisors $B_k$ on $Y$ such that, for every $k$, the intersection $B_k\cap Y_0$ is nonempty and the line bundle
\[
r D_{c_{k+1}, 1} - (r D_{c_k, 1} + B_k)
\]
admits a nonzero section.
By Proposition~\ref{prop_more-psef-det}, for every $N\in\ZZ_{+}$, we obtain
\[
c_{1}(rG + d E)
\geq c_{1}(r\phi^{*}D_{c_N, 1})
\geq c_{1}(r\phi^{*}D_{c_1, 1}) + \sum_{k=1}^N c_{1}(\phi^{*}B_k).
\]

Here and below, for Bott--Chern classes $\alpha$ and $\beta$, we write
$\alpha\geq\beta$ if the class $\alpha-\beta$ is pseudo-effective.

On the other hand, by Lemma~\ref{lem-Qfact}, we can write
\begin{align}\label{eq-dec2}
\phi^{*} D_{c,m}
\sim_{\mathbb{Q}} \pi^{*} \Lambda_{c,m} + F_{c,m},
\quad\text{and}\quad
\phi^{*} B_k
\sim_{\QQ} \pi^{*} \Lambda_k + F_k,
\end{align}
where $F_{c,m}$ and $F_k$ are $\pi$-exceptional divisors.

Note that $\Lambda_k$ is an effective Weil $\ZZ$-divisor on $X$.
Indeed, as in the proof of Lemma~\ref{lem-Qfact}, 
the divisor $\Lambda_k$ can be  defined
as the pushforward $\Lambda_k:=\pi_*(\phi^*B_k)$ as Weil divisors. 
Moreover, the strong $\mathbb{Q}$-factoriality of $X$ implies that
$\Lambda_k$ is $\mathbb{Q}$-Cartier.

Since $G=\pi^{*}A+F$ by Lemma~\ref{lemm-good}\textup{(1)}, 
after pushing forward by $\pi$ and taking intersection numbers, we obtain
\[
c_{1}(rA) \cdot \{\omega^{n-1}\}
\geq
\left(c_{1}(r\Lambda_{c_{1},1})+\sum_{k=1}^N c_{1}(\Lambda_k)\right)\cdot 
\{\omega^{n-1}\}
\]
for every $N\in\ZZ_{+}$, where $\omega$ is a K\"ahler form on $X$.
The left-hand side is independent of $N$.
Each current $[\Lambda_k]$ is semipositive 
and the sequence of partial sums is bounded. 
Hence, the series
$\sum_{k=1}^{\infty}[\Lambda_k]\cdot \{\omega^{n-1}\}$ converges.
In particular, we see that 
$c_{1}(\Lambda_k)\cdot\{\omega^{n-1}\}\to 0$ as $k\to\infty$.

After passing to a subsequence, we may assume that the integration currents
$[\Lambda_k]$ converge weakly to a semipositive $(1,1)$-current $T$.
The above argument shows that 
$$
T\cdot\omega^{n-1}
=\lim_{k\to\infty}[\Lambda_k]\cdot\omega^{n-1}=0,$$
and hence $T=0$.
Since each $\Lambda_k$ is an effective Weil $\ZZ$-divisor, 
by Lemma~\ref{lemm-uniform}, there is a uniform lower bound for
$c_{1}(\Lambda_k)\cdot\omega^{n-1}$ whenever $\Lambda_k\neq 0$.
Therefore, we can conclude that $\Lambda_k=0$ for all sufficiently large $k$.
Consequently, for all sufficiently large $k$, we have
$\phi^{*}B_k\sim_{\mathbb{Q}}F_k$.
Since $F_k$ is $\pi$-exceptional, we have 
$$
\int_M c_1(\phi^{*}B_k) \wedge \pi^*\omega^{n-1}
=
\int_M c_1(F_k) \wedge \pi^*\omega^{n-1}
=
0.
$$
This contradicts the fact that $B_k\cap Y_0$ is nonempty.
\end{proof}

\begin{lemm}\label{lemm-uniform}
Let $X$ be a normal compact analytic variety, and let $\omega$ be a Hermitian form on $X$.
Then, there exists a constant $\delta>0$, depending only on $(X,\omega)$, such that
\[
\int_Z \omega^{\dim Z}\geq \delta
\]
for every subvariety $Z\subset X$.
\end{lemm}

\begin{rem}
The lemma is immediate when $X$ is projective.
Indeed, in this case, we may take $\omega$ to be a K\"ahler form in the first Chern class of a very ample line bundle.
Then, the above integral is the number of intersection points in
$Z \cap H_{1} \cap \cdots \cap H_{p}$, which is at least one, 
where $p:=\dim Z$ and $\{H_i\}_{i=1}^p$ are general hypersurfaces.
In the non-projective case, we use the monotonicity in the Lelong--Jensen formula, 
which also naturally appears in the definition of the Lelong number.
\end{rem}

\begin{proof}
Set $n:=\dim X$.
We first choose finitely many local charts with uniform geometric estimates.
Since $X$ is compact, there exist finitely many open subsets
\[
W_i\Subset V_i\Subset U_i\Subset X
\qquad (1\leq i\leq N)
\]
such that $\{W_i\}_{i=1}^{N}$ is an open cover of $X$ and each $U_i$ admits a closed holomorphic embedding
$\iota_i\colon U_i\hookrightarrow\Omega_i\subset\mathbb{C}^{N_i}$.
We may assume that there exists a constant $r_i>0$ such that  
\[
B(\iota_i(x),r_i)\cap\iota_i(U_i)\subset\iota_i(V_i)
\]
for every $x\in W_i$. 
Let $\beta_i$ be a sufficiently small positive multiple of the pullback via
$\iota_i$ of the standard K\"ahler form on $\mathbb{C}^{N_i}$.
We may assume that $\omega|_{V_i}\geq\beta_i$.

Let $Z\subset X$ be a subvariety of dimension $p$.
Choose a point $x\in Z$ and an index $i$ such that $x\in W_i$.
The Lelong--Jensen formula shows that
\[
r^{-2p}
\int_{\iota_i(Z\cap U_i)\cap B(\iota_i(x),r)}\beta_i^p
\]
is nondecreasing in $r$, and that its limit as $r\to 0$ is equal 
to the Lelong number of the integration current $[Z]$ at $x$, up to a positive constant $a_p$ depending only on $p$.
Since $\multiplicity_x Z\geq 1$, we obtain
\[
\int_{\iota_i(Z\cap U_i)\cap B(\iota_i(x),r_i)}\beta_i^p
\geq a_p r_i^{2p}.
\]
Using $\omega|_{V_i}\geq \beta_i$, we conclude that
\[
\int_Z\omega^p
\geq
\int_{Z\cap V_i}\omega^p
\geq
\int_{\iota_i(Z\cap U_i)\cap B(\iota_i(x),r_i)}\beta_i^p
\geq
 a_p r_i^{2p}.
\]
The constant 
\[
\delta_p:=
\min_{1\leq i\leq N}
\left\{ a_p r_i^{2p}\right\},
\]
satisfies the desired conclusion for subvarieties of dimension $p$.
\end{proof}

\begin{prop}
\label{prop_iso-direct-image-Y0}
Let everything be as above.
Let $c_0$ be the integer given by Proposition~\ref{prop_det-stablize}, and let $c \geqslant c_0$.
Then, the following assertions hold:
\begin{itemize}
\item[\rm(a)]
The sheaf $\phi_* \mathcal{O}_M(G + cE)$ is isomorphic to
$\phi_* \mathcal{O}_M(G + kE + \lceil pE'\rceil)$ over $Y_0$ for every
$k \geqslant c$ and every $p \in \ZZ_{\geqslant 0}$;

\item[\rm(b)]
The sheaf $\phi_* \mathcal{O}_M(G + cE)$ is
$(1/m) D_{c,m}$-weakly positively curved over $Y_0$ for every
$m\in\ZZ_{+}$.
\end{itemize}
\end{prop}

\begin{proof}

Recall that $E' := F_{+} - \phi^{*} N_Y$ by \eqref{eq-e'}.
The divisor $F_{+}$ is an effective $\pi$-exceptional divisor by \eqref{eq-gpair}, 
whereas $N_Y$ is effective and supported in $Y \setminus Y_{0}$ by Proposition~\ref{prop_bir-geometry-psi}\textup{(a)} and the definition of $Y_{0}$.
The latter property yields
\[
\phi_* \mathcal{O}_M(G + kE + \lceil pE'\rceil)
\cong
\phi_* \mathcal{O}_M(G + kE + \lceil pF_{+}\rceil)
\]
over $Y_0$.
Thus, by \eqref{eq-rank}, we see that
$\rank\!\phi_*\scrO_M(G+kE+\lceil pE'\rceil)=r_1$ 
and by Proposition~\ref{prop_det-stablize} , we see 
that the natural injection over $Y_{0}$
\[
\scrO_Y(rD_{c,1})
\cong
\det\!\phi_*\scrO_M(G+cE)
\hookrightarrow
\det\!\phi_*\scrO_M(G+kE+\lceil pE'\rceil)
\]
is an isomorphism.
Therefore, by \cite[Lemma~1.20]{DPS94}, the natural inclusion on $Y$
\begin{equation}
\label{eq_incl-direct-image}
\phi_*\scrO_M(G+cE)
\hookrightarrow
\phi_*\scrO_M(G+kE+\lceil pE'\rceil)
\end{equation}
is an isomorphism over the locally free locus of
$\phi_*\scrO_M(G+kE+\lceil pE'\rceil)\big|_{Y_0}$.
Since $\phi$ is flat over $Y_0$, both
$\phi_*\scrO_M(G+cE)$ and
$\phi_*\scrO_M(G+kE+\lceil pE'\rceil)$ are reflexive over $Y_0$.
Hence, the morphism \eqref{eq_incl-direct-image} is an isomorphism over $Y_0$.
This proves \textup{(a)}.

\smallskip

We now prove \textup{(b)}.
To apply Proposition~\ref{prop-psef}, we verify the assumptions for
$G+dE+\lceil pE'\rceil$ for sufficiently large integers $p \gg d \gg 1$.

By Proposition \ref{prop_more-psef-det}, we can take $d \gg 1$ such that 
$$
G-\frac{1}{m}\phi^\ast\!D_{c,m}+dE
$$ is pseudo-effective  on $M$.

We consider 
\begin{align*}
&G+dE+\lceil pE'\rceil \\
=&p(K_{M/Y}+\Delta_M)+(G+dE)+\lceil pE'\rceil-pE' 
+\bigl(-p(K_{M/Y}+\Delta_M)+pE'\bigr).
\end{align*}
The line bundle $G+dE$ is $\phi$-big, the divisor
$-(K_{M/Y}+\Delta_M)+E'$ is nef, and 
$L:=\Delta_M+(1/p)(\lceil pE'\rceil-pE')$ 
satisfies the klt condition for $p\gg 1$.
Moreover, since $G+dE$ is
$(1/m)\phi^{*}D_{c,m}$-weakly positively curved, 
Proposition~\ref{prop-psef} (2) shows that the direct image sheaf 
$$\phi_*\scrO_M(G+dE+\lceil pE'\rceil)$$ is
$(1/m)D_{c,m}$-weakly positively curved.
Combining this with \textup{(a)}, we conclude that
$\phi_*\scrO_M(G+cE)$ is
$(1/m)D_{c,m}$-weakly positively curved on $Y_0$.
This proves \textup{(b)}.
\end{proof}

We complete the proof of Theorem~\ref{thm-faked-flat} by proving the following two propositions.
The initial construction of $G$ (in particular, property (3)) plays a crucial role in this step (see Lemma~\ref{lemm-good}).

\begin{prop}
\label{prop_positivity-Um-Vm}

For $1 \ll m\in\ZZ_{+}$ divisible by $r$, we set
\begin{align*}
\calU_{c,m}
&:=\Sym^m\!\phi_*\mathcal{O}_M(G+cE)
\otimes\mathcal{O}_Y(-mD_{c,1}), \\
\calV_{c,m}
&:=\phi_*\mathcal{O}_M(m(G+cE))
\otimes\mathcal{O}_Y(-mD_{c,1})
=\phi_*\mathcal{O}_M(L_m).
\end{align*}
Then both $\calU_{c,m}$ and $\calV_{c,m}$ are weakly positively curved on $Y_0$.
In particular, the sheaf
$\mathcal{E}_{m}=\pi_{[*]}(\phi^{*}\calV_{c,m})$
is weakly positively curved on $X$.
\end{prop}

\begin{proof}
By Proposition~\ref{prop_iso-direct-image-Y0} \textup{(b)},
the sheaf $\phi_*\mathcal{O}_M(G+cE)$ is
$D_{c,1}$-weakly positively curved on $Y_0$.
This shows that
\[
\phi_*\mathcal{O}_M(G+cE)
\otimes\mathcal{O}_Y(-D_{c,1})
\]
is weakly positively curved on $Y_0$,
and hence so is its symmetric power $\calU_{c,m}$.
 
By Lemma~\ref{lemm-good}\textup{(3)}, the natural morphism
$\calU_{c,m}\to\calV_{c,m}$ is generically surjective.
Hence, the sheaf $\calV_{c,m}$ is also weakly positively curved on $Y_0$ 
with respect to a K\"ahler form on $Y$.

The assumptions of Lemma~\ref{lem-pull}\textup{(2)} are satisfied by
Proposition~\ref{prop_bir-geometry-psi}\textup{(c)}.
Hence, the latter conclusion follows by applying
Lemma~\ref{lem-pull}\textup{(2)}.
\end{proof}

It remains to prove that $c_{1}(\mathcal{E}_{m})=0$.
Let $C\subset Y$ be a Zariski-closed subset of codimension at least two such that
$\calV_{c,m}$ is locally free on $Y\setminus C$.
Then, on $M\setminus\phi^{-1}(C)$, we have
\[
\det\phi^{*}\calV_{c,m}
=
\phi^{*}\det\calV_{c,m}
=
\phi^{*}\mathcal{O}_{Y}\bigl(r_m(D_{c,m}-mD_{c,1})\bigr).
\]
Since $\phi\colon M\to Y$ is flat over $Y_0$, the inverse image of
$C\cap Y_0$ has codimension at least two in $\phi^{-1}(Y_0)$.
On the other hand, we have 
$\codim \pi(\phi^{-1}(Y\setminus Y_0))\geq 2$
by Proposition~\ref{prop_bir-geometry-psi} \textup{(c)}.
Therefore, by reflexivity, we obtain
\[
\det\pi_{[*]}\phi^{*}\calV_{c,m}
=
\pi_{[*]}\phi^{*}
\mathcal{O}_{Y}\bigl(r_m(D_{c,m}-mD_{c,1})\bigr)
\sim_{\bb{Q}}
\mathcal{O}_{X}\bigl(r_m(\Lambda_{c,m}-m\Lambda_{c,1})\bigr),
\]
where $\Lambda_{c,m}$ is defined as in \eqref{eq-dec2}.

Hence, to complete the proof of Theorem~\ref{thm-faked-flat}, it remains to prove the following proposition.

\begin{prop}[{\cite[Proposition~3.6]{CH19}}]
\label{prop_det-direct-image-power}
For every $m\in\ZZ_{+}$ divisible by $r:=r_1$, we have
\[
\Lambda_{c,m}\equiv m\Lambda_{c,1}.
\]
\end{prop}

\begin{proof}
Since
$\mathcal{E}_{m}=\pi_{[*]}(\phi^{*}\calV_{c,m})$
is weakly positively curved, so is its determinant. Moreover, we have 
\[
\det\pi_{[*]}(\phi^{*}\calV_{c,m})
=
\pi_{[*]}\phi^{*}
\mathcal{O}_Y\bigl(r_m(D_{c,m}-mD_{c,1})\bigr)
\sim_{\bb{Q}}
\mathcal{O}_X\bigl(r_m(\Lambda_{c,m}-m\Lambda_{c,1})\bigr).
\]
It follows that the class
$\Lambda_{c,m}-m\Lambda_{c,1}$ is pseudo-effective.

On the other hand, since
$\phi_*\mathcal{O}_M(G+cE)$ is
$(1/m)D_{c,m}$-weakly positively curved over $Y_0$
by Proposition~\ref{prop_iso-direct-image-Y0}\textup{(b)}, the determinant sheaf 
\[
\det\phi_*\mathcal{O}_M(G+cE)-\frac{r}{m}D_{c,m}
\]
is weakly positively curved on $Y_0$.
Thus, applying Lemma~\ref{lem-pull} (2) again, we conclude that
\[
\pi_{[*]}\phi^{*}
\mathcal{O}_Y(mD_{c,1}-D_{c,m})
\sim_{\bb{Q}}
\mathcal{O}_X(m\Lambda_{c,1}-\Lambda_{c,m})
\]
is weakly positively curved on $X$.
Hence, the class $m\Lambda_{c,1}-\Lambda_{c,m}$ is also pseudo-effective.
Therefore, both
$\Lambda_{c,m}-m\Lambda_{c,1}$ and its dual are pseudo-effective, 
showing that $\Lambda_{c,m}\equiv m\Lambda_{c,1}$.
\end{proof}

We conclude this subsection by deducing Corollary~\ref{thm-flat} from Theorem~\ref{thm-faked-flat}.

\begin{proof}[Proof of Corollary~\ref{thm-flat}]
The proof is the same as that of \cite[Theorem~3.2]{MW}.
For the reader's convenience, we provide the complete details.
For simplicity of notation, we omit the subscript $m$ from $\mathcal{V}_{m}$ and $\mathcal{E}_{m}$ throughout the proof. 

Using the slightly different notation introduced in Lemma~\ref{lemm-good}, we consider the following diagram:
\begin{align}\label{graph:factorisation-diagram}
\xymatrix@C=4.5em@R=3.2em{
M
  \ar@/^2.6pc/[drrr]^{\pi}
  \ar[dr]^{\gamma}
  \ar@/_2.4pc/[ddr]_{\phi}
&&& \\
& \Gamma
    \ar[rr]^{\mu}
    \ar[d]^{\varphi}
&& X
   \ar@{-->}@/^1.8pc/[dll]^{\psi} \\
& Y. && 
}
\end{align}
Here $\Gamma$ is the normalization of the graph of $\psi$, equipped with the natural morphisms
$\varphi \colon \Gamma \to Y$ and $\mu \colon \Gamma \to X$, and $M$ is a log resolution of the generalized pair as in Setting~\ref{setting}.
As in the proof of Lemma~\ref{lemm-good}, although the $\pi$-exceptional locus
$E=\Exc(\pi)$ may dominate $Y$, the $\mu$-exceptional locus
$E_\mu:=\Exc(\mu)$ does not dominate $Y$.
This is an advantage of working with $\Gamma$.

Let $B \subset Y$ be the union of the non-locally-free locus of $\mathcal{V}$ and the non-flat locus of
$\varphi \colon \Gamma \to Y$.
Note that $B$ has codimension at least two.
We first verify the following claim.

\begin{claim}\label{claim2}
Let everything be as above. 
Then, we have 
\[
\pi^*\mathcal{E}
=
\phi^{*}\mathcal{V}
\quad\text{on}\quad
M \setminus \bigl(\gamma^{-1}(E_\mu) \cup \phi^{-1}(B)\bigr).
\]
\end{claim}

\begin{proof}
The sheaf $\varphi^*\mathcal{V}$ is locally free on
$\Gamma \setminus \varphi^{-1}(B)$.
Hence, applying the projection formula to the restriction
\[
\gamma \colon M \setminus \phi^{-1}(B)
\to
\Gamma \setminus \varphi^{-1}(B),
\]
we obtain
\[
\gamma_*\gamma^*\varphi^*\mathcal{V}
=
\varphi^*\mathcal{V}
\quad\text{on}\quad
\Gamma \setminus \varphi^{-1}(B).
\]
Using the diagram~\eqref{graph:factorisation-diagram} and the isomorphism
$\mu\colon \Gamma\setminus E_\mu\cong X\setminus\mu(E_\mu)$, we obtain
\[
\mu^*\pi_*\phi^*\mathcal{V}
=
\mu^*\mu_*\gamma_*\gamma^*\varphi^*\mathcal{V}
=
\mu^*\mu_*\varphi^*\mathcal{V}
=
\varphi^*\mathcal{V}
\]
on $\Gamma \setminus (E_\mu \cup \varphi^{-1}(B))$.
This implies that 
\begin{align*}
\pi^*\mathcal{E}
=
\gamma^*\mu^*\bigl((\pi_*\phi^*\mathcal{V})^{**}\bigr) 
=
\gamma^*\bigl((\mu^*\pi_*\phi^*\mathcal{V})^{**}\bigr) 
=
\phi^*\mathcal{V}
\end{align*}
on $M \setminus (\gamma^{-1}(E_\mu) \cup \phi^{-1}(B))$.
Here, we used the definition of $\mathcal{E}$ and the reflexivity of
$\varphi^*\mathcal{V}$ on $\Gamma \setminus \varphi^{-1}(B)$.
\end{proof}

Let us return to the proof of Corollary~\ref{thm-flat}.
Under the additional assumption that $X$ is maximally quasi-\'etale, it follows from \cite{CDM26} that $\mathcal{E}$ is a numerically flat locally free sheaf. 
The pullback $\pi^*\mathcal{E}$ is also a numerically flat locally free sheaf on $M$,
and hence $\pi^*\mathcal{E}$ admits a flat connection $\nabla$ by
\cite[Corollary~3.10 and the subsequent discussion]{Sim92} and
\cite[Theorem~1.18]{DPS94}.

By the definition of $B$ and  $\codim B\geqslant 2$, 
it is enough to show that there is a Zariski-closed subset $C\subset Y_0$ of codimension $\geq 2$ 
such that $\mathcal{V}$ admits a flat connection on $Y_0\setminus(B\cup C)$. 
For simplicity of notation, after replacing $Y$ with $Y\setminus B$, 
we may assume that $B$ is empty, 
that is, $\mathcal{V}$ is locally free and
$\varphi\colon\Gamma\to Y$ is flat.
Then, The base $Y$ is  noncompact, but the compactness will not be used in the argument below.

We construct the desired flat connection on $\mathcal{V}$ over $Y_0\setminus C$ according to the following strategy:
\begin{itemize}
\item[(1)]
Using the equality $\pi^*\mathcal{E}=\phi^*\mathcal{V}$ away from
$\gamma^{-1}(E_\mu)$, we show that the connection $\nabla$ on
$\pi^*\mathcal{E}$ descends to a flat connection on $\mathcal{V}$ over a Zariski-open subset of $Y$.

\item[(2)]
We show that the flat connection on $\mathcal{V}$ constructed in \textup{(1)} extends over $Y_0\setminus C$ for some Zariski-closed subset
$C\subset Y_0$ of codimension $\geq 2$.
\end{itemize}

\smallskip

\noindent
\textup{(1)}
We first work in a neighborhood of a given point of $Y$ in the analytic topology.
Thus, after replacing $Y$ with a sufficiently small open ball, we fix a local frame of $\mathcal{V}$ on $Y$.
By Claim~\ref{claim2}, recalling that $B$ is empty, 
this frame induces a frame of $\pi^*\mathcal{E}$ on
$M\setminus\gamma^{-1}(E_\mu)$.
With respect to this frame, the flat connection $\nabla$ on
$\pi^*\mathcal{E}$ can be written as $$\nabla=d+\Xi,$$ where $d$ is the exterior derivative and $\Xi$ is a $(1,0)$-form with values in
$\End(\pi^*\mathcal{E})$.
We regard $\Xi$ as a matrix of $(1,0)$-forms with respect to the chosen frame.
Since $\nabla$ is flat, the matrix $\Xi$ satisfies
$$d\Xi+\Xi\wedge\Xi=0.$$
The $(1,1)$-part of the left-hand side is $\dbar\Xi$, and hence
$\dbar\Xi=0$.
This implies that $\Xi$ is a matrix of holomorphic one-forms on
$M\setminus\gamma^{-1}(E_\mu)$.

We claim that $\Xi$ descends to a holomorphic $(1,0)$-form $\Xi'$ with values in $\End(\mathcal{V})$ on $Y_1$, where 
\begin{equation*}
Y_1:=
\left\{
y\in Y\setminus\varphi(E_\mu)
\ \middle|\
M_y:=\phi^{-1}(y)\text{ is smooth and rationally connected}
\right\}.
\end{equation*}
Since $\psi\colon X\dashrightarrow Y$ is an MRC fibration, the fiber
$X_y:=\psi^{-1}(y)$ at a general point $y\in Y$ is rationally connected, and so is $M_y$.
Thus, the subset $Y_1$ is a nonempty Zariski-open subset of $Y$.
Although it is possible that $\varphi(E)=Y$, we have
$\varphi(E_\mu)$ is properly contained in $Y$; 
this is the reason for considering $E_\mu$ rather than $E$.

For every $y\in Y_1$, the fiber $M_y$ has no nonzero holomorphic one-forms.
Hence, the restrictions to $M_y$ of the entries of $\Xi$ vanish.
Moreover, by Claim~\ref*{claim2}, we have
$\End(\pi^*\mathcal{E})=\phi^*\End(\mathcal{V})$ over 
$Y\setminus\varphi(E_\mu)$.
Therefore, there exists a holomorphic section
$\Xi'\in H^0(Y_1,\Omega^1_Y\otimes\End(\mathcal{V}))$ such that
$\Xi=\phi^*\Xi'$.

\smallskip

\noindent
\textup{(2)}
The complement $Y \setminus Y_{1}$ may contain a divisorial component, 
so we need to show that $\Xi'$ extends to a section over $Y_0\setminus C$.
To this end, let $P\subset Y_0$ be a divisorial component of $Y\setminus Y_1$.
It is enough to extend $\Xi'$ across a general point $y\in P$.
By the definition of $Y_0$, the divisor $\phi^*P$ is not $\pi$-exceptional.
Hence, by Conclusion~\textup{(e)} of
Proposition~\ref{prop_bir-geometry-psi},
there exists a reduced component $P_0$ of the divisor $\phi^*P$ 
such that $P_{0} $ is not $\pi$-exceptional. 
In particular, we also see that 
$$P_0\not\subset\gamma^{-1}(E_\mu),$$ since
$\gamma^{-1}(E_\mu)$ is $\pi$-exceptional.
Moreover, the induced morphism
$\phi|_{P_0}\colon P_0\to P$ is surjective since 
$\phi\colon M\to Y$ is flat over $Y_0$ and $P\subset Y_0$.
Therefore, for a general point $y\in P$, there exists a point
$x\in P_0\cap M_y$ such that 
$$\text{
$x\notin\gamma^{-1}(E_\mu)$ \quad and \quad $x$ is a smooth point of $\phi \colon M \to Y$.
}
$$

The first property implies that all entries of the matrix
$\Xi=\phi^*\Xi'$ are bounded in a neighborhood of $x$, 
since $\Xi$ is defined not only on $\phi^{-1}(Y_1)$ but also on
$M\setminus\gamma^{-1}(E_\mu)$.
The second property provides a local holomorphic section $s \colon Y \to M$ of
$\phi\colon M\to Y$ passing through $x$ over a neighborhood of $y$.
Pulling back the equality $\Xi=\phi^*\Xi'$ by this local section $s$, 
all entries of $\Xi'$ are bounded near $y$.
Therefore, the Riemann extension theorem implies that $\Xi'$ extends across a general point of $P$.

The locally extended connections on $\mathcal{V}$ satisfy the gluing condition and remain flat,
since these properties hold on the dense open subset $Y_1$.
Consequently, they glue to a flat connection on $\mathcal{V}$ over
$Y_0\setminus C$.
Since $C$ has codimension at least two, the associated local system extends
uniquely across $C$. Its holomorphic vector bundle agrees with $\mathcal{V}$
by reflexivity. Thus, the sheaf $\mathcal{V}|_{Y_0}$ is locally free  and admits a
flat connection.
\end{proof}

\section{The structure theorem in the strongly \texorpdfstring{$\mathbb{Q}$}{text}-factorial case}
\label{sec-mainthm}

The purpose of Sections~\ref{sec-mainthm} and~\ref{sec-red}
is to complete the proof of Theorem~\ref{thm-main}.
Before giving the details, we briefly outline the strategy of the proof.

Let $(X,\Delta,\boldsymbol{\beta})$ be a generalized klt pair of Calabi--Yau type.
By Theorem~\ref{thm-maxqetale}, there exists a maximally quasi-\'etale cover
$X'\to X$ such that $X'$ admits the structure of a generalized klt pair
of Calabi--Yau type.
Thus, after replacing $X$ with $X'$, we may assume that $X$ is maximally
quasi-\'etale.
By Theorem~\ref{thm-Q-factorialization}, there then exists a small
$\mathbb{Q}$-factorialization
$q\colon X^{\qf}\to X$
such that $X^{\qf}$ also admits the structure of a generalized klt pair
of Calabi--Yau type.
A key point is that $X^{\qf}$ remains maximally quasi-\'etale by
Theorem~\ref{thm-maxqetale}.
Note that a maximally quasi-\'etale cover of a strongly
$\mathbb{Q}$-factorial variety need not be strongly
$\mathbb{Q}$-factorial.
Therefore, it is important to first take a maximally quasi-\'etale cover
and then a $\mathbb{Q}$-factorialization.

The proof of Theorem~\ref{thm-main} proceeds as follows:
\begin{itemize}
\item[$\bullet$]
In Section~\ref{sec-mainthm}, we prove Theorem~\ref{thm-main}
under the additional technical assumptions that $X$ is strongly
$\mathbb{Q}$-factorial and maximally quasi-\'etale.
More precisely, we prove that $X^{\qf}$ admits the desired fibration
$X^{\qf}\to Y^{\qf}$.

\item[$\bullet$]
In Section~\ref{sec-red}, we prove that the fibration
$X^{\qf}\to Y^{\qf}$ constructed above descends to a fibration
$X\to Y$ with the desired properties.
\end{itemize}

\subsection{Splitting of tangent sheaves}\label{subsec-splitting}

This subsection is devoted to the proof of the following splitting theorem.

\begin{thm}
\label{thm_splitting}
Let $(X,\Delta,\boldsymbol{\beta})$ be a K\"ahler generalized klt pair
of Calabi--Yau type, as in Theorem~\ref{thm-main}.
Assume further that $X$ is strongly $\mathbb{Q}$-factorial and
maximally quasi-\'etale.
Then, the MRC fibration of $X$ induces a decomposition of the tangent
sheaf $T_X$ 
\[
T_X\cong\calF\oplus\calG
\]
such that
\begin{itemize}
\item[$\bullet$]
$\calF$ is an algebraically integrable foliation, and its general leaf
is a fiber of the MRC fibration of $X$;

\item[$\bullet$]
$\calG$ is a foliation whose canonical sheaf $K_{\calG}:=\det (\mathcal{G}^{*})$ satisfies $K_{\calG}\sim_{\QQ}0$. 
\end{itemize}
\end{thm}

\begin{proof}
We use the notation of Setting~\ref{setting}. 
The reflexive sheaf $\mathcal{E}_m$ defined in Theorem~\ref{thm-faked-flat} is locally free on $X$ by Corollary~\ref{thm-flat}, 
since $X$ is assumed to be maximally quasi-\'etale. 

The proof essentially follows the arguments in \cite[\S 3.C, Proof of Theorem~1.2]{CH19}. 
For the reader's convenience, we recall the proof here with some simplifications and clarifications. 
The main idea is to find a bimeromorphic model $X_G$ of $X$ that admits a locally constant fibration over $Y_0$. 
This induces a splitting of the tangent sheaf of $X_G$ over $Y_0$. 
If the bimeromorphic model $X_G$ is carefully chosen, this splitting descends to a splitting of $T_X$. 
We divide the proof into four steps:

\setcounter{step}{0}

\begin{step}[Construction of the bimeromorphic model $X_G$]
Let everything be as in Setting~\ref{setting}.
By \cite[Theorem~1.2]{CH24}, the morphism
$\phi\colon M\to Y$ is projective.
Hence, we may choose a $\phi$-big line bundle $G$ as in
Lemma~\ref{lemm-good}.
Let $c\geqslant c_0$ be a positive integer as in
Proposition~\ref{prop_det-stablize}.

After possibly blowing up $M$, we may assume that the
$\phi$-relative base locus of $G+cE$ is divisorial.
We write
\[
G+cE=G_{c,\free}+G_{c,\base},
\]
where $G_{c,\base}$ is the fixed part of the $\phi$-relative linear
system defined by $G+cE$, and
\[
G_{c,\free}:=(G+cE)-G_{c,\base}
\]
is its movable part.
The line bundle $G_{c,\free}$ is relatively generated over $Y$.
The relative linear system associated with $G_{c,\free}$ induces a
morphism over $Y$
\[
\pi_G\colon
M\to\PP\bigl(\phi_*\OX_M(G_{c,\free})\bigr)
\]
such that
\[
\OX_M(G_{c,\free})
=
\pi_G^\ast
\OX_{\PP(\phi_*\OX_M(G_{c,\free}))}(1).
\]

After replacing $G$ and $c$ by sufficiently large and divisible
multiples, we may assume that $\pi_G$ coincides with the
$\phi$-relative Iitaka fibration associated with $G+cE$.
In particular, the morphism $\pi_G$ is bimeromorphic, and its image, denoted by $X_G$, is normal.
Let $\psi_G\colon X_G\to Y$ be the induced morphism.
We then obtain the following commutative diagram:
\[
\xymatrix@C=4.5em@R=3.2em{
X \ar@{-->}[dr]_{\psi}
&
M \ar[l]^{\pi}
  \ar[d]_{\phi}
  \ar[r]^{\pi_G}
&
X_G \ar@{^{(}->}[r]
    \ar[dl]_{\psi_G}
&
{\PP\bigl(\phi_*\OX_M(G_{c,\free})\bigr)}
    \ar[dll]^{p}
\\
& Y. & & \\
& & &
}
\]

\end{step}

\begin{step}[The local constancy of $X_G$ over $Y_0$]

By Theorem~\ref{thm-faked-flat}, the sheaf $\mathcal{E}_m$ is numerically flat and locally free,
and by Theorem~\ref{thm-flat}, the sheaf $\calV_m$ is a flat vector bundle over $Y_0$ for every $m$.
Moreover, by our construction, 
the flat connection on $\mathcal{E}_m$ descends to a flat connection
on $\calV_m$.
Therefore, the sheaf $\calV_m$ satisfies condition~($\bullet$) in
Proposition~\ref{prop_bir-geometry-psi}(d).

Set
\[
\calW_m
:=
\psi_{G\ast}\mathcal{O}_{X_G}(m)
\otimes\mathcal{O}_Y(-mD_{c,1}),
\]
where
\[
\mathcal{O}_{X_G}(1)
:=
\mathcal{O}_{\PP(\phi_*\mathcal{O}_M(G_{c,\free}))}(1)
\big|_{X_G}.
\]
We claim that $\calW_m\cong\calV_m$ over $Y_0$.

\begin{proof}[Proof of the claim]
By Lemma~\ref{lemm-good}(3), there exists a generically surjective
morphism $\calU_m\to\calV_m$, which corresponds to a global section
\[
s\in H^0(Y,\calU_m^*\otimes\calV_m).
\]
Since both $\calU_m|_{Y_0}$ and $\calV_m|_{Y_0}$ are flat vector
bundles satisfying condition~($\bullet$) in
Proposition~\ref{prop_bir-geometry-psi}(d), so is
$(\calU_m^*\otimes\calV_m)|_{Y_0}$.
Therefore, Proposition~\ref{prop_bir-geometry-psi}(d) implies that
$s|_{Y_0}$ is parallel with respect to the flat connection. 
Consequently, the morphism $\calU_m\to\calV_m$ has constant rank over
$Y_0$ and is therefore surjective over $Y_0$.

Consider the natural inclusion
\[
\phi_*\mathcal{O}_M(mG_{c,\free})
\hookrightarrow
\phi_*\mathcal{O}_M(mG).
\]
After tensoring with $\mathcal{O}_Y(-mD_{c,1})$, we obtain the
commutative diagram
\[
\xymatrix@C=2.5cm@R=1.2cm{
\Sym^{m}\!\phi_*\mathcal{O}_M(G_{c,\free})
\otimes\mathcal{O}_Y(-mD_{c,1})
\ar[r]^-{\cong} \ar[d]
&
\calU_m \ar[d]
\\
\phi_*\mathcal{O}_M(mG_{c,\free})
\otimes\mathcal{O}_Y(-mD_{c,1})
\ar@{^{(}->}[r]
&
\calV_m.
}
\]
The right vertical arrow is surjective over $Y_0$, as shown above.
Since its composition with the top horizontal isomorphism factors
through the bottom horizontal inclusion, the latter is also
surjective over $Y_0$.
Hence, the bottom horizontal arrow is an isomorphism over $Y_0$.

Indeed, the composition map 
$$\Sym^{m}\!\phi_*\mathcal{O}_M(G_{c,\free})\otimes \mathcal{O}_Y(-mD_{c,1}) \to \mathscr{V}_m$$ is surjective.

Now, consider the short exact sequence
\[
0
\to\mathcal{I}_{X_G}
\to\mathcal{O}_{\PP(\phi_*\mathcal{O}_M(G_{c,\free}))}
\to\mathcal{O}_{X_G}
\to 0,
\]
where $\mathcal{I}_{X_G}$ is the ideal sheaf of $X_G$ in
$\PP(\phi_*\mathcal{O}_M(G_{c,\free}))$.
By relative Serre vanishing, for $m\gg 0$, the natural morphism
\[
\begin{aligned}
p_*\mathcal{O}_{\PP(\phi_*\mathcal{O}_M(G_{c,\free}))}(m)
\cong
\Sym^m\!\phi_*\mathcal{O}_M(G_{c,\free}) 
\longrightarrow
\psi_{G*}\mathcal{O}_{X_G}(m)
\end{aligned}
\]
is surjective.
Since $\pi_G\colon M\to X_G$ is bimeromorphic and $X_G$ is normal,
we have
$
\mathcal{O}_{X_G}\cong\pi_{G*}\mathcal{O}_M.
$
The projection formula yields
\[
\psi_{G*}\mathcal{O}_{X_G}(m)
\cong
\phi_*\mathcal{O}_M(mG_{c,\free}).
\]
After tensoring with $\mathcal{O}_Y(-mD_{c,1})$, we obtain an inclusion
\[
\begin{aligned}
\calW_m
&=
\psi_{G*}\mathcal{O}_{X_G}(m)
\otimes\mathcal{O}_Y(-mD_{c,1}) \\
&\cong
\phi_*\mathcal{O}_M(mG_{c,\free})
\otimes\mathcal{O}_Y(-mD_{c,1})
\hookrightarrow\calV_m.
\end{aligned}
\]
Consider the composition
\[
\calU_m
\cong
\Sym^m\!\phi_*\mathcal{O}_M(G_{c,\free})
\otimes\mathcal{O}_Y(-mD_{c,1})
\twoheadrightarrow
\calW_m
\hookrightarrow
\calV_m,
\]
which is surjective over $Y_0$.
It follows that the inclusion $\calW_m\hookrightarrow\calV_m$ is an
isomorphism over $Y_0$.
This proves the claim.
\end{proof}

As a consequence of the claim, for every $m$, the sheaf $\calW_m$ is
a flat vector bundle over $Y_0$ and satisfies condition~($\bullet$) in
Proposition~\ref{prop_bir-geometry-psi}(d).
Therefore, the sheaf 
\[
\SheafHom_{\mathcal{O}_Y}(\Sym^m\calW_1,\calW_m)
\]
is also a flat vector bundle over $Y_0$ satisfying the same condition.
Hence, each of its global sections is parallel.
In particular, condition~(2) in Proposition~\ref{prop_flat-lcf} is
satisfied.
Therefore, by Proposition~\ref{prop_flat-lcf}, we conclude that the restriction of
$\psi_G\colon X_G\to Y$ over $Y_0$ is a locally constant fibration.

Finally, we claim that $(X_G)_y\cong F$ for every $y\in Y_0$, where
$F$ denotes a general fiber of $\psi$; recall that $\psi$ is almost
holomorphic.
Since $\pi_G$ coincides with the $\phi$-relative Iitaka fibration
associated with $G+cE$, its restriction
\[
\pi_G|_{M_y}\colon M_y\to (X_G)_y
\]
is bimeromorphic for general $y\in Y$.
On the other hand, by the construction of $G$ in
Lemma~\ref{lemm-good}, we have
\[
G+cE
=
\pi^*A+\text{an effective }\pi\text{-exceptional divisor},
\]
where $A|_F$ is very ample.
By \eqref{eq-rank}, the linear system
$\left|\pi^*A|_{M_y}\right|$ is the free part of
$\left|(G+cE)|_{M_y}\right|$.
Hence, the morphism $\pi_G|_{M_y}$ factors through
\[
\pi|_{M_y}\colon M_y\to X_y\cong F.
\]

The induced morphism $F\to (X_G)_y$ is bimeromorphic and defined by
the very ample line bundle $A|_F$.
Therefore, it is an isomorphism onto its image, and hence
$(X_G)_y\cong F$ for general $y\in Y_0$.
Since $\psi_G$ is locally constant over $Y_0$, the same conclusion
holds for every $y\in Y_0$.
\end{step}

\begin{step}[Every divisorial component of the exceptional locus of $\pi_G|_{\phi^{-1}(Y_0)}\colon \phi^{-1}(Y_0)\to\psi_G^{-1}(Y_0)$ is $\pi$-exceptional]

Let $D$ be a divisorial component of the exceptional locus of
$\pi_G|_{\phi^{-1}(Y_0)}$.
We distinguish the following two cases.

\begin{itemize}
\item[\textup{Case 1:}]
Assume that $D$ is $\phi$-horizontal.
For general $y\in Y_0$, the restriction $D|_{M_y}$ is exceptional for
$\pi_G|_{M_y}$.
As shown in Step 2, after identifying $(X_G)_y$ with
$X_y\cong F$, the morphism $\pi_G|_{M_y}$ coincides with
$\pi|_{M_y}$.
Therefore $D|_{M_y}$ is contained in $E|_{M_y}$, and hence
$D\subset E$.

\item[\textup{Case 2:}]
Assume that $D$ is $\phi$-vertical.
Since $\phi$ is flat over $Y_0$, the image $\phi(D)$ is a divisor.
For general $y\in\phi(D)$, the fiber $D_y$ of
$\phi|_D\colon D\to\phi(D)$ is an irreducible component of the fiber
$M_y$ of $\phi$.
Since $\pi_G$ is a morphism over $Y$ and $D$ is
$\pi_G$-exceptional, the variety $D_y$ is mapped by $\pi_G$ to a
variety of dimension strictly smaller than $\dim D_y$.
On the other hand, by Lemma~\ref{lemm-good}(2), if $D_y$ is not
contained in $E$, then $(G+cE)|_{D_y}$ is big.
In this case, we see that $D_y$ cannot be mapped by $\pi_G$ to a variety of
smaller dimension.
Therefore $D_y\subset E$, and consequently $D\subset E$.
\end{itemize}
\end{step}

\begin{step}[Proof of the splitting of tangent sheaves]
By Step 2, the restriction of $\psi_G\colon X_G\to Y$ over $Y_0$ is a locally constant fibration.
In particular, by \cite[Remark~2.2]{Wang20}, we have a splitting
\[
T_{\psi_G^{-1}(Y_0)}
\cong
T_{X_G/Y}|_{\psi_G^{-1}(Y_0)}
\oplus
\psi_G^*T_{Y_0}.
\]
By Step 3, every divisorial component of the exceptional locus of
\[
\pi_G|_{\phi^{-1}(Y_0)}
\colon
\phi^{-1}(Y_0)\to\psi_G^{-1}(Y_0)
\]
is contained in the $\pi$-exceptional locus $E$.
Thus, the above splitting induces a splitting on an open subset of
$\phi^{-1}(Y_0)\setminus E$ whose complement has codimension at
least~2.
Taking reflexive extensions and pushing forward by $\pi$, we obtain a
splitting
\[
T_X\cong\mathcal{F}\oplus\mathcal{G},
\]
since both $\pi(E)$ and
$\pi\bigl(\phi^{-1}(Y\setminus Y_0)\bigr)$ have codimension at least two
in $X$.

By Step 2, a general fiber of $\psi_G$ is isomorphic to $F$.
Therefore, $\mathcal{F}$ agrees, on a Zariski-open subset of $X$, with
the foliation induced by the rational map
$\psi\colon X\dashrightarrow Y$, and hence coincides with this
foliation.
In particular, the foliation $\mathcal{F}$ is algebraically integrable, and its
general leaf is a fiber of the MRC fibration
$\psi\colon X\dashrightarrow Y$ and is therefore rationally connected.

By construction, the restriction of $\mathcal{G}$ to the Zariski-open set $\pi(\phi^{-1}(Y_0) \setminus E)$ is $\pi_* \phi^* T_{Y_0}$, which implies that $\mathcal{G}=\pi_{[*]} \phi^* T_{Y_0}$.
By Proposition~\ref{prop_bir-geometry-psi}(a), we have
\[
K_{\mathcal{G}}
\big|_{\pi(\phi^{-1}(Y_0)\setminus E)}
\sim_{\QQ}0.
\]
Since both $\pi(E)$ and
$\pi\bigl(\phi^{-1}(Y\setminus Y_0)\bigr)$ have codimension at least two
in $X$, we conclude that 
$
K_{\mathcal{G}}\sim_{\QQ}0.
$
This completes the proof of Theorem~\ref{thm_splitting}.
\end{step}
The proof of Theorem \ref{thm_splitting} is completed. 
\end{proof}

\subsection{Structure of MRC fibrations}\label{subsec-str}

In this subsection, we prove Theorem~\ref{thm-main}
under the additional technical assumption that $X$ is strongly
$\mathbb{Q}$-factorial and maximally quasi-\'etale. 
Our proof follows the strategy of previous work \cite{MW}, 
but several modifications are required due to the lack of projectivity. 
For Steps~1 and~2, it suffices to modify a few technical points. 
By contrast, the latter half of Step~3 requires a new argument even in the projective case.

\begin{proof}[Proof of Theorem~\ref{thm-main} under the additional technical assumption. ]

Let $(X,\Delta,\boldsymbol{\beta})$ be a K\"ahler generalized klt pair
of Calabi--Yau type, as in Theorem~\ref{thm-main}.
Suppose that $X$ is strongly $\mathbb{Q}$-factorial and maximally quasi-\'etale.
The proof can be divided into three steps:
\setcounter{step}{0}
\begin{step}[{\bf Holomorphicity of the MRC fibration}]
A key observation is that both $\calF$ and $\calG$ are weakly regular
and that
$\calF\subset T_X$ is algebraically integrable.
This allows us to apply an analytic version of
\cite[Theorem~17]{DGP24} and conclude that $\calF$ is induced by an
equidimensional fibration.
The validity of an analytic analogue of
\cite[Theorem~17]{DGP24} is highly nontrivial. 
Indeed, its proof involves several delicate arguments.
We provide a rigorous proof of the required analytic version in Section~\ref{sec-druel}.

As explained above, we can apply an analytic analogue of \cite[Theorem~17]{DGP24} 
to obtain an equidimensional fibration $f \colon X\to Y$ onto a K\"ahler variety $Y$. 

In this step, we verify the following properties of $f \colon X\to Y$. 
\begin{enumerate}
\item[$(1)$] $\calF$ coincides with the foliation induced by $f \colon X\to Y$;
\item[$(2)$] $f \colon X\to Y$ is an equidimensional projective MRC fibration; 
\item[$(3)$] $Y$ is strongly $\mathbb{Q}$-factorial and maximally quasi-\'etale. 
\item[$(4)$] $Y$ has canonical singularities and satisfies 
 $K_Y\sim_{\QQ}0$. 
\end{enumerate}

\smallskip
\noindent \textup{(1), (2)}
Since $\mathcal{F}$ is obtained from an MRC fibration, 
properties~(1) and~(2), except for projectivity, follow directly from the construction of $f \colon X\to Y$. 
The projectivity of $f \colon X\to Y$ follows from the projectivity criterion in \cite{CH24} 
since  $X$ is strongly $\mathbb{Q}$-factorial 
and the general fibers are rationally connected. 

\smallskip
\noindent \textup{(3)}
The strong $\mathbb{Q}$-factoriality of $Y$ follows from the equidimensionality of $f$ and the strong $\mathbb{Q}$-factoriality of $X$ by Lemma~\ref{lem-base}(1) (cf.~\cite[Lemma~4.1]{IMZ23}). 
Let $Y' \to Y$ be a quasi-\'etale cover of $Y$, 
and let $X'$ be the normalization of $X \times_{Y} Y'$. 
Then, the natural morphism $X' \to X$ is also quasi-\'etale by the equidimensionality of $f \colon X\to Y$. 
Since $X$ is maximally quasi-\'etale, the morphism $X'\to X$ is \'etale. 
It follows that $Y' \to Y$ is \'etale, so $Y$ is maximally quasi-\'etale.

\smallskip
\noindent \textup{(4)}
The fibration $f \colon X \to Y$ is semistable in codimension one by Proposition~\ref{prop_bir-geometry-psi}(e). 
Hence, the ramification divisor of $f$ is zero (cf.\,\cite[Definition 2.16]{CKT16}). 
Then \cite[Lemma~2.31]{CKT16} and the condition $K_{\calG}\sim_{\QQ}0$ yield 
\[
K_{X/Y}\sim K_{\calF}\sim_{\QQ} K_X. 
\] 
This implies that $K_Y\sim_{\QQ}0$, 
and hence $Y$ has canonical singularities by Proposition~\ref{prop_bir-geometry-psi}(a).

Indeed, take a log resolution $\alpha\colon\widetilde{Y}\to Y$. We can write
\[
\alpha^*K_Y+E=K_{\widetilde{Y}}+F,
\]
where $E$ and $F$ are effective $\alpha$-exceptional divisors with no
common components. It suffices to prove that $F=0$.
The variety $\widetilde{Y}$ appears as the base of the MRC fibration
considered in Proposition~\ref{prop_bir-geometry-psi}. Hence, by
Proposition~\ref{prop_bir-geometry-psi}\textup{(a)}, there exist an integer
$m>0$ and an effective divisor $D$ such that
$m(E-F)\sim D$.
Since $m(E-F)$ is $\alpha$-exceptional, so is $D$. Therefore,
$m(E-F)-D$ is an $\alpha$-exceptional divisor that is numerically trivial
over $Y$. Applying the negativity lemma to this divisor and its negative,
we obtain
\[
m(E-F)=D.
\]
In particular, $E-F\geq 0$, and hence $E\geq F$. Since $E$ and $F$ have no
common components, it follows that $F=0$.

\end{step}

\begin{step}[{\bf Direct image sheaves of relatively ample line bundles}]
In this step, using the equidimensionality of $f$, 
we construct an $f$-relatively ample $\mathbb{Q}$-line bundle $L_X$ on $X$ 
such that the direct image sheaves $f_* \mathcal{O}_{X}(mL_X)$ are flat vector bundles, 
which enables us to apply the criterion for locally constant 
fibrations in Proposition~\ref{prop_flat-lcf}.

To this end, let $\mu\colon Y'\to Y$ be a desingularization that is
an isomorphism over $Y_{\reg}$, let $X'$ be the normalization of the
fiber product $X\times_Y Y'$, and let $M$ be a desingularization of
$X'$, as in the following diagram:
\begin{center}
\begin{tikzpicture}[scale=2.0]
\node (A) at (0,0) {$Y$.};
\node (B) at (0,1) {$X$};
\node (C) at (-1,0) {$Y'$};
\node (D) at (-1,1) {$X'$};
\node (E) at (-1.6,1.6) {$M$};
\path[->,font=\scriptsize,>=angle 90]
(B) edge node[right]{$f$} (A)
(D) edge node[left]{$f'$} (C)
(D) edge node[above]{$\mu_X$} (B)
(C) edge node[below]{$\mu$} (A)
(E) edge node[above right]{$\pi'$} (D)
(E) edge [bend left=20] node[above right]{$\pi$} (B)
(E) edge [bend right=20] node[below left]{$\phi$} (C);
\end{tikzpicture}
\end{center}
We return to the situation of Setting~\ref{setting}.
However, there is a slight difference in notation:
the smooth variety $Y'$ in the present situation plays the role of $Y$ in Setting~\ref{setting}.

We recall the construction of $G$ from the proof of Lemma~\ref{lemm-good}. 
The morphism $\beta \colon X' \to X$ in that lemma is an isomorphism in the present situation. 
Since $f$ is projective, we can take an $f$-relatively ample line bundle $A$ on $X$ 
and apply the argument to $G':=A$. 
Then $F'=F=0$ in that lemma. 
Thus, the desired line bundle $G$ coincides with $\pi^{*}A$.

Let $m \in \ZZ_{+}$ be a sufficiently large and divisible integer.
Since $X$ is strongly $\QQ$-factorial, 
we can find an effective Cartier divisor $E_X$ on $X'$ 
such that $-E_X$ is $\mu_X$-relatively ample and $\Supp(E_X)=\Exc(\mu_X)$

(for example, see the proof of \cite[Complement~2.10, p.~74]{Kollar07}).
Then, we have an isomorphism 
\begin{align}\label{eq-isompre}
\notag 
\calU_{c,m}:&=f'_\ast\OX_{X'}\bigl(m(\mu_X^\ast A+cE_X)\bigr)\\
&\cong\phi_*\OX_M\bigl(m(\pi^\ast A+c(\pi'^{\ast}E_X+E'))\bigr),
\end{align}
where $E'$ denotes the exceptional divisor of $\pi'$. 
As in Section~\ref{sec-main}, there exists $c_0\in\ZZ_{+}$ 
such that both $\rank(\calU_{c,m})$ and $\det\calU_{c,m}|_{\mu^{-1}(Y_{\reg})}$ are 
independent of $c\geqslant c_0$. 
Thus, in what follows, we fix such a $c\geqslant c_0$ and omit the subscript $c$. 
The variety $X$ is maximally quasi-\'etale. 
Thus, by Corollary~\ref{thm-flat},  we see that 
\[
\mathcal{E}_{m}:=\pi_{[\ast]}\phi^\ast\mathcal{V}_{m}
\]
is a numerically flat locally free sheaf on $X$, 
where $\mathcal{V}_{m}$ is defined by \eqref{eq-line2} 
after replacing $E$ in \eqref{eq-line} with $\pi'^{\ast}E_X+E'$. 
Using the isomorphism \eqref{eq-isompre}, we may write 
\[
\calV_m = f'_\ast\scrO_{X'}(mL_{X'}), 
\quad\text{where}\quad
L_{X'}:=(\mu_X^\ast A+cE_X)-\frac{1}{r}\bigl(f'^{\ast}\det \calU_1\bigr),
\]
and $r:=\rank\calU_{1}$. 
Define
\[
L_X:=\mu_{X[*]}L_{X'} 
\quad\text{and}\quad 
\calW_m:=f_\ast\scrO_X\bigl(mL_X\bigr).
\]

Since $\mu\colon Y' \to Y$ is an isomorphism over $Y_{\reg}$, we have $\mu_{X[*]}f'^{\ast}\det \calU_1=f^\ast\mu_{[*]}\det \calU_1$ on $X$. 
This shows that $L_X$ is $f$-relatively ample.

By \cite[Corollary~1.7]{Har80}, 
the sheaf $\calW_m=f_\ast\scrO_X(mL_X)$ is reflexive since $f$ is equidimensional. 
Since $\mu_X \colon X' \to X$ is an isomorphism over $f^{-1}(Y_{\reg})$, we have 
\[
\calW_m \cong\calW_m^{\ast\ast}
\cong \mu_{[\ast]}\calV_m\bigr.
\]
Moreover, the equidimensionality of $f$ gives 
\[
\mathcal{E}_{m}
=\pi_{[\ast]}\phi^\ast\mathcal{V}_{m}
\cong \mu_{[\ast]} f'^{\ast}\mathcal{V}_m 
\cong f^{[\ast]}\mathcal{W}_m
\quad\text{on }X.
\]
There exists a Zariski-closed subset $C_{m} \subset Y$ with $\codim C_{m} \geqslant 2$ such that $\calW_m$ is locally free on $Y_{\reg} \setminus C_{m}$. 
Since $\mathcal{E}_{m}$ is numerically flat and locally free, we have $\mathcal{E}_{m}\cong f^{*}\calW_m$ over $Y_{\reg} \setminus C_{m}$. 
After possibly replacing $C_{m}$, the argument in the proof of Theorem~\ref{thm-flat} shows that $\calW_m$ is flat on $Y_{\reg} \setminus C_{m}$. 
Since $\codim C_m\geqslant 2$, the inclusion $Y_{\reg}\backslash C_m\hookrightarrow Y_{\reg}$ induces an isomorphism $\pi_1(Y_{\reg}\backslash C_m)\cong\pi_1(Y_{\reg})$, so $\calW_m|_{Y_{\reg}}$ is also flat. 
Since $Y$ is maximally quasi-\'etale, the reflexivity of $\calW_m$ and Theorem \ref{thm-maxqetale}~(2) imply that $\calW_m$ is a flat vector bundle.

\end{step}

\begin{step}[{\bf Conclusion}]
In this final step, we establish that $f$ is locally constant. 
By replacing $L_X$ with $m_0 L_X$ for some sufficiently large and divisible $m_0 \in \ZZ_{+}$, we may assume that $\calW_m$ is a flat vector bundle for every $m\in\ZZ_{+}$. 
In particular, these direct images are locally free, and the flatness criterion implies that $f \colon X \to Y$ is flat (see, for example, \cite[Vol.~II, Proposition~2.5, p.~8]{ACG11}). 
Applying Proposition~\ref{prop_bir-geometry-psi}(d) to the flat vector bundle \(\SheafHom_{\OX_Y}\bigl(\Sym^m\calW_1,\calW_m\bigr)\), we see that the natural morphism 
\[
\Sym^m\calW_1 \longrightarrow \calW_m
\]
is a morphism of local systems over $Y_{\reg}$
(i.e.\ it is induced by a morphism of representations of $\pi_1(Y_{\reg})$). 
Since $Y$ is maximally quasi-\'etale, 
it is also induced by a morphism of representations of $\pi_1(Y)$ by Theorem \ref{thm-maxqetale}~(2). 
Therefore, condition~\textup{(2)} of Proposition~\ref{prop_flat-lcf} is satisfied, and hence $f$ is a locally constant fibration.

\smallskip
We show that the fibration $f\colon X\to Y$ is locally constant with
respect to the boundary $\Delta$ as well.
To this end, consider the representation
$\rho\colon\pi_1(Y)\to\Aut(F)$ defining the locally constant fibration
$f\colon X\to Y$.
It remains to show that $\Delta_F:=\Delta|_F$ is invariant under the
action of $\pi_1(Y)$ and that there is an isomorphism of pairs over $Y$
\[
(X,\Delta)
\cong
\bigl(\Unv{Y}\times(F,\Delta_F)\bigr)/\pi_1(Y).
\]
It suffices to prove the corresponding assertion for
$\lfloor p\Delta\rfloor$ for every $p\in\mathbb{Z}_{+}$.
Indeed, since
$
(1/p)\lfloor p\Delta_F\rfloor
\to
\Delta_F
$
as $p\to\infty$, the corresponding assertions for
$\lfloor p\Delta\rfloor$ imply the desired conclusion for
$(X,\Delta)$.

From now on, let $m$ be a sufficiently large and divisible integer.
By Lemma~\ref{lemma_lc_line}, there exists an $f$-relatively ample
line bundle $A$ on $X$ such that $f_*\mathcal{O}_X(mA)$ is a flat
locally free sheaf.
Note that $\Delta$ is horizontal by
Proposition~\ref{prop_bir-geometry-psi}\textup{(b)}.
By the criterion for locally constant fibrations of pairs
(see Proposition~\ref{prop_flat-lcf}), it suffices to prove that, for
every fixed $p\in\mathbb{Z}_{+}$ and $m\gg p$, the direct image sheaf
\[
f_*\mathcal{O}_X\bigl(mA-\lfloor p\Delta\rfloor\bigr)
\]
is locally free and numerically flat.

We first show that $f_*\mathcal{O}_X(mA)$ is numerically flat and that
$A$ is pseudo-effective.
Let
$
\rho_m\colon\pi_1(Y)\to\GL(r,\CC)
$
be the linear representation associated with the flat locally free
sheaf $f_*\mathcal{O}_X(mA)$.
Since $Y$ has canonical singularities and $K_Y\sim_{\QQ}0$, we have
$\kappa(Y')=0$.
Hence $Y'$ is special in the sense of Campana.
Applying \cite[Theorem~7.8]{Cam04} to the representation
\[
\rho_m\circ\mu_*\colon\pi_1(Y')\to\GL(r,\CC),
\]
we find that its image is virtually abelian.
Since the natural homomorphism
$\mu_*\colon\pi_1(Y')\to\pi_1(Y)$ is surjective, the image
$\Image (\rho_m)$ is also virtually abelian.
Thus, there exists a finite \'etale cover
$\nu\colon\tilde{Y}\to Y$ such that the representation
corresponding to
$\nu^*(f_*\mathcal{O}_X(mA))$ has abelian image.

By simultaneous triangularization, the flat vector bundle
$\nu^*(f_*\mathcal{O}_X(mA))$ admits a filtration by flat subbundles
whose successive quotients are flat line bundles.
Since every flat line bundle is numerically flat and extensions of
numerically flat vector bundles are numerically flat, it follows from
\cite[Proposition~1.15(ii)]{DPS94} that
$\nu^*(f_*\mathcal{O}_X(mA))$ is numerically flat.
Numerical flatness descends under finite \'etale covers; hence,
$f_*\mathcal{O}_X(mA)$ is numerically flat.

The relative global generation of $\mathcal{O}_X(mA)$ gives a
surjective evaluation morphism
\[
f^*f_*\mathcal{O}_X(mA)
\longrightarrow
\mathcal{O}_X(mA).
\]
The vector bundle on the left-hand side is nef, and hence the line
bundle $mA$ is nef.
In particular, $A$ is pseudo-effective.

We work with the following commutative diagram:
\[
\xymatrix@C=5em@R=4em{
M \ar[r]^{\pi} \ar[d]_{\phi}
&
X \ar[d]^{f}
\\
Y' \ar[r]_{\mu}
&
Y.
}
\]
For a sufficiently large and divisible integer $k$ satisfying
$k\gg m\gg p$, we consider the following direct image sheaf on $Y'$:
\begin{align*}
&\phi_*\left(
\mathcal{O}_M\left(
kK_{M/Y'}
+m\pi^*A
+\bigl(
k\Delta_M-\lfloor p\pi^{-1}_*\Delta\rfloor
+\lceil kF_+\rceil-kF_+
\bigr)
+k(F_+-K_M-\Delta_M)
\right)
\right)\\
&\qquad =
\phi_*\left(
\mathcal{O}_M
\bigl(
\lceil kF_+\rceil
+m\pi^*A
-\lfloor p\pi^{-1}_*\Delta\rfloor
\bigr)
\right)
\otimes\mathcal{O}_{Y'}(-kK_{Y'}).
\end{align*}
The line bundle $m\pi^*A$ is pseudo-effective and relatively
$\phi$-big, while $F_+-K_M-\Delta_M$ is nef.
Moreover, since $\pi^{-1}_*\Delta\leq\Delta_M$, the divisor
\[
\frac{1}{k}
\bigl(
k\Delta_M-\lfloor p\pi^{-1}_*\Delta\rfloor
+\lceil kF_+\rceil-kF_+
\bigr)
\]
is effective and satisfies the klt condition for
$k\gg 1$.
Thus, Proposition~\ref{prop-psef} shows that the above direct image
sheaf is weakly positively curved on $Y'$.

The divisor
$
\pi^*(\lfloor p\Delta\rfloor)
-\lfloor p\pi^{-1}_*\Delta\rfloor
$
is effective and $\pi$-exceptional.
Moreover, the divisor
$
K_{Y'}-\mu^*K_Y
$
is effective and $\mu$-exceptional since $Y$ has canonical
singularities.
Using also the fact that $\lceil kF_+\rceil$ is $\pi$-exceptional, the
projection formula yields
\begin{align}\label{eq-push}
&\mu_{[*]}\left(
\phi_*\mathcal{O}_M
\bigl(
\lceil kF_+\rceil
+m\pi^*A
-\lfloor p\pi^{-1}_*\Delta\rfloor
\bigr)
\otimes\mathcal{O}_{Y'}(-kK_{Y'})
\right)\notag\\
\cong{}&
f_*\mathcal{O}_X\bigl(mA-\lfloor p\Delta\rfloor\bigr)
\otimes\mathcal{O}_Y(-kK_Y).
\end{align}
By Lemma~\ref{lem-push}, the sheaf on the right-hand side of
\eqref{eq-push} is weakly positively curved.
Since $K_Y\sim_{\QQ}0$ and $k$ is sufficiently divisible, it follows
that
$
f_*\mathcal{O}_X\bigl(mA-\lfloor p\Delta\rfloor\bigr)
$
is weakly positively curved and that its determinant
$
\det f_*\mathcal{O}_X\bigl(mA-\lfloor p\Delta\rfloor\bigr)
$
is pseudo-effective.

We next show that the dual of this determinant sheaf is also
pseudo-effective.
Applying Proposition~\ref{prop_more-psef-det} to
$m\pi^*A-\lfloor p\pi^{-1}_*\Delta\rfloor$, with this divisor playing
the role of $mG$, we obtain an effective $\pi$-exceptional divisor
$d E$ on $M$ such that
\[
m\pi^*A-\lfloor p\pi^{-1}_*\Delta\rfloor+dE-\phi^*\det\phi_*\mathcal{O}_M
\bigl(
m\pi^*A-\lfloor p\pi^{-1}_*\Delta\rfloor+cE
\bigr)
\]
is pseudo-effective.
Let $\theta$ be a smooth representative of the class of the
determinant sheaf appearing above.
It follows that
$m\pi^*A+dE$ is $\phi^*\theta$-weakly positively curved.

Applying the same direct-image argument to $m\pi^*A+F$, we find that
\[
\phi_*\mathcal{O}_M
\bigl(
\lceil kF_+\rceil+m\pi^*A+dE
\bigr)
\otimes\mathcal{O}_{Y'}(-kK_{Y'})
\]
is $\theta$-weakly positively curved.
Taking determinants, we conclude that
\[
c_1(\phi_*\mathcal{O}_M
\bigl(
\lceil kF_+\rceil+m\pi^*A+F
\bigr)
\otimes\mathcal{O}_{Y'}(-kK_{Y'}))
-
r_m c_1(m\pi^*A-\lfloor p\pi^{-1}_*\Delta\rfloor+cE)
\]
is pseudo-effective, where $r_m$ denotes the rank of the direct image
sheaf.
Since $E$, and $F_+$ are exceptional divisors, we have
\begin{align*}
&\mu_{[*]}\left(
\phi_*\mathcal{O}_M
\bigl(
\lceil kF_+\rceil+m\pi^*A+F
\bigr)
\otimes\mathcal{O}_{Y'}(-kK_{Y'})
\right)
\cong
f_*\mathcal{O}_X(mA)\otimes\mathcal{O}_Y(-kK_Y)
\quad\text{and}\\
&\mu_{[*]}\phi_*\mathcal{O}_M
\bigl(
m\pi^*A-\lfloor p\pi^{-1}_*\Delta\rfloor+cE
\bigr)
\cong
f_*\mathcal{O}_X\bigl(mA-\lfloor p\Delta\rfloor\bigr).
\end{align*}
The first Chern class of the first one is numerically flat 
by the numerical flatness of $f_*\mathcal{O}_X(mA)$ together with $K_Y\sim_{\QQ}0$. 
Thus, by pushing the preceding pseudo-effective class forward by $\mu$, 
we conclude that
\[
r_m
c_1\left(
\det f_*\mathcal{O}_X
\bigl(mA-\lfloor p\Delta\rfloor\bigr)
\right)
\]
is pseudo-effective.

Therefore, we have proved that
\[
c_1\left(
f_*\mathcal{O}_X
\bigl(mA-\lfloor p\Delta\rfloor\bigr)
\right)
=
c_1\left(
\det f_*\mathcal{O}_X
\bigl(mA-\lfloor p\Delta\rfloor\bigr)
\right)
=
0.
\]
The flatness criterion for pseudo-effective sheaves with vanishing
first Chern class shows that
$
f_*\mathcal{O}_X\bigl(mA-\lfloor p\Delta\rfloor\bigr)
$
is locally free and numerically flat.
\end{step}
\vspace{-0.3cm}
\end{proof}

We conclude this subsection by proving the lemmas used in Step 1.
Although they are presumably well known to experts, we include the proofs for completeness.
\begin{lemm}\label{lem-Weil}
Let $X$ be a normal analytic variety and $\mathcal{L}$ a reflexive
sheaf of rank one on $X$.
Then, locally on $X$, there exists an effective Weil divisor $D$ such
that
$
\mathcal{L}\cong\mathcal{O}_X(D).
$
\end{lemm}

\begin{proof}
We may assume that $X$ is Stein.
By Cartan's theorem A, we can choose a nonzero section
\[
s\in H^0(X,\mathcal{L})
\]
that generates $\mathcal{L}$.
This section determines an effective Weil divisor $D$.
Since $\mathcal{L}$ is reflexive, the associated divisorial sheaf
$\mathcal{O}_X(D)$ is isomorphic to $\mathcal{L}$.
\end{proof}

\begin{lemm}\label{lem-base}
Let $\phi\colon X\to Y$ be an equidimensional projective fibration
between normal analytic varieties.
If $X$ is strongly $\mathbb{Q}$-factorial, then so is $Y$.
\end{lemm}

\begin{proof}
Let $\mathcal{L}$ be a reflexive sheaf of rank one on $Y$.
Since $X$ is strongly $\mathbb{Q}$-factorial, the reflexive pullback
$\phi^{[*]}\mathcal{L}$ is $\mathbb{Q}$-Cartier.
Since the assertion is local on $Y$, Lemma~\ref{lem-Weil} allows us
to write
$
\mathcal{L}\cong\mathcal{O}_Y(D)
$
for some Weil divisor $D$ on $Y$.

Choose a relative embedding
$
X\hookrightarrow\mathbb{P}^N\times Y
$
over $Y$, and set $r:=\dim X-\dim Y$.
By taking sufficiently general relative hyperplane sections
$H_1,\ldots,H_r$ of $X$, we obtain
\[
Z:=H_1\cap\cdots\cap H_r
\]
such that the induced morphism
$\phi_Z:=\phi|_Z\colon Z\to Y$ is finite and surjective.
Since $\phi^{[*]}\mathcal{L}$ is $\mathbb{Q}$-Cartier, there exists
$m_0\in\mathbb{Z}_{+}$ such that $m_0\phi^*D$ is Cartier.

Since $\phi$ is equidimensional, the pullbacks $\phi^*D$ and
$\phi_Z^*D$ of the Weil divisor $D$ are well-defined and satisfy
$
\phi_Z^*D=(\phi^*D)|_Z.
$
Thus, $m_0\phi_Z^*D$ is Cartier.
Since $\phi_Z\colon Z\to Y$ is finite and surjective, the norm
argument in \cite[Lemma~5.16]{KM98} implies that $D$ is
$\mathbb{Q}$-Cartier.
\end{proof}

\section{Descent of locally constant fibrations} \label{sec-red}

\subsection{Descent: the product case} \label{subsec-baby}

In this subsection, we explain how to descend the conclusion of Theorem~\ref{thm-main}
from a strongly $\mathbb{Q}$-factorialization.
This descent was established in \cite[Section~4]{MW} in the projective
case, but the argument does not extend directly to the K\"ahler setting.
One of the main obstacles
(see \cite[Theorem~4.12]{MW}) is the need for a transcendental
base-point-free theorem for Calabi--Yau varieties, asserting that every
nef and big Bott--Chern class on a Calabi--Yau variety is semiample
(see \cite[Conjecture~6.1]{Tos09} and
\cite[Conjecture~2.7]{Tos26}).
This particular obstacle can now be overcome by the recent resolution
of the transcendental base-point-free conjecture in \cite{HX26}.
Nevertheless, further difficulties specific to the K\"ahler setting remain,
and the projective argument still cannot be applied.
Our argument provides a complementary perspective on the role of the
transcendental base-point-free theorem in this descent problem:
it bypasses this theorem and simultaneously overcomes the additional
difficulties through a new approach based on the following idea.

Let $(X', \Delta')$ be a K\"ahler klt pair such that $X$ is maximally quasi-\'etale, 
and take a strongly $\mathbb{Q}$-factorialization
$\pi \colon (X, \Delta) \to (X', \Delta')$. 
Note that the notation is different from Subsection~\ref{subsec-str}. 
By the results of Subsection~\ref{subsec-str}, the pair $(X, \Delta)$ admits
a locally constant fibration $(X, \Delta) \to Y$ onto a Calabi--Yau variety,
with rationally connected fibers.
Our goal is to descend this locally constant fibration to $(X', \Delta')$.
We first consider the special case in which
$(X, \Delta) \cong Y \times (F, \Delta_{F})$.
This case is treated in Proposition~\ref{prop-baby}, which generalizes
\cite[Lemma~7.5]{BGL22}.
In general, the fibration $(X, \Delta) \to Y$ need not be a product,
but the fiber product $(X, \Delta)\times_{Y} \Unv{Y}$ is a product.
By adapting the argument for the product case, we can complete
the descent argument (see Theorem~\ref{thm-reduction}).
The definition of locally constant fibrations plays a crucial role in this argument.
We emphasize that this approach is new even in the projective case.

We begin with the product case.

\begin{prop}\label{prop-baby}
Let $Z$ be a $($not necessarily compact$)$ K\"ahler variety with rational singularities 
and let $F$ be a compact K\"ahler variety with rational singularities and irregularity zero. 
Consider the product $W:=Z \times F$ 
and a fibration $q \colon W \to W'$ onto a K\"ahler variety $W'$. 
Then, the following assertions hold:

$(1)$ There exist fibrations $q_{1} \colon Z \to Z'$ and $q_{2} \colon F \to F'$ 
onto K\"ahler varieties $Z'$ and $F'$, respectively, and an isomorphism $W' \cong Z' \times F'$ 
that fits into the following commutative diagram:
\begin{align*}
\xymatrix{
& W=Z \times F \ar[ld]_{q}  \ar[rd]^{\quad q':=q_{1} \times q_{2}} & \\
W'  \ar[rr]^{\quad\cong \quad }  &  & Z' \times F'.
}
\end{align*}

$(2)$ Suppose that $q \colon W \to W'$ is a bimeromorphic morphism. 
Then, the fibrations $q_{1} \colon Z \to Z'$ and $q_{2} \colon F \to F'$ obtained in $(1)$ 
are bimeromorphic. 
Furthermore, for a Weil divisor $\Delta_{F}$ on $F$, set
\[
\Delta' := q_{*}\pr_{2}^{*} \Delta_{F} 
\quad \text{ and } \quad 
\Delta_{F}' := q_{2*} \Delta_{F}, 
\]
where $\pr_{2} \colon W=Z \times F \to F$ is the second projection. 
Then, the isomorphism $W' \cong Z' \times F'$ obtained in $(1)$ 
is compatible with the boundary divisors; that is,
\[
(W', \Delta')\cong Z' \times (F', \Delta_{F}'). 
\]
\end{prop}

\begin{cor}\label{cor-baby}
Let $W$, $W'$, and $Z$ be normal analytic varieties. 
Assume that $Z$ has rational singularities. 
Let $W' \to Z$ be a fibration and let $(W, \Delta)\to Z$ be a locally trivial fibration 
whose fiber $F$ is a compact K\"ahler variety with rational singularities and irregularity zero. 

$(1)$
Assume that there exists a fibration $q \colon W \to W'$ over $Z$. 
Then $W' \to Z$ is a locally trivial fibration.

$(2)$
Assume that there exists a bimeromorphic morphism $q \colon W \to W'$ over $Z$. 
Then $(W', \Delta') \to Z$ is a locally trivial fibration, where $\Delta':=q_{*}\Delta$.
\end{cor}

\begin{proof}[Proof of Proposition~\ref{prop-baby}]
We first consider the case where $Z$ is compact. 
The noncompact case requires an additional technical argument, 
which we give at the end of the proof.

Fix a K\"ahler class $\alpha$ on $W'$. 
We show that $W'$ splits as a product 
by decomposing the pullback $q^{*}\alpha$ on $W$. 
Since $Z$ and $F$ may be singular, 
we first pass to resolutions and decompose the pullback of $\alpha$ on the resulting smooth varieties. 
Let $\pi_{1} \colon \tilde{Z} \to Z$ and $\pi_{2} \colon \tilde{F} \to F$ 
be resolutions of singularities of $Z$ and $F$. 
\begin{rem}
Since $Z$ and $F$ are compact K\"ahler, 
the natural morphisms 
\begin{align}\label{eq-inj}
\tilde j_{1} \colon H^{1,1}_{\BottChern}(\tilde{Z}) \to H^{2}(\tilde{Z}, \mathbb{C})
\quad  
\text{ and }  \quad   
\tilde j_{2} \colon H^{1,1}_{\BottChern}(\tilde{F}) \to H^{2}(\tilde{F}, \mathbb{C}) 
\end{align}
are injective. 
Note that 
\[
\tilde j \colon H^{1,1}_{\BottChern}(\tilde{Z}\times \tilde{F}) \to H^{2}(\tilde{Z}\times \tilde{F}, \mathbb{C})
\]
is also injective. 
From now on, we identify the Bott--Chern classes with their images. 
However, this is not the case when $Z$ is noncompact; 
we will resolve this issue at the end of the proof.
\end{rem}

Set $\tilde{W}:=\tilde{Z} \times \tilde{F}$ and 
consider the following diagram:
\[
\xymatrix{
\tilde{F} & \tilde{W}:=\tilde{Z}\times \tilde{F}
\ar[l]_{\tilde\pr_{2}\qquad}
\ar[d]^{\pi:=\pi_{1}\times \pi_{2}}
\ar[r]^{\qquad \tilde\pr_{1}} & \tilde{Z} \\
F & W=Z\times F
\ar[l]_{\pr_{2}\qquad}
\ar[r]^{\qquad \pr_{1}} & Z}.
\]
We first establish the following claim.
\begin{claim}\label{claim-baby1}
The class $\pi^{*}(q^{*}\alpha)$ decomposes in the second cohomology group as the sum of classes pulled back from the two factors. 
\end{claim}
\begin{proof}
Since $\tilde{F}$ is a compact K\"ahler manifold, 
the Hodge decomposition shows that 
\[
h^{1}(\tilde{F}, \mathbb{C})
=
2h^{1}(\tilde{F}, \mathcal{O}_{\tilde{F}})
=
2h^{1}(F, \mathcal{O}_{F})
=
0.
\]
Here, the last two equalities follow from the facts that $F$ has rational singularities and irregularity zero. 
The K\"unneth formula then yields the decomposition 
\begin{align*}
H^{2}(\tilde{W}, \mathbb{C})
&=
H^{2}(\tilde{Z}, \mathbb{C})
\oplus
\bigl(H^{1}(\tilde{Z}, \mathbb{C}) \otimes H^{1}(\tilde{F}, \mathbb{C})\bigr)
\oplus
H^{2}(\tilde{F}, \mathbb{C}) \\
&=
H^{2}(\tilde{Z}, \mathbb{C}) \oplus H^{2}(\tilde{F}, \mathbb{C}).
\end{align*}
This decomposition allows us to write
\begin{align}\label{eq-dec}
\pi^{*} (q^{*}\alpha)
=
\tilde \pr_{1}^{*} \tilde \beta+\tilde\pr_{2}^{*}\tilde\gamma \quad \text{in }H^{2}(\tilde{W}, \mathbb{C}),
\end{align}
where $\tilde \beta$ and $\tilde \gamma$ are degree-two cohomology classes on $\tilde{Z}$ and $\tilde{F}$, respectively. 
\end{proof}

We need a decomposition of $q^{*}\alpha$ into Bott--Chern cohomology classes, rather than merely a decomposition in the second cohomology group.
However, a priori $\tilde \beta$ and $\tilde \gamma$ need not be Bott--Chern cohomology classes, 
and the pushforwards $\pi_{1*}\tilde \beta$ and $\pi_{2*}\tilde \gamma$ may not be defined as Bott--Chern classes because of the singularities. 
Therefore, we prove the following claim. 

\begin{claim}\label{claim-baby2}
The class $q^{*}\alpha$ decomposes as the sum of Bott--Chern classes pulled back from the two factors. 
\end{claim}
\begin{proof}
Choose general points of $Z$ and $F$, and consider their lifts to $\tilde{Z}$ and $\tilde{F}$, respectively. 
For simplicity of notation, we denote these points and their lifts by the same symbol $*$. 

Let $F'' \subset W'$ be the image of $\{*\}\times F\subset W$ under $q \colon W \to W'$, 
and let $\nu_{2}\colon F'\to F''$ be the finite morphism obtained 
from the Stein factorization of $\{*\}\times F \to F''$. 
Applying the same construction to $Z\times \{*\}$, we obtain 
$Z'' \subset W'$ and $\nu_{1}\colon Z'\to Z''$. 
We have the following commutative diagram:
\[
\begin{tikzpicture}[baseline=(m-2-2.base)]
  \matrix (m) [
    matrix of math nodes,
    column sep=4.5em,
    nodes={inner sep=1pt}
  ] {
    \{*\}\times \tilde{F} &
    \tilde{W}=\tilde{Z}\times \tilde{F} &
    \tilde{Z}\times \{*\} \\[3.2em]
    \{*\}\times F &
    W=Z\times F &
    Z\times \{*\} \\[6.0em]
    F'' &
    W' &
    Z'' \\
  };

  \draw[right hook->] (m-1-1) -- (m-1-2);
  \draw[right hook->] (m-2-1) -- (m-2-2);
  \draw[right hook->] (m-3-1) -- (m-3-2);

  \draw[left hook->] (m-1-3) -- (m-1-2);
  \draw[left hook->] (m-2-3) -- (m-2-2);
  \draw[left hook->] (m-3-3) -- (m-3-2);

  \draw[->] (m-1-2) to[bend right=18] node[above] {$\tilde\pr_{2}$} (m-1-1);
  \draw[->] (m-1-2) to[bend left=18] node[above] {$\tilde\pr_{1}$} (m-1-3);

  \draw[->] (m-2-2) to[bend left=18] node[below] {$\pr_{2}$} (m-2-1);
  \draw[->] (m-2-2) to[bend right=18] node[below] {$\pr_{1}$} (m-2-3);

  \draw[->] (m-1-1) -- node[right] {$\pi_2$} (m-2-1);
  \draw[->] (m-1-3) -- node[right] {$\pi_1$} (m-2-3);
  \draw[->] (m-1-2) -- node[right] {$\pi:=\pi_{1}\times \pi_{2}$} (m-2-2);
  \draw[->] (m-2-2) -- node[right] {$q$} (m-3-2);

  \node (Fp) at ([xshift=-1.8em]$(m-2-1)!0.55!(m-3-1)$) {$F'$};
  \node (Zp) at ([xshift=1.8em]$(m-2-3)!0.55!(m-3-3)$) {$Z'$};

  \draw[->] (m-2-1) to[bend right=15] node[left] {$q_{2}$} (Fp);
  \draw[->] (Fp) to[bend right=15] node[left] {$\nu_{2}$} (m-3-1);

  \draw[->] (m-2-3) to[bend left=15] node[right] {$q_{1}$} (Zp);
  \draw[->] (Zp) to[bend left=15] node[right] {$\nu_{1}$} (m-3-3);
\end{tikzpicture}
\]
Define the Bott--Chern classes $\beta'$ and $\gamma'$ by 
\[
\beta':=\nu_{1}^{*} (\alpha|_{Z''}) \quad \text{ and } \quad 
\gamma':=\nu_{2}^{*} (\alpha|_{F''}). 
\]
Note that both $\beta'$ and $\gamma'$ are K\"ahler classes since each $\nu_{i}$ is a finite morphism. 
From the diagram, we obtain 
\begin{equation}\label{eq-inj2}
\begin{aligned}
\pi_1^{*} (q_{1}^{*} \beta')
&=
\pi_1^{*} (q_{1}^{*}(\nu_{1}^{*} (\alpha|_{Z''} )))
= (\pi^{*} (q^{*} \alpha))|_{ \tilde{Z} \times \{*\} }
= \tilde \beta, \\
\pi_2^{*} (q_{2}^{*}\gamma')
&=
\pi_2^{*} (q_{2}^{*}(\nu_{2}^{*} (\alpha|_{F''})))
= (\pi^{*} (q^{*} \alpha))|_{ \{*\}\times \tilde{F} }
= \tilde \gamma.
\end{aligned}
\end{equation}
Here, the last equalities follow from \eqref{eq-dec}. 
Consequently, we deduce that $\tilde{\beta}$ and $\tilde{\gamma}$ are Bott--Chern classes. 
Moreover, by the commutative identities 
\[
\pi_2  \circ \tilde \pr_{2} = \pr_{2} \circ \pi
\text{ \quad and \quad }
\pi_1  \circ \tilde \pr_{1} = \pr_{1} \circ \pi, 
\]
we have 
\begin{align}\label{eq-inj3}
\pi^{*}(q^{*}\alpha)
&= \tilde\pr_{1}^{*}\tilde \beta+\tilde\pr_{2}^{*}\tilde\gamma \\
&=(\pi_1\circ \tilde\pr_{1})^{*} (q_{1}^{*}\beta') +(\pi_2\circ \tilde\pr_{2})^{*} (q_{2}^{*}\gamma') \notag \\
&=\pi^{*}
\bigl(\pr_{1}^{*} (q_{1}^{*}\beta') +
\pr_{2}^{*} (q_{2}^{*}\gamma') \bigr).  \notag
\end{align}
Since $W$ has rational singularities, 
the natural morphism 
\[
\pi^* \colon 
H^{2}(W, \mathbb{C}) \to H^{2}(\tilde{W}, \mathbb{C})
\]
is injective (for example, see \cite[Proposition B.2.7]{Kir15}). 
Hence, we can conclude that 
\begin{align}\label{eq-1}
q^{*}\alpha = \pr_{1}^{*} (q_{1}^{*}\beta') +
\pr_{2}^{*} (q_{2}^{*}\gamma'). 
\end{align}
\end{proof}

Consider the following fibrations:
\begin{align*}\label{eq-comm}
\xymatrix{
& W=Z \times F \ar[ld]_{q}  \ar[rd]^{\quad q':=q_{1} \times q_{2}} & \\
W' &  & Z' \times F'.
}
\end{align*}
We prove that there is an isomorphism 
$W' \cong Z' \times F'$ fitting into the desired commutative diagram.
To this end, we compare the subvarieties contracted by $q \colon W \to W'$ 
with those contracted by $q' \colon W \to Z' \times F'$ using the following claim.

\begin{claim}\label{claim-baby3}
Let $C \subset W$ be a subvariety. 
Then $C$ is contracted to a point by $q$ 
if and only if it is contracted to a point by $q'$. 
In particular, $W'$ is isomorphic to $Z' \times F'$ in a way compatible with the desired commutative diagram.
\end{claim}
\begin{proof}
Since $\alpha$ is a K\"ahler class, 
the subvariety $C$ is contracted to a point by $q \colon W \to W'$ 
if and only if $(q^{*}\alpha)|_{C}\equiv 0$. 
More precisely, the condition $(q^{*}\alpha)|_{C}\equiv 0$ is interpreted 
after pulling back to a resolution of singularities of the normalization of $C$, 
but we use this notation for simplicity.

On the other hand, by \eqref{eq-1} and the commutative diagram
\[
\xymatrix{
\{*\} \times F \ar[d]^{q_{2}}&  \ar[l]_{\pr_{2}}  W=Z \times F  \ar[r]^{\pr_{1}}   \ar[d]^{q':=q_{1} \times q_{2}}   
& Z \times \{*\}  \ar[d]^{q_{1}}\\
F' &  \ar[l]^{\pr'_{2}}  Z' \times F'  \ar[r]_{\pr'_{1}}     & Z',
}
\]
we obtain 
\begin{align}\label{eq-2}
q^{*}\alpha
&=
q'^{*}\Bigl(
\pr'^{*}_{1}\beta'
+
\pr'^{*}_{2}\gamma'
\Bigr). 
\end{align}
Since $q^{*}\alpha$ can be written as the pullback of the K\"ahler class 
$\pr'^{*}_{1}\beta'+\pr'^{*}_{2}\gamma'$ under $q'$, 
the condition $(q^{*}\alpha)|_{C}\equiv 0$ is equivalent to the condition that 
$C$ is contracted to a point by $q' \colon W \to Z' \times F'$. 
This proves Claim \ref{claim-baby3} and $(1)$ of Proposition \ref{prop-baby} when $Z$ is compact. 
\end{proof}

We now prove~$(2)$ of Proposition \ref{prop-baby} when $Z$ is compact. 
Suppose that $q \colon W \to W'$ is a bimeromorphic morphism. 
By construction, 
$q_{1}\colon Z \to Z'$ and 
$q_{2}\colon F \to F'$ are bimeromorphic morphisms. 
The compatibility with the boundary divisors follows directly from the commutative diagram in~$(1)$.

\medskip
We now consider the case where $Z$ is noncompact. 
Since the natural morphism $\tilde{j}_{1}$ in \eqref{eq-inj} 
need not be injective, the preceding argument requires some care. 
More precisely, the decomposition \eqref{eq-dec} must be replaced by  
\begin{align}\label{eq-dec'}
\tilde{j} \big( \pi^{*} (q^{*}\alpha) \big)
=
\tilde{\pr}_{1}^{*} \tilde{\beta}+\tilde{\pr}_{2}^{*}\tilde{\gamma} \quad \text{in }H^{2}(\tilde{W}, \mathbb{C}),
\end{align}
and \eqref{eq-inj2} must be replaced by  
\begin{equation}\label{eq-inj2'}
\begin{aligned}
\tilde{j}_{1} \big( \pi_1^{*} (q_{1}^{*} \beta') \big)
=
\tilde{j}_{1} \big(  \pi_1^{*} (q_{1}^{*}(\nu_{1}^{*} (\alpha|_{Z''} ))) \big)
&=\tilde{j}_{1} \big(    (\pi^{*} (q^{*} \alpha))|_{ \tilde{Z} \times \{*\} }\big)= \tilde{\beta}, \\
\tilde{j}_{2} \big( \pi_2^{*} (q_{2}^{*}\gamma') \big)
=
\tilde{j}_{2} \big(  \pi_2^{*} (q_{2}^{*}(\nu_{2}^{*} (\alpha|_{F''}))) \big)
&= \tilde{j}_{2}\big( (\pi^{*} (q^{*} \alpha))|_{ \{*\}\times \tilde{F} }\big)= \tilde{\gamma}. 
\end{aligned}
\end{equation}
Then, the decomposition \eqref{eq-inj3} is replaced by 
\begin{align}\label{eq-inj3'}
\tilde{j} \big( \pi^{*}(q^{*}\alpha) \big)
=\tilde{j} \big( \pi^{*}
\bigl(\pr_{1}^{*} (q_{1}^{*}\beta') +
\pr_{2}^{*} (q_{2}^{*}\gamma') \bigr) \big). 
\end{align}
Let $j\colon H^{1,1}_{\BottChern}(W)\to H^2(W,\mathbb{C})$ be the natural morphism. 
Note that $\tilde{j} \circ \pi^{*}= \pi^{*}\circ j$ 
and that
\[
\pi^{*}\colon H^{2}(Z\times F, \mathbb{C}) \to 
H^{2}(\tilde{Z}\times \tilde{F}, \mathbb{C})
\]
is injective since $Z \times F$ has rational singularities. 
Hence, \eqref{eq-inj3'} implies that there exists $K \in \Ker(j)$ such that
\begin{align}\label{eq-1'}
q^{*}\alpha = \pr_{1}^{*} (q_{1}^{*}\beta') +
\pr_{2}^{*} (q_{2}^{*}\gamma') + K.
\end{align}
Consequently, as in \eqref{eq-2}, we obtain 
\begin{align}\label{eq-2'}
q^{*}\alpha
&=
q'^{*}\Bigl(
\pr'^{*}_{1}\beta'
+
\pr'^{*}_{2}\gamma'
\Bigr) + K. 
\end{align}

In the proof of Claim~\ref{claim-baby3}, 
we use \eqref{eq-2'} after restricting it to compact subvarieties $C$,  
so $K$ does not affect the argument. 
Indeed, let $C \subset W$ be a subvariety contracted to a point 
by $q \colon W \to W'$ or by $q' \colon W \to Z' \times F'$. 
Let $\pi_{C} \colon \tilde{C}\to C$ be a resolution of the normalization of $C$, 
factoring through the normalization $\nu_C\colon C_n\to C$. 
We may choose a resolution of singularities $\bar \pi \colon \bar W \to W$ such that
$\tilde{C}$ is a subvariety of $\bar W$. 
The pullback $\pi_{C}^* (K|_C)=(\bar \pi^{*} K)|_{\tilde{C}}$ 
lies in the kernel of 
$j_C \colon H^{1,1}_{\BottChern}(\tilde{C}, \mathbb{C}) \to H^2(\tilde{C}, \mathbb{C})$. 
The map $j_C$ is injective since $\tilde{C}$ is a compact K\"ahler manifold. 
Thus, $\pi_{C}^*(K|_C)=0$. 
The same argument as in the proof of Claim~\ref{claim-baby3} applies.
\end{proof}

Corollary~\ref{cor-baby} follows from a slight modification of the proof of Proposition~\ref{prop-baby}. 

\begin{proof}[Proof of Corollary~\ref{cor-baby}]
Since $(W, \Delta)\to Z$ is a locally trivial fibration, after replacing $Z$ by a sufficiently small analytic open subset, 
we may assume that 
\[
(W, \Delta)=Z \times (F, \Delta_{F})
\]
and that the following diagram commutes:
\[
\xymatrix{
(W, \Delta)=Z \times (F, \Delta_{F}) \ar[rr]^{\qquad q} \ar[rd]_{\pr_{1}} & & W' \ar[ld] \\
& Z. &
}
\]
The additional point, compared with Proposition~\ref{prop-baby}, is that 
$q \colon W\to W'$ is a fibration over $Z$. 
Repeating the proof of Proposition~\ref{prop-baby}, we obtain fibrations
$q_{1}\colon Z\to Z'$ and $q_{2}\colon F\to F'$ together with an isomorphism
\[
W'\cong Z'\times F'
\]
under which $q$ is identified with $q_{1}\times q_{2}$.

Let $g=(g_1,g_2)\colon Z' \times F' \simeq W'\to Z$ denote the given structure morphism. 
Since $g\circ q=\pr_{1}$, $g_1\circ q_1=\id$ which implies that $q_1$ is injective. 
As $q_{1}$ is a fibration between normal varieties, it follows that $q_{1}$ is an isomorphism. 
After identifying $Z'$ with $Z$, the above isomorphism becomes an isomorphism over $Z$:
\[
W'\cong Z\times F'.
\]
This proves that $W'\to Z$ is a locally trivial fibration. 
If $q$ is bimeromorphic, the compatibility with the boundary divisor follows from Proposition~\ref{prop-baby}~$(2)$, yielding
\[
(W',\Delta')\cong Z\times(F',\Delta_F')
\]
locally over $Z$.
\end{proof}

\subsection{Descent: the general case} \label{subsec-gen}

We prove Theorem~\ref{thm-reduction}, 
which may be regarded as a generalization of the product case treated in Proposition~\ref{prop-baby} 
to the case of locally constant fibrations.

\begin{thm}\label{thm-reduction}
Let $q \colon (X, \Delta) \to (X', \Delta')$ 
be a small bimeromorphic morphism, 
where $q_{*}\Delta=\Delta'$, 
$(X, \Delta)$ is a compact K\"ahler klt pair,  
and $X'$ is a compact K\"ahler variety that is locally potentially klt. 
Assume that $X$ admits a fibration $\phi \colon X \to Y$ with the following properties: 
\begin{itemize}
\item The fibration $\phi \colon (X, \Delta) \to Y$ is projective and locally constant. 
\item The fiber $F$ of $\phi \colon X \to Y$ is rationally connected. 
\item The base variety $Y$ is Calabi--Yau $($i.e.,\,$K_{Y}$ is numerically trivial$)$.
\end{itemize}
Then after replacing $(X', \Delta')$ with a finite quasi-\'etale cover, 
the pair $(X', \Delta')$ admits 
a fibration $\phi' \colon X' \to Y'$ with the same properties as above. 
Moreover, there exists a small bimeromorphic morphism
$Y\to Y'$ fitting into the following commutative diagram:
\begin{equation*}\label{eq-comm}
\xymatrix{
X \ar[d]^{\phi} \ar[r]^{q}
&
X' \ar[d]^{\phi'}
\\
Y \ar[r]
&
Y'.
}
\end{equation*}
\end{thm}

\begin{proof}
We first make several reductions by proving Claims~\ref{claim1} and~\ref{claim2'}.
 
\begin{claim}\label{claim1}
By replacing $(X, \Delta)$ and $(X', \Delta')$ with finite quasi-\'etale covers, 
we may assume that $Y \cong Z \times T$, 
where $T$ is a complex torus and 
$Z$ is a klt Calabi--Yau variety with vanishing augmented irregularity $\hat{q}(Z)=0$. 
\end{claim}

\begin{proof}[Proof of Claim~\ref{claim1}]
By the singular Beauville--Bogomolov--Yau decomposition \cite[Theorem A]{BGL22}, 
there exists a finite quasi-\'etale cover $Y_{\rm{\text{q\'e}}} \to Y$ 
such that $Y_{\rm{\text{q\'e}}} \cong Z \times T$, 
where $T$ and $Z$ are as in the claim. 
Let $X_{\rm{\text{q\'e}}}$ be the normalization of the fiber product of 
$Y_{\rm{\text{q\'e}}} \to Y$ and $\phi \colon X \to Y$. 

The induced morphism $\sigma \colon X_{\rm{\text{q\'e}}} \to X$ is a finite quasi-\'etale cover, 
since $\phi \colon X \to Y$ is a locally constant fibration (hence equidimensional), 
and $Y_{\rm{\text{q\'e}}} \to Y$ is a finite quasi-\'etale cover. 
Let $X_{\rm{\text{q\'e}}} \to X'_{\rm{\text{q\'e}}} \to X'$ be the Stein factorization of 
$X_{\rm{\text{q\'e}}} \to X'$, 
as in the following diagram: 
\[
\xymatrix{
& X'_{\rm{\text{q\'e}}} \ar@/_1pc/[dl]_{\tau}& \\
X' & X \ar[l]_{q} \ar[d]_{\phi} & X_{\rm{\text{q\'e}}} \ar[l]_{\sigma} \ar[d] \ar@/_1pc/[ul]_{q_{\rm{\text{q\'e}}}} \ar@{}[dl]|{\square} \\
& Y & Y_{\rm{\text{q\'e}}}. \ar[l]
}
\]
The morphism $\tau \colon X'_{\rm{\text{q\'e}}} \to X'$ is a finite quasi-\'etale cover, 
and $q_{\rm{\text{q\'e}}} \colon X_{\rm{\text{q\'e}}} \to X'_{\rm{\text{q\'e}}}$ is a small bimeromorphic morphism, 
since $q \colon X \to X'$ is a small bimeromorphic morphism 
and $\sigma \colon X_{\rm{\text{q\'e}}} \to X$ is a finite quasi-\'etale cover. 

Define the Weil divisors $\Delta'_{\rm{\text{q\'e}}}$ and $\Delta_{\rm{\text{q\'e}}}$ by
$\Delta'_{\rm{\text{q\'e}}}:=\tau^{*}\Delta'$ and $\Delta_{\rm{\text{q\'e}}}:=\sigma^{*}\Delta$. 
These pullbacks are well-defined since $\tau$ and $\sigma$ are finite quasi-\'etale morphisms. 
By construction, we have $(q_{\rm{\text{q\'e}}})_{*}\Delta_{\rm{\text{q\'e}}}=\Delta'_{\rm{\text{q\'e}}}$. 
Furthermore, the induced morphism 
$(X_{\rm{\text{q\'e}}}, \Delta_{\rm{\text{q\'e}}})\to Y_{\rm{\text{q\'e}}}$ is 
projective and locally constant. 
Thus, the reduction in the claim is achieved by replacing all the data with the corresponding data carrying the subscript $\rm{\text{q\'e}}$. 
\end{proof}
\begin{lemm}
\label{lem-local-const-base-change}
Let $\phi\colon X\to Y$ be a fibration between normal analytic varieties,
and let $\Delta$ be a Weil $\RR$-divisor on $X$.
Assume that the morphism $\phi\colon (X,\Delta)\to Y$ is locally constant.
Let $f\colon Y'\to Y$ be a morphism between normal analytic varieties.
Then the base change
\[
\begin{tikzcd}
X':=X\times_Y Y' \arrow[r, "f'"] \arrow[d, "\phi'"] &
X \arrow[d, "\phi"] \\
Y' \arrow[r, "f"'] &
Y
\end{tikzcd}
\]
induces a locally constant fibration $(X',f'^*\Delta)\to Y'$.
\end{lemm}

\begin{proof}
By construction, $\phi'$ is an analytic fiber bundle with fiber $F$,
and every component of $f'^*\Delta$ is horizontal.
Since $\phi$ is locally constant, there exist a Weil $\RR$-divisor
$\Delta_F$ on $F$ and a representation
\[
\rho\colon\pi_1(Y)\to\Aut(F)
\]
such that:
\begin{itemize}
\item
the divisor $\Delta_F$ is invariant under the action of $\pi_1(Y)$;

\item
the pair $(X,\Delta)$ is isomorphic over $Y$ to the quotient
\[
\bigl(\Unv{Y}\times(F,\Delta_F)\bigr)/\pi_1(Y).
\]
Here $\Unv{Y}$ denotes the universal cover of $Y$, and
$\gamma\in\pi_1(Y)$ acts diagonally on $\Unv{Y}\times F$ by
\[
\gamma\cdot(y,z)=\bigl(\gamma\cdot y,\rho(\gamma)(z)\bigr).
\]
\end{itemize}
Let
$q\colon\Unv{Y}\times(F,\Delta_F)\to(X,\Delta)$ denote the quotient map.
Let $p'\colon\Unv{Y'}\to Y'$ be the universal covering. The holomorphic map
$f\circ p'\colon\Unv{Y'}\to Y$ admits a lift
$\widetilde f\colon\Unv{Y'}\to\Unv{Y}$. After choosing compatible base points,
$\widetilde f$ is $f_*$-equivariant, namely
\[
\widetilde f(\gamma'\cdot y')=f_*(\gamma')\cdot\widetilde f(y')
\]
for every $\gamma'\in\pi_1(Y')$. Set
$\rho':=\rho\circ f_*\colon\pi_1(Y')\to\Aut(F)$. Since $\Delta_F$ is invariant
under the action of $\pi_1(Y)$, it is also invariant under the action of
$\pi_1(Y')$ via $\rho'$.

Consider the holomorphic map
\[
\widetilde f'\colon\Unv{Y'}\times(F,\Delta_F)\to(X',f'^*\Delta),\qquad
(y',z)\longmapsto\bigl(p'(y'),q(\widetilde f(y'),z)\bigr).
\]
This is well-defined since
$\phi(q(\widetilde f(y'),z))=f(p'(y'))$. Moreover, it is
$\pi_1(Y')$-equivariant, where $\gamma'\in\pi_1(Y')$ acts on
$\Unv{Y'}\times F$ by
\[
\gamma'\cdot(y',z)=\bigl(\gamma'\cdot y',\rho'(\gamma')(z)\bigr).
\]
Indeed, this follows from the $f_*$-equivariance of $\widetilde f$ and the
$\pi_1(Y)$-invariance of the quotient map $q$. Thus $\widetilde f'$ induces
a morphism
\[
\bigl(\Unv{Y'}\times(F,\Delta_F)\bigr)/\pi_1(Y')
\longrightarrow(X',f'^*\Delta).
\]
The induced map is a fiberwise biholomorphism, and locally on $Y'$
 both sides are the same locally trivial bundle with fiber $(F,\Delta_F
)$. Hence it is an isomorphism. Hence
\[
(X',f'^*\Delta)\simeq
\bigl(\Unv{Y'}\times(F,\Delta_F)\bigr)/\pi_1(Y'),
\]
where the action of $\pi_1(Y')$ on $F$ is given by $\rho'=\rho\circ f_*$. Therefore
$(X',f'^*\Delta)\to Y'$ is locally constant.
\end{proof}
\begin{rem}\label{rem-replace}
(1) The proof of Claim~\ref{claim1} shows 
that, given a finite quasi-\'etale cover of $Y$, 
we may freely replace $(X, \Delta)$ and $(X', \Delta')$ with the associated finite quasi-\'etale covers. 
This will also be used in the proof of Claim~\ref{claim2'}. 

(2) $X_{\rm{\text{q\'e}}}$ cannot be expected to be strongly $\mathbb{Q}$-factorial 
even if $X$ is so. 
This is why we formulate Theorem~\ref{thm-reduction} without assuming that $X$ is strongly $\mathbb{Q}$-factorial. 
Instead of imposing $\mathbb{Q}$-factoriality, we assume that the fibration $X \to Y$ is projective, 
a property preserved by the induced fibration $X_{\rm{\text{q\'e}}} \to Y_{\rm{\text{q\'e}}}$.
\end{rem}

We now consider the representation associated with the locally constant fibration $\phi \colon (X, \Delta) \to Y$. 
\begin{claim}\label{claim2'}
Let $\rho \colon \pi_{1}(Y) \to \Aut(F, \Delta_{F})$ 
be the representation associated with the locally constant fibration $\phi \colon (X, \Delta) \to Y$, 
where $(F, \Delta_{F})$ is a fiber of $\phi \colon (X, \Delta) \to Y$. 
After replacing $(X, \Delta)$ and $(X', \Delta')$ with finite quasi-\'etale covers, 
we may assume that 
$\rho \colon \pi_{1}(Y) \to \Aut(F, \Delta_{F})$ factors through $\pi_{1}(T)$:
\[
\rho \colon \pi_{1}(Y) \xrightarrow{\pr_{2*}} \pi_{1}(T) \longrightarrow \Aut(F, \Delta_{F}),
\]
where $\pr_{2}\colon Y \cong Z \times T \to T$ denotes the second projection. 
\end{claim}

\begin{proof}[Proof of Claim~\ref{claim2'}]
We first show that $\Image(\rho)$ is contained in a linear algebraic group. 
Since $\phi \colon X \to Y$ is a projective fibration, there exists a $\phi$-relatively ample line bundle $A$. 
By Lemma~\ref{lemma_lc_line} applied to $A$, 
there exist an ample line bundle $A_F$ on $F$ and a positive integer $m$ 
such that $mA_F$ is $\rho$-equivariant and $\rho$-linearizable. 
Hence, the image $\Image(\rho)$ is contained in $\Aut(F, mA_F)$. 
The latter is a linear algebraic group of automorphisms preserving $mA_F$. 

Using the decomposition $Y \cong Z \times T$, 
consider the induced representation 
\[
\rho_{Z} \colon \pi_{1}(Z) \times \{1\} \subset \pi_{1}(Y) 
\xrightarrow{\rho} \Aut(F, \Delta_{F}).
\]
As shown above, the image $\Image(\rho_Z)$ is also contained in a linear algebraic group. 
Since the Calabi--Yau variety $Z$ is special in the sense of Campana, 
by \cite[Theorem~7.8]{Cam04}, we see that $\Image(\rho_{Z})$ is virtually abelian. 
After replacing $Y$ with a suitable finite \'etale cover, 
we may assume that $\Image(\rho_{Z})$ is abelian 
(see Remark~\ref{rem-replace}(1)). 

By the universal property of abelianization, 
the representation $\rho_{Z} \colon \pi_{1}(Z)\times \{1\} \to \Aut(F, \Delta_{F})$ 
factors through the abelianization of $\pi_{1}(Z)$. 
The abelianization of $\pi_{1}(Z)$ is $H_{1}(Z, \mathbb{Z})$. 
This group is finitely generated and torsion by the assumption $\hat{q}(Z)=0$, 
and is therefore finite. 
Thus, after replacing $Y$ by a further finite \'etale cover if necessary, 
we may assume that 
\[
\rho_{Z} \colon \pi_{1}(Z)\times \{1\} \to \Aut(F, \Delta_{F})
\]
is trivial. 
It follows that $\rho \colon \pi_{1}(Y) \to \Aut(F, \Delta_{F})$ factors through $\pi_{1}(T)$. 
\end{proof}

From now on, we assume that the representation 
$\rho \colon \pi_{1}(Y) \to \Aut(F, \Delta_{F})$ 
satisfies the conclusion of Claim~\ref{claim2'}. 
Consider the composition 
\[
\varphi:=\pr_{2} \circ \phi\colon X \to Y=Z \times T \to T.
\]
The rational map $\varphi\circ q^{-1}\colon X' \dashrightarrow T$ 
extends to a fibration 
\(
\psi \colon X' \to T
\)
by \cite[Corollary~1.4]{Fuj26} (see also \cite[Corollary~1.7]{HM07}), 
since $T$ contains no rational curves and $X'$ is locally potentially klt.

We now have the following diagram:
\[
\xymatrix{
X' \ar@/_1pc/[rdd]_{\psi} & X \ar[d]^{\phi} \ar[l]^{q}\ar@/^5pc/[dd]_{\varphi}& \\
& Y=Z \times T \ar[d]^{\pr_{2}} & \\
& T. &
}
\]
Define 
\[
W:= Z \times F 
\quad \text{ and }\quad 
\Delta_{W}:=\pr_{2}^{*}\Delta_{F}=\pr_{2}^{*}(\Delta|_{F}).
\]
Let $\rho_T\colon\pi_1(T)\to\Aut(F,\Delta_F)$ be the representation through which $\rho$ factors. 
Consider the representation
\begin{align}\label{eq-rep}
\tau \colon \pi_{1}(T) \xrightarrow{\rho_T} \Aut(F, \Delta_F) 
\xrightarrow{\id_{Z}\times \bullet} 
\Aut(W, \Delta_{W}).
\end{align} 
By Claim~\ref{claim2'}, the fibration
\(
\varphi \colon (X, \Delta)\to T
\)
is locally constant, its fiber is $(W=Z \times F, \Delta_W)$, 
and its associated representation is $\tau$.

We now base-change the above diagram under the universal cover 
$p_{T}\colon \Unv{T} \to T$. 
Let $\tilde{\bullet}$ denote the base change of $\bullet$ under 
$p_{T}\colon \Unv{T} \to T$. 
The definition of locally constant fibrations gives the following commutative diagram:
\begin{equation}\label{eq-big}
 \begin{aligned}
\begin{tikzpicture}[
  >=stealth,
  font=\normalsize,
  every node/.style={inner sep=2pt}
]

  \node (B) at (3.8,2.2) {$\tilde{X}\cong \Unv{T}\times W$};
  \node (C) at (7.2,2.2) {$W$};

  \node (E) at (-1.4,0) {$X'$};
  \node (F) at (1.4,0) {$X$};
  \node (D) at (3.8,0) {$\Unv{T}$};

  \node (G) at (1.4,-2.0) {$T$.};

  \draw[->] (B) -- node[above] {$\pr_{2}$} (C);
  \draw[->] (B) -- node[right] {$\pr_{1}$} (D);

  \draw[->] (F) -- node[above] {$q$} (E);
  \draw[->] (B) -- node[below left] {$p_{X}$} (F);

  \draw[->] (E) -- node[left] {$\psi$} (G);
  \draw[->] (F) -- node[right] {$\varphi$} (G);
  \draw[->] (D) -- node[below right] {$p_{T}$} (G);
\end{tikzpicture}
 \end{aligned}
\end{equation}

Fix a K\"ahler class $\alpha$ on $X'$. 
We have $H^{1}(W, \mathbb{C})=0$, 
since both $Z$ and $F$ have rational singularities and vanishing irregularity. 
Hence, by the K\"unneth theorem, there exist a Bott--Chern class $\beta$ on $\Unv{T}$ 
and a Bott--Chern class $\gamma$ on $W$ such that 
\begin{align}\label{eq-decom_uni}
p_{X}^{*} q^{*} \alpha = \pr_{1}^{*}\beta + \pr_{2}^{*}\gamma 
\end{align} 
in the second cohomology group. 
More precisely, since $Z$ may have singularities, 
this decomposition follows from the same argument as in Proposition~\ref{prop-baby}
after passing to resolutions of singularities (see the proofs of Claims~\ref{claim-baby1} and~\ref{claim-baby2}). 
As explained in the proof of Proposition~\ref{prop-baby}, 
the decomposition \eqref{eq-decom_uni} might not hold in Bott--Chern cohomology, 
but the difference $K$ is in the kernel of 
$j \colon H_{\BottChern}^{1,1}(\Unv{T} \times W) \to  H^{2}(\Unv{T} \times W, \mathbb{C})$. 
We use \eqref{eq-decom_uni} only after restricting it to a compact K\"ahler subvariety of $\Unv{T} \times W$. 
After pulling back to a resolution, the restriction of $K$ to any such subvariety vanishes. 
Thus, we can effectively use \eqref{eq-decom_uni} in Bott--Chern cohomology.

Fix a general point $* \in T$ and consider the base change of \eqref{eq-big} over $\{*\}$. 
The subscript $*$ denotes this base change; for example, $X_{*}$ denotes the fiber of $\varphi \colon X \to T$ over $*$. 
Writing $\Unv{T}_{*}=p_{T}^{-1}(*)$ as a discrete set of points $\coprod_{i\in I} \{t_{i}\}$, 
we have
$
\tilde{X}_{*} \cong \coprod_{i\in I} \{t_{i}\} \times W.
$
The base change of \eqref{eq-big} to $\{*\}$ is summarized by the following diagram, 
where we omit the symbols for the restrictions of morphisms for simplicity of notation.

\begin{equation}\label{eq-big-base}
\begin{aligned}
\begin{tikzpicture}[
  >=stealth,
  font=\normalsize,
  every node/.style={inner sep=2pt}
]

  \node (B) at (3.8,2.2) {$\tilde{X}_{*} \cong \coprod_{i\in I} \{t_{i}\} \times W$};
  \node (C) at (7.2,2.2) {$W$};

  \node (E) at (-1.4,0) {$X'_{*}$};
  \node (F) at (1.4,0) {$X_{*}$};
  \node (D) at (3.8,0) {$\Unv{T}_{*}=\coprod_{i\in I} \{t_{i}\}$};

  \node (G) at (1.4,-2.0) {$\{*\}$.};

  \draw[->] (B) -- node[above] {$\pr_{2}$} (C);
  \draw[->] (B) -- node[right] {$\pr_{1}$} (D);

  \draw[->] (F) -- node[above] {$q$} (E);
  \draw[->] (B) -- node[below left] {$p_{X}$} (F);

  \draw[->] (E) -- node[left] {$\psi$} (G);
  \draw[->] (F) -- node[right] {$\varphi$} (G);
  \draw[->] (D) -- node[below right] {$p_{T}$} (G);
\end{tikzpicture}
 \end{aligned}
\end{equation}
By restricting the decomposition \eqref{eq-decom_uni} to $\tilde{X}_{*}$, 
we obtain
\begin{align}\label{eq-ind}
p_{X}^{*} q^{*}(\alpha|_{X'_{*}})
=
(p_{X}^{*} q^{*} \alpha)|_{\tilde{X}_{*}}  =
(\pr_{2}^{*} \gamma)|_{\coprod_{i\in I} \{t_{i}\} \times W}.
\end{align}
Indeed, the term $\pr_{1}^{*}\beta$ vanishes after restriction to $\tilde{X}_{*}$, 
since $\tilde{X}_{*}$ lies over the discrete set 
$\Unv{T}_{*}=p_{T}^{-1}(*)$. 
For each $i$, the induced morphism 
$\pr_{2} \colon \{t_{i}\} \times W \to W$ is an isomorphism. 
Hence, the restrictions
\[
\{(\pr_{2}^{*} \gamma)|_{\{t_{i}\} \times W}\}_{i \in I}
\]
are all identified with $\gamma$ via these isomorphisms. 

From now on, we take $*$ as the base point for the fundamental group $\pi_{1}(T)$. 
Then, we have the following claim.
\begin{claim}\label{claim3}
The Bott--Chern class $\gamma$ on $W$ is $\tau$-invariant; that is, 
$g^{*}\gamma=\gamma$ for every $g\in \Image(\tau)$. 
\end{claim}

\begin{proof}[Proof of Claim~\ref{claim3}]
Fix an arbitrary automorphism $g \in \Image(\tau) \subset \Aut(W, \Delta_{W})$. 
Choose $\ell \in \pi_{1}(T,*)$ such that 
$g=\tau(\ell) \in \Aut(W, \Delta_{W})$.
The loop $\ell$ acts naturally on the universal cover $\Unv{T}$. 
Thus, for a given point $\tilde{t} \in p_{T}^{-1}(*) \subset \Unv{T}$, 
we obtain another point $\ell \cdot \tilde{t}$. 
By the definition of locally constant fibrations,
the isomorphism $g=\tau(\ell) \in \Aut(W, \Delta_{W})$ coincides with the composition 
\begin{align}\label{eq-isom}
W \cong \{\tilde{t}\}\times W
\xrightarrow[\cong]{\;p_X|_{\{ \tilde{t}\}\times W}\;}
X_{*}
\xleftarrow[\cong]{\;p_X|_{\{\ell\cdot \tilde{t}\}\times W}\;}
\{\ell\cdot \tilde{t}\}\times W \cong W.
\end{align}
Here the isomorphisms 
$W \cong \{\tilde{t}\}\times W$ and 
$\{\ell\cdot \tilde{t}\}\times W \cong W$ 
are induced by the second projection $\pr_{2}$. 
By \eqref{eq-ind}, under the left-hand part 
\[
W \cong \{\tilde{t}\} \times W \cong X_{*}
\]
of this isomorphism, the Bott--Chern class $\gamma$ on the left-hand copy of $W$ is identified with 
$q^{*}(\alpha|_{X'_{*}})$ on $X_{*}$. 
The same argument applies to the right-hand part
\[
X_{*} \cong \{\ell \cdot \tilde{t}\} \times W \cong W
\]
of the isomorphism. 
Hence, the class $\gamma$ on the left-hand copy of $W$ is identified with the same class on the right-hand copy of $W$ through the composition \eqref{eq-isom}. 
Since this composition coincides with $g$, we conclude that $g^{*}\gamma=\gamma$. 
\end{proof}

Fix an isomorphism $X_{*} \cong W$. 
We now apply Proposition~\ref{prop-baby} to 
\[
\xymatrix{
W':=X'_{*} \ar[dr]_{\psi} & & \ar[ll]_{q_W:=q |_{X_{*}}\quad} X_{*} \cong W  = Z \times F \ar[dl]^{\varphi}  \\
 & \{*\}.  &\\
}
\]
Then, there exist bimeromorphic morphisms $q_{1} \colon Z \to Z'$ and $q_{2} \colon (F, \Delta_{F}) \to (F', \Delta'_{F})$ 
onto normal compact K\"ahler varieties $Z'$ and $F'$ 
such that  
\[
(W', \Delta'_{W})\cong Z' \times (F', \Delta'_{F}),
\] 
where $\Delta'_{F}:=q_{2*}\Delta_{F}$ and $\Delta'_{W}:=(q_{*} \Delta)|_{X'_{*}}=\Delta'|_{X'_*}$. 
Note that $q_{1} \colon Z \to Z'$ and $q_{2} \colon (F, \Delta_{F}) \to (F', \Delta'_{F})$ 
are small bimeromorphic morphisms 
by the proof of Proposition~\ref{prop-baby}.  
We now prove the following claim.

\begin{claim}\label{claim4}
For every automorphism $g \in \Image(\tau) \subset \Aut(W, \Delta_{W})$, 
there exists a unique automorphism $g' \in \Aut(W', \Delta'_{W})$  
satisfying the following properties:
\begin{itemize}
\item[$(1)$]
$g'$ fits into the following commutative diagram:
 \[
\xymatrix{
W \ar[r]^g \ar[d]_{q_W}& W \ar[d]^{q_W}\\
W' \ar[r]_{g'}  & W'.
}
\] 
\item[$(2)$] $g'$ can be written as $\id_{Z'} \times g''$ for some $g''\in \Aut(F',\Delta'_{F})$. 
\end{itemize}
Consequently, we obtain a representation 
\[
\tau' \colon \pi_{1}(T) \to \Aut(F', \Delta'_{F}) 
\xrightarrow{\id_{Z'}\times \bullet} 
\Aut(W', \Delta'_{W}). 
\]
\end{claim}
\begin{proof}[Proof of Claim~\ref{claim4}]
Let $C \subset W$ be a subvariety contracted by $q_W$ to a point. 
To descend $g$ to $g'$, it is enough to show that the image $g(C)$ is also contracted by $q_W$ to a point. 
\eqref{eq-ind} shows that the Bott--Chern class $\gamma$ can be identified with 
\begin{align}\label{eq-gamma}
\gamma=(q|_{X_*})^{*}(\alpha|_{X'_{*}})=q_W^{*}(\alpha|_{X'_{*}})
\end{align}
under the identification \eqref{eq-isom}. 
Thus, since $\alpha|_{X'_{*}}$ is a K\"ahler class and $q_W(C)$ is a point, 
we have $\gamma|_{C}\equiv 0$. 
On the other hand, 
since $\gamma$ is $\tau$-invariant by Claim~\ref{claim3}, 
we obtain 
\[
(g|_C)^*(\gamma|_{g(C)})
=
(g^{*} \gamma)|_{C}
=
\gamma |_{C} \equiv 0. 
\]
This shows that $g(C)$ is also contracted by $q_W$ to a point. 
Since $q_W$ has connected fibers,  we see that $q_W\circ g$ factors uniquely through $q_W$. 
Applying the same argument to $g^{-1}$, we obtain an automorphism $g'$ making the above diagram commute. 

The automorphism $g'$ preserves the boundary divisor $\Delta'_{W}$ 
since $g$ preserves $\Delta_{W}$ and $(q_W)_{*}\Delta_{W} = \Delta'_{W}$. 
Furthermore, by \eqref{eq-rep}, the automorphism $g$ has the form $\id_Z\times g_F$ on $W \cong Z \times F$. 
Since $q_W=q_1\times q_2$, its descent has the form $g'=\id_{Z'}\times g''$ for some $g''\in\Aut(F',\Delta'_{F})$.
This proves the claim.
\end{proof}

The representation 
\(
\tau' \colon \pi_{1}(T) \to \Aut(W', \Delta'_{W}) 
\)
obtained from Claim~\ref{claim4} determines a locally constant fibration 
\[
X^{\circ}:=(\Unv{T}\times W')/{\pi_{1}(T)} \to T. 
\]
The representation $\tau'$ takes values in $\{\id_{Z'}\}\times\Aut(F',\Delta'_{F})$ by Claim~\ref{claim4}(2), 
and hence we obtain a locally constant fibration $X^{\circ} \to T \times Z'$. 
It remains to prove the following claim.
\begin{claim}\label{claim5}
$X'$ is isomorphic to $X^{\circ}$. 
\end{claim}

\begin{proof}[Proof of Claim~\ref{claim5}]
Consider the bimeromorphic morphism
\[
\Unv{T} \times (W, \Delta_{W}) 
\xrightarrow{\ \id_{\Unv{T}} \times q_W\ } 
\Unv{T} \times (W', \Delta'_{W})
\]
between the covering spaces of $X$ and $X^{\circ}$. 
Take a closed loop $\ell \in \pi_{1}(T)$. 
Put $g:=\tau(\ell)$ and $g':=\tau'(\ell)$, where $g'$ is the automorphism obtained in Claim~\ref{claim4}. 
By Claim~\ref{claim4}, the actions given by $(\ell,g)$ and $(\ell,g')$ are compatible with the bimeromorphic morphism 
$\id_{\Unv{T}} \times q_W$. 
Thus, the morphism descends to a bimeromorphic morphism 
$\sigma \colon X \to X^{\circ}$ 
between the quotients by the $\pi_{1}(T)$-actions, and we obtain the following commutative diagram:
\[
\xymatrix{
X \ar[r]^q \ar[d]_\sigma  & X' \ar[d]\\
X^{\circ} \ar[r]  & T.
}
\]

We prove that the morphisms $q$ and $\sigma$ contract precisely the same subvarieties, as in Claim~\ref{claim-baby3}. 
Let $C \subset X$ be a subvariety contracted to a point by $q \colon X \to X'$. 
Since $\alpha$ is a K\"ahler class on $X'$, we have $(q^{*}\alpha)|_{C}\equiv 0$. 
Let $*\in T$ be the image of $C$ under the morphism $X \to T$. 

Then $C$ is contained in the fiber $X_{*}$ of $X \to T$, 
and hence it can be regarded as a subvariety of $W$ via the identification of $X_{*}$ with $W$. 
Under this identification, \eqref{eq-ind} gives $\gamma|_{C}\equiv 0$. 
By \eqref{eq-gamma}, it follows that $C\subset W$ is contracted to a point by $q_W \colon W \to W'$. 
Therefore, by the construction of $X^{\circ}$, the subvariety $C$ is contracted to a point by 
$\sigma \colon X \to X^{\circ}$. 

Conversely, let $C \subset X$ be a subvariety contracted to a point by 
$\sigma \colon X \to X^{\circ}$. 
The corresponding subvariety of $W$ is contracted to a point by $q_W \colon W\to W'$, 
and hence $\gamma|_{C}\equiv 0$. 
On the other hand, the subvariety $C$ is contained in a fiber $X_{*}$ of $X\to T$. 
This shows that $(q^{*}\alpha)|_{C}\equiv 0$. 
Since $\alpha$ is a K\"ahler class on $X'$, this implies that $C$ is contracted to a point by 
$q \colon X \to X'$. 

Thus $q$ and $\sigma$ contract precisely the same subvarieties. 
Consequently, the induced bimeromorphic map $X' \dashrightarrow X^{\circ}$ is an isomorphism. 
This completes the proof of the claim.
\end{proof}

This completes the proof of Theorem~\ref{thm-reduction}. 
\end{proof}

As explained at the beginning of Section~\ref{sec-mainthm}, 
the proof of Theorem~\ref{thm-main} is completed by combining the results of 
Subsections~\ref{subsec-str} and~\ref{subsec-gen}.

\section{An analytic version of Druel's theorem}
\label{sec-druel}

In this section, we generalize Druel's theorem
\cite[Theorem~17]{DGP24} to the analytic setting. More precisely, we prove
the following theorem.
\begin{thm}
\label{thm-Druel}
Let $X$ be a compact normal analytic variety that is strongly
$\QQ$-factorial and has klt singularities, and let $\calF$ be a weakly
regular algebraically integrable foliation on $X$. See
Subsubsection~\ref{subsub-foliations-invariant}\textup{(A2)} for the
definition of weakly regular foliations.
\begin{itemize}
\item[\textup{(a)}] There exists a surjective equidimensional morphism
$p:X\to Y$ onto a compact normal analytic variety $Y$ such that
$\calF=T_{X/Y}$.
\item[\textup{(b)}] Moreover, there exists a Zariski open subset
$Y^\circ\subseteq Y$ whose complement has codimension $\geq 2$ such that
$X_y:=p^{-1}(y)$ is irreducible for every $y\in Y^\circ$.
\end{itemize}
\end{thm}

Recall that a foliation $\calF$ is called {\itshape algebraically
integrable} if every leaf $L^\circ$ of $\calF$ is open in its Zariski
closure. Here, the Zariski closure of $L^\circ$ means the smallest analytic
subvariety containing $L^\circ$ and is not necessarily its closure in the
analytic topology.

The main idea of the proof is similar to that of
\cite[Theorem 17]{DGP24}. Nevertheless, we need to check carefully that
each step in \cite[Proof of Theorem 17]{DGP24} remains valid in the analytic
setting and modify the argument when difficulties arise. Much of the
necessary preparatory work was carried out in
\S\,\ref{subsec-foliations}.

Before proving Theorem~\ref{thm-Druel}, we introduce some notation. As in
\cite[\S 4.2]{DGP24}, we consider a setting slightly more general than that
of Theorem~\ref{thm-Druel}.
\begin{setup}[see {\cite[Setup 18]{DGP24}}]
\label{setup-Druel}
Let $X$ and $Y$ be normal analytic varieties, and let
$p':X\dashrightarrow Y$ be a dominant meromorphic map
(see \cite[Definition 2.2, p.~13]{Ueno75}) with
$r:=\dim X-\dim Y\in\ZZ_{>0}$. Here, dominant means that the image of $p'$
contains a Zariski open subset of $Y$.

Let $Z_0$ be the closure of the graph of $p'$ in $Y\times X$, let
$\nu:Z\to Z_0$ be its normalization
of the main component, and let $p:Z\to Y$ and $q:Z\to X$ be the
natural morphisms. Let $\calF$ be the foliation induced by $p'$, and set
$\calF_Z:=q^{-1}\calF$.
\end{setup}

We first make the following observation.
\begin{prop}[see {\cite[Proposition 19(1) \& Corollary 20(1)]{DGP24}}]
\label{prop-alg_foliation_canonical}
Assume the notation and setting of Setup~\ref{setup-Druel}, and suppose
that $mK_{\calF}$ is a line bundle for some $m\in\ZZ_{>0}$. Then the Pfaff
field
$\eta_{\calF}:\Omega_X^r\to\OX(K_{\calF})$ associated with $\calF$ induces
a morphism
\[
(\Omega_Z^r)^{\otimes m}\to q^\ast\OX_X(mK_{\calF})
\]
that factors through the morphism
\[
(\Omega_Z^r)^{\otimes m}\to\OX_Z(mK_{\calF_Z})
\]
induced by the Pfaff field associated with $\calF_Z$. In particular, there
exists an effective $q$-exceptional $\QQ$-divisor $B\geq 0$ on $Z$ such
that
\[
K_{\calF_Z}+B\sim_{\QQ}q^\ast K_{\calF}.
\]
Moreover, if $X$ has klt singularities and $\calF$ is weakly regular, then
$\calF_Z$ is weakly regular and
$K_{\calF_Z}\sim_{\QQ}q^\ast K_{\calF}$.
\end{prop}

\begin{proof}
We first note that \cite{DGP24} proves a similar result under the
assumptions that $X$ is quasi-projective and $K_{\calF}$ is Cartier. Thus,
our result slightly generalizes
\cite[Proposition 19(1) \& Corollary 20(1)]{DGP24} to the setting where
$X$ is analytic and $K_{\calF}$ is $\QQ$-Cartier.

Consider the foliation
\[
\calG
:=
\pr_2^\ast\calF
\subseteq
\pr_2^\ast T_X
\subseteq
\pr_1^\ast T_Y\oplus\pr_2^\ast T_X
=
T_{Y\times X}
\]
on $Y\times X$. It gives rise to a Pfaff field
\[
\Omega_{Y\times X}^r
\simeq
\bigwedge^r(\pr_1^\ast\Omega_Y\oplus\pr_2^\ast\Omega_X)
\xrightarrow{\text{projection}}
\pr_2^\ast\Omega_X^r
\xrightarrow{\pr_2^\ast\eta_{\calF}}
\OX_{Y\times X}(\pr_2^\ast K_{\calF})
=
\OX_{Y\times X}(K_{\calG}),
\]
where $K_{\calG}=\pr_2^\ast K_{\calF}$. Clearly, $Z_0$ is
$\calG$-invariant, and hence we obtain a factorization
\[
(\Omega_{Y\times X}^r)^{\otimes m}|_{Z_0}
\twoheadrightarrow
(\Omega_{Z_0}^r)^{\otimes m}
\to
\OX_{Y\times X}(mK_{\calG})|_{Z_0}.
\]
By Seidenberg's theorem
(see \cite[Lemma 2.4.16]{Wang-thesis}), the Pfaff field
\[
(\Omega_{Z_0}^r)^{\otimes m}
\to
\OX_{Z_0}(mK_{\calG}|_{Z_0})
\]
extends to the normalization and gives rise to a generically surjective
morphism
\[
(\Omega_Z^r)^{\otimes m}
\to
\nu^\ast\OX_{Z_0}(m\pr_2^\ast K_{\calF}|_{Z_0})
=
q^\ast\OX_X(mK_{\calF}).
\]
Over $Z\backslash\Exc(q)$, this morphism factors through
\[
(\Omega_Z^r)^{\otimes m}\to\OX_Z(mK_{\calF_Z}),
\]
which is induced by the Pfaff field of $\calF_Z$. By
\cite[Lemma 2.4.14]{Wang-thesis}, the factorization holds over the locus
$Z^\circ$ where
$(\Omega_Z^r)^{\otimes m}\to\OX_Z(mK_{\calF_Z})$ is surjective. Since
$Z^\circ$ contains the locus where $\calF_Z$ is a subbundle of
$T_{Z_{\reg}}$, the complement $Z\backslash Z^\circ$ has codimension
$\geq 2$. The factorization therefore extends to $Z$ by reflexivity. In
particular, we obtain a natural morphism of rank-one reflexive sheaves
\[
\OX_Z(mK_{\calF_Z})
\to
q^\ast\OX_X(mK_{\calF}).
\]
This morphism is an isomorphism away from $\Exc(q)$. Consequently, there
exists an effective $q$-exceptional $\QQ$-Weil divisor $B\geq 0$ such that
\[
K_{\calF_Z}+B\sim_{\QQ}q^\ast K_{\calF}
\]
as $\QQ$-Weil divisors.

We now prove the final assertion. Assume that $X$ has klt singularities and
that $\calF$ is weakly regular. The discussion above, together with the
reflexive pullback from Theorem~\ref{thm-pullback-reflexive} and
\cite[Lemma 2.4.16]{Wang-thesis}, 
gives the following commutative
diagram:
\begin{center}
\begin{tikzpicture}[scale=1.8]
\node (A) at (0,0) {$(\Omega_Z^{[r]})^{[\otimes] m}$};
\node (A1) at (0,1) {$(q^\ast\Omega_X^{[r]})^{[\otimes] m}$};
\node (B) at (2,0) {$\OX_Z(mK_{\calF_Z})$};
\node (B1) at (2,1) {$q^\ast\OX_X(mK_{\calF})$};
\path[->,font=\scriptsize,>=angle 90]
(A1) edge node[left]{$(d_{\text{refl}}q)^{[\otimes] m}$} (A)
(B) edge (B1)
(A1) edge (B1)
(A) edge (B);
\end{tikzpicture}
\end{center}
Since $\calF$ is weakly regular, the twisted Pfaff field
\[
\Omega_X^r[\otimes]\OX_X(-K_{\calF})\to\OX_X
\]
is surjective. Taking its $m$-th tensor power and applying
\cite[Lemma 2.4.14]{Wang-thesis}, we obtain a surjective morphism
\[
(\Omega_X^{[r]})^{[\otimes]m}
\otimes\OX_X(-mK_{\calF})
\twoheadrightarrow
\OX_X.
\]
It follows that the right vertical morphism in the diagram is surjective
and hence is an isomorphism. In particular,
$K_{\calF_Z}\sim_{\QQ}q^\ast K_{\calF}$. Lemma
\ref{lemma-foliation_bimero}\textup{(b)} then implies that
$\calF_Z=q^{-1}\calF$ is weakly regular.
\end{proof}

The proof of Theorem~\ref{thm-Druel} is based on the following key lemma.
\begin{lemm}
\label{lemma-key-Druel}
Assume the notation and setting of Setup~\ref{setup-Druel}, and let
$B\geq 0$ be the effective $q$-exceptional $\QQ$-divisor such that
\[
K_{q^{-1}\calF}+B\sim_{\QQ}q^\ast K_{\calF},
\]
as obtained in Proposition~\ref{prop-alg_foliation_canonical}. If $E$ is a
$q$-exceptional prime divisor such that $p(E)=Y$, then
$E\subseteq\Supp(B)$. In other words, $\calF$ is not canonical along
$q(E)$.
\end{lemm}

\begin{proof}
Since the problem is local over $X$, we may replace $X$ by a small
neighbourhood of a general point of $q(E)$. We may thus assume that $X$ is
a locally analytic Stein subset of $\CC^N$ for some $N\in\ZZ_{>0}$ and
that
$\OX_X(mK_{\calF})\simeq\OX_X$ for some $m\in\ZZ_{>0}$. We then replace
$Z$ by $q^{-1}(X)$. Although the image of $Z$ in $Y$ need not be open,
the assumption $p(E)=Y$ allows us to shrink $Y$ further so that
$p:Z\to Y$ is flat, by Frisch's generic flatness theorem
(Theorem~\ref{thm-generic-flatness}), and hence open. We argue by induction
on $r$, which is the rank of $\calF$ by Setup~\ref{setup-Druel}.

\paragraph{Case $r=1$}
Let $f:X_1\to X$ be the associated cyclic cover, which is quasi-\'etale,
and let $Z_1$ be the normalization of the main component of
$X_1\times_XZ$. We obtain the following commutative
diagram:
\begin{center}
\begin{tikzpicture}[scale=2.0]
\node (A) at (0,0) {$Z$};
\node (A1) at (0,1) {$Z_1$};
\node (B) at (1.2,0) {$X$};
\node (B1) at (1.2,1) {$X_1$};
\node (C) at (0,-1) {$Y$};
\path[->,font=\scriptsize,>=angle 90]
(A1) edge node[left]{$g$} (A)
(B1) edge node[right]{$f$} (B)
(A1) edge node[above]{$q_1$} (B1)
(A) edge node[below]{$q$} (B)
(A) edge node[left]{$p$} (C);
\end{tikzpicture}
\end{center}
By the proof of Proposition~\ref{prop-alg_foliation_canonical}, there
exists an effective $q_1$-exceptional divisor $B_1\geq 0$ such that
\[
K_{q_1^{-1}(f^{-1}\calF)}+B_1
\sim_{\QQ}
q_1^\ast K_{f^{-1}\calF}.
\]
By our choice of $X$, we have $K_{\calF}\sim_{\QQ}0$. Since $f$ is
quasi-\'etale, Lemma~\ref{lemma-foliation_finite} gives
$K_{f^{-1}\calF}\sim_{\QQ}f^\ast K_{\calF}\sim_{\QQ}0$. Thus,
\begin{align*}
K_{q^{-1}\calF}+B &\sim_{\QQ}0,\\
K_{g^{-1}(q^{-1}\calF)}+B_1 &\sim_{\QQ}0.
\end{align*}
Moreover, Lemma~\ref{lemma-inv_divisor} gives
\[
K_{g^{-1}(q^{-1}\calF)}
\sim
g^\ast K_{q^{-1}\calF}
+
\Ramification(g)^{{\rm non}-q^{-1}\calF-{\rm inv}},
\]
where
$\Ramification(g)^{{\rm non}-q^{-1}\calF-{\rm inv}}$ denotes the
non-$q^{-1}\calF$-invariant part of the ramification divisor of $g$.
Consequently,
\[
g^\ast B
\sim_{\QQ}
B_1+\Ramification(g)^{{\rm non}-q^{-1}\calF-{\rm inv}}.
\]
The difference between the two sides is $q_1$-exceptional and numerically
trivial over $X_1$. Applying the negativity lemma to this difference and
its negative, we obtain
\[
g^\ast B
=
B_1+\Ramification(g)^{{\rm non}-q^{-1}\calF-{\rm inv}}.
\]
In particular,
$\Supp(B_1)\subseteq g^{-1}(\Supp(B))$.

Since $p(E)=Y$, the divisor $E$ is not $q^{-1}\calF$-invariant. Let $E_1$
be a prime divisor on $Z_1$ such that $g(E_1)=E$. Then $E_1$ is not
$g^{-1}(q^{-1}\calF)$-invariant. In view of the preceding inclusion, it
suffices to prove that $E_1\subseteq\Supp(B_1)$. We distinguish two cases:
\begin{itemize}
\item Suppose that $q_1(E_1)\subseteq(X_1)_{\sing}$. By Seidenberg's
theorem, Theorem~\ref{thm-Seidenberg}, $(X_1)_{\sing}$ is invariant under
every local derivation on $X_1$ and hence under $f^{-1}\calF$. Since $E_1$
is not invariant, Lemma~\ref{lemma-exc_invariant_rank1} implies that
$E_1\subseteq\Supp(B_1)$.
\item Suppose that $q_1(E_1)\cap(X_1)_{\reg}\neq \emptyset$. After
restricting to a neighbourhood of a general point of $q_1(E_1)$, we may
assume that $X_1$ is smooth. The subvariety $q_1(E_1)$ must then be
contained in $\Sing(f^{-1}\calF)$: otherwise, if $f^{-1}\calF$ were regular
at a general point $x\in q_1(E_1)$, then $q_1$ would be an isomorphism near
$x$. By Lemma~\ref{lemma-singular_invariant},
$\Sing(f^{-1}\calF)$ is $f^{-1}\calF$-invariant. Therefore,
Lemma~\ref{lemma-exc_invariant_rank1} again gives
$E_1\subseteq\Supp(B_1)$.
\end{itemize}
It follows that $E\subseteq\Supp(B)$.

\paragraph{Case $r\geq 2$}
Recall that we have reduced to the case where $X\subseteq\CC^N$ is a
locally analytic Stein subset. Let $V$ be the linear system on $X$
consisting of the restrictions of affine-linear functions on $\CC^N$.
Thus, the corresponding divisors are affine hyperplane sections, and $V$
is basepoint-free. Let $H$ be a component of $H_v$ for an analytically
sufficiently general $v\in V$, and set $H_Z:=q^{-1}(H)$. Then:
\begin{itemize}
\item by the Bertini-type theorem, Theorem~\ref{thm-Bertini}, both $H$ and
$H_Z$ are normal;
\item by the Bertini-type theorem for foliations,
Theorem~\ref{thm-Bertini-foliation},
$\calG:=\calF|_H\cap T_H$ is a foliation of rank $r-1$ on $H$, and
\[
q_H^{-1}\calG
=
T_{H_Z}\cap(\calF_Z|_{H_Z}),
\]
where $q_H:=q|_{H_Z}:H_Z\to H$.
\end{itemize}
Since the restriction of $q$ to a general fibre of $p$ is finite, we have
$\dim q(E)\geq r-1\geq 1$. Thus, $H$ meets $q(E)$, and
$E|_{H_Z}$ is a nonzero divisor on $H_Z$.

By Theorem~\ref{thm-Bertini-foliation}\textup{(b)} and \textup{(c)}, we
have
\[
K_{\calG}\sim(K_{\calF}+H)|_H.
\]
By Theorem~\ref{thm-Bertini-foliation}\textup{(a)}, we also have
\[
K_{q_H^{-1}\calG}
\sim
(K_{q^{-1}\calF}+H_Z)|_{H_Z}-\Delta_{H_Z}
\]
for some effective $q_H$-exceptional divisor
$\Delta_{H_Z}\geq 0$ on $H_Z$. Since
$K_{q^{-1}\calF}+B\sim_{\QQ}q^\ast K_{\calF}$, it follows that
\[
K_{q_H^{-1}\calG}
+B|_{H_Z}
+\Delta_{H_Z}
\sim_{\QQ}
q_H^\ast K_{\calG}.
\]
By the induction hypothesis,
\[
\Supp(E|_{H_Z})
\subseteq
\Supp(B|_{H_Z}+\Delta_{H_Z}).
\]

Moreover, by Theorem~\ref{thm-Bertini-foliation}\textup{(a)}, $\calF_Z$ is
tangent to $H_Z$ at a general point of every component of
$\Delta_{H_Z}$. If a component of $\Delta_{H_Z}$ dominated $Y$, then
$H_Z$ would contain a component of a general fibre of $p$, which is
impossible by the dimension formula. Thus, $\Delta_{H_Z}$ is
$p$-vertical. On the other hand, by the general choice of $H_v$, every
component of $E|_{H_Z}$ dominates $Y$. Hence,
\[
\Supp(E|_{H_Z})
\subseteq
\Supp(B|_{H_Z}),
\]
which implies that $E\subseteq\Supp(B)$.
\end{proof}

The remainder of the proof of Theorem~\ref{thm-Druel} is similar to that of
\cite[Theorem 17]{DGP24}. It follows from the next two lemmas, which are
analytic analogues of \cite[Lemma 21]{DGP24} and
\cite[Lemma 22]{DGP24}, respectively.

\begin{lemm}
\label{lemma-foliation_wr_irred}
Assume the notation and setting of Setup~\ref{setup-Druel}. Suppose that
$X$ has klt singularities and that $\calF$ is weakly regular. Then there
exists an open subset $Y^\circ\subseteq Y$ such that
$Y\backslash Y^\circ$ is analytically meagre of codimension $\geq 2$ and,
for every $y\in Y^\circ$, either $Z_y:=p^{-1}(y)$ is empty or every
connected component of $Z_y$ is irreducible. If $p$ is proper, then
$Y^\circ$ can be chosen to be Zariski open.
\end{lemm}

\begin{proof}
Suppose, for contradiction, that there exists a locally analytic subset
$D\subseteq Y$ of codimension $1$ such that, for a general $y\in D$,
$Z_y\neq\emptyset$ and some connected component of $Z_y$ is reducible.
Then we can find a locally analytic subvariety $S$ of $p^{-1}(D)$ such
that:
\begin{itemize}
\item $S$ is analytic in $p^{-1}(D)$ after removing an analytically meagre
subset;
\item for a general point $z\in S$, at least two irreducible components of
$Z_{p(z)}$ pass through $z$.
\end{itemize}
We choose $S$ to have maximal dimension among all such subvarieties. We
briefly explain its construction. Consider the normalization
\[
\nu_{p^{-1}(D)}:(p^{-1}(D))^\nu\to p^{-1}(D).
\]
After removing from $D$ an analytically meagre subset and removing its
$p$-preimage, we may assume that $\nu_{p^{-1}(D)}$ is also the fibrewise
normalization. Let $S'$ be the analytic subset consisting of the points
$z\in p^{-1}(D)$ at which $\nu_{p^{-1}(D)}$ is not an isomorphism. By
assumption, $S'$ has an irreducible component dominating $D$, and we take
$S$ to be such a component.

We proceed in three steps.

\paragraph{Step 1: Reduction and construction}
After removing a suitable locally analytic subset of codimension $\geq 2$
from $Y$, generic flatness, together with the irreducibility of $Z$ and
normality of $Y$, allows us to assume that $p$ is flat and $Y$ is smooth.
After replacing $X$ by a small open neighbourhood of $q(S)$, we may also
assume that
$\OX_X(mK_{\calF})\simeq\OX_X$ for some $m\in\ZZ_{>0}$. Let
$f:X_1\to X$ be the associated cyclic cover, which is quasi-\'etale, and
let $Z_1$ be the normalization of $X_1\times_XZ$, with natural morphisms
$g:Z_1\to Z$ and $q_1:Z_1\to X_1$.

We also need an auxiliary construction in Step 2. Since the problem is
local over $Y$, we may shrink $Y$ and construct a semistable reduction of
$p_1:=p\circ g$ in codimension $1$, denoted by $Y_2\to Y$. Let $Z_2$ be
the normalization of $Y_2\times_YZ_1$. We obtain the following commutative
diagram:
\begin{center}
\begin{tikzpicture}[scale=2.0]
\node (A) at (0,0) {$Z$};
\node (A1) at (0,1) {$Z_1$};
\node (B) at (1.2,0) {$X$};
\node (B1) at (1.2,1) {$X_1$};
\node (A2) at (-1.2,1) {$Z_2$};
\node (C) at (0,-1) {$Y$};
\node (C2) at (-1.2,-1) {$Y_2$};
\path[->,font=\scriptsize,>=angle 90]
(A1) edge node[right]{$g$} (A)
(B1) edge node[right]{$f$} (B)
(A1) edge node[above]{$q_1$} (B1)
(A) edge node[above]{$q$} (B)
(A) edge node[right]{$p$} (C)
(A2) edge node[above]{$g_1$} (A1)
(A2) edge node[left]{$p_2$} (C2)
(C2) edge (C)
(A1) edge[bend right] node[left]{$p_1$} (C);
\end{tikzpicture}
\end{center}
The morphism $Y_2\to Y$ is a cyclic cover ramified along the branch locus
of $p_1$ and is therefore Galois. Consequently, $g_1$ and
$g\circ g_1$ are also Galois.

\medskip

\paragraph{Step 2: Some regularity results}
Under the reductions made in Step 1, we prove the following:
\begin{itemize}
\item[(2A)] $\calF_Z$ is weakly regular.
\item[(2B)] $Z$ has klt singularities and therefore has quotient
singularities away from an analytic subset of codimension $\geq 3$.
\item[(2C)] $S$ has codimension $2$ in $Z$.
\end{itemize}
Assertion (2A) follows directly from
Proposition~\ref{prop-alg_foliation_canonical}, which also gives
\[
K_{\calF_Z}\sim_{\QQ}q^\ast K_{\calF}.
\]
Since $f$ is quasi-\'etale, $X_1$ has klt singularities by
\cite[Proposition 5.20, p.~160]{KM98}, and
$\calF_{X_1}:=f^{-1}\calF$ is weakly regular by
Lemma~\ref{lemma-foliation_finite}\textup{(b)}. The foliation
$\calF_{Z_1}:=q_1^{-1}\calF_{X_1}$ is induced by $p_1$, and $Z_1$ can be
identified with the normalization of the closure of the graph of the
meromorphic map $p_1\circ q_1^{-1}$. Proposition
\ref{prop-alg_foliation_canonical} again shows that $\calF_{Z_1}$ is
weakly regular and that
\[
K_{\calF_{Z_1}}\sim q_1^\ast K_{\calF_{X_1}}.
\]
By Lemma~\ref{lemma-foliation_Gorenstein_wr_canonical},
$\calF_{X_1}$ is canonical. The preceding equality and
Lemma~\ref{lemma-foliation_bimero}\textup{(a)} then imply that
$\calF_{Z_1}$ is canonical.

Every divisorial component of the branch locus of $g_1$ comes from a
divisor on $Y$ and is therefore $\calF_{Z_1}$-invariant. Hence,
Lemma~\ref{lemma-foliation_finite}\textup{(a)} implies that
$\calF_{Z_2}$ is canonical. Since $p_2$ has reduced fibres in codimension
$1$ and $Y_2$ is smooth, Lemma~\ref{lemma-canonical_foliation_variety}
implies that $Z_2$ has canonical singularities. Since
$g\circ g_1:Z_2\to Z$ is a finite Galois cover, the Riemann--Hurwitz
formula for Galois covers yields an effective $\QQ$-divisor $\Delta_Z$ on
$Z$ such that
\[
K_{Z_2}
\sim_{\QQ}
(g\circ g_1)^\ast(K_Z+\Delta_Z).
\]
It follows from \cite[Proposition 5.20, p.~160]{KM98} that
$(Z,\Delta_Z)$ is klt, and hence $Z$ has klt singularities. The final
assertion of (2B) now follows from
Proposition~\ref{prop-local_structure_klt_codim2}.

Since $Z$ is Cohen--Macaulay by (2B), Hartshorne's connectedness theorem,
Theorem~\ref{thm-connectedness}, implies that $S$ has codimension $2$ in
$Z$. This proves (2C).

\medskip

\paragraph{Step 3: Conclusion}
Let $z\in S$ be a general point. By (2B) and (2C), there exist an open
neighbourhood $U\subseteq Z$ of $z$, an open neighbourhood $W$ of the
origin in $\CC^{\dim Z}$, and a finite subgroup
$G\subseteq\GL_{\dim Z}(\CC)$ without quasi-reflections such that
$U\simeq W/G$. The quotient map
\[
\pi_U:W\to W/G\simeq U
\]
is \'etale over $U_{\reg}$ and hence is quasi-\'etale. Therefore,
\[
\calF_W:=\pi_U^{-1}(\calF_Z|_U)
\]
is a regular foliation on $W$ by
Lemma~\ref{lemma-foliation_finite}\textup{(b)}.

Let $F_1$ and $F_2$ be two distinct irreducible components of
$Z_{p(z)}$ passing through $z$. Since $S$ dominates $D$, neither $F_1$ nor
$F_2$ is contained in the singular locus of $\calF_Z$. Hence, both
$\pi_U^{-1}(F_1\cap U)$ and $\pi_U^{-1}(F_2\cap U)$ are disjoint unions of
leaves of $\calF_W$. Moreover,
\[
\pi_U^{-1}(F_1\cap U)\cap
\pi_U^{-1}(F_2\cap U)
\neq
\emptyset.
\]
A leaf of $\calF_W$ passing through a point of this intersection is a
connected component of both preimages. It follows that
$F_1\cap U=F_2\cap U$, contradicting the choice of $F_1$ and $F_2$.
\end{proof}

\begin{lemm}
\label{lemma-foliation_wr_projectable}
Assume the notation and setting of Setup~\ref{setup-Druel}. Suppose that
$X$ has klt singularities and that $\calF$ is weakly regular. Let $E$ be a
prime $q$-exceptional divisor on $Z$ that is not $p$-exceptional; that is,
$\dim p(E)\geq\dim Y-1$.
\begin{itemize}
\item[\textup{(a)}] Then $\dim p(E)=\dim Y-1$. In particular, $E$ is
$\calF_Z$-invariant.
\item[\textup{(b)}] If $z\in E$ is a general point, then there exists a
local curve $T\subseteq E$ passing through $z$ such that $\dim p(T)=1$ and
\[
q\bigl(E_{p(t_1)}(t_1)\bigr)
=
q\bigl(E_{p(t_2)}(t_2)\bigr)
\]
for every $t_1,t_2\in T$, where $E_{p(t)}(t)$ denotes the irreducible
component of
$E_{p(t)}:=E\cap Z_{p(t)}$ passing through $t\in T\subseteq E$.
\end{itemize}
\end{lemm}

\begin{proof}
Proposition~\ref{prop-alg_foliation_canonical} shows that $\calF_Z$ is
weakly regular and that
$K_{\calF_Z}\sim_{\QQ}q^\ast K_{\calF}$. By
Lemma~\ref{lemma-key-Druel}, $p(E)\subsetneqq Y$. Since $E$ is not
$p$-exceptional, it follows that $\dim p(E)=\dim Y-1$. This proves
\textup{(a)}.

Set $\Gamma:=q(E)$, and choose a suitable Zariski open subset
$E^\circ\subseteq E_{\reg}\cap Z_{\reg}$. After shrinking $X$, we may
assume that $\Gamma$ is smooth. To prove \textup{(b)}, it suffices to show
that the foliation $\calF_{E^\circ}$ induced by $\calF_Z$ on $E^\circ$ is
projectable under
\[
a:=q|_{E^\circ}:E^\circ\to\Gamma.
\]
We prove this claim below.

We may assume that $X$ is an open neighbourhood of a general point of
$q(E)$ such that
$\OX_X(mK_{\calF})\simeq\OX_X$ for some $m\in\ZZ_{>0}$. After passing to
the associated cyclic cover, we may further assume that $K_{\calF}$ is
Cartier. Indeed, let $f:X_1\to X$ be the cyclic cover, and let $Z_1$ be
the normalization of $Z\times_XX_1$, with natural morphisms
$g:Z_1\to Z$ and $q_1:Z_1\to X_1$. Then $X_1$ is klt by
\cite[Proposition 5.20, p.~160]{KM98}, and $f^{-1}\calF$ is weakly regular
by Lemma~\ref{lemma-foliation_finite}. Moreover, if
$p_1:=p\circ g$, then $Z_1$ can be identified with the normalization of
the closure of the graph of the meromorphic map
$p_1\circ q_1^{-1}:X_1\dashrightarrow Y$. For every prime divisor $E_1$
on $Z_1$ such that $g(E_1)=E$, the divisor $E_1$ is $q_1$-exceptional and
\[
\dim p(E)=\dim p_1(E_1).
\]
After this reduction, we continue to use the original notation.

Consider the following commutative diagram:
\begin{center}
\begin{tikzpicture}[scale=2.0,>=angle 90]
\node (A) at (0,0) {$\Omega_Z^r|_{E^\circ}=\Omega_Z^{[r]}|_{E^\circ}$};
\node (A1) at (0,1) {$a^\ast(\Omega_X^{[r]}|_\Gamma)\simeq q^\ast\Omega_X^{[r]}|_{E^\circ}$};
\node (B) at (2,0) {$\OX_Z(K_{\calF_Z})|_{E^\circ}$};
\node (B1) at (2,1) {$q^\ast\OX_X(K_{\calF})|_{E^\circ}$};
\node (C) at (-2,-1) {$\Omega_{E^\circ}^r$};
\node (C1) at (-2,2) {$a^\ast\Omega_\Gamma^r$};
\path[->,font=\scriptsize]
(A1) edge node[left]{$d_{\text{refl}}q|_{E^\circ}$} (A)
(B) edge node[right]{$\simeq$} (B1)
(C) edge[bend right=20] (B)
(C1) edge[bend left=20] (B1)
(C1) edge node[left]{$da$} (C);
\path[->>,font=\scriptsize]
(A1) edge (B1)
(A) edge (B)
(A1) edge node[below left]{$a^\ast(d_{\text{refl}}i)$} (C1)
(A) edge node[above left]{$d_{\text{refl}}j$} (C);
\end{tikzpicture}
\end{center}
Here:
\begin{itemize}
\item the square on the right is the diagram appearing in the proof of
Proposition~\ref{prop-alg_foliation_canonical}, in the case $m=1$;
\item the trapezoid on the left consists of the natural pullback morphisms
for reflexive differentials constructed in
Theorem~\ref{thm-pullback-reflexive};
\item the lower triangle arises from the fact that $E^\circ$ is
$\calF_Z$-invariant and from Lemma~\ref{lemma-invariant=Pfaff}. In
particular, we obtain a Pfaff field
\[
\Omega_{E^\circ}^r
\to
\OX_Z(K_{\calF_Z})|_{E^\circ}
\]
that induces a foliation $\calF_{E^\circ}$ on $E^\circ$;
\item the morphism
\[
a^\ast\Omega_\Gamma^r
\to
q^\ast\OX_X(K_{\calF})|_{E^\circ}
=
a^\ast(\OX_X(K_{\calF})|_\Gamma)
\]
is the composition
\[
a^\ast\Omega_\Gamma^r
\xrightarrow{da}
\Omega_{E^\circ}^r
\to
\OX_Z(K_{\calF_Z})|_{E^\circ}
\xrightarrow{\simeq}
q^\ast\OX_X(K_{\calF})|_{E^\circ}.
\]
\end{itemize}
We therefore obtain the following triangle:
\begin{center}
\begin{tikzpicture}[scale=2.0,>=angle 90]
\node (A) at (0,0) {$\Omega_\Gamma^r$};
\node (A1) at (0,1) {$\Omega_X^{[r]}|_\Gamma$};
\node (B1) at (2,1) {$\OX_X(K_{\calF})|_\Gamma$};
\path[->,font=\scriptsize]
(A) edge (B1);
\path[->>,font=\scriptsize]
(A1) edge (B1)
(A1) edge node[left]{$d_{\text{refl}}i$} (A);
\end{tikzpicture}
\end{center}
Its pullback by $a^\ast$ is the upper triangle in the preceding diagram.
In particular, we obtain a surjective Pfaff field
\[
\Omega_\Gamma^r
\twoheadrightarrow
\OX_X(K_{\calF})|_\Gamma,
\]
which induces a regular foliation $\calF_\Gamma$ on $\Gamma$, and
\[
da(\calF_{E^\circ})=a^\ast\calF_\Gamma.
\]
Thus, $\calF_{E^\circ}$ is projectable under $a$. The characterization of
projectable foliations in Definition~\ref{def-projectable} now yields
\textup{(b)}.
\end{proof}

Finally, we deduce Theorem~\ref{thm-Druel} from
Lemmas~\ref{lemma-foliation_wr_irred} and
\ref{lemma-foliation_wr_projectable}.
\begin{proof}[Proof of Theorem~\ref{thm-Druel}]
Apply Lemma~\ref{lemma-foliation_wr_irred} to the family of leaves
$p:Z\to Y$ of $\calF$. It remains to prove that the proper modification
$q:Z\to X$ is an isomorphism, or equivalently, that
$\Exc(q)=\emptyset$.

Suppose otherwise, and let $E$ be an irreducible component of $\Exc(q)$.
Since $X$ is strongly $\QQ$-factorial, $E$ is a divisor. By
Lemma~\ref{lemma-foliation_wr_irred}, $Z_y$ is irreducible for a general
$y\in p(E)$. Lemma~\ref{lemma-foliation_wr_projectable}\textup{(a)} then
implies that
\[
E=p^{-1}(p(E)).
\]
Moreover, Lemma~\ref{lemma-foliation_wr_projectable}\textup{(b)} yields a
local curve $S\subseteq p(E)$ passing through a general point of $p(E)$
such that
\[
q(Z_{s_1})=q(Z_{s_2})
\]
for every $s_1,s_2\in S$. Thus, the natural morphism
\[
Y\to\text{Chow}(X)
\]
contracts a positive-dimensional subvariety, contradicting the
construction of the family of leaves. This proves the theorem.
\end{proof}

\section{Applications of the Beauville--Bogomolov--Yau decomposition}  \label{sec-application}

\subsection{Splitting and uniformization}
\label{subsec-bby}

This subsection is devoted to the proof of Corollaries~\ref{cor-main1} and \ref{cor-bby}.

\begin{proof}[Proof of Corollary~\ref{cor-main1}]
After replacing $X$ with a finite quasi-\'etale cover,
we may assume that $X$ is maximally quasi-\'etale.
By Theorem~\ref{thm-Q-factorialization}, we take a small
$\mathbb{Q}$-factorialization
$q\colon X^{\qf}\to X$.

We first prove that $(X^{\qf},\Delta^{\qf})$ admits the desired
uniformization; that is, there exist a finite quasi-\'etale cover
$\tau\colon X'\to X^{\qf}$ and a possibly infinite \'etale cover
$p\colon X''\to X'$ such that $(X'',\Delta'')$ admits the desired
decomposition.
By Theorem~\ref{thm-main}, the pair
$(X^{\qf},\Delta^{\qf})$ admits a locally constant fibration
$f_{\qf}\colon (X^{\qf},\Delta^{\qf})\to Y^{\qf}$.
By the singular Beauville--Bogomolov--Yau decomposition \cite[Theorem A]{BGL22},
there exists a finite quasi-\'etale cover
$T\times Z\to Y^{\qf}$,
where $T$ is a complex torus and $Z$ is a product of strict
Calabi--Yau varieties and irreducible holomorphic symplectic varieties.
In particular, $Z$ has vanishing augmented irregularity.

Let $\tau\colon X'\to X^{\qf}$ be the finite quasi-\'etale cover
obtained by base change via $T\times Z\to Y^{\qf}$.
Since $X^{\qf}$ is strongly $\mathbb{Q}$-factorial,
the fibration $f_{\qf}$ admits a relatively $f_{\qf}$-ample line bundle.
The pullback of this line bundle is relatively ample for the induced
locally constant fibration $f'\colon X'\to T\times Z$.

Let
\[
\rho\colon
\pi_1(T\times Z)
\longrightarrow
\Aut(F',\Delta'_{F'})
\]
be the monodromy representation defining $f'\colon X'\to T\times Z$.
After replacing $X'$ with a further finite quasi-\'etale cover,
Claim~\ref{claim2'} shows that the composition
\[
\pi_1(Z)
\longrightarrow
\pi_1(T\times Z)
\xrightarrow{\rho}
\Aut(F',\Delta'_{F'})
\]
is trivial.
Therefore, the pair $(X',\Delta')$ admits a product decomposition
\[
(X',\Delta')
\cong
(W,\Delta_W)\times Z
\]
over $T\times Z$. Moreover, the pair $(W,\Delta_W)$ admits a locally constant
fibration $g\colon (W,\Delta_W)\to T$ such that 
$(X',\Delta')\to T\times Z$ coincides with $g\times\id_Z$.
Thus, after pulling back the fibration by the universal cover
$\mathbb{C}^d\to T$, we obtain a product decomposition
\[
(X'',\Delta'')
\cong
(F',\Delta'_{F'})\times Z\times\mathbb{C}^d.
\]
This proves the desired uniformization of
$(X^{\qf},\Delta^{\qf})$.

Let $X'\to W\to X$ be the Stein factorization of the composition
$X'\to X^{\qf}\to X$.
We have the following commutative diagram:
\[
\xymatrix@C=40pt@R=30pt{
X' \ar[r] \ar[d]_{\tau}
&
W \ar[d]
\\
X^{\qf} \ar[r]_{q}
&
X.
}
\]
Then $W\to X$ is a finite quasi-\'etale cover, whereas
$X'\to W$ is a small bimeromorphic morphism.
Since $X$ is maximally quasi-\'etale, the morphism $W\to X$ is
\'etale.

By the proof of Theorem~\ref{thm-reduction}, after replacing $W$
with a finite quasi-\'etale cover,
the pair $(W,\Delta_W)$ admits a locally constant fibration
$f_W\colon (W,\Delta_W)\to Y_W=T\times Z_W$ with fiber
$(F_W,\Delta_{F_W})$.
Moreover, if
$\rho_W\colon\pi_1(T\times Z_W)\to
\Aut(F_W,\Delta_{F_W})$
denotes the corresponding monodromy representation, then the composition
\[
\pi_1(Z_W)
\longrightarrow
\pi_1(T\times Z_W)
\xrightarrow{\rho_W}
\Aut(F_W,\Delta_{F_W})
\]
is trivial since $X' \to T\times Z$ has the same property.
Consequently, after pulling back the fibration by the universal cover
$\mathbb{C}^d\to T$, where $d:=\dim T$, we obtain a possibly infinite
\'etale cover $p\colon W''\to W$ and a product decomposition
\[
(W'',\Delta_{W''})
\cong
(F_W,\Delta_{F_W})\times Z_W\times\mathbb{C}^d.
\]
This finishes the proof. 
\end{proof}

\begin{proof}[Proof of Corollary~\ref{cor-bby}]
By replacing $(X, \Delta)$ with a maximally quasi-\'etale cover,
we may assume that $(X, \Delta)$ is maximally quasi-\'etale
and still satisfies the assumptions of Corollary~\ref{cor-bby}.
Consider a strongly $\mathbb{Q}$-factorialization
$
(X^{\qf}, \Delta^{\qf}) \to (X, \Delta).
$
Since $(X^{\qf}, \Delta^{\qf}) \to (X, \Delta)$ is a small bimeromorphic
morphism, the log canonical divisor
$K_{X^{\qf}}+\Delta^{\qf}$ is also numerically trivial.
By Proposition \ref{prop-baby},
it suffices to show that $(X^{\qf}, \Delta^{\qf})$ admits a product
structure after a finite quasi-\'etale base change of $Y^{\qf}$
(see also Remark~\ref{rem-replace}).

By Theorem~\ref{thm-main}, there exists a locally constant fibration
\[
f_{\qf}\colon
(X^{\qf}, \Delta^{\qf}) \to Y^{\qf}
\]
with the properties stated in Theorem~\ref{thm-main}.
Since $X^{\qf}$ is strongly $\mathbb{Q}$-factorial, this fibration is
projective by \cite[Corollary~4.2]{CH24}.
Thus, by \cite[Proposition~4.4, Theorem~4.7]{Amb05}, after a finite
quasi-\'etale base change, the fibration admits a splitting
\[
(X^{\qf}, \Delta^{\qf})|_U
\cong
(F^{\qf}, \Delta_{F^{\qf}}) \times U
\]
over a Zariski-open subset $U \subset Y^{\qf}$, where we retain the
same notation for the base-changed spaces.
Here the strongly $\mathbb{Q}$-factorial model is used to ensure that
the fibration is projective.

We claim that this splitting extends over the whole base $Y^{\qf}$
by the local triviality of
$(X^{\qf}, \Delta^{\qf}) \to Y^{\qf}$.
Indeed, since the fibration is locally trivial, for every sufficiently
small open subset $B\subset Y^{\qf}$, we may choose an isomorphism
\[
\Phi_B\colon
(X^{\qf},\Delta^{\qf})|_B
\cong
(F^{\qf},\Delta_{F^{\qf}})\times B
\]
over $B$.
Over $U\cap B$, this trivialization and the product decomposition above
differ by a meromorphic map
\[
\gamma_B\colon
B\dashrightarrow
\operatorname{Aut}
(F^{\qf},\Delta_{F^{\qf}}).
\]
After shrinking $B$ and composing $\Phi_B$ with a fixed automorphism
of the fiber, we may assume that $\gamma_B$ takes values in the identity
component
$\operatorname{Aut}^{0}(F^{\qf},\Delta_{F^{\qf}})$.
By \cite[Theorem~4.2]{Amb05}, we have that $\kappa(K_{F^{\qf}}+\Delta_{F^{\qf}})=0$. 
By \cite[Proposition~4.6]{Amb05}, this connected automorphism group is
an abelian variety and therefore contains no rational curves.
Since $B$ has klt singularities, the meromorphic map $\gamma_B$
extends to a morphism on $B$ by
\cite[Corollary~1.4]{Fuj26}
(see also \cite[Corollary~1.7]{HM07}).
Modifying $\Phi_B$ by this extended automorphism, we obtain a local
trivialization that agrees with the given product decomposition over
$U\cap B$.
These modified local trivializations agree on overlaps since they agree
over a dense open subset. They therefore glue to an isomorphism
\[
(X^{\qf},\Delta^{\qf})
\cong
(F^{\qf},\Delta_{F^{\qf}})\times Y^{\qf}.
\]
The conclusion follows from 
the singular Beauville--Bogomolov--Yau decomposition  \cite[Theorem A]{BGL22} applied to
$Y^{\qf}$.
\end{proof}

\subsection{Hacon--\texorpdfstring{M\textsuperscript c}{text}Kernan-type inequalities} 
\label{subsec-HM}

This subsection is devoted to the proof of Corollary~\ref{cor-HM}.
We first prepare the following proposition, 
which enables us to compare the Kodaira dimension
on the original variety with that on a finite \'etale cover.

\begin{prop}\label{prop-red}
Let $f\colon (X,\Delta)\to Y$ be a locally trivial fibration between
normal analytic varieties, and let
$\nu\colon X'\to X$ be a finite \'etale cover.
Consider the Stein factorization of the composition
$f\circ\nu\colon X'\to Y$:
\[
\xymatrix{
X'\ar[r]^{\nu}\ar[d]_{f'} & X\ar[d]^{f}\\
Y'\ar[r]_{\tau} & Y.
}
\]
Set $\Delta':=\nu^*\Delta$. Then the following assertions hold:
\begin{enumerate}
\item[$(1)$]
The morphism $\tau\colon Y'\to Y$ is a finite \'etale cover.

\item[$(2)$]
Assume that a fiber of $f$ is simply connected.
Then the natural morphism
\[
(\nu,f')\colon (X',\Delta')
\longrightarrow
(X,\Delta)\times_Y Y'
\]
is an isomorphism of pairs.
In particular, if we further assume that
$f\colon (X,\Delta)\to Y$ is a locally constant fibration, then
$f'\colon (X',\Delta')\to Y'$ is also a locally constant fibration.
\end{enumerate}
\end{prop}
\begin{proof}
Fix a point $y\in Y$. Take a sufficiently small contractible open
neighborhood $U$ of $y$ such that 
\[
(X,\Delta)|_U\simeq (F,\Delta_F)\times U.
\]
The induced morphsim 
\[
X'\times_Y U\longrightarrow F\times U
\]
is a finite \'etale cover. Since $U$ is contractible, every connected
finite covering of $F\times U$ is pulled back from a connected finite
covering of $F$. Hence, we have 
\[
X'\times_Y U
\simeq
\coprod_{\alpha=1}^r(F'_\alpha\times U),
\]
where each $F'_\alpha\to F$ is a connected finite \'etale cover.
Then, the Stein
factorization of $F'_\alpha\times U\to U$ is the morphism itself.
Consequently, we have 
\[
Y'\times_Y U\simeq\coprod_{\alpha=1}^rU.
\]
Thus $\tau\colon Y'\to Y$ is finite \'etale.

Assume now that $F$ is simply connected. Then every
$F'_\alpha\to F$ is an isomorphism. 
Thus, by the above local description, the natural morphism
\[
(\nu,f')\colon X'\longrightarrow X\times_Y Y'
\]
is an isomorphism. Since $\Delta'=\nu^*\Delta$, it is an isomorphism
of pairs. 
The final assertion follows by pulling back the local product presentation
of $f$ along the finite \'etale morphism $Y'\to Y$.
\end{proof}

We prove the following slight generalization of Theorem~\ref{thm-main}~\textup{(1)}.

\begin{cor}\label{cor-HM2}
In the setting of Theorem~\ref{thm-main},
after replacing $X$ with a finite quasi-\'etale cover, 
let $f\colon X\to Y$ be the resulting locally constant fibration,
and let $F$ be a fiber of $f$.
Then
\[
\kappa(-(K_{X}+\Delta))
\leq
\kappa(-(K_{F}+\Delta|_{F})).
\]
\end{cor}

\begin{proof}[Proof of Corollary~\ref{cor-HM2}]

Let $f\colon X\to Y$ be the locally constant fibration obtained in Theorem~\ref{thm-main}. 
We first observe that we can replace $X$ with a finite \'etale cover. 
Indeed, let $\nu \colon X' \to X$ be a finite \'etale cover, 
and consider the situation as in Proposition \ref{prop-red}. 
The Kodaira dimension (as well as the numerical dimension) of $-(K_X+\Delta)$ 
is unchanged under the finite \'etale cover by $K_{X'}+\Delta'=\nu^*(K_X+\Delta)$.
Note that the fiber of $f\colon X\to Y$ is rationally connected, and thus simply connected. 
Proposition \ref{prop-red} shows that $(F, \Delta_F) \cong (F', \Delta_{F'})$, 
where $F$ and $F'$ are the fibers of $f \colon X \to Y$ and $f' \colon X' \to Y'$, respectively. 
The same argument using Theorem \ref{thm-reduction} shows that 
we can replace $X$ with a small $\mathbb{Q}$-factorialization $X^{\qf} \to X$. 
Therefore, we may assume that $X$ is maximally quasi-\'etale and strongly $\mathbb{Q}$-facotiral.

If $\kappa(-(K_X+\Delta))=-\infty$, there is nothing to prove.
Fix a sufficiently divisible positive integer $m$ and an arbitrary
nonzero section
\[
s\in H^0(X,-m(K_X+\Delta)),
\]
and set
$
D:=\frac{1}{m}\operatorname{div}(s).
$
Then $D$ is an effective $\mathbb{Q}$-divisor such that
$
D\sim_{\bb{Q}}-(K_X+\Delta).
$
For $0<\delta\ll 1$, the generalized pair
$
(X,\Delta+\delta D,(1-\delta)\boldsymbol{\beta})
$
is still klt and of Calabi--Yau type.
Indeed, let $\pi\colon M\to X$ be a log resolution of
$(X,\Delta+D,\boldsymbol{\beta})$, and consider the formula
\eqref{eq-gpair} for the canonical divisors.
By the (classical) negativity lemma, there exists an effective
$\pi$-exceptional divisor $E$ such that
\[\boldsymbol{\beta}_{M}
-\pi^{*}\boldsymbol{\beta}_{X}
=-E.
\]
Writing, as in \eqref{eq-gpair},
\begin{align*}
\pi^*(K_X+(\Delta+\delta D)+(1-\delta)\boldsymbol{\beta}_{X})
+F_{+,\delta}
\equiv
K_M+\Delta_{M,\delta}
+(1-\delta)\boldsymbol{\beta}_{M},
\end{align*}
we obtain
\[
\Delta_{M,\delta}
=
\Delta_M+\delta\pi^*D-\delta E.
\]
It follows that the new generalized pair remains klt for all
sufficiently small $\delta>0$.
Moreover,
\[
K_X+\Delta+\delta D+(1-\delta)\boldsymbol{\beta}_X
\equiv
\delta(D-\boldsymbol{\beta}_X)
\equiv 0.
\]

Applying Theorem~\ref{thm-main} to this new generalized pair, we see
that $f\colon (X, \Delta + \delta D)\to Y$ remains a locally constant fibration.
Note that the locally constant MRC fibration does not depend on the boundary
(see Subsection~\ref{subsec-str}).
In particular, the local product description shows that $D$ contains
no fiber of $f$.
Consequently, the restriction $s|_F$ is nonzero.
Hence, the restriction
map
\[
H^{0}(X,-m(K_X+\Delta))
\longrightarrow
H^{0}(F,-m(K_F+\Delta|_F))
\]
is injective.
Comparing the asymptotic growth of the dimensions, we obtain
\[
\kappa(-(K_X+\Delta))
\leq
\kappa(-(K_F+\Delta|_F)).
\]
\end{proof}

\begin{proof}[Proof of Theorem~\ref{cor-HM}~\textup{(2)}]
As explained at the beginning of the proof of
Corollary~\ref{cor-HM2}, we may replace $(X,\Delta)$ with a finite
quasi-\'etale cover and a small $\mathbb{Q}$-factorialization.
Therefore, we may assume that $X$ is maximally quasi-\'etale and
strongly $\mathbb{Q}$-factorial.

Let $f\colon (X,\Delta)\to Y$ be the locally constant MRC fibration.
Set
\[
L:=-(K_X+\Delta),\qquad
L_F:=L|_F=-(K_F+\Delta|_F),
\]
and let $n:=\dim X$, $m:=\dim Y$, and $d:=n-m=\dim F$.
Recall the following diagram:
\begin{equation*}\label{locally-constant2}
\begin{gathered}
\xymatrix@C=40pt@R=30pt{
X \ar[d]_{f}
&
X\times_Y\Unv{Y}\cong\Unv{Y}\times F
\ar[d]_{\pr_1}
\ar[l]^{p_X\qquad\qquad}
\ar[r]^{\qquad\qquad\pr_2}
&
F
\\
Y
&
\ar[l]^{p_Y}\Unv{Y}.
}
\end{gathered}
\end{equation*}

The fibration admits a relatively $f$-ample line bundle.
By Lemma~\ref{lemma_lc_line}, after replacing it with a sufficiently
divisible multiple and twisting it by the pullback of a line bundle
on $Y$, we obtain a relatively $f$-ample $\mathbb{Q}$-line bundle
$A_f$ on $X$ such that
\[
p_X^*A_f\sim_{\mathbb{Q}}\pr_2^*A_F,
\qquad
A_F:=A_f|_F.
\]
Moreover, by the local constancy of $f$ with respect to $\Delta$, together
with $K_Y\sim_{\mathbb{Q}}0$, we see that 
$$
p_X^*L\sim_{\mathbb{Q}}\pr_2^*L_F.
$$
After taking sufficiently divisible multiples, we may assume that $L$, $A_f$, $L_F$, and $A_F$ are line
bundles. Both linear equivalences are $\pi_1(Y)$-equivariant. 

We first show that the desired conclusion follows from the following
intersection-theoretic vanishing:
\begin{equation}\label{eq-pushforward-vanishing}
\int_X
c_1(L)^a c_1(A_f)^b
\wedge f^*\omega_Y^{n-a-b}
=0
\qquad
\text{if }d<a+b\leq n,
\end{equation}
for some K\"ahler form $\omega_Y$ on $Y$.
All intersection numbers on singular spaces are understood after
pulling back to resolutions.

Fix a K\"ahler form $\omega_Y$ for which
\eqref{eq-pushforward-vanishing} holds.
Since $A_f$ is relatively ample, after replacing $\omega_Y$ with a
sufficiently large positive multiple, the class
\[
\{\omega_X\}:=c_1(A_f)+f^*\{\omega_Y\}
\]
is a K\"ahler class on $X$.
Set $r:=\nu(L_F)$.
Since $L_F$ is nef and $A_F$ is ample, we have
\[
\int_F c_1(L_F)^r c_1(A_F)^{d-r}>0.
\]
Expanding $c_1(L)^r\cdot\{\omega_X\}^{n-r}$, every term containing
$(f^*\{\omega_Y\})^k$ with $k>m$ vanishes for dimensional reasons.
If $k<m$, then $r+(n-r-k)=n-k>d$, so the corresponding term vanishes
by \eqref{eq-pushforward-vanishing}.
Consequently, we obtain 
\begin{align*}
c_1(L)^r\cdot\{\omega_X\}^{n-r}
&=
\binom{n-r}{m}
c_1(L)^r c_1(A_f)^{d-r}
\cdot f^*\{\omega_Y\}^m\\
&=
\binom{n-r}{m}
\left(
\int_F c_1(L_F)^r c_1(A_F)^{d-r}
\right)
\left(
\int_Y\omega_Y^m
\right)>0.
\end{align*}
Here the second equality follows from the projection formula.
It follows that $\nu(L)\geq r=\nu(L_F)$.

Conversely, set $s:=\nu(L)$.
Then $c_1(L)^s\cdot\{\omega_X\}^{n-s}>0$.
Suppose that $s>d$.
Since $n-s<m$, every term in the binomial expansion of this
intersection number has the form
\[
c_1(L)^s c_1(A_f)^{n-s-k}
\cdot f^*\{\omega_Y\}^k
\]
with $k<m$.
Moreover, $s+(n-s-k)=n-k>d$, so every such term vanishes by
\eqref{eq-pushforward-vanishing}, a contradiction.
Thus $s\leq d$.

As above, in the expansion of
$c_1(L)^s\cdot\{\omega_X\}^{n-s}$, all terms vanish except the term
containing $f^*\{\omega_Y\}^m$.
Therefore, we obtain 
\begin{align*}
0
<
c_1(L)^s\cdot\{\omega_X\}^{n-s}
=
\binom{n-s}{m}
\left(
\int_F c_1(L_F)^s c_1(A_F)^{d-s}
\right)
\left(
\int_Y\omega_Y^m
\right).
\end{align*}
It follows that
$\int_F c_1(L_F)^s c_1(A_F)^{d-s}>0$, and hence
$s\leq\nu(L_F)$.
This finishes the proof. 

It remains to prove \eqref{eq-pushforward-vanishing}.
To this end, we reduce the problem to the smooth case and apply the
 Grothendieck--Riemann--Roch formula.
The vanishing in \eqref{eq-pushforward-vanishing} and the numerical
dimensions under consideration are unchanged after passing to a
finite quasi-\'etale cover.
Thus, by the argument in the proof of Corollary~\ref{cor-main1}, after
passing to a further finite quasi-\'etale cover, we may assume that
\[
Y\cong T\times Z
\quad \text{ and }\quad 
(X,\Delta)\cong(W,\Delta_W)\times Z
\]
over $T\times Z$.
Here $T$ is a complex torus, $K_Z$ is trivial, and
$(W,\Delta_W)$ admits a projective locally constant MRC fibration
$f_W\colon(W,\Delta_W)\to T$ such that the fiber is $(F,\Delta|_F)$ 
and $f=f_W\times\id_Z$.
There exist line bundles $L_W$ and $A_{f,W}$ on $W$
such that
\[
L=\pr_W^*L_W,
\qquad
A_f=\pr_W^*A_{f,W},
\]
where $\pr_W\colon W\times Z\to W$ is the first projection.
Set $w:=\dim W$ and $t:=\dim T$, so that $w=d+t$.
Choose K\"ahler forms $\omega_T$ and $\omega_Z$ on $T$ and $Z$,
respectively, and set
$\omega_Y:=\pr_T^*\omega_T+\pr_Z^*\omega_Z$.

If $a+b>w$, then the left-hand side of
\eqref{eq-pushforward-vanishing} vanishes for dimensional reasons.
For $a+b\leq w$, by the K\"unneth projection formula, we have 
\eqref{eq-pushforward-vanishing} to
\begin{equation}\label{eq-pushforward-vanishing-W}
\int_W
c_1(L_W)^a c_1(A_{f,W})^b
\wedge f_W^*\omega_T^{w-a-b}
=0
\qquad
\text{if }d<a+b\leq w.
\end{equation}

The base $T$ is smooth, whereas $W$ may be singular.
Since $f_W$ is locally constant, the singularities of $W$ arise from
those of the fiber $F$.
Let $\mu_F\colon\tilde F\to F$ be a functorial resolution of
singularities \cite[Theorem~5.3.2]{Tem12}.
In particular, since every automorphism of $F$ lifts functorially to
$\tilde F$, we obtain 
\[
\Aut(F)\longrightarrow\Aut(\tilde F).
\]
Composing the monodromy representation of $f_W$ with this
homomorphism, we obtain a smooth locally constant fibration
$f_{\tilde W}\colon\tilde W\to T$ with fiber
$\tilde F$, together with a resolution
$\pi_{\tilde W}\colon\tilde W\to W$.
Set
\[
\tilde L:=\pi_{\tilde W}^*L_W,
\qquad
\tilde A:=\pi_{\tilde W}^*A_{f,W}.
\]
By the definition of the intersection numbers on $W$, it is enough
to prove
\begin{equation}\label{eq-pushforward-vanishing-tildeW}
\int_{\tilde W}
c_1(\tilde L)^a c_1(\tilde A)^b
\wedge f_{\tilde W}^*\omega_T^{w-a-b}
=0
\qquad
\text{if }d<a+b\leq w.
\end{equation}

To prove \eqref{eq-pushforward-vanishing-tildeW}, 
we first check that all higher direct images
\[
R^qf_{\tilde W,*}
\mathcal{O}_{\tilde W}(u\tilde L+v\tilde A)
\]
are flat vector bundles on $T$.
Let $p_T\colon\mathbb{C}^t\to T$ be the universal cover.
The locally constant structure gives
\[
\tilde W\times_T\mathbb{C}^t
\cong
\mathbb{C}^t\times\tilde F.
\]
Put
$\tilde L_F:=\tilde L|_{\tilde F}$ and
$\tilde A_F:=\tilde A|_{\tilde F}$.
The equivariant identifications coming from the local constancy and
Lemma~\ref{lemma_lc_line} give
\[
p_{\tilde W}^*
\mathcal{O}_{\tilde W}(u\tilde L+v\tilde A)
\cong
\pr_2^*
\mathcal{O}_{\tilde F}
(u\tilde L_F+v\tilde A_F)
\]
as $\pi_1(T)$-linearized line bundles, where
$p_{\tilde W}\colon\mathbb{C}^t\times\tilde F
\to\tilde W$ is the natural covering map.
Flat base change then gives
\begin{align*}
p_T^*
R^qf_{\tilde W,*}
\mathcal{O}_{\tilde W}(u\tilde L+v\tilde A)
&\cong
H^q\bigl(
\tilde F,
\mathcal{O}_{\tilde F}
(u\tilde L_F+v\tilde A_F)
\bigr)
\otimes_{\mathbb{C}}\mathcal{O}_{\mathbb{C}^t}.
\end{align*}
The monodromy linearization induces a representation
\[
\rho_{q,u,v}\colon
\pi_1(T)
\longrightarrow
\operatorname{GL}\left(
H^q\bigl(
\tilde F,
\mathcal{O}_{\tilde F}
(u\tilde L_F+v\tilde A_F)
\bigr)
\right),
\]
and the above isomorphism is $\pi_1(T)$-equivariant.
Consequently,
\[
R^qf_{\tilde W,*}
\mathcal{O}_{\tilde W}(u\tilde L+v\tilde A)
\cong
\left(
\mathbb{C}^t\times
H^q\bigl(
\tilde F,
\mathcal{O}_{\tilde F}
(u\tilde L_F+v\tilde A_F)
\bigr)
\right)\big/\pi_1(T).
\]
Thus, it is a flat vector bundle on $T$.

Since both $\tilde W$ and $T$ are smooth, 
by the  Grothendieck--Riemann--Roch formula (for example, see \cite{OTT85, Wu23c}), we 

\xr{reduce}
\begin{align}
f_{\tilde W,*}\left(
e^{u c_1(\tilde L)+v c_1(\tilde A)}
\operatorname{td}(T_{\tilde W/T})
\right)
&=
\operatorname{ch}\left(
\sum_q(-1)^q
\left[
R^qf_{\tilde W,*}
\mathcal{O}_{\tilde W}
(u\tilde L+v\tilde A)
\right]
\right).
\label{eq-singular-grr}
\end{align}
Since the higher direct images on the right-hand side are flat, their
rational Chern classes of positive degree vanish.
Consequently, the right-hand side equals
\[
\chi\bigl(
\tilde F,
\mathcal{O}_{\tilde F}
(u\tilde L_F+v\tilde A_F)
\bigr).
\]
By the classical Hirzebruch--Riemann--Roch theorem on the smooth
variety $\tilde F$, this Euler characteristic is a polynomial in
$u$ and $v$ of total degree at most $d$.

Fix $a,b\geq0$ and set $c:=a+b$, where $d<c\leq w$.
Pair both sides of \eqref{eq-singular-grr} with
$\omega_T^{w-c}$ and integrate over $T$.
The coefficient of $u^av^b$ on the right-hand side vanishes because
the above Euler characteristic has total degree at most $d$.
On the left-hand side, only the degree-zero component of the Todd class
$\operatorname{td}(T_{\tilde W/T})$ contributes for dimensional
reasons.
Thus, we obtain
\[
\int_{\tilde W}
c_1(\tilde L)^a c_1(\tilde A)^b
\wedge f_{\tilde W}^*\omega_T^{w-c}
=0.
\]
This proves \eqref{eq-pushforward-vanishing-tildeW}, and hence
\eqref{eq-pushforward-vanishing-W}.

\end{proof}

\subsection{On generic nefness of tangent sheaves}\label{subsec-gen-new}

In this subsection, we prove Theorems~\ref{nefness} and~\ref{thm-second},
which extend the results of \cite{Ou23} and \cite{IM22}, respectively,
to K\"ahler generalized pairs.
The proofs rely on the following result of Cao--P\u{a}un \cite[Theorem~1.2]{CP25}.

\begin{thm}[{cf.~\cite[Theorem~1.2]{CP25}}]
\label{new positivity*}
Let $f\colon X\to Y$ be a fibration between compact K\"ahler manifolds.
Let $\theta$ be a $d$-closed smooth $(1,1)$-form on $Y$, and let
$L$ be a $\QQ$-line bundle on $X$.
Assume that the following conditions hold:
\begin{enumerate}[label=\textbullet]
    \item The $\QQ$-line bundle $L$ admits a singular Hermitian metric
    $h_L$ such that
    $\sqrt{-1}\Theta_{h_L}\geq f^*\theta$ and
    $\mathcal{I}(h_L)=\mathcal{O}_X$.
    \item The restriction $(K_X+L)|_{X_y}$ is pseudo-effective for a
    general fiber $X_y$.
    \item There exist a K\"ahler form $\omega_X$ and a holomorphic
    $2$-form $\sigma$ on $X$ such that
    $\omega_X+\sigma+\bar{\sigma}$ represents a rational cohomology
    class and $\sigma|_{X_y}=0$ for a general fiber $X_y$.
\end{enumerate}
Then, the $\QQ$-line bundle $K_{X/Y}+L-R(f)$ is
$f^*\theta$-positively curved; that is, there exists a $(1,1)$-current
$T\in c_1(K_{X/Y}+L-R(f))$ such that $T\geq f^*\theta$.
\end{thm}
Here $R(f)$ denotes the ramification divisor
\[
R(f)=\sum_{i=1}^k(m_i-1)W_i,
\]
where the $W_i$ are the irreducible divisors on $X$ whose images
$Z_i:=f(W_i)$ are divisors on $Y$, and $m_i$ is the multiplicity of
$f^*Z_i$ along $W_i$.

\begin{proof}
The proof is essentially the same as that of
\cite[Theorem~1.2]{CP25}, so we explain only the necessary
modifications.
Let $y$ be a general point of $Y$, and choose a sufficiently small
neighborhood $U$ of $y$ such that $\theta=dd^c\varphi$ on $U$.
The metric $h_Le^{f^*\varphi}$ has semipositive curvature on
$f^{-1}(U)$, and its multiplier ideal remains trivial.
Applying \cite[Lemma~3.3 and Theorem~3.4]{CP25}, we obtain a
semipositively curved metric $h_{\min,\varphi}$ on
$(K_{X/Y}+L)|_{f^{-1}(U)}$.
These results also apply to $\QQ$-line bundles by passing to a
sufficiently divisible multiple and rescaling the resulting weights.

In the present setting, the local metrics $h_{\min,\varphi}$ with
minimal singularities need not glue to a global metric on
$K_{X/Y}+L$.
Nevertheless, the same computation as in the proof of
\cite[Theorem~1.2]{CP25} shows that the twisted metrics
$h_{\min,\varphi}e^{-f^*\varphi}$ are compatible on overlaps.
They glue on a Zariski-open subset of $X$, and their
curvature currents are bounded from below by $f^*\theta$.
The Riemann type extension theorem produces a global curvature
current $T_0\in c_1(K_{X/Y}+L)$ satisfying $T_0\geq f^*\theta$.

It remains to account for the ramification divisor.
The computation at a general point of each $W_i$ in the proof of
\cite[Theorem~1.2]{CP25} shows that the generic Lelong number of
$T_0-f^*\theta$ along $W_i$ is at least $m_i-1$.
Consequently, the Siu decomposition gives
\[
T_0-f^*\theta-[R(f)]\geq 0.
\]
Thus $T:=T_0-[R(f)]$ belongs to
$c_1(K_{X/Y}+L-R(f))$ and satisfies $T\geq f^*\theta$.
\end{proof}

We now prove the generic nefness for K\"ahler generalized pairs.
The following result extends to polarizations defined by arbitrary
movable classes, but the treatment of singular varieties requires
some care (see \cite[Setting~6.2]{IJZ25}).
Therefore, we consider only polarizations arising from K\"ahler forms.

\begin{thm}[{cf.~\cite[Theorem~1.4]{Ou23}}]
\label{nefness}
Let $(X,\Delta,\boldsymbol{\beta})$ be a K\"ahler generalized klt pair
of Calabi--Yau type.
Let $\mathcal{Q}$ be a torsion-free quotient of the tangent sheaf
$T_X$:
\[
T_X \longrightarrow \mathcal{Q} \longrightarrow 0.
\]
Let $\alpha=\{\omega_1\wedge\cdots\wedge\omega_{n-1}\}$, where the
$\omega_i$ are K\"ahler forms on $X$.
Then, we have 
\[
c_1(\mathcal{Q})\cdot\alpha\geq 0.
\]
\end{thm}

\begin{rem}\label{rem-inter}
When $X$ is not strongly $\mathbb{Q}$-factorial, the first Chern class
$c_1(\mathcal{Q})$ may not be defined as a Bott--Chern class.
Nevertheless, it can be defined as a functional on classes represented
by $d$-closed forms as follows.
Choose a modification $\pi\colon M\to X$ such that $M$ is smooth and
$(\pi^*\mathcal{Q})/\tor$ is locally free, where $\tor$ denotes the
maximal torsion subsheaf of $\pi^*\mathcal{Q}$.
For a class represented by a $d$-closed $(n-1,n-1)$-form $\gamma$, we define
\[
\int_X c_1(\mathcal{Q})\wedge\gamma
:=
\int_M
c_1\bigl((\pi^*\mathcal{Q})/\tor\bigr)\wedge\pi^*\gamma.
\]
This definition is independent of the choice of the modification
$\pi\colon M\to X$.
\end{rem}

\begin{proof}
We first reduce the problem to the case where $X$ is strongly
$\mathbb{Q}$-factorial.
To this end, take a strongly $\mathbb{Q}$-factorialization
$q\colon X^{\qf}\to X$ and assume that the assertion has been proved
for $X^{\qf}$.
Note that $X^{\qf}$ is also the underlying space of a K\"ahler
generalized klt pair of Calabi--Yau type.

Let $\omega^{\qf}$ be a K\"ahler form on $X^{\qf}$.
For $\delta>0$, consider the class
\[
\beta_{\delta}
:=
\bigl\{
(q^*\omega_1+\delta\omega^{\qf})
\wedge\cdots\wedge
(q^*\omega_{n-1}+\delta\omega^{\qf})
\bigr\},
\]
which is represented by a product of K\"ahler forms on $X^{\qf}$.
Taking reflexive pullbacks, we obtain a generically surjective morphism
\[
\gamma\colon
T_{X^{\qf}}
\cong q^{[*]}T_X
\longrightarrow q^{[*]}\mathcal{Q}.
\]
The sheaf $\Image (\gamma)$ need not coincide with
$q^{[*]}\mathcal{Q}$.
However, since $q$ is small, their quotient is supported in codimension
at least two, and hence
$c_1(\Image (\gamma))=c_1(q^{[*]}\mathcal{Q})$ (see Proposition \ref{det-tor} and Proposition \ref{det-exact}).
Applying the assertion on $X^{\qf}$ to the torsion-free quotient
$\Image (\gamma)$, we obtain
$
c_1(q^{[*]}\mathcal{Q})\cdot\beta_{\delta}\geq 0.
$
Letting $\delta\to0$, we conclude that
\[
c_1(\mathcal{Q})\cdot\alpha
=
c_1(q^{[*]}\mathcal{Q})\cdot q^*\alpha
\geq0.
\]

For the remainder of the proof, we assume that $X$ is strongly
$\mathbb{Q}$-factorial.
Suppose, to the contrary, that
$\mu_{\alpha,\min}(T_X)<0$, where $\mu_{\alpha,\min}(T_X)$ denotes
the minimal slope, defined as the infimum of
$\mu_{\alpha}(\mathcal{Q})$ over all nonzero torsion-free quotients
$T_X\to\mathcal{Q}\to0$.
Together with
\[
-K_X\cdot\alpha
=
(\Delta+\boldsymbol{\beta}_X)\cdot\alpha
\geq0,
\]
this shows that $T_X$ is not $\alpha$-semistable.
Hence, it admits a nontrivial Harder--Narasimhan filtration
\[
0=\mathcal{E}_0\subsetneqq\mathcal{E}_1\subsetneqq\cdots
\subsetneqq\mathcal{E}_r=T_X
\]
with $r\geq2$.
We have
$\mu_{\alpha}(T_X/\mathcal{E}_{r-1})
=\mu_{\alpha,\min}(T_X)<0$, whereas
$\mu_{\alpha}(T_X)\geq0$.
Let $k\in\{1,\ldots,r-1\}$ be the smallest index such that
$\mu_{\alpha}(T_X/\mathcal{E}_k)<0$, and set
$\mathcal{F}:=\mathcal{E}_k$.
Then, we have 
\[
K_X\cdot\alpha>K_{\mathcal{F}}\cdot\alpha.
\]

By the minimality of $k$, we have
$\mu_{\alpha}(T_X/\mathcal{E}_{k-1})\geq0$.
It follows that
\[
\mu_{\alpha,\min}(\mathcal{F})
=
\mu_{\alpha}(\mathcal{E}_k/\mathcal{E}_{k-1})
>0.
\]
A standard slope argument 
 (see e.g. \cite[Lemma 4.10]{CP19}) then shows that $\mathcal{F}$ is closed
under the Lie bracket. Thus, the subsheaf $\mathcal{F}\subset T_X$ is a foliation on $X$.

Let $\mu\colon X'\to X$ be a resolution of singularities.
By \cite[Corollary~1.5]{CP25}, the pullback foliation
$\mu^{-1}\mathcal{F}$ is algebraically integrable, and the closures of
its leaves are rationally connected.
It follows that $\mathcal{F}$ is algebraically integrable and that the
closures of its leaves are rationally connected.
Let $f\colon X\dashrightarrow Y'$ be the dominant meromorphic map
associated with $\mathcal{F}$, where $Y'$ is smooth.
By the universal property of the MRC quotient, there exists a dominant
rational map $g\colon Y'\dashrightarrow Y$ such that
$\psi=g\circ f$, where $\psi\colon X\dashrightarrow Y$ is the MRC
quotient in Setting~\ref{setting}.

After replacing $Y'$ with a suitable compact K\"ahler manifold, we may assume that $g$ is holomorphic.
We may also take a modification $\pi\colon M\to X$ that resolves the
indeterminacies of $f$ and $\psi$, as in the following diagram:
\[
\begin{tikzcd}[column sep=7em, row sep=5em]
M
  \arrow[r, "\pi"]
  \arrow[d, "\phi"']
  \arrow[dr, "\tilde f"]
&
X
  \arrow[d, dashed, "f"]
\\
Y
&
Y'.
  \arrow[l, "g"']
\end{tikzcd}
\]

By $K_X\cdot\alpha>K_{\mathcal{F}}\cdot\alpha$, we obtain
\begin{align}\label{b}
K_M\cdot\pi^*\alpha
=
K_X\cdot\alpha
>
K_{\tilde{\mathcal{F}}}\cdot\pi^*\alpha
=
K_{\mathcal{F}}\cdot\alpha,
\end{align}
where $\tilde{\mathcal{F}}:=T_{M/Y'}\subset T_M$ is the foliation
induced by $\mathcal{F}$ via $\pi$.

As in \cite[Lemma~5.12]{CP19}, we have
\[
K_{\tilde{\mathcal{F}}}
=
K_{M/Y'}-R(\tilde f)-E
\]
for some $\tilde f$-exceptional divisor $E$, where
$R(\tilde f)$ is the ramification divisor of $\tilde f$.
By Proposition~\ref{prop_bir-geometry-psi}~\textup{(b)}, the divisor
$E$ is also $\pi$-exceptional.
By \cite[Theorem~8.5]{Ou25}, there exist a K\"ahler form $\omega_M$
and a holomorphic $2$-form $\sigma$ on $M$ such that
$[\omega_M+\sigma+\bar{\sigma}]\in H^2(M,\mathbb{Q})$.
Since a general fiber $F$ of $\tilde f$ is rationally connected, we
have $\sigma|_F=0$.
The projectivity criterion \cite[Corollary~4.2]{CH24} 
implies that $\tilde f$ is projective.
Hence, there exists a relatively $\tilde f$-ample line bundle $L$
on $M$.

Recall that
\[
\boldsymbol{\beta}_M
\equiv
F_+-K_M-\Delta_M
\]
for some effective $\pi$-exceptional $\RR$-divisor $F_+$, and that
$\boldsymbol{\beta}_M$ is nef.
After a further modification, we may assume that
$\Delta_M+F_+$ has simple normal crossing support.
For $m\gg1$, set
\[
\mathcal{L}_m
:=
-K_M+\frac{1}{m}\lceil mF_+\rceil+\frac{1}{m}L.
\]
We apply Theorem~\ref{new positivity*} to $\mathcal{L}_m$.
To verify its curvature assumption, write
\[
\mathcal{L}_m
=
(F_+-K_M-\Delta_M)+\frac{1}{m}L
+
\left(
\Delta_M+\frac{1}{m}\lceil mF_+\rceil-F_+
\right).
\]
Choose a K\"ahler form $\omega_{Y'}$ on $Y'$.
Since $L$ is relatively $\tilde f$-ample, we can choose a smooth
metric $h_L$ on $L$ and a constant $C>0$ such that
\[
\sqrt{-1}\Theta_{h_L}
+C\tilde f^*\omega_{Y'}
\geq
\omega_M.
\]
The nefness of $F_+-K_M-\Delta_M$ allows us to choose a metric whose
curvature is bounded from below by $-(1/m)\omega_M$.
It follows that the first two terms in the above decomposition admit
a metric whose curvature is bounded from below by
$-(C/m)\tilde f^*\omega_{Y'}$.

For sufficiently large $m$, the divisor
\[
B_m
:=
\Delta_M+\frac{1}{m}\lceil mF_+\rceil-F_+
\]
has coefficients strictly smaller than one.
Its canonical singular metric has semipositive curvature and trivial
multiplier ideal.
Thus, $\mathcal{L}_m$ admits a singular Hermitian metric $h_m$ such
that
\[
\sqrt{-1}\Theta_{h_m}
\geq
-\frac{C}{m}\tilde f^*\omega_{Y'},
\qquad
\mathcal{I}(h_m)=\mathcal{O}_M.
\]
Moreover, the restriction of $K_M+\mathcal{L}_m$ to a general fiber
of $\tilde f$ is pseudo-effective.
The remaining assumption of Theorem~\ref{new positivity*} follows
from the choice of $\omega_M$ and $\sigma$ above.

Applying Theorem~\ref{new positivity*} with
$\theta=-(C/m)\omega_{Y'}$, we obtain a current
\[
T_m
\in
c_1\left(
K_{M/Y'}-K_M+\frac{1}{m}\lceil mF_+\rceil
+\frac{1}{m}L-R(\tilde f)
\right)
\]
such that $T_m\geq-(C/m)\tilde f^*\omega_{Y'}$.
After passing to a subsequence and letting $m\to\infty$, we conclude
that
\[
-\tilde f^*K_{Y'}+F_+-R(\tilde f)
\]
is pseudo-effective.

Consequently, by standard result (see, for example,
\cite[ Lemma 2.31]{CKT16}
), the class
$
K_{\tilde{\mathcal{F}}}-K_M+F_++E
$
is pseudo-effective.
Since $F_++E$ is $\pi$-exceptional, we obtain
\[
0
\leq
\bigl(K_{\tilde{\mathcal{F}}}-K_M+F_++E\bigr)
\cdot\pi^*\alpha
=
\bigl(K_{\tilde{\mathcal{F}}}-K_M\bigr)\cdot\pi^*\alpha,
\]
contradicting \eqref{b}.
\end{proof}

We now consider the case where $X$ is smooth.
The first result gives an affirmative solution to the Kodaira problem
for compact K\"ahler manifolds with nef anti-canonical bundle,
in the strong sense that properties~(1) and~(3) below hold.

\begin{proof}[Proof of Corollary~\ref{cor-kodaira}]
By combining \cite[Theorem~D]{NW23} with Theorem~\ref{thm-main},
we obtain an algebraic approximation of $X$ along which the
anti-canonical bundles of the fibers remain nef.
Although sufficiently small deformations of compact K\"ahler manifolds
remain K\"ahler, the associated deformation morphism need not be a
K\"ahler morphism.
However, in our setting, the construction in \cite[Theorem~D]{NW23}
yields an algebraic approximation that is a K\"ahler morphism.
Indeed, consider the MRC fibration $f \colon X \to Y$, which is shown
to be a locally constant fibration by Theorem~\ref{thm-main}.
We construct an algebraic approximation of $X$ by deforming its complex
structure through the coarse moduli space of polarized Calabi--Yau
manifolds constructed in \cite[Theorem~5.4]{FS90}, as follows.
There exists an algebraic approximation
$\mathfrak{Y} \to \mathbb{D}$ of $Y$, where $\mathbb{D}$ is a small
disc in $\mathbb{C}$, that is a K\"ahler morphism by \cite{FS90}.

By the proof of Proposition~\ref{prop_flat-lcf}, there exists a flat vector
bundle $E$ over $Y$ together with a relative embedding
$X \hookrightarrow \mathbb{P}(E)$ over $Y$.
Since the inclusion $Y=\mathfrak{Y}_0 \hookrightarrow \mathfrak{Y}$
induces an isomorphism
$\pi_1(Y) \cong \pi_1(\mathfrak{Y})$, the bundle $E$ extends to a flat
vector bundle $\mathfrak{E}$ over $\mathfrak{Y}$.
Moreover, the defining equations of $X$ on the universal cover of
$\mathbb{P}(E)$ do not depend on the point of the universal cover of
$Y$ (see the proof of Proposition~\ref{prop_flat-lcf}).
Thus, the same equations define a projective submersion
$\mathfrak{X} \to \mathfrak{Y}$ via a relative embedding
$\mathfrak{X} \hookrightarrow \mathbb{P}(\mathfrak{E})$.
Since $\mathfrak{Y} \to \mathbb{D}$ is a K\"ahler morphism and
$\mathfrak{X} \to \mathfrak{Y}$ is projective, the composition
$\mathfrak{X} \to \mathbb{D}$ is a K\"ahler morphism.
\end{proof}

The following result is a slight generalization of the structure theorems
in \cite[Corollary~1.7]{IM22} and \cite[Theorem~1.8]{Ou23}.
We note that, after the first version of this paper was posted on arXiv,
a similar result was obtained by a different method in \cite{IJZ25}.

\begin{thm}[{cf.~\cite[Corollary~1.7]{IM22}, \cite[Theorem~1.8]{Ou23}}]
\label{thm-second}
Let $X$ be a compact K\"ahler manifold of dimension $n \geq 2$
with nef anti-canonical bundle.
Then, the Miyaoka--Yau-type inequality
\[
c_2(T_X) \cdot \{\omega_1\} \cdot \dots \cdot \{\omega_{n-2}\} \geq 0
\]
holds for any K\"ahler forms $\omega_1, \dots, \omega_{n-2}$.

Moreover, the following conditions are equivalent:
\begin{enumerate}
\item[$(1)$] There exist K\"ahler forms
$\omega_1, \dots, \omega_{n-2}$ such that
\[
c_2(T_X) \cdot \{\omega_1\} \cdot \dots
\cdot \{\omega_{n-2}\} = 0.
\]
\item[$(2)$] There exists a finite \'etale cover
$\tilde{X} \to X$ such that $\tilde{X}$ is isomorphic either to a complex torus
or to a $\mathbb{P}^1$-bundle over a complex torus.
In the latter case, after possibly passing to a further double
\'etale cover, this is equivalent to
$\tilde{X} \cong \PP(V)$ for some vector bundle $V$ of rank $2$
$($see \cite[Lemma~7.4]{CP91}$)$.
\end{enumerate}
\end{thm}

\begin{proof}

In the projective case, the first assertion was proved in
\cite[Corollary~1.5]{Ou23}, based on \cite[Theorem~6.1]{Miy87}.
In the compact K\"ahler case, it follows from
\cite[Remark~5.4.11]{Cao13} by combining the generic nefness result
in Theorem~\ref{nefness} with the Gieseker--Bogomolov inequality and
the Khovanskii--Teissier inequality.
For a single polarization, the required inequalities are classical:
they follow from the results of Uhlenbeck--Yau and the Hodge index
theorem, respectively.
In our setting, however, we need their multi-polarized versions.
The required Gieseker--Bogomolov inequality is established in
\cite[Corollary~3.7]{CW24}, while the required
Khovanskii--Teissier inequality follows from
\cite[Proposition~5.2]{Dem93}.
The latter can also be proved by the argument of
\cite[Proposition~6.2.1]{Cao13}, using
\cite[Theorem~2.4.B]{Gro90} or \cite[Theorems~A and~C]{DN06}.

The implication from \textup{(2)} to \textup{(1)} is immediate.
Therefore, it suffices to prove that \textup{(1)} implies \textup{(2)}.
We will show that, under Condition~\textup{(1)}, the augmented
irregularity of $X$ is at least $n-1$.
This yields \textup{(2)}.
Our argument uses an algebraic approximation together with
\cite[Theorem~1.8]{Ou23}, which proves the assertion in the
projective case.

First, by \cite[Lemma~5.1]{Ou23}, if the equality in
Condition~\textup{(1)} holds for one collection of K\"ahler classes,
then it holds for every collection of K\"ahler classes.
Although \cite[Lemma~5.1]{Ou23} is stated for ample classes,
its proof applies without change to arbitrary K\"ahler classes.

By Corollary~\ref{cor-kodaira}, there exists an algebraic approximation
$p \colon \mathfrak{X} \to \mathbb{D}$ of $X$ that is a K\"ahler
fibration and contains a sequence of projective fibers converging to
the central fiber $\mathfrak{X}_0 \cong X$.
Choose a relative K\"ahler class on $\mathfrak{X}$, and denote its
restriction to $\mathfrak{X}_t$ by $\alpha_t$.
By the preceding observation, Condition~\textup{(1)} gives
\[
c_2(T_{\mathfrak{X}_0}) \cdot \alpha_0^{n-2}=0.
\]
This intersection number is invariant under deformation, and hence
\[
c_2(T_{\mathfrak{X}_t}) \cdot \alpha_t^{n-2}=0
\]
for every $t \in \mathbb{D}$.
Applying \cite[Theorem~1.8]{Ou23} to any projective fiber of this
algebraic approximation, we conclude that its augmented irregularity
is at least $n-1$.
Since the augmented irregularity is invariant under this deformation,
the augmented irregularity of $X$ is also at least $n-1$.
This proves \textup{(2)}.
\end{proof}

\appendix
\section{Supplementary results of independent interest} \label{sec-app}

We begin by explaining the contents of the appendix.
Subsections~\ref{subsec-psef} and~\ref{subsec-segre}
are devoted to topics related to the flatness criterion established in~\cite{CDM26},
which plays a crucial role in applying the result of Section~\ref{sec-main}.
In Subsection~\ref{subsec-psef}, we establish a characterization of
pseudo-effectivity for sheaves that is needed in the proof of~\cite{CDM26}.
In Subsection~\ref{subsec-segre}, assuming that the ambient space $X$ is smooth,
we give an alternative proof of the flatness criterion using the second Chern
class and Segre currents.
Although this argument is no longer needed once the result of~\cite{CDM26} is established, 
it is of independent interest, 
since an approach based on the second Chern class naturally generalizes the previous work \cite{HP16}.

In Subsection~\ref{subsec-app}, we collect some basic facts about determinant
sheaves and the Chern classes of coherent sheaves.

The material in these subsections may be familiar to experts, with the exception of Subsections~\ref{subsec-psef} and~\ref{subsec-segre}.
However, to the best of our knowledge, these results have not been explicitly
formulated and proved in the existing literature.
We include complete proofs for future reference.

\subsection{A characterization of pseudo-effective sheaves}\label{subsec-psef}

The following flatness criterion was generalized to (potentially klt) singular K\"ahler varieties in \cite{CDM26}. 
The proof there uses the characterization of pseudo-effective sheaves given in Lemma~\ref{lemm-extension}. 
In this subsection, we prove Lemma~\ref{lemm-extension}. 

\begin{thm}\label{thm-flatness}
Let $\mathcal{E}$ be a torsion-free sheaf on a compact K\"ahler manifold $X$. 
If $\mathcal{E}$ is a pseudo-effective sheaf with $c_{1}(\mathcal{E})=0$, 
then the reflexive hull $\mathcal{E}^{**}$ is locally free and numerically flat on $X$. 
\end{thm}

In this subsection, 
we generalize the basic theory of pseudo-effective sheaves, 
as developed in \cite{Mat23, BKK}, to the complex analytic setting.  
We work with the following setup throughout this subsection.

\begin{setup}\label{setup}
Let $\cal{E}$ be a torsion-free sheaf on a compact K\"ahler variety $X$. 
Let $\pi_{\mathcal{E}} \colon \bb{P}(\cal{E}) \to X$ denote the normalization of the irreducible component dominating $X$ of  
the projectivization ${\bf{Proj}}(\bigoplus_{m=0}^{\infty}\operatorname{Sym}^{m} \cal{E})$ 
of the graded sheaf $\bigoplus_{m=0}^{\infty}\operatorname{Sym}^{m} \cal{E}$. 
(Note that this notation slightly differs from that in \cite{Mat23}, where $\bb{P}(\cal{E})$ was defined before normalization.)
Let $\pi \colon P \to \bb{P}(\cal{E})$ 
be a projective log resolution of $\bb{P}(\cal{E})$, 
equipped with the commutative diagram:
\[
\xymatrix{
\bb{P}(\cal{E}) \ar[d]_{\pi_{\mathcal{E}}}&  & P\ar[ll]_{\hspace{0cm} \pi}  \ar@/^16pt/[dll]^{p} \\
X. &  & 
}
\]
Set $X_{0} := X_{\reg} \cap X_{\cal{E}}$ and $P_{0} := p^{-1}(X_{0})$. 
Note that $\codim (X \setminus X_{0}) \geq 2$. 
We may assume that $\pi \colon P \to \bb{P}(\cal{E})$ is isomorphic over $X_{0}$, 
so that $p \colon P \to X$ can be identified with the projective bundle 
$\bb{P}(\cal{E}|_{X_{0}}) \to X_{0}$ associated to the locally free sheaf $\mathcal{E}|_{X_{0}}$. 

We use the following notation:
\begin{itemize}
  \item[$\bullet$] $L := \pi^{*} \mathcal{O}_{\bb{P}(\cal{E})}(1)$, 
    where $\mathcal{O}_{\bb{P}(\cal{E})}(1)$ is the pullback of the hyperplane bundle on 
    ${\bf{Proj}}(\bigoplus_{m=0}^{\infty} \operatorname{Sym}^{m} \cal{E})$;
  \item[$\bullet$] $\omega_{P}$ is a K\"ahler form on $P$;
  \item[$\bullet$] $\omega_{X}$ is a K\"ahler form on $X$;
  \item[$\bullet$] $\Lambda$ is an effective $p$-exceptional divisor such that
    \[
    p_{*}\mathcal{O}_{P}(m(L + \Lambda)) = \operatorname{Sym}^{[m]}\mathcal{E}.
    \]
    See \cite[\S III.5, Lemma 5.10, pp.~107--108]{Nak04} for the existence of such $\Lambda$.
\end{itemize}
\end{setup}

We can easily construct singular Hermitian metrics $\{h_{\varepsilon}\}_{\varepsilon > 0}$  
on the restriction $L|_{P_{0}}$ satisfying properties (1) and (2) in Lemma~\ref{lemm-metric},  
which can then be extended to metrics on $L + \Lambda$ by Lemma~\ref{lemm-extension}.

\begin{lemm}\label{lemm-metric}
Under Setup~\ref{setup}, assume that $\mathcal{E}$ is pseudo-effective 
in the sense of Definition \ref{def-weak}. 
Then, the line bundle $L|_{P_{0}}$ admits 
singular Hermitian metrics $\{h_{\e}\}_{\e>0}$ such that 
\begin{itemize}
\item[$(1)$] $\sqrt{-1}\Theta_{h_{\e}} + \e p^{*}\omega_{X} \geq \delta_{\e} \omega_{P}$ holds on $P_{0}$ 
for some $\delta_{\e}>0$; 
\item[$(2)$] $\{x \in P_{0}\,|\, \nu(h_{\e}, x) >0 \}$ is not dominant over $X_{0}$, 
where $\nu(h_{\e}, x)$ denotes the Lelong number of the local weight of $h_{\e}$. 
\end{itemize}
\end{lemm}
\begin{proof}
We first confirm that there exists a singular Hermitian metric $h$ on $L$
such that $h$ is smooth on $P_0$ and satisfies
\[
\sqrt{-1}\Theta_h + C\, p^{*} \omega_X \geq a\, \omega_P
\]
for some constants $C \gg 1$ and $1\gg a >0$.
The hyperplane bundle associated with ${\bf{Proj}}(\bigoplus_{m=0}^{\infty} \mathrm{S}^m \mathcal{E})$
is relatively ample with respect to ${\bf{Proj}}(\bigoplus_{m=0}^{\infty} \mathrm{S}^m \mathcal{E}) \to X$, 
and the structure sheaf $\mathcal{O}_{\mathbb{P}(\mathcal{E})}$
is relatively ample with respect to the normalization $\mathbb{P}(\mathcal{E}) \to {\bf{Proj}}(\bigoplus_{m=0}^{\infty} \mathrm{S}^m \mathcal{E})$.
This shows that the line bundle $\mathcal{O}_{\mathbb{P}(\mathcal{E})}(1)$ is relatively $\pi_{\mathcal{E}}$-ample.
Hence, there exists a smooth Hermitian metric $g_1$ on $\mathcal{O}_{\mathbb{P}(\mathcal{E})}(1)$ such that
$
\sqrt{-1}\Theta_{g_1} + C\, \pi_{\mathcal{E}}^{*} \omega_X > 0 
$
for some constant $C > 0$.
Let $E$ be an effective $\pi$-exceptional divisor such that $-E$ is relatively $\pi$-ample.
Then $E$ admits a smooth Hermitian metric $g_{2}$ such that  
\begin{align*}
&\sqrt{-1}\Theta_{\pi^{*}g_{1}}(L) + C p^{*} \omega_{X} - \delta \sqrt{-1}\Theta_{g_{2}}(E) \\
= {}&\pi^{*}\big( \sqrt{-1}\Theta_{g_{1}}(L) + C \pi_{\mathcal{E}}^{*} \omega_{X} \big) - \delta \sqrt{-1}\Theta_{g_{2}}(E) >0
\end{align*}
for $1 \gg \delta > 0$.
We define the singular Hermitian metric on $L = (L - \delta E) + \delta E$ by  
\[
h := (\pi^{*}g_{1}) \cdot g_{2}^{-\delta} \cdot g_{[\delta E]},
\]
where $g_{[\delta E]}$ denotes the singular Hermitian metric associated with the effective divisor $\delta E$.  
Then $h$ satisfies the required properties.

Next, we construct singular Hermitian metrics $\{g_{\e}\}_{\e > 0}$ on $L|_{P_{0}}$ satisfying the following properties: 
\begin{itemize}
\item[$\bullet$] $\sqrt{-1}\Theta_{g_{\e}} + \e p^{*}\omega_{X} \geq 0$ on $P_{0}$;
\item[$\bullet$] $\{x \in P_{0} \mid \nu(g_{\e}, x) > 0\}$ is not dominant over $X_{0}$.
\end{itemize}
By the definition of pseudo-effective sheaves,  
there exist singular Hermitian metrics $\{h_{m}\}_{m \in \mathbb{Z}_{+}}$ on ${\rm{Sym}}^{[m]} \mathcal{E}$  
such that $\sqrt{-1}\Theta_{h_{m}} \geq -\omega_{X} \otimes \mathrm{id}_{{\rm{Sym}}^{[m]} \mathcal{E}}$.  
The pullbacks $p^{*} h_{m}$ define singular Hermitian metrics on $p^{*} \Sym^{m} \mathcal{E}|_{P_{0}}$  
such that $\sqrt{-1}\Theta_{p^{*} h_{m}} \geq -p^{*} \omega_{X} \otimes \mathrm{id}$ on $P_{0}$.
Let $H_{m}$ be the induced singular Hermitian metric on $mL$ via the natural surjective morphism:
\[
p^{*} \Sym^{m} \mathcal{E}|_{P_{0}} \to mL|_{P_{0}}.
\]
By construction, we see that $\sqrt{-1}\Theta_{H_{m}} \geq -p^{*} \omega_{X}$  
and $\{x \in P_{0} \mid \nu(H_{m}, x) > 0\}$ is not dominant over $X_{0}$.  
The metrics $H_{m}^{1/m}$ on $L|_{P_{0}}$ yield the desired family $\{g_{\e}\}_{\e > 0}$ as follows.
For each $\e>0$, choose a positive integer $m_\e$ such that $1/m_\e<\e$, and define
$g_\e:=H_{m_\e}^{1/m_\e}$. Then, we have 
\[
\sqrt{-1}\Theta_{g_\e}(L)
= \sqrt{-1}\Theta_{H_{m_\e}^{1/m_\e}}(L)
\geq -\frac{1}{m_\e}p^{*}\omega_{X}
\geq -\e p^*\omega_X.
\]

Finally, we define singular Hermitian metrics on $L|_{P_{0}}$ by  
\[
h_{\e} := g_{\e/2}^{1 - \eta} \cdot h^{\eta}.
\]
From the curvature properties of $h$ and $g_{\e}$, we obtain
\[
\sqrt{-1}\Theta_{h_{\e}} \geq -\left( \frac{\e}{2} + \eta C \right) p^{*} \omega_{X} + \eta a \omega_{P}.
\]
Setting $\eta := \min\{1/2,\e/(2C)\}$, we conclude that $h_{\e}$ satisfies the required properties.
\end{proof}

In general, the metric $h_{\varepsilon}$ on $L|_{P_{0}}$ \textit{cannot} 
extend to a singular Hermitian metric on $L$ while preserving the curvature condition. 
For instance, let $\mathcal{E}$ be the ideal sheaf of a submanifold $Z \subset X$ of codimension $\ge 2$.  
Then $\mathcal{E}$ is pseudo-effective, and we may take $P$ to be the blow-up of $X$ along $Z$ 
with the exceptional divisor $E$, so that $L = -E$.  
Then, the line bundle $L = -E$ never admits metrics $h_{\varepsilon}$ with $\sqrt{-1}\Theta_{h_{\e}} \geq - \e \omega_{P}$. 

As noted above, extending a metric $h_{\varepsilon}$ defined only on $P_{0}$ to $P$ is a delicate problem, 
since $P \setminus P_{0}$ can be divisorial.  
Nevertheless, by considering $L + \Lambda$ instead of $L$,  
we can extend $h_{\e}$ to a metric on $L + \Lambda$  
using the identification $L|_{P_{0}} \cong (L + \Lambda)|_{P_{0}}$,  
which is a key step in the proof of Theorem~\ref{thm-flatness}.

\begin{lemm}\label{lemm-extension}
Under the assumptions of Lemma~\ref{lemm-metric},  
the metric $h_{\e}$ constructed in Lemma~\ref{lemm-metric} 
extends to a singular Hermitian metric on $L+\Lambda$ satisfying the same properties~$(1)$ and~$(2)$ over $X$.
\end{lemm}
\begin{proof}
We omit the subscript $\e$ in $h_{\e}$ and regard the metric $h$ on $L|_{P_{0}}$ as a metric on $(L+\Lambda)|_{P_{0}}$ via the identification $L|_{P_{0}} \cong (L+\Lambda)|_{P_{0}}$.  
It suffices to show that the local weight $\varphi$ of $h$ is bounded from above in a neighborhood of $P\setminus P_{0}$.  
The problem is local on $X$, so we may assume that $X$ is Stein.  
Then, replacing $h$ by $h\,e^{-\psi}$, where $\psi$ is a local potential of $\e\,p^{*}\!\omega_{X}$, 
we may assume that $\sqrt{-1}\Theta_{h}$ is a K\"ahler current by property~$(1)$.  

For a positive integer $m\in\mathbb{Z}_{+}$, we consider the $L^{2}$-space  
\[
H^{0}_{L^{2}}(P_{0},m(L+\Lambda))_{h^{m},\omega_{P}}
\]
of holomorphic sections $\sigma$ on $P_{0}$ with finite $L^{2}$-norm $\|\sigma\|_{h^{m},\omega_{P}}< \infty$.  
Let $B_{m}$ denote the Bergman kernel of this Hilbert space,  
which is the singular Hermitian metric on $(L+\Lambda)|_{P_{0}}$ whose local weight $b_{m}$ is defined to be 
\[
b_{m}(p)
          :=\frac{1}{m}\sup\Bigl\{\log|\sigma|(p) \;\Big|\;\|\sigma\|_{h^{m},\omega_{P}}\le 1\Bigr\} \quad \text{ for } p\in P_{0}. 
\]

Then, for every relatively compact set $K\Subset X$, there exists a constant $C=C(K)$ (independent of $m$) such that
\begin{equation}\label{eq-dem}
\varphi(p)\le b_{m}(p)+\frac{C}{m}\qquad\text{for }p\in P_{0}\cap p^{-1}(K).
\end{equation}
This can be verified by modifying Demailly's approximation technique (cf. \cite[(13.21) Theorem, p.~138]{Dem10}). 
We briefly explain the proof for the reader's convenience.

Fix a point $p_{0} \in P_{0}\cap p^{-1}(K)$. 
Take an open ball $B \Subset P$ centered at $p_{0}$ and 
a hypersurface $H\subset X$ such that $X\setminus X_{0} \subset H$ and $p_{0} \not \in p^{-1}(H)$.  
Note that $P\setminus p^{-1}(H)$ is a weakly pseudoconvex K\"ahler manifold and $B\setminus p^{-1}(H)$ is Stein. 
This allows us to solve $\dbar$-equations on $P\setminus p^{-1}(H)$ and 
to apply the $L^{2}$-extension theorem on $B\setminus p^{-1}(H)$.
If $\varphi(p_{0})=-\infty$, the inequality is trivial, so we may assume $\varphi(p_{0})>-\infty$.  
For a given $a\in\mathbb{C}$, 
by the $L^{2}$-extension theorem, there exists a holomorphic function $f$ on $B\setminus p^{-1}(H)$ such that
\[
f(p_{0})=a \quad \text{ and } \quad 
\int_{B\setminus p^{-1}(H)}|f|^{2}e^{-m\varphi}\,dV_{\omega_{P}}\le C_{1}|a|^{2}e^{-m\varphi(p_{0})}.
\]
Fix a relatively compact subset $K_{P}\Subset B$ and consider a cut-off function $\theta$ supported in $B$ and equal to~$1$ on $K_{P}$.  
Since $\sqrt{-1}\Theta_{h}$ is a K\"ahler current, 
by the standard technique of $\dbar$-equations, 
we can solve $\overline{\partial}g=\overline{\partial}(\theta f)$ on $P\setminus p^{-1}(H)$ with the $L^{2}$-estimate
\[
\int_{P\setminus p^{-1}(H)}|g|_{h^{m}}^{2}\,e^{-2n\theta\log|z-z(p_{0})|}\,dV_{\omega_{P}}
\le C_{2}\int_{B\setminus p^{-1}(H)}|f|^{2}e^{-m\varphi}\,dV_{\omega_{P}}, 
\]
where $z$ is a coordinate of $B$ and $n:=\dim P$. 
Then, the section $\sigma:=\theta f-g$ satisfies  
\begin{align}\label{eq-norm}
\sigma(p_{0})= f(p_{0})=a \quad \text{ and } \quad  
\|\sigma\|_{h^{m},\omega_P}^{2}\le C_{3}|a|^{2}e^{-m\varphi(p_{0})}, 
\end{align}
where we identify $\sigma$ with a holomorphic function via a fixed local frame.
Since $h$ is defined near every point of $p^{-1}(H)\cap P_{0}$, 
any $L^{2}$-sections on $P\setminus p^{-1}(H)$ extend to $P_{0}$, giving the isomorphism
\[
H^{0}_{L^{2}}(P_{0},m(L+\Lambda))_{h^{m},\omega_{P}}
\cong H^{0}_{L^{2}}(P\setminus p^{-1}(H),m(L+\Lambda))_{h^{m},\omega_{P}}.
\]
This shows that $\sigma$ lies in the left-hand space. 
We can choose $a$ so that the right-hand side in \eqref{eq-norm} equals $1$, and thus obtain \eqref{eq-dem}.

\smallskip 
Henceforth, to simplify notation, we assume $m = 1$ and set $L_{\Lambda} := L + \Lambda$.
We will show that the weight $b:=b_{1}$ of the Bergman kernel is locally bounded from above near $P \setminus P_{0}$.  
The crucial point is that every section of $L_{\Lambda}$ on $P_{0}$ extends holomorphically to $P$ 
because $p_{*}(L_{\Lambda})$ is reflexive by the choice of $\Lambda$. 

Fix a trivialization of $L_{\Lambda}$ on an open ball $B \subset P$.  
Then, for each point $p \in B_{0}:=B \cap P_{0}$, we consider the evaluation map: 
\[
\mathrm{ev}_{p} \colon 
H^{0}_{L^{2}}(P_{0},L_{\Lambda})_{h,\omega_{P}}
\longrightarrow
H^{0}(B_{0},L_{\Lambda}) \cong H^{0}(B_{0},\mathcal{O}_{B})
\longrightarrow
\mathbb{C}. 
\]

Varying the point $p$ induces a map
\[
\mathrm{ev}^{*}: B_{0}\xrightarrow{ }
\bigl(H^{0}_{L^{2}}(P_{0},L_{\Lambda})_{h,\omega_{P}}\bigr)^{*}.
\]
Here the dual space is equipped with the strong topology.
By Montel's theorem, the map $\mathrm{ev}^{*}$ is continuous.
Taking the operator norm, 
we can interpret the weight $b$ as the composition 
\[
b=b_{1}\colon\,
B_{0}\xrightarrow{\ \mathrm{ev}^{*}\ }
\bigl(H^{0}_{L^{2}}(P_{0},L_{\Lambda})_{h,\omega_{P}}\bigr)^{*}
\xrightarrow{\ \|\cdot\|_{\mathrm{op}}\ }
\mathbb{R}. 
\]

Hence, it is enough to extend $\mathrm{ev}^{*}$ continuously to $B$; 
once this is done, $b=b_{1}$ is bounded on compact subsets of $B$.

The restriction map
\[
H^{0}(P,L_{\Lambda}) \longrightarrow H^{0}(P_{0},L_{\Lambda})
\]
is an isomorphism of Fr\'echet spaces with the topology induced by local sup-norms. 
Indeed, this map is clearly continuous, injective, and linear. 
Moreover, it is surjective 
since we have
\[
H^{0}(P,L_{\Lambda})
\cong H^{0}(X,p_{*}(L_{\Lambda}))
\cong H^{0}(X_{0},p_{*}(L_{\Lambda}))
\cong H^{0}(P_{0},L_{\Lambda})
\]
by reflexivity. 
Hence, the map is a homeomorphism by the open mapping theorem.
This induces a continuous map
\[
H^{0}(P,L_{\Lambda})^{*}
\cong H^{0}(P_{0},L_{\Lambda})^{*}
\to
\bigl(H^{0}_{L^{2}}(P_{0},L_{\Lambda})_{h,\omega_{P}}\bigr)^{*}. 
\]
On the other hand, there exists a continuous map
\[
B \to H^{0}(B, L_{\Lambda})^{*} \to H^{0}(P,L_{\Lambda})^{*}.
\]
The composition of these maps extends $\mathrm{ev}^{*}$ continuously to $B$. 
Consequently, $b=b_1$ is locally bounded from above on $B$.
\end{proof}

\subsection{Segre currents and second Chern classes}\label{subsec-segre}

The purpose of this subsection is to establish a flatness criterion 
for pseudo-effective sheaves on compact K\"ahler manifolds (see Theorem~\ref{thm-flatness}). 
The strategy of our proof follows that of \cite{Wu22} , but requires substantial improvements.  

The generalization of Theorem~\ref{thm-flatness} to singular K\"ahler varieties 
was proved without using second Chern classes. 
However, as in the projective case, 
it is of interest to consider another proof using second Chern classes. 
In this subsection, we give such a proof 
by generalizing the results of \cite{Wu22}, which were originally formulated for projective space bundles, 
to a broader setting. 

\begin{thm}\label{Segre}
Let $f\colon X \to Y$ be a fibration of relative dimension $r-1$ 
between compact K\"ahler manifolds, smooth outside a Zariski-closed subset $Z \subset Y$ with  
$\codim Z \geq k+1$.
Let $T$ be a semipositive $(1,1)$-current on $X$ with analytic singularities such that  
\begin{itemize}
\item[$\bullet$] $T=\alpha+\ddbar\varphi$, where $\alpha$ is a $d$-closed smooth $(1,1)$-form on $X$;
\item[$\bullet$] $T$ is smooth outside a Zariski-closed set  
      $A\subset X$ with $A=f^{-1}\!\bigl(f(A)\bigr)$ and $\operatorname{codim}f(A)\ge k$.
\end{itemize}
Let $c\colon Y\setminus Z \to \RR_{>0}$ be a continuous function. 
Assume that for every point $x\in f^{-1}(Y\setminus Z)$, 
there exists a quasi-psh function $\psi$ on $X$ (depending on $x$) with analytic singularities such that 
\begin{itemize}
\item[$\bullet$] $\alpha+\ddbar\psi$ is semipositive on a neighborhood of the fiber $X_{f(x)}$;
\item[$\bullet$] $\psi$ is smooth outside a Zariski-closed set $B \subset X$ with $\codim B\geq k+r$;
\item[$\bullet$] $\nu(\psi,x)\ge c(f(x))$.
\end{itemize}

Then, there exists a semipositive $(k,k)$-current $S$ in the class $f_{*}\{\alpha\}^{\,r+k-1}$ such that for every $y\in Y\setminus Z$,
\[
\nu(S,y)\;\ge\;
\min\!\Bigl\{\bigl(\sup_{x\in X_{y}}\nu(\varphi,x)\bigr)^{\,r+k-1},\;c(y)\Bigr\}.
\]
\end{thm}

\begin{proof}
The theorem follows by suitably adapting the arguments of \cite{Wu22}.  
To avoid repetition, we explain only the essential modifications here.

Following \cite[Section~4]{Wu22}, we define a $(k,k)$-current $S$ on $Y$ by
\[
S := \wlim_{\delta \to 0} f_* \big( ( \alpha + \ddbar \log (e^{\varphi} + \delta) )^{r+k-1} \big)
\in f_* \{\alpha\}^{r+k-1}.
\]
The regularized potential is locally bounded for every $\delta>0$, so the current is well-defined. 
Fix a point $y \in Y \setminus Z$, and take $x \in f^{-1}(y)$.
Using the function $\psi$ in the assumptions, we consider
\[
\wlim_{\delta \to 0} f_* \big( ( \alpha + \ddbar \log (e^{\varphi} + \delta e^{\psi}) )^{r+k-1} \big)
\in f_* \{\alpha\}^{r+k-1}, 
\]
which is clearly a semipositive current near $y$. 
Moreover, this coincides with $S$, showing that $S$ is semipositive on $Y \setminus Z$.
This follows from the argument in \cite[Lemma~14]{Wu22}. 
Indeed, under the codimension condition on $f(A)$, the difference between two weak limits arising from different choices of $\psi$ is a normal current supported on $f(A)$.  
By Demailly's support theorem, this difference is a finite linear combination of integration currents over irreducible components of $f(A)$, with coefficients given by the generic Lelong numbers along those components.  
Thanks to the codimension assumption on $Z$, we can compute these Lelong numbers at generic points outside $Z$, 
and hence the local argument in \cite[Lemma~14]{Wu22} applies without modification.

By a version of the Skoda--El Mir extension theorem (see \cite[pp.~178, Cor.~3.2]{Sib85}), the $(k,k)$-current $S$ extends across the analytic subset $Z \subset Y$ of codimension $\geq k+1$.  
By
\[
H^{2k}(Y \setminus Z, \RR) \cong H^{2k}(Y, \RR),
\]
the extended current represents the same cohomology class.
The Lelong number estimate is obtained by a local computation 
(see \cite[pp.~52--53]{Wu22} for details).
\end{proof}

We now apply Theorem~\ref{Segre} to prove the seminegativity of the second Chern class 
$c_{2}(\cal{E}^{**})$, which is a key step in the proof of Theorem~\ref{thm-flatness}. 

\begin{prop}\label{prop-c2}
Under Setup~\ref{setup}, assume that $\mathcal{E}$ is a pseudo-effective reflexive sheaf 
such that $c_{1}(\mathcal{E})=0$.

Then, the second Chern class $c_{2}(\cal{E})$ is represented by a seminegative $(2,2)$-current on $X$.
\end{prop}
\begin{proof}

We may assume that the metric $h_{\varepsilon}$ in Lemma~\ref{lemm-extension} has analytic singularities by Demailly's approximation theorem (see \cite[(13.21) Theorem, p.~138]{Dem10}).  
Let $\varphi_{\varepsilon}$ be the quasi-psh function associated with $h_{\varepsilon}$, so that
\[
\alpha_{\varepsilon} + \ddbar \varphi_{\varepsilon}
:= (\beta + \varepsilon p^{*}\omega_{X}) + \ddbar \varphi_{\varepsilon} 
\geq \delta_{\varepsilon} \omega_{P} \geq 0,
\]
where $\beta$ is a smooth representative of $c_{1}(L + \Lambda)$.
The nonsmooth locus $Z:=X \setminus X_{0}$ of the fibration $p \colon P \to X$ 
satisfies the assumptions of Theorem~\ref{Segre} for $k=2$
since $\mathcal{E}$ is reflexive. 
The singular locus $A$ of $\varphi_{\e}$ may not satisfy 
the assumption of Theorem~\ref{Segre} for $k=2$, 
but satisfies this assumption for $k=1$ by property~(2). 
For every $x\in p^{-1}(X_0)$, a function $\psi$ satisfying the assumptions of Theorem~\ref{Segre} exists (see \cite[p.~52]{Wu22}).
Thus, applying Theorem~\ref{Segre} to each $\varepsilon$ in the case $k=1$, 
we obtain a semipositive $(1,1)$-current $S_{\varepsilon}$ on $X$ such that
\begin{align}\label{eq-lelong}
S_{\varepsilon} \in p_* \{\alpha_{\varepsilon}\}^{r} \quad\text{and}\quad
\nu(S_{\varepsilon}, y) \geq  
\min\!\Bigl\{\bigl(\sup_{x\in P_{y}}\nu(\varphi_{\varepsilon},x)\bigr)^{\,r},\;c(y)\Bigr\}
\end{align}
for every $y \in X_0 = X \setminus Z$.  
Note that this estimate may not hold for $y \in Z$.

We now compute the cohomology class 
\[
p_*\left(c_1(L + \Lambda)\right)^{r} = c_1(\mathcal{E}) = 0
\quad \text{on } X_0.
\]
The same equality holds on $X$ by Proposition~\ref{topo} and the isomorphism $H^2(X, \mathbb{R}) \cong H^2(X_0, \mathbb{R})$. 
Thus, the classes containing $S_{\varepsilon}$ satisfy
\[
S_{\varepsilon} \in 
p_* \{\alpha_{\varepsilon}\}^{r}
= p_*\left( c_1(L + \Lambda) + \varepsilon p^*\omega_X \right)^r
\to 0
\quad \text{as } \varepsilon \to 0.
\]
This shows that $\sup_{y \in X} \nu(S_{\varepsilon}, y) \to 0$ (see, for example, \cite[Lemma~15]{Wu22}).  
By \eqref{eq-lelong}, 
this implies that $\nu(\varphi_{\varepsilon}, x) \to 0$ locally uniformly on $p^{-1}(X \setminus Z)$.

Let $\{K_j\}_{j \in \mathbb{Z}_{+}}$ be an increasing exhaustion of $X \setminus Z$ by relatively compact open subsets
such that $K_j \Subset K_{j+1} \Subset X \setminus Z$. 
Fix a decreasing sequence $\varepsilon_j \to 0$.  
By Demailly's regularization theorem with smooth forms (see \cite[Th\'eor\`eme 9.1]{Dem82}, \cite[Theorem 1.1]{Dem94}), 
we find smooth quasi-psh functions $\tilde{\varphi}_j$ on $P$ such that
\[
\tilde{T}_j := \alpha_{\varepsilon_j} + \ddbar \tilde{\varphi}_j \geq -\varepsilon_j \omega_P
\quad \text{on } p^{-1}(K_j).
\]
This estimate follows from the local uniform convergence of $\nu(\varphi_{\varepsilon}, x)$ on $p^{-1}(X \setminus Z)$.

We can define $(2,2)$-currents on $X$ by
\[
S_j := p_* \big( \tilde{T}_j + \varepsilon_j \omega_P \big)^{r+1}
\in p_*\left\{ \alpha_{\varepsilon_j} + \varepsilon_j \omega_P \right\}^{r+1}.
\]

Note that $S_j$ is not semipositive on $X$, but it is smooth on $X \setminus Z$ and belongs to a bounded family of cohomology classes.  
By compactness, after passing to a subsequence, we may assume that $S_j$ converges weakly to a current $S$.
For a fixed $j_0$, the family $\{S_j\}_{j \ge j_0}$ is semipositive on $K_{j_0}$.
Hence, the limit $S$ is semipositive on each $K_{j_0}$ and therefore on $X \setminus Z$.
The current $S$ extends to a semipositive current on $X$ by the Skoda--El Mir extension theorem, since $Z$ has codimension at least three. 

The current $S$ represents the class 
\[
\lim_{j \to \infty} p_*\left\{ \alpha_{\varepsilon_j} + \varepsilon_j \omega_P \right\}^{r+1}
=p_*\left( c_1(L + \Lambda) \right)^{r+1}.
\]
Observe that 
\[
p_*\left( c_1(L + \Lambda) \right)^{r+1}
= c_1^2(\mathcal{E}) - c_2(\mathcal{E}) = -c_2(\mathcal{E})
\]
on $X_0$. 
The same equality holds globally on $X$ 
by Proposition~\ref{topo} and the isomorphism $H^4(X, \mathbb{R}) \cong H^4(X_0, \mathbb{R})$. 
This completes the proof.

\end{proof}

With the above preparations, Theorem~\ref{thm-flatness} follows easily from the arguments in \cite{CCM19, Wu22}. 

\begin{proof}[Proof of Theorem~\ref{thm-flatness}]
It is sufficient to show that $\mathcal{E}^{**}$ is locally free. 
Replacing $\mathcal{E}$ with $\mathcal{E}^{**}$, which remains pseudo-effective and has vanishing first Chern class, 
we may assume that $\mathcal{E}$ is reflexive. Once its local freeness is established, 
its numerical flatness follows directly from \cite[Main Theorem]{Wu22}. 

Take a reflexive subsheaf $\mathcal{S} \subset \mathcal{E}$ 
of minimal rank among the subsheaves satisfying 
\begin{align}\label{eq-slope}
\int_X c_{1} (\mathcal{S}) \cdot \omega_X^{n-1} \geq 0.
\end{align}
Consider the exact sequence
\begin{align}\label{exact}
0 \rightarrow \mathcal{S} \rightarrow \mathcal{E} \rightarrow 
\mathcal{Q} := \mathcal{E}/\mathcal{S} \rightarrow 0.
\end{align}
The quotient sheaf $\mathcal{Q}/\mathcal{Q}_{\tor}$ modulo torsion is pseudo-effective, 
since it is also a quotient of the pseudo-effective sheaf $\mathcal{E}$. 
In particular, its first Chern class $c_{1}(\mathcal{Q})$ is also pseudo-effective, 
since the torsion part contributes only a pseudo-effective first Chern class by Proposition~\ref{det-tor}.
We can then deduce that $c_{1}(\mathcal{S}) = c_{1}(\mathcal{Q}) = 0$. 
Indeed, we have $c_{1}(\mathcal{S}) = -c_{1}(\mathcal{Q})$ from 
$c_{1}(\mathcal{E}) = 0$, and the equalities follow from \eqref{eq-slope}.

If $\mathcal{S} = \mathcal{E}$, then $\mathcal{S}$ is stable by its minimality. 
Hence, Bogomolov's inequality and the condition $c_1(\mathcal{S}) = 0$ give
\[
\int_{X} c_{2}(\mathcal{S}) \cdot \omega_{X}^{n-2} \geq 0.
\]
Combined with Proposition~\ref{prop-c2}, this gives 
$\int_X c_2(\mathcal{S})\cdot\omega_X^{n-2} = 0$, 
and hence $\mathcal{S} = \mathcal{E}$ is locally free by \cite[Corollary~3]{BS94}. 

It remains to consider the case $\rank \mathcal{S} < \rank \mathcal{E}$. 
Applying \cite[Lemma 2.9]{IMZ23} (cf.\,\cite[Lemma~1.20]{DPS94}) to the morphism $\mathcal{Q}^{*} \to \mathcal{E}^*$, 
we deduce that both $\mathcal{Q}$ and $\mathcal{S}$ are vector bundles on $X_{\mathcal{E}}$, 
and that the morphism is a bundle morphism over $X_{\mathcal{E}}$. 
Consider the surjective bundle morphism 
\[
\Lambda^{m+1} \mathcal{E} \otimes \det \mathcal{Q}^* \to \mathcal{S}
\]
on $X_{\mathcal{E}}$, where $m = \rank \mathcal{Q}$. 
This implies that the reflexive sheaf $\mathcal{S}$ is pseudo-effective 
by $c_1(\mathcal{Q}) = 0$. 
By the induction hypothesis on the rank, 
we conclude that $\mathcal{S}$ is a numerically flat locally free sheaf on $X$. 
Similarly, the reflexive hull $\mathcal{Q}^{**}$ is a numerically flat locally free sheaf on $X$ by the same hypothesis.

By \cite[Lemma~16]{Wu22}, the extension class of \eqref{exact}, defined on $X_{\mathcal{E}}$, 
extends to an extension class on $X$ between $\mathcal{S}$ and $\mathcal{Q}^{**}$. 
This extension determines a vector bundle whose restriction to $X_{\mathcal{E}}$ coincides with $\mathcal{E}$. 
Hence, by reflexivity, we conclude that $\mathcal{E}$, and therefore the reflexive hull of the original sheaf, is locally free. 
\end{proof}

\subsection{Chern classes of coherent sheaves}\label{subsec-app}

In this subsection, 
we collect some basic facts about determinant sheaves and Chern classes of sheaves.  
Although these results may be known to experts,  
we include the proofs for the reader's convenience.

We begin by defining the determinant sheaf of a sheaf $\mathcal{E}$ on a normal analytic variety $X$.  
Our definition agrees with the standard one given in Section~\ref{subsec-notation} when $\mathcal{E}$ is torsion-free. 
For the corresponding definition on a smooth complex manifold, we refer the reader to \cite[Section 6, Chap.~V]{Kob87}.

\begin{defi}
\label{det-sheaf}
Let $\mathcal{E}$ be a coherent sheaf on a normal analytic variety $X$. 
Consider the short exact sequence
\[
0 \to \mathcal{E}_{\mathrm{tor}} \to \mathcal{E} \to \mathcal{E}/\mathcal{E}_{\mathrm{tor}} \to 0, 
\]
where $\mathcal{E}_{\mathrm{tor}}$ denotes the maximal torsion subsheaf of $\mathcal{E}$.
Since $\mathcal{E}/\mathcal{E}_{\mathrm{tor}}$ is torsion-free, its determinant sheaf is defined as
\[
\det(\mathcal{E}/\mathcal{E}_{\mathrm{tor}}) := \wedge^{[\mathrm{rk}\,\mathcal{E}]}(\mathcal{E}/\mathcal{E}_{\mathrm{tor}}).
\]
By \cite[Proposition 6.14, Chap.~V]{Kob87}, the determinant line bundle $\det \mathcal{E}_{\mathrm{tor}}|_{X_{\mathrm{reg}}}$ (defined only on the regular locus $X_{\mathrm{reg}}$) admits a non-trivial section whose zero locus is contained in the support of $\mathcal{E}_{\mathrm{tor}}$.
By the Remmert--Stein theorem (for example, see \cite[Theorem 8.7, Chap.~II]{agbook}) and the normality of $X$, 
this zero locus extends to a Weil divisor across the singular locus $X_{\sing}$.
Using this extended Weil divisor $D_{\mathrm{tor}}$, we define
\[
\det(\mathcal{E}_{\mathrm{tor}}) := \mathcal{O}_X(D_{\mathrm{tor}}).
\]
Then, we define the determinant sheaf of $\mathcal{E}$ as
\[
\det(\mathcal{E}) := \det(\mathcal{E}_{\mathrm{tor}}) [\otimes] \det(\mathcal{E}/\mathcal{E}_{\mathrm{tor}}).
\]
Note that $c_{1}(\mathcal{E}_{\mathrm{tor}})$ is represented by an effective Weil divisor by definition. 
\end{defi}

This definition agrees with the usual one when $X$ is smooth 
(see \cite[Lemma~6.3, Proposition~6.10, Proposition~6.14, Chap.~V]{Kob87}).  
By Definition~\ref{det-sheaf}, we immediately obtain the following propositions: 

\begin{prop}
\label{det-tor}
Let $\mathcal{E}$ be a torsion coherent sheaf on a strongly $\mathbb{Q}$-factorial analytic variety $X$.
Then, the determinant sheaf $\det \mathcal{E}$ is represented by an effective $\mathbb{Q}$-Cartier Weil divisor whose support is contained in the support of $\mathcal{E}$.
\end{prop}

\begin{prop}
\label{det-open}
Let $\mathcal{E}$ be a coherent sheaf on a normal analytic variety $X$.  
Let $U \subset X$ be an open subset with the open immersion $i_U \colon U \to X$.  
Then, we have  
\[
i_U^*\mathrm{det}(\mathcal{E}) = \mathrm{det}(\mathcal{E}|_{U}).
\]
\end{prop}

The determinant line bundle is additive with respect to short exact sequences of 
locally free sheaves on a smooth complex manifold (for example, see \cite[Proposition 6.9, Chap.~V]{Kob87}). 
This additivity extends to determinant sheaves on normal analytic varieties as well. 

\begin{prop}\label{det-exact}
Let
\[
0 \to \mathcal{S} \to \mathcal{E} \to \mathcal{Q} \to 0
\]
be a short exact sequence of coherent sheaves on a normal analytic variety $X$.  
Then, we have
\[
\det(\mathcal{E}) = \det(\mathcal{S}) [\otimes] \det(\mathcal{Q}).
\]
\end{prop}

\begin{proof}
Since $X$ is normal, it suffices to verify the equality on the regular locus of $X$.  
By Proposition~\ref{det-open}, the desired equality then follows from \cite[Proposition 6.9, Chap.~V]{Kob87}.
\end{proof}

We now check a compatibility property for Chern characters.

\begin{prop}\label{topo}
Let $\mathcal{E}$ be a coherent sheaf on a compact complex manifold $X$,  
and let $U \subset X$ be a Zariski-open subset with the open immersion $i_U \colon U \to X$.  
Then, we have 
\[
i_{U}^{*}\!\operatorname{ch}(\mathcal{E}) \;=\; \operatorname{ch}\bigl(\mathcal{E}|_{U}\bigr)
\]
in de Rham cohomology.
\end{prop}

\begin{proof}
The statement follows directly from the construction of Chern classes via locally free resolutions;  
see \cite[Remark 1]{Wu23b} for details.
\end{proof}

\begin{lemm}
\label{topo-lemma}
Let $\pi\colon M\to X$ be a proper bimeromorphic morphism between normal
analytic varieties.
Let $S$ and $T$ be $(1,1)$-currents admitting local potentials on $M$ such that
\[
[S]=[T]\in H_{\BottChern}^{1,1}(M)
\quad\text{and}\quad
\pi_*S=\pi_*T.
\]
Assume that $S-T$ can be written as the difference of semipositive closed
$(1,1)$-currents.
Then $S=T$ holds.
\end{lemm}

\begin{proof}
By the equality $[S]=[T]$ in $H_{\BottChern}^{1,1}(M)$, 
there exists a $(0,0)$-current $g$ such that
$
T-S=\ddbar g.
$
By assumption, we can locally write
$
T-S=\ddbar(v-u)
$
for some psh functions $u$ and $v$.
It follows that $g-(v-u)$ is a pluriharmonic distribution, and the real part of a
holomorphic function.
Thus, the $(0,0)$-current $g\in L^1_{\mathrm{loc}}(M)$ is locally a difference of psh functions.

Since the pushforward of currents commutes with $\ddbar$, we have
\[
\ddbar(\pi_*g)
=
\pi_*(T-S)
=
0.
\]
This shows that the $(0,0)$-current $\pi_*g$ is represented by a pluriharmonic
function on $X$.
Over the open subset on which $\pi$ is an isomorphism, we have
\[
g=\pi^*(\pi_*g).
\]
Both sides are locally integrable on $M$, and the complement of this
open subset is an analytic subset of measure zero.
Hence, the equality holds on the whole variety $M$ as an equality of currents.
Consequently, we can conclude that 
\[
T-S
=
\ddbar g
=
\ddbar\pi^*(\pi_*g)
=
\pi^*\bigl(\ddbar(\pi_*g)\bigr)
=
0.
\]
\end{proof}

\subsection{The generic principle for analytic varieties}
\label{subsec-connectedness}
The aim of this appendix is to provide a proof of the generic principle (Proposition \ref{prop-generic}) that is the crucial ingredient for Bertini-type Theorem \ref{thm-Bertini}.

First recall Remmert's mapping theorem 
\begin{thm}[{\cite[Satz 19]{Rem57}, \cite[\S (0.2), Lemma, p.267]{Man82}}]
\label{thm-Remmert-mapping}
Let $f:X\to Y$ be morphism between analytic varieties, then $f(X)$ is a countable union of locally analytic subset of $Y$.
\end{thm}

\begin{proof}
This result is essentially established in \cite[Satz 19]{Rem57}, although not written explicitly in the form we state here. More importantly, the notion of analytic varieties in Remmert's time is a bit too restricted: they are assumed to be locally irreducible (see \cite[\S 1.4, Def. 3, p.337]{Rem57}). For these reasons, we will try to provide a proof here based on Remmert's original argument, for the convenience of the readers, although this is definitely well known to experts. 

First let us do some reductions: up to normalizing $X$ we can assume $X$ is normal, in particular, locally irreducible; moreover, since by our convention, all analytic varieties are assumed to be $\sigma$-compact, the result is local over $Y$, then we may assume $Y$ is a locally analytic subset of $\CC^N$, and hence assume $Y=\CC^N$. By these reduction process, we can apply the results of \cite{Rem57} to $f:X\to Y$. In particular, by \cite[\S 3.2, Satz 17, p.348]{Rem57} we obtain a finite chain of analytic subsets 
\[
\emptyset=X(-1)\subseteq X(0)\subseteq\cdots\subseteq X(r)=X,
\] 
where $X(i)$ denotes the subset of $X$ that consists of all points $x\in X$ such that $\rank_f(x)\le i$, where $\rank_f(x)$ denotes the rank of $f$ at $x$ and is defined to be the codimension of $f^{-1}(f(x))$ in $X$, and $r:=\sup_{x\in X}\rank_f(x)$ (see \cite[\S 3.1, Def. 13, p.346]{Rem57}). Now apply \cite[\S 4.1, Satz 19, p.351]{Rem57} to each $f|_{X(i)\backslash X(i-1)}$, we see that $f(X(i)\backslash X(i-1))$ is a countable union of locally analytic subset in $Y$, thus so is $f(X)$.
\end{proof}

\newpage

\bibliographystyle{alpha}
\bibliography{anti-nef_Kah_corrected}

\end{document}